\documentclass[12pt,twoside,a4paper]{amsart}
\usepackage{a4wide}
\usepackage{amssymb,amsthm,amsmath,xy,lscape}
\xyoption{all}
\usepackage[swedish,english]{babel}
\usepackage[latin1]{inputenc}
\usepackage{thumbpdf}
\usepackage[hypertex,
colorlinks=true, urlcolor=blue, pagecolor=black,
linkcolor=black, citecolor=black, filecolor=black, menucolor=black, bookmarks=true, bookmarksnumbered=true,
pdfpagelayout=TwoColumnRight, pdfstartview=FitH, plainpages=false,
pdfpagelabels]{hyperref}
\usepackage{helvet}

\font\fivesy=cmsy5 
\font\sevensy=cmsy7 \font\eightsy=cmsy8 \font\ninesy=cmsy9
\font\tensy=cmsy10 \font\fivei=cmmi5 \font\sixi=cmmi6
\font\seveni=cmmi7 \font\eighti=cmmi8 \font\ninei=cmmi9
\font\teni=cmmi10  \font\fivebf=cmbx5
\font\sixbf=cmbx6 \font\sevenbf=cmbx7 \font\eightbf=cmbx8
 \font\tenbf=cmbx10 \font\twelvbf=cmbx12
\font\fiverm=cmr5 \font\sixrm=cmr6 \font\sevenrm=cmr7
\font\eightrm=cmr8 \font\ninerm=cmr9 \font\tenrm=cmr10

\font\eightmsbm=msbm8

\font\tenib=cmmib10

\makeatletter
\newsavebox{\rotbox}

\theoremstyle{plain}
\newtheorem{theorem}{Theorem}[section]
\newtheorem{proposition}[theorem]{Proposition}
\newtheorem{lemma}[theorem]{Lemma}
\newtheorem{corollary}[theorem]{Corollary}
\newtheorem*{varth}{Theorem}
\newtheorem*{varcor}{Corollary}
\newtheorem*{varprop}{Proposition}
\newtheorem*{varlemma}{Lemma}
\theoremstyle{definition}
\newtheorem{definition}[theorem]{Definition}
\newtheorem*{vardef}{Definition}
\newtheorem*{varex}{Example}
\newtheorem{note}[theorem]{Note}
\newtheorem*{varnote}{Note}

\theoremstyle{remark}

\newcommand {\Cbf}{{\hbox{\eightbf C}}}
\newcommand {\Xit}{{\hbox{\eighti  X}}}
\newcommand {\Yit}{{\hbox{\eighti  Y}}}

\newcommand {\Deight}{\hbox{\scriptsize\bf$\Delta$}}
\newcommand {\oadd}{\hbox{\scriptsize\bf$\oplus$}}

\newenvironment{varnorem}{\noindent  \rm}{}

\makeatother

\begin{document}
\title{On the Foundation of Algebraic Topology}
\author{G\"oran Fors}
\address{Department of Mathematics, University of Stockholm, Se-104 05
Stockholm, Sweden. goranf@math.su.se}

\begin{abstract}{
In the 70:th, combinatorialists begun to systematically relate simplicial complexes and polynomial algebras, named Stanley-Reisner rings or face rings. This demanded an algebraization of the simplicial complexes, that turned the
empty simplicial complex into a zero object w.r.t. to simplicial join, losing its former role as join-unit - a role taken over by a new (-1)-dimensional simplicial complex containing only the empty simplex.

There can be no realization functor targeting the classical category of topological spaces that turns the contemporary simplicial join into topological join unless a (-1)-dimensional space is introduced as a topological join-unit.
This algebraization of general topology enables a homology theory that unifies the classical relative and reduced homology functors and allows a Künneth Theorem for simplicial resp. topological pair-joins.
}\end{abstract}

\maketitle

\tableofcontents
\def\myl{{\vcenter{\hrule width 0.13 true in height 0.005 true in}}}
\def\myline{$\myl$}
\def\mmyline{\myl}

\def\mysqcuplus{{\rlap{$\sqcup$}{\raise1.5pt\hbox{\sevenbf +}}}}
\def\mysqcupplus{\mysqcuplus}
\def\sqcupP{\ \!{{\rlap{$_{_{\!}}\sqcup
$}{\vbox{\moveleft2.8pt\hbox{{\ }\raise1.0pt\hbox{\fivebf P}}}}} } }
\def\mysqcupP{\sqcupP}
\def\Psqcup#1{\ \!{{\rlap{$_{_{\!}}\sqcup
$}{\vbox{\moveleft2.0pt\hbox{{\ }\raise1.0pt\hbox{\fivebf
\hbox{#1}}}}}} } }
\def\myPsqcup#1{\Psqcup#1~}

\def\mmysubsetneqq{\subsetneqq}
\def\mmysupsetneqq{\supsetneqq}

\def\nabl#1
{\font\mysyfont=cmsy#1
{
{\!_{\!}}^{_{_{\hbox{\mysyfont \char"72}}}}\!}}
\def\mynabla#1{$\nabl#1 $}
\def\mmynabla#1{\nabl#1 }

\def\barnabl#1{\!\bar{\ \!\nabl#1 }}
\def\mybarnabla#1{\barnabl#1 }

\def\thickl{{\vcenter{\hrule width 0.1 true in height 0.015 true in }}}
\def\thickline{$\thickl$}
\def\mthickline{\thickl}

\def\DEFhatH{{
\hat{\hbox{\bf H}}}}
\def\hatH{$\DEFhatH$}
\def\mhatH{\DEFhatH}

\def\eDEFhatH{{
\hat{\hbox{\eightbf H}}}}
\def\ehatH{$\eDEFhatH$}
\def\mehatH{\eDEFhatH}

\def\DEFsmallhatH{{
\hat{\hbox{\eightbf H}}}}
\def\smallhatH{$\DEFsmallhatH$}
\def\msmallhatH{\DEFsmallhatH}

\def\DEFminihatH{{
\hat{\hbox{\sixbf H}}}}
\def\minihatH{$\DEFminihatH$}
\def\mminihatH{\DEFminihatH}

\def\FFrame#1#2#3{                   %
     \vbox{\hrule height#2           
          \hbox{\vrule width#2       %
                 \hskip#1            
                 \vbox{\vskip#1{}    %
                       #3            
                       \vskip#1}     %
                 \hskip#1
                 \vrule width#2}
           \hrule height#2}}

\def\DEFhatbar{\hbox{$\mathbf{\hat{\vert}}$}}
\def\hatbar{$\DEFhatbar$}
\def\mhatbar{\DEFhatbar}

\def\DEFhatSigma{\hbox{$\mathbf{\hat{\Sigma}}$}}
\def\hatSigma{$\DEFhatSigma$}
\def\mhatSigma{\DEFhatSigma}

\def\DEFHatAst{{
{
{\hbox{$^{{\land}}$}}{\vbox{\hbox{$\!\!\!\ast$}}}}
}}
\def\hatast{$\DEFHatAst$}
\def\mhatast{\DEFHatAst}

\def\DEFsmallHatAst{{\ \!\!{\rlap{{\lower3.5pt\hbox
{\vbox{\moveright-0.5pt\hbox{$^{^{{\hbox{\fivesy \char"5E}}}}$}}}}}
{{{\hbox{\fivesy \char"03}}}}}\ \! }}
\def\topast{$\DEFsmallHatAst$}
\def\mtopast{\DEFsmallHatAst}

\def\DEFsmallHattAst{{\ \!\!{\rlap{{\lower2.0pt\hbox
{\vbox{\moveright-0.5pt\hbox{$^{^{{\hbox{\fivesy \char"5E}}}}$}}}}}
{{{\hbox{\sevensy \char"03}}}}}\ \! }}
\def\toppast{$\DEFsmallHattAst$}
\def\mtoppast{\DEFsmallHattAst}

\def\DEFCircast{{{\rlap{\lower4.5pt\hbox{$^{^{_{{_{\!}}\bigcirc}}}$}}
{\vbox{\moveright-0.5pt\hbox{$\ast$}}}}}}
\def\circast{$\DEFCircast$}
\def\mcircast{\DEFCircast}

\def\DEFsimpline{\bullet\!\!\!-\!-\!-\!\!\!\bullet}
\def\simplline{$\DEFsimpline$}
\def\msimplline{\DEFsimpline}

\def\DEFsimpsquare{{\hbox{\hsize1.23cm\hskip0.1cm
\raise0.7pt\vbox{\hrule width 1.0pt height 0.3cm depth
0.32pt}\hskip-0.48cm
\vbox{$_{\hbox{$\bullet\!\raise2.8pt\vbox{\hrule width 0.25cm height 0.7pt depth 0.32pt}\!\bullet$}}
^{^{\hbox{$\bullet\!\raise2.8pt\vbox{\hrule width 0.25cm height
0.7pt depth 0.32pt}\!\bullet$}}}$}
\hskip-0.5cm\raise0.7pt\vbox{\hrule width 1.0pt height 0.3cm depth
0.32pt}\hskip0.1cm
}
}}
\def\simplsquare{$\DEFsimpsquare$}
\def\msimplsquare{\DEFsimpsquare}

\def\DEFringhom#1{
{\lower4.0pt\hbox{$^{_{\hskip0.05cm\hbox{{\sevenbf #1}}
\atop{\hbox{$\cong$}}}}$}}}
\def\ringHom#1{$\DEFringhom#1$}
\def\mringHom#1{\DEFringhom#1}

\def\DEFringhombb#1{
{\hbox{$^{_{{\mathbb{#1}}
\atop{\hbox{$\cong$}}}}$}}}
\def\ringHombb#1{$\DEFringhombb#1$}
\def\mringHombb#1{\DEFringhombb#1}

\def\DEFsmallRinghom#1{
{\lower4.0pt\hbox{$^{_{\hbox{{\fivebf #1}}
\atop{{\raise1.4pt\hbox{{\eightbf =}}}
\hskip-0.25cm{\raise3.8pt\hbox{{\eightsy {\char"18}}}} }}}$}}}
\def\smallRingHom#1{$\DEFsmallRinghom#1$}
\def\msmallRingHom#1{\DEFsmallRinghom#1}

\def\DEFbfdel#1{\hbox{\tenbf#1$\hskip-3.0pt$#1$\hskip-3.0pt$#1}}
\def\bfDel#1{$\DEFbfdel#1$}
\def\mbfDel#1{\DEFbfdel#1}

\def\DEFbigbfdel#1{\hbox{\tenbf\big#1$\hskip-4.5pt$\big#1$\hskip-4.5pt$\big#1}}
\def\BfDel#1{$\DEFbigbfdel#1$}
\def\mBfDel#1{\DEFbigbfdel#1}

\def\DEFbfbrace#1{\hbox{$#1$\hskip-5.8pt$#1$\hskip-5.8pt$#1$}}
\def\bfBrace#1{$\DEFbfbrace#1$}
\def\mbfBrace#1{\DEFbfbrace#1}

\def\DEFbigbfbrace#1{\hbox{$\big#1$\hskip-6.75pt$\big#1$\hskip-6.75pt$\big#1$}}
\def\BfBrace#1{$\DEFbigbfbrace#1$}
\def\mBfBrace#1{\DEFbigbfbrace#1}

\def\DEFbfCupCap#1#2{\hbox{$#1$\hskip-0.33cm$\lower3pt\hbox{$_{_{#2}}$}$}}
\def\bfcupcap#1{$\DEFbfCupCap#1$}
\def\mbfcupcap#1{\DEFbfCupCap#1}

\def\mytinynabla{\lower1.0pt\hbox{$^{_{_{\hbox{\sevensy \char"72} }}}$}}

\def\LHS{{\bf L\hatH S}\hskip2pt}

\part{Introduction}\label{SecP1:1}
\section{The Algebraic Topological Foundation is in Need of a Makeover}
\subsection{Simplicial complexes within Combinatorics}

In Gräbe's paper \cite{12} \S 4.1 p.~171, the reduced homology groups, with field coefficients, of the simplicial join of two complexes were given as
\htmladdnormallink{{\it direct sums} }{http://en.wikipedia.org/wiki/Direct_sum}
of \htmladdnormallink{{\it tensor products,}}{http://en.wikipedia.org/wiki/Direct_sum}
the relative analogue of which seemed non-existing, cf. p.~\pageref{UnderavdP5:KynnethFormForJoins}:
\vskip-0.2cm
$$\
\widetilde{H}\!\!_{_{\!q+1\!}}
        (\Sigma_{1}\ast \Sigma_{2};{\bf k})
\cong
{\rlap{$_{_{_{_{i+j=q}}}}$}} {\ \raise0.0pt\hbox{$\oplus$}}\
\hbox{\bf(}\widetilde{H}_{i}
           (\Sigma_{1};{\bf k})
\otimes\!_{_{\bf{k}}}\widetilde{H}_{\!j}
(\Sigma_{2};{\bf k})\hbox{\bf)}.
$$
\cite{12} \S 4.2 p.~171 gives the boundary of the simplicial join of quasi-manifolds:
$${\rm Bd} (\Sigma_{1}
\ast
\Sigma_{2})= (({\rm Bd} \Sigma_{1})
\ast
\Sigma_{2})\cup (\Sigma_{1}
\ast
({\rm Bd} \Sigma_{2})).
$$
This formula only works under Definition~1 below, partly due to the fact that manifolds now can have ${\emptyset}$ as well as the ``new" $(-1)$-dimensional join unit $\{\emptyset\}$ as its boundary.

\goodbreak
\begin{vardef}{1)} \label{DefP1:1}
\noindent An
\htmladdnormallink{(augmented abstract){\it simplicial complex}}
{http://en.wikipedia.org/wiki/Abstract_simplicial_complex}
$\Sigma$ on a vertex set {\bf V}$_{^{_{\!\!\Sigma}}}\!$ is a
collection of finite subsets $\sigma$, the {\it simplices}, of {\bf V}$_{^{_{\!\!\Sigma}}}\!$ satisfying:
$$ {\rm If}\ v\in\ {\bf V}_{^{_{\!\!\Sigma}}}, {\rm then}\ \{v\}\in
\Sigma.\leqno{\bf (a)}$$
\vskip-0.8cm
$$ {\rm If}\ \sigma\in \Sigma\ {\rm and}\
 \tau\subset \sigma,\ {\rm then}\ \tau\in \Sigma.\leqno{\bf (b)}$$
\smallskip
\noindent{\rm Define the} {\it simplicial join} ``$\Sigma_1\ast
\Sigma_2$'', of two complexes $\Sigma_1$ and
$\Sigma_2$ with
$V_{\Sigma_{1}}\cap V_{\Sigma_{2}}=\emptyset$
to be:
{$$\Sigma_1\ast\Sigma_2:=\{\sigma_1\cup\sigma_2\mid\
\sigma_i\in\Sigma_i \ (i=1,2) \}.$$}
\end{vardef}
Definition\ 1 has three levels: a {\it complex} (level~1)
is a set of {\it simplices} (level~2), which are {\it finite} sets (empty or non-empty) of {\it
vertices} (level~3).
The empty set, $\emptyset$, plays a dual role; it is both a $(-1)$-dimensional
simplex and a $(-\infty)$-dimensional complex.
See p.~\pageref{DefP7:SimpDim12}
for the definition of ``dimension''.

\begin{vardef}{2)} \label{DefP2:2}
{\rm(\cite{30}\ p.~108-9)}
\noindent A
\htmladdnormallink{\underbar{classical} $($abstract$)$ {\it simplicial complex}
}{http://planetmath.org/encyclopedia/SimplicialComplex.html}
$\Sigma$ on a vertex set {\bf V}$\!_{_{^{\!\!\Sigma\!}}}$
is a collection of \underbar{non-empty} finite subsets $\sigma$, the {\it simplices}, of {\bf V}$\!_{_{^{\!\!\Sigma\!}}}$ satisfying:
$$ {\rm If}\ v\in\ {\bf V}_{^{_{\!\!\Sigma}}}, {\rm then}\ \{v\}\in
\Sigma.\leqno{\bf (a)}$$
$$ {\rm If}\ \sigma\in \Sigma\ {\rm and}\
 \tau\subset \sigma,\ {\rm then}\ \tau\in \Sigma.\leqno{\bf (b)}$$

The \underbar{\it classical} {\it simplicial join} ``$\Sigma_1\ast \Sigma_2$'', of two disjoint \underbar{\it classical}
complexes is defined to be:
$$\Sigma_1\ast\Sigma_2:=\Sigma_1\cup\Sigma_2\cup\{\sigma_1\cup\sigma_2\mid\
\sigma_i\in\Sigma_i \ (i=1,2) \}.$$
\end{vardef}

If $
{\rm card}(\sigma)=q+1$ for a simplex $\sigma$ then $\dim\sigma:=q$
\label{DefP7:SimpDim12}and $\sigma$ is said to be a
$q$-{\it face} or a
$q$-{\it simplex} of $\Sigma$ and
$\dim\Sigma:=$sup$\{\dim (\sigma)\vert \sigma\in
\Sigma\}$.
Writing\ $\emptyset{_{^{\!_o}}}\!$ when using $\emptyset$ as a
simplex, we get dim$(\emptyset)\!=\!-\infty$ and
$\dim(\{\emptyset{_{^{\!_o}}}\!\})\!=\!\dim(\emptyset{_{^{\!_o}}}\!)\!=\!-1.$

So, the dimension of any non-empty simplicial complex is a well-defined integer $\ge-1$,which for the join of two
simplicial complexes with disjoint vertex sets results in the following dimension formula:
$$
\dim(\Sigma_1\ast\Sigma_2)=\dim\Sigma_1 + \dim\Sigma_2 +1
$$
implying that the join unit must be of dimension $-1$,which by definition is the dimension of the augmented abstract simplicial complex $\{\emptyset\}$.

The difference between Definition~1 and Definition~2, is
that the former allows the empty set $\emptyset$ as a {\it simplex} that moreover is contained in every nonempty simplicial complex,
implying that the simplicial complex generated by any vertex
$\hbox{\it v}$, called a (combinatorial) $0$-ball, is the set
$\bullet:=\{\emptyset, \{v\}\}$, resp. $\{\{v\}\}$.
Similarly, the (combinatorial) $0$-sphere is the set
$\bullet\bullet:=\{\emptyset, \{v\}, \{w\}\}$,
resp. $\{\{v\}, \{w\}\}$.
In dimension $1$ and below, the boundary Bd of any reasonably
defined simplicial manifold of dimension $n$ is generated by the simplices of dimension $n-1$ that is a face of only one simplex of dimension $n$.
So, in particular, ${\rm Bd}(\bullet)=\{\emptyset\}\ne\emptyset$ according to Definition~1 but ${\rm Bd}(\bullet)=\emptyset$ according to Definition~2, while ${\rm Bd}(\bullet\bullet):=\emptyset$ with respect to either of the two definitions.

\noindent
\framebox
{\vbox{
The classical algebraic topological structure needed two separate homology functors, the relative and the ad-hoc invented reduced homology functor.
The change of combinatorial foundation from Definition~2 to Definition~1 took place ($\approx$1970) the minute combinatorialists began to use (``clashed into'', cf. \cite{MacLane})
commutative algebra to solve combinatorial problems.
The only homology apparatus invented by the combinatorialists
to handle this ingenious algebraization of the category of simplicial complexes were a jargon like -``We will use reduced homology with $\widetilde{\rm H}_{-1}(\{\emptyset\})={\mathbb{Z}}$'', which added even more inefficiency to algebraic topology than it already suffered from as a result of that unfortunate but necessary classical use of the \underbar{two} complementing homology functors just mentioned.

The category of topological spaces and continuous functions were left unattended. So, suddenly there were two vertex free simplicial complexes, $\emptyset$ and $\{\emptyset\}$, but still only one point-free topological space, i.e. $\emptyset$, which made it impossible to define, in any reasonable way, a faithful realization functor.

This unfortunate state of affairs is rectified in this article by performing an analogous algebraization of general topology by simply introducing a $(-1)$-dimensional topological join-unit (denoted "$\{\wp$\}") that also is the realization of the $(-1)$-dimensional simplicial complex $\{\emptyset\}$.
Now {\it the topological join} has
\htmladdnormallink{{\it the tensor product} }{http://en.wikipedia.org/wiki/Tensor_product} as its functorial counterpart within the tensor categories, which underlines its central role within algebraic topology.

The resulting homology theory, i.e. {\it the augmental homology theory}, is complete in the sense that no {\it reduced} homology functor is needed.}}

\subsection{Looking elsewhere for ``simplicial complexes''}\label{UnderavdP3:CheckingLiterature}

In \cite{11} Ch.~3
p.~110, Fritsch and Piccinini give the following definition and example.
\begin{quote}
\noindent
{\bf Definition\ 3}
{\rm(Equivalent to our Definition\ 1.)}
\noindent A  {\it Simplicial Complex} is a set $\Sigma$ of finite
sets closed under the formation of subsets, i.e., any subset of a
member of $\Sigma$ is also a member of $\Sigma$; more formally:
$$(\sigma\!\in \!\Sigma)\ \!\land\ \!
 (\tau \!\subset\! \sigma)\ \Longrightarrow\ \tau \!\in\! \Sigma.$$
\noindent
{\bf Example}
(\cite{11} p.~110)
Let $\{U_\lambda: \lambda\in\Lambda\}$ be a family of arbitrary sets; then the set $K(\Lambda)$ of all finite subsets of $\Lambda$ such that
$\bigcap_{\lambda\in x}{U_\lambda} \neq \emptyset $
is a simplicial complex. Note that in general the vertex set of this simplicial complex $K(\Lambda)$ is not the index set $\Lambda$ itself, but only its subsets consisting of the indices $\lambda$ with ${U_\lambda} \neq \emptyset$.
Now, if $Z$ is a space and $\{U_\lambda: \lambda\in\Lambda\}$ is a covering of $Z$ (see Section A.3), then the simplicial complex $K(\Lambda)$ obtained in this way is called the {\it nerve} of the covering $\{U_\lambda: \lambda\in\Lambda\}$.
\end{quote}
\begin{varnote}
$\bigcap_{\lambda\in \emptyset}{U_\lambda} = $
Universe $\!\neq\emptyset\Rightarrow \!x\!=\!\emptyset\!\in\!
K(\Lambda)$ (cf.~{\rm \spaceskip2pt\cite{19}, Kelley: {\it General Topology}, p.~256).}

So, we note that any
\htmladdnormallink{{\it nerv} }{http://en.wikipedia.org/wiki/Nerve_of_an_open_covering}
$K(\Lambda)$ constructed on a nonempty
index set $\Lambda$ within any nonempty universe,
is nonempty, since it contains at least the {\it simplex}
$\emptyset$.
\end{varnote}
\cite{11} is an excellent book and due to its detailed excellence we have a chance to dwell somewhat on a difference between its chapter~3 and chapter~4.

In \cite{11} Ch.~3
the empty simplex is allowed, as seen above, but the
realization of the complex $\{\emptyset\}$ is nowhere explicit. This new simplicial complex $\{\emptyset\}$ has no natural realization within any Euclidian $n$-space within the old classical category of topological spaces and continuous maps.

In \cite{11} Ch.~4
however, Fritsch and Piccinini use the ``topologist's $\Delta$'', here denoted $\hbox{\tenib\char"01}\!^\circ$, which they call the {\it category of finite ordinals}. This $\hbox{\tenib\char"01}\!^\circ$ actually excludes the empty simplex and therefore also the complex $\{\emptyset\}$, which makes the realization functor identical to the classical and therefore non-problematic.

Next definition, equivalent to our Definition~2, is given in
p.~2 in \cite{JandT}.
\begin{quote}
\noindent {\bf Definition\ 4.} A {\it Simplicial Complex} $K$ is a
collection of {\it non-empty}, finite subsets (called {\it simplices}) of a given set $V$ (of {\it vertices}) such that any {non-empty} subset of a simplex is a simplex.
\end{quote}
Deleting the word ``{\it non-empty}'' makes it equivalent to our Definition~1, i.e.:

\smallskip
\noindent {\bf Definition\ 5.} An (augmented abstract) {\it
Simplicial Complex} $K$ is a collection of finite subsets (called {\it simplices}) of a given set $V$ (of {\it vertices}) such that any subset of a simplex is a simplex.

\section{Remaking the Topological Foundation}

\subsection{Synchronizing Logic, Combinatorics, General
Topology and Algebra}\label{UnderavdP6:RelGenTopToComb}\hfill\break\indent
The ``$\sigma_1\cup\sigma_2$''-operation in Definition~1 p.~\pageref{DefP1:1}
is strongly ``ordinal addition''-related, see p.~\pageref{P4:OrdinalAddition},
and we learn from Note~iii p.~\pageref{NoteP21:iii}
that the Stanley-Reisner (St-Re) ring of a join of two simplicial complexes is the tensor product of their respective St-Re rings.
This simple observation reveals a very strong relation between logics and algebra that also involves combinatorics and general topology.
Still, the join operations is very much ignored in the literature, except of course for
\cite{2},
\cite{Ehlers&Porter}\
and \cite{Fritsch&Golasinski}.

A moment of reflexion on ``realization functor''-candidates
$|\cdot|$, reveals the need for a real topological join-unit,
$\{\wp\}=|\{\emptyset\}|\neq|\emptyset| =\nolinebreak\emptyset$,
which also becomes our new $(-1)$-dimensional geometrical
\htmladdnormallink{standard simplex.}
{http://en.wikipedia.org/wiki/Simplex}

We are obviously dealing with a non-classical situation, which we
must not try to squeeze into a classical framework.
The
\htmladdnormallink{category}{http://plato.stanford.edu/entries/category-theory/} of
\htmladdnormallink{topological spaces}{http://en.wikipedia.org/wiki/Topological_space}
and continuous maps has been
unchanged ever since the publication of
\htmladdnormallink{Felix Hausdorff's}{http://en.wikipedia.org/wiki/Felix_Hausdorff}
{\it Grundzüge der Mengenlehre} (1914),
but now, modern algebra-based techniques makes the introduction of an $(-1)$-dimensional topological space $\{\wp\}$ unavoidable.
The result of such a $\{\wp\}$-introduction is an algebraization of
\htmladdnormallink{{\it general}}
{http://en.wikipedia.org/wiki/General_topology}
and
\htmladdnormallink{{\it algebraic topology}}
{http://en.wikipedia.org/wiki/Algebraic_topology}, i.e. a synchronization of these categories to that of the
\htmladdnormallink{{\it tensor categories}
}{http://en.wikipedia.org/wiki/Monoidal_category}
within commutative algebra that is
completely compatible with the switch from Definition~2 to Definition~1 within combinatorics in the 70th.

The pure ``Def.~2''-based
\htmladdnormallink{{\it homology theory},}
{http://en.wikipedia.org/wiki/Homology_theory}
i.e. classical homology theory, uses a combination of the
\htmladdnormallink{{\it relative homology theory}}
{http://en.wikipedia.org/wiki/Relative_homology} \
and the ad-hoc invented
\htmladdnormallink{{\it reduced homology functor},}
{http://en.wikipedia.org/wiki/Reduced_homology} the mere existence of which is a proof of the insufficiency of classical relative homology theory.
These two functors are simply unified into one single homology functor $\mhatH$ through the external adjunction, to the classical category of topological spaces, of the non-empty $-1$-dimensional object $\{\wp\}$ - unit element with respect to topological join. This is, trough $\{\wp\}:=|\{\emptyset\}|$, fully compatible with the $\{\emptyset\}$-introduction in Definition~1 above: ``\underbar{Externally} adjoined'' means that $\wp$ is assumed not to have been an element of any \underbar{classical} topological space.

The category of sets mustn't be tempered with, but add, using topological sum, to each classical topological space $X\in \mathcal{D}$ an external element $\wp$, resulting in $X_{\!_{^{\wp}}}\!:=X \!+ \{\wp\}\in\mathcal{D}\!_{\wp}$.
Finally we add the universal initial object $\emptyset$.
This resembles a familiar routine from Homotopy Theory providing all free spaces with a common base point, cf. \cite{35} p.~103.
Working with $X_{\!_{^{\wp}}}\!$ instead of $X$ makes it possible, e.g., to deduce the Künneth Formula for joins of topological \underbar{pair}-spaces, as in Theorem~\ref{TheoremP14:4} p.~\pageref{TheoremP14:4}.

The simplicial as well as the topological join-operation are (modulo realization, equivalent) cases of \label{FunctorialityOfJoin}
\htmladdnormallink{colimits,}{http://en.wikipedia.org/wiki/Limit_(category_theory)}
 as being simplicial resp. topological attachments.
Restricting to topological $k$-spaces, \cite{11} p.~157ff remains true also with our new (cf. p.~\pageref{DefP8:alpha0}) realization functor, i.e. {\it geometric realization preserves finite limits and all colimits.}

Let ${{\mathcal K}}$ be the classical category of simplicial
complexes and simplicial maps resulting from Definition~2 p.~\pageref{DefP2:2}, and let
${{\mathcal K}\!_{o}}$ be the category of simplices resulting from Definition~1 p.~\pageref{DefP1:1}.
Define ${\mathcal{E}}_{o}\!:\!{\mathcal K}\rightarrow{{\mathcal
K}\!_{o}}$\label{DefP2:E0} to be the functor adjoining $\emptyset$,
as a simplex, to each classical simplicial complex.
$\!{\mathcal{E}}_{o}\!$ has an inverse ${{\mathcal{E}}}\!:\!{\mathcal
K}\!_{o}\!\!\rightarrow\!{\mathcal K}\ \hbox{\rm deleting}\
{\emptyset}$.
This idea is older than category theory itself. It was suggested by
S. Eilenberg and S. MacLane in \cite{8} (1942),
where they in page 820, when exploring the 0-dimensional homology of
a simplicial complex $K$ using the {\it incidence numbers} $[\ :\ ]$,\nolinebreak\ wrote:
\begin{quote}
An alternative procedure would be to consider $K$ ``augmented'' by a single
$(-1)$-cell $\sigma^{_{_{-1}}}$ such that
$[\sigma_i^{_{_{0}}}:\sigma^{_{_{-1}}}]=1$ for all
$\sigma_i^{_{_{0}}}$.
\end{quote}
This approach is explored by Hilton and Wylie in \cite{17} p.~16~ff.

As for ${\mathcal{E}}_{o}\!$ above,
let ${{\mathcal F}\!_{\wp}}\!:\!{\mathcal D}\rightarrow{{\mathcal
D}\!_{\wp}}$ be the functor adjoining to each classical topological space a
\htmladdnormallink{{\it new non-final element}}
{http://en.wikipedia.org/wiki/Initial_object}
denoted $\wp$.
Put ${X}_{\!\wp}:= X + \{\wp\}$, which makes ${X}_{\!\wp}$
resemble Maunder's modified $K^+$-object defined in \cite{22}
p.~317 and further explored in pp.~340-341.

$\!{\mathcal F}\!_{\wp}\!$ has an inverse ${\mathcal F}:{\mathcal D}\!_{\wp}\!\longrightarrow\!{\mathcal D}$ deleting ${\wp}$.

\smallskip
Our Proposition~\ref{PropP15:1} p.~\pageref{PropP15:1},
partially quoted below, reveals a connection between the combinatorial and the topological local structures and is in itself a strong motivation for introducing a topological (-1)-dimensional object $\{\wp\}$, imposing the following definition p.~\pageref{DefP8:PointSetMinus0}
of a ``{\it point-setminus}'' ``$\setminus{_{_{^{\!o\!}}}}$'', in ${{\mathcal D}\!_{\wp}}$ using the convention $x\leftrightarrow\{x,\wp\}$, as opposed to the classical $x\leftrightarrow\{x\}$.
``$\smallsetminus$" denotes the classical ``{\it setminus}'', also known in set theory as the \htmladdnormallink{``{\it relative}}
{http://en.wikipedia.org/wiki/Complement_(set_theory)}
\htmladdnormallink{{\it complement}''}
{http://en.wikipedia.org/wiki/Complement_(set_theory)}
while ``$\emptyset{_{^{\!_o}}}\!$'' denotes the {\it empty set}, $\emptyset$, when regarded as a \underbar{simplex} and \underbar{not} as a simplicial complex.

We can not risk dropping out of our new category so for
${X}\!_{\wp}\neq\emptyset$;
$$
{X}_{\!\wp}\setminus\!{_{_{^{\!o\!}}}}{x}:=
\left\{\begin{array}{ll}
\emptyset  & \text{if}\ \ x=\wp \\
{{\mathcal F}\!_{\wp}}\bigl({X}
\smallsetminus x\bigl)& \text{otherwise}
\end{array}\hskip-0.1cm\right.
$$

In Part {\bf \ref{PartIII}} we will find constant use for the
{\it link-construction}, denoted ${{\rm Lk}}_{_{^{\Sigma\!\!}}}\sigma$, of
a simplex $\sigma$ with respect to a simplicial complex $\Sigma$ and defined by:
$$
{{\rm Lk}}_{_{^{\Sigma\!\!}}}\sigma\!:=\! \{\tau\in \Sigma\mid
[\sigma\cap \tau =\emptyset]\land [\sigma\cup \tau \!\in\!
\Sigma]\}\indent\indent(\Rightarrow{\rm
Lk}_{_{\!\Sigma}}\emptyset\!_{o}\!\!=\Sigma).$$

We will also often use:
$$``\textrm{The {\it contrastar} of}\
\sigma\!\in\!\Sigma\hbox{''}\!=\hbox{\rm cost}_{_{\!{\Sigma}}}\!\sigma\!:=\!
\{\tau\!\in\!\Sigma\mid\ \tau\!\not\supseteq\sigma\},$$
implying that:
$\hbox{\rm
cost}_{_{{\!\Sigma}}}\!\emptyset\!_{o}\!\!=\!\emptyset\!$
\ and\
${\hbox{\rm cost}_{_{{\!\Sigma}}}\!\sigma\!=\!\Sigma\
\underline{\hbox{\rm iff}}\ \sigma\!\not\in\!\Sigma}.$
In classical literature, the {\it contrastar}, or {\it costar}, of $\Sigma$ w.r.t. $\sigma$ is known as the {\it complement of} $\Sigma$ w.r.t. $\sigma$, e.g., see \cite{9} p.~74 exercise E. See also \cite{Mi}.

In p.~\pageref{DefP8:alpha0} we define {\it the realization} $\vert\Sigma\vert$ of the complex $\Sigma$ by slightly modifying Spanier's {\it realisation-functor}, which, relative to the classical definition, in the range space exhibits an additional {\it coordinate function} $\alpha_0$ with $\alpha_0(v)\equiv 0\ \forall\ v\in {\bf V}\!\!_{_{^{\Sigma}}}$.

\begin{varprop}\hskip-0.1cm{\bf\ref{PropP15:1}}
{\rm (p.~\pageref{PropP15:1})}
Let ${\bf G}$ be any
\htmladdnormallink{module}{http://en.wikipedia.org/wiki/Module_(mathematics)}
over a \htmladdnormallink{commutative ring}{http://en.wikibooks.org/wiki/Abstract_algebra/Rings,_fields_and_modules}
{\bf A} with unit. With $\alpha\in\{\alpha\in|\Sigma|\mid[{\bf
v}\in\sigma ] \Longleftrightarrow [\alpha({\bf v})\neq0]
\}$ and $\ {\!}\alpha_{\!}=_{\!}
\alpha_{_{^{\!0}}}\!$ {\underbar{\hbox{\rm iff}}}\
$\sigma\!=\!\emptyset\!_{_{^{o}}}{\!}$,
the following {\bf A}-module isomorphisms are all induced by chain equivalences:
\begin{eqnarray}
\mhatH
\hbox{$_{i-\#\sigma}$}
\hbox{$({\rm Lk}_{{\Sigma}}\sigma;{\bf
G})$}
\mringHom{A}
\mhatH_{i}(\Sigma,\hbox{\rm cost}_{{\Sigma}}\sigma;{\bf G})
\mringHom{A}
\mhatH_{i}(\vert\Sigma\vert,\vert\hbox{\rm
cost}_{{\Sigma}}\sigma\vert;{\bf G})
\mringHom{A}
\mhatH_{i}(|\Sigma|,|\Sigma|\setminus\!{_{_{^{\!o\!}}}}\alpha;{\bf
G}).\nonumber
\end{eqnarray}
{\rm (The $\mhatH$ in the first two $\mhatH$omology groups denotes the simplicial $\mhatH$omology functor, while the two last represent the singular.
``\#'' stands for ``{\it cardinality}''.
The proof depends on \cite{26} p.~279 Th.\ 46.2 quoted here in p.~\pageref{TheoremP13:Munkres}.)}
\end{varprop}
\begin{varnorem}
{\it Astro-Physical inspiration}:
Currently, astronomers hold it likely that every galaxy,
including the Milky Way, possesses a single
\htmladdnormallink{{\it supermassive black hole},}{http://en.wikipedia.org/wiki/Supermassive_black_hole}
active or inactive,
here interpreted as $\{\wp\}$, and so, eliminating it, would dispose off the whole galaxy, hinting at the above
$\!``\!{X}_{\!\wp}\!\setminus\!{_{_{^{\!o\!}}}}{\wp}\!:=\emptyset\hbox{''}\!$.
Actually, the last proposition leaves us with no other option, since ${\rm Lk}_{_{\!\Sigma}}\emptyset\!_{o}\!\!=\Sigma$ by definition.
\end{varnorem}

\subsection{Exploring the boundary formula}\label{UnderavdP7:ExplBundFormula}

The absence of {\it a unit-space with respect to topological join} cripples classical algebraic topology and this shortcoming can not be eliminated by the use of any technique from any mathematical field. The only remedy is to invent such a unit space $\{\wp\}$.
We use the `Definition~1'-generated object $\{\emptyset_o\}$ from p. \pageref{DefP1:1}
as a model-object and define $\{\wp\}$ to be its realization, i.e. set: $\{\wp\}:=\vert\{\emptyset_o\}\vert$.

Any link of a simplicial homology manifold is either a simplicial homology ball or a sphere.
Its boundary is the set of simplices having ball-links, i.e.
$${\rm Bd}_{_{^{\!\ }}}\!\Sigma:=
\{\sigma\in \Sigma \mid\
\mhatH\hskip-0.3cm_{_{^{n-\#\sigma}}}\!(\hbox{\rm
Lk}_{_{^{\!\Sigma\!\!}}}\sigma)\!=0\}\indent (\approx\{\sigma\in \Sigma \mid \widetilde{\hbox{H}}\!\!\!\!\!\!_{_{^{n-\#\sigma}}}\!(\hbox{\rm
Lk}_{_{^{\!\Sigma\!\!}}}\sigma)\!=0\}\ \text{classically}).$$
\label{DefP2:BoundaryOfSimplManif}
Since, ${\rm Lk}_{_{\!\Sigma}}\emptyset\!_{o}\!\!=\Sigma$
(``{\it the missing link}''),
the
\htmladdnormallink{{\it real projective plane}}
{http://en.wikipedia.org/wiki/Real_projective_plane}
\ ${{\mathbb{RP}}}\!\!_{_{^{{\bf \
}}}}^{^{_{\bf \ 2}}}\!$
has the boundary
${\rm Bd}{{\mathbb{RP}}}^{^{_{\fivebf 2}}}=\{\emptyset_o\}\ne \emptyset$,
which also holds for any one-point space $\bullet$, i.e. ${\rm Bd}^{\!}(\bullet)\!=\!{\rm
Bd}^{\!}(\{\emptyset_o, \{{{v}}\}\})\!=\!\{\emptyset_o\}\!=$``the
join-unit''. So, the {boundary of a $0$-ball is the $-1$-}
sphere.
The one-point complex $\bullet:=\{\emptyset_o, \{{{v}}\}\}$ is the only finite
\htmladdnormallink{{\it orientable manifold}}
{http://en.wikipedia.org/wiki/Orientability}
having $\{\emptyset_o\}$ as its boundary.
For the $0$-sphere ${{\rm Bd}}(\bullet\bullet)\!=\!\emptyset^{\!}$ -
the join-zero with respect to Definition~1 p.~\pageref{DefP1:1}. The \underbar{classical} setting gives;
$$\hbox{${{\rm Bd}}(\{\{{v}\}\})\!\!=\emptyset=\!{{\rm Bd}}(\{\{{v}\}, \{{w}\}\})\!=\!\emptyset=$``the join-unit\hbox{''} with respect to Definition~2,}$$
implying that the boundary formula below, for the join of two
homology manifolds, can not hold classically but, as will be shown, will always hold in the ordinary $\wp$-augmented categories.
Moreover, for any $n$-dimensional locally orientable simplicial manifold, the boundary is generated by those $(n-1)$-dimensional simplices whose link contains exactly one vertex.
$${{{\mathcal{E}}}}({{\rm Bd}}\Sigma)={{\rm Bd}}({{{\mathcal{E}}}}(\Sigma)),$$
i.e., what remains after deleting $\emptyset_o$ is exactly the
classical boundary - always!
The new simplicial complex $\{\emptyset_o\}$ will serve as the ({\it abstract simplicial}) $(-1)$-{\it standard simplex}.

We will see that the join-operation in the category of simplicial complexes and simplicial sets in their augmented form, will turn into (graded) tensor product when applied to either our homology functor {\hatH} or the {\it Stanley-Reisner ring functor} p.~\pageref{DefP21}, which is made explicit here in p.~\pageref{Ex:HomologyVsSt-Re}.

We give a few low-dimensional examples how the non-classical approach enriches algebraic topology through the following manifold  boundary formula for topological as well as simplicial join, as given on p.~\pageref{TheoremP17:7} resp. p.~\pageref{theorem:TheoremP31:12}:\label{BundaryIntro}
$${{\rm Bd}} ({\hbox{{\it M}}}_{1}{\!}\ast_{\!}
{\it M}_{2}\!){\!}
={\!} (({{\rm Bd}} {\hbox{{\it M}}}_{1}\!)\ast {{\it
M}}_{2}{\!})\ {\!}\cup\ {\!} (_{\!}{\hbox{{\it
M}}}_{1}\!\ast ({{\rm Bd}} {\hbox{{\it M}}}_{2}\!)).
$$

$\{\emptyset_o,\{v_1\},\{v_2\}\{v_1,v_2\}\}
=\bullet\ast\bullet$ is a figurative way of saying that
the simplicial $1$-ball, generated by a $1$-simplex equals the join of two \underbar{distinct} simplicial $0$-balls,  each generated by a $0$-simplex, so:
$${\rm Bd}(1\text{-ball})
={\rm Bd}(\bullet\ast\bullet)=(\{\emptyset_o\}\ast\bullet)\cup(\bullet\ast\{\emptyset_o\})=
\bullet\cup\bullet=:\bullet\bullet=0\text{-sphere}.$$

$
{(\bullet\bullet)\ast(\bullet\bullet)}=$ (the
simplicial) 1-sphere. So, since ${\rm
Bd}(\bullet\bullet)=\emptyset$:
$${\rm Bd}(1\text{-sphere})
={\rm Bd}((\bullet\bullet)\ast(\bullet\bullet))
=(\emptyset\ast(\bullet\bullet))\cup((\bullet\bullet)\ast\emptyset)=
\emptyset\cup\emptyset=\emptyset.$$

The topological join might be thought of as the two factor-spaces along with arcs going from each point in one factor to each point in the other.

Our formula makes it trivial to calculate the boundaries of the following three low-dimensional joins, where
$\dim(X_1\ast X_2)=\dim X_1+\dim X_2\ +1,$
which, by the way, implies that a {\it unit space w.r.t. join} has to be of dimension $-1$.

\begin{enumerate}
\item $\mathcal{M} \ast \bullet=$
the cone of the M\"obius band $\mathcal{M}_{\!},$\ \
\smallskip
\item $\mathbb{RP}^2\ast\bullet=$
the cone of the real projective plane $\mathbb{RP}^2$,
\smallskip
\item $\mathbb{S}^2\ast \bullet=$ ``the cone of the 2-sphere
$\mathbb{S}^2\hbox{''}={\mathbb{B}}^3$,
the solid $3$-dimensional ball or $3$-ball.
\end{enumerate}
\goodbreak

The real projective plane $\mathbb{RP}^2$ frequently appears in different mathematical fields. For example, within cobordism-theory it is known as the ``first'' manifold (with respect to dimension), that is not the boundary of any manifold, when the type of ``manifold'' is specified to mean ``topological manifold'', which in dimension $2$ is equivalent to both ``homology manifold'', and to ``{\it quasi-manifold}''. A simplicial complex is a quasi-manifold if and only if its links are all {\it pseudo-manifolds}, see p.~\pageref{DefP27:2} for definitions.

Our intuition easily gets hold of a useful picture of ``the cone of
an up-right standing cylinder'' and also of a correct picture of its ``boundary'' if the cone vertex, in a realization, is placed in the point of inertia of the cylinder.\label{BundaryCylCone} Here ``boundary'' should be interpreted as those
points where the local singular homology is trivial in the highest
dimension, i.e. in dimension $3$.

Now, make a straight cut in the cylinder, from the top to the bottom and glue these cut-ends back together again, after twisting one of
them $180^\circ$. The result is a Möbius band. With the cone-point in the point of inertia of the Möbius band we create ${\mathcal M} \ast \bullet$ ``the cone of a Möbius band''. We engage the same set of points as in the cylinder case, but differently arranged, and there are reasons to believe that no sharp intuitive picture is revealing itself. ${\mathcal M} \ast \bullet$ is $3$-dimensional but it can not be inscribed in our ordinary room.

Regarding the factors as quasi-manifolds and using, say, the 3-element prime field
${{\mathbb{Z}_{_{{\bf 3}}}}}\!$
as the coefficient module for the homology groups,
we get:\\
($\mathbb{B}^n$ is the $n$-dimensional disk or
``$n$-ball'' and $\mathbb{S}^n$ is the ``$n$-sphere''.)

\begin{enumerate}
\item ${\rm Bd}_{_{{\mathbb{Z}_{\bf 3}}}}\!\!({\mathcal M} \ast \bullet)=$
$(({\rm Bd}\!_{_{{\mathbb{Z}_{\bf 3}}}}\!\!
{\mathcal M})\ast\bullet) \cup ({\mathcal M} \ast
{\rm Bd}_{_{{\mathbb{Z}_{\bf 3}}}}\!\!
\bullet)
=
(({\mathbb{S}}^1)\ast\bullet)\cup ({\mathcal M}
\ast
{\rm Bd}_{_{{\mathbb{Z}_{\bf 3}}}}\!\! \bullet)
=
\mathbb{B}^2 \cup
({\mathcal M}\ast \{\emptyset_o\})
=\mathbb{B}^2 \cup {\mathcal M}
= {{\mathbb{RP}}}^2\!\!:=$
\htmladdnormallink{{\it real projective plane}}
{http://en.wikipedia.org/wiki/Real_projective_plane}.

\medskip\noindent
\item ${\rm Bd}_{_{{\mathbb{Z}_{\bf 3}}}}\!\!(
{{\mathbb{RP}}}^2 \ast \bullet)=$
$ (({\rm Bd}_{_{{\mathbb{Z}_{\bf 3}}}}\!\!
{\mathbb{RP}}^2)
\ast\bullet) \cup ({\mathbb{RP}}^2
\ast {\rm Bd}_{_{{\mathbb{Z}_{\bf 3}}}}\!\! \bullet)
=
(\{\emptyset_o\}\ast\bullet) \cup ({\mathbb{RP}}^2
\ast \{\emptyset_o\})
= \bullet \cup {\mathbb{RP}}^2.
$

\medskip\noindent
\item ${\rm Bd}_{_{{\mathbb{Z}_{\bf3}}}}\!\!({\mathbb{S}}^{2}
\ast
\bullet)=$
$ (({\rm Bd}_{_{{\mathbb{Z}_{\bf 3}}}}\!\!
{\mathbb{S}}^2)\ast\bullet)
\cup
({\mathbb{S}}^2
\ast
{\rm Bd}_{_{{\mathbb{Z}_{\bf 3}}}}\!\!
\bullet)
=
(\emptyset\ast\bullet)
\cup
({\mathbb{S}}_{\!}^{^{_{2}}}
\ast
\{\emptyset_o\})
={\mathbb{S}}^2.
$
\end{enumerate}

\medskip
$\!{\mathbb{RP}}^2\!\ast
\bullet$ is
a homology$_{_{{\mathbb{Z}_{\bf p}}}}$
$3$-manifold  \underbar{with} boundary
$\bullet\cup{\mathbb{RP}}^2$
if  ${\bf p}\neq2$, by \cite{13} p.~36 and Th.~\ref{TheoremP31:12}
p.~\pageref{TheoremP31:12}, while ${\mathcal M} \ast \bullet$ fails
to be one due to the cone-point.
However, ${\mathcal M} \ast \bullet$ is a quasi-$3$-manifold, see
p.~\pageref{DefP27:2}, which by the boundary definition p.~\pageref{DefP2:BoundaryOfSimplManif}
has the real projective plane ${{\mathbb{RP}}}^{^{_{2}}}\!$ as its well-defined boundary, as seen above.
The topological realisation, cf.
p.~\pageref{DefP8:RealizationOfSimpCompl},
of
${\mathbb{RP}}^2\!\ast\bullet$
is not a
\htmladdnormallink{{\it topological manifold},}
{http://en.wikipedia.org/wiki/Topological_manifold}
due to the properties at the cone-point.
This is also true classically.
Topological $n$-manifolds can only have $n-1$-, $-1$- or
$-\infty$-dimensional boundaries as homology manifolds, due to their local orientability.
The \underbar{global} non-orientability of ${\mathcal M}$ and
${\mathbb{RP}}^2$ becomes \underbar{local} non-orientability at the cone-point in
${{\mathcal M}}\ast\bullet$,
resp.\nolinebreak\ ${\mathbb{RP}}^2\!\ast\bullet$.
It is obviously rewarding to study combinatorial manifolds such as pseudo- and quasi-manifolds, also when the main object of study is topological manifolds.
Our section~\ref{section:SecP27:III} is therefore devoted to these  generalized notions of
\htmladdnormallink{{\it simplicial manifolds}}
{http://en.wikipedia.org/wiki/Simplicial_manifold}.

\bigskip \
\goodbreak
\subsection{\spaceskip2.2pt Why classical search for
``The \underline{Relative} Künneth Formula for Joins'' failed }\label{UnderavdP5:KynnethFormForJoins}\hfill\break
\indent
The short answer is that the objects in the classical categories
simply do not contain any $(-1)$-dimensional object that at the same time is a sub-object contained in every non-empty object and a unit element with respect to join.
With respect to the category of simplicial complexes, Definition~1  p.~\pageref{DefP1:1} provide such a $(-1)$-dimensional {\it simplicial complex} through $\{\emptyset_o\}$, which is generated by the new empty {\it simplex} $\emptyset_o$. Note that this new complex $\{\emptyset_o\}$ also is the simplicial join-unit.
Maybe this answer becomes more enlightening after a thorough check of the proof of the Eilenberg-Zilber theorem for simplicial join in
p.~\pageref{SimplEilenbergZilber}.
Our externally adjoined element $\wp$ generates the join-unit
$\{\wp\}$ that allows an Eilenberg-Zilber theorem for {\it topological} join, resembling the classical one for products.
The rest is not really category specific and is covered by
\htmladdnormallink{{\it homological algebra,}
}{http://en.wikipedia.org/wiki/Homological_algebra}
as founded in the {\it Cartan-Eilenberg} masterpiece \cite{CartanandEilenberg},
for which the initial impetus was precisely the Künneth Formula.

That the {\it join}-{\it operation} is at the heart of algebraic
topology, becomes apparent, for instance, by Milnor's
construction of the universal principal fiber bundle in \cite{23}
where he in Lemma~2.1 also reaches the limit for classical algebraic topology as he formulates a reduced version of the non-relative
Künneth formula for joins of general topological spaces as (no
\underbar{classical} \underbar{relative} version can exist due to the missing
$(-1)$-dimensional join-unit):
$$\
\widetilde{H}\!\!_{_{\!q+1\!}}
        (\!X\!_{\!_1}\!\ast Y\!_{_1})
\ \! {{{_{\mathbb{Z}}}} \atop{\raise3.5pt\hbox{$\cong$}}}
{\rlap{$_{_{_{_{i+j=q}}}}$}} {\ \raise0.0pt\hbox{$\oplus$}}\
\hbox{\bf(}\widetilde{H}_{i}
           (X\!_{_1})
\otimes\!_{_{\mathbb{Z}}}\widetilde{H}_{\!j}
(Y\!\!_{_1})\hbox{\bf)}\ \
{\rlap{$\!\!_{_{_{\!i+j=q-\!1}}}$}{\ \raise0.2pt\hbox{$\oplus$}}}\ \
 \hbox{\rm Tor}_1^{^{_{\mathbb{Z}}}\!}\!\hbox{\bf(}\widetilde{H}_{i}(X\!_{_1}),\widetilde{H}_{\!j}
      (Y\!\!_{_1})\hbox{\bf)},\!
$$ \label{Milnor'sReducedFormula}
i.e. the $``X_2=Y_2=\emptyset\hbox{''}$-case in our Th.~\ref{TheoremP14:4}
p.~\pageref{TheoremP14:4}.
Milnor's results, apparently, inspired\ \  G.W. Whitehead to
introduce the {\it augmental total chain complex}
$\widetilde{\bf S}(\circ)$ and {\it augmental homology},
$\widetilde{H}_{^{_{\star\!}}}(\circ),$ in \cite{36}.
$\widetilde{\bf S}(\circ)$ comprise an
in-built dimension shift, which in our version becomes an
application of the suspension operator on the underlying singular chain.
G.W. Whitehead gives the empty space, $\emptyset,$ the status of $(-1)$-dimensional standard simplex, but in his pair space theory, $\widetilde{H}_{^{_{\star\!}}}(X_1,X_2:{\bf G})$, he never took into account that the join-unit $\emptyset$ then would have the identity map, Id$_{_{^{\emptyset^{\!}}}}$, as a generator for its $(-1)$-dimensional singular augmental chain group, which, correctly interpreted, actually makes his pair space theory equivalent to the ordinary relative homology functor.

G.W. Whitehead states that $\widetilde{\bf S} (X\!\ast Y) $ and $\widetilde{\bf S}(X)\otimes\widetilde{\bf S}(Y) $ \hbox{\rm are chain equivalent,} {$\lower3.6pt\hbox{$^{\approx}$}\ \!$},
which, modulo the introduction of a $-1$-dimensional topological object, is confirmed by our Theorem\ \ref{TheoremP13:3}
p.~\pageref{TheoremP13:3},
and then in a footnote Whitehead points out, referring to
\cite{23}
p.~431 Lemma 2.1, that;
\begin{quote}
This fact does not seem to be stated explicitly in the literature but is not difficult to deduce from Milnor's proof of the ``Künneth theorem'' for the \hbox{homology groups of the join.}
\end{quote}

On the chain level it is indeed ``{\it not\hskip0.2cm
difficult}''\label{EasyIso} to see what is needed to achieve
$\widetilde{\bf S}(X\ast Y) \approx\ \widetilde{\bf
S}(X) \otimes\widetilde{\bf S}(Y)$, since the right hand
side is well-known, but then to actually do it for the classical category of topological pair-spaces and within the frames of the
\htmladdnormallink{``{\it Eilenberg-}}
{http://en.wikipedia.org/wiki/Eilenberg-Steenrod_axioms}
\htmladdnormallink{{\it Steenrod axioms}''}
{http://en.wikipedia.org/wiki/Eilenberg-Steenrod_axioms}
for relative homology
is, unfortunately, impossible, since
the need for a $(-1)$-dimensional standard simplex is
indisputable and the initial object $\emptyset$ just
won't do, since it already plays another incompatible
role within this axiomatization, cf. \cite{9}
p. 3-4. Indeed, the ``convention''
$\!\hbox{\rm H}_{{{\hbox{\seveni i}}}}(\cdot)=
\hbox{\rm H}_{{{\hbox{\seveni i}}}}(\cdot,\emptyset;{\mathbb{Z}})$,
cf. \cite{9}
pp.~3 + 273,
is more than a mere ``convention'' in that it connects the single space and the pair space theories and thereby determine the chain-groups of $\emptyset$ to be trivial in all degrees.
This was also observed in \cite{22}
p.~108, leading up to a refutation of Whitehead's approach, sketched above.

Our  Theorem~\ref{TheoremP13:3} p.~\pageref{TheoremP13:3},
answers G.W. Whitehead's half a century old quest above for a proof.
For arbitrary ``Definition~1''-simplicial complexes,
G.W. Whitehead's ``chain equi\-valent''-claim above (=the non-relative Eilenberg-Zilber theorem for joins), is easily
proven as shown here in p.~\pageref{SimplEilenbergZilber},
where we quote the proof of this fact, taken from
\cite{Fors2},
which was the first \TeX-version of the present article.

Our K\"unneth Theorem for topological joins, i.e.
Theorem~\ref{TheoremP14:4} p.~\pageref{TheoremP14:4},
holds for arbitrary topological spaces.

The generalization, i.e.
Theorem~\ref{TheoremP17:7}
p.~\pageref{TheoremP17:7},
of the
\htmladdnormallink{``{\it derivation}''}
{http://en.wikipedia.org/wiki/Derivation_(abstract_algebra)}
-like boundary formula p.~\pageref{BundaryIntro}
above, also holds for non-triangulable
homology manifolds,
which we have restricted to be locally compact Hausdorff, but
{\it weak Hausdorff \ $k$-spaces}, cf. \cite{11}
p.~$\!$243ff, would also do, since we only need the
\htmladdnormallink{``${\bf T}_1$-{\it separation axiom}''
}{http://en.wikipedia.org/wiki/Separation_axiom}
($\Leftrightarrow$ ``points are closed'').

\bigskip
Since classical Combinatorics and classical General Topology cannot separate the ``join zero'' from the ``join-unit''
or the boundary of the $0$-ball from that of the $0$-sphere,
any theory constructed thereon --- classical Algebraic Topology in particular --- will need ad-hoc definitions/reasoning.
The use of a relative and a separate reduced homology functor
instead of one single homology functor, is only one example.
These are probably the reasons for the ongoing marginalization of \underbar{\it classical} algebraic topology, as witnessed by M.\ Hovey
at \htmladdnormallink{http://math.wesleyan.edu/
\lower4pt\hbox{$\tilde{\ }$}mhovey/problems/.}
{http://math.wesleyan.edu/~mhovey/problems/}

\medskip
We quote Paul Dirac
as he comments on the increasing abstraction within the realm of mathematics and physics.
(P.A.M. Dirac, Proc.~Roy. Soc. A 133, 60 (1931).
\htmladdnormallink{{\it Quantised}}
{http://rspa.royalsocietypublishing.org/content/133/821/60.citation}
\htmladdnormallink{{\it Singularities in the Electromagnetic Field})}
{http://rspa.royalsocietypublishing.org/content/133/821/60.citation}

\begin{quote}
It seems likely that this process of increasing abstraction will continue in the future and that advance in physics is to be associated with a continual modification and generalization of the axioms at the base of the mathematics rather than with a logical development of any one mathematical scheme on a fixed foundation.
\end{quote}

\part{Augmental Homology Theory}
\section{Simplicial and singular augmental homology theory}\label{Ch3:AugmHomology}
\subsection{Definition of underlying augmental categories}
\S 3.1 presents a firm formal foundation for an overall augmented environment, suitable for (augmented)
\htmladdnormallink{{\it homology theories}}
{http://en.wikipedia.org/wiki/Homology_theory}.

The combination of the complex-part of Definition~2 and the join-part of Definition~1 p.~\pageref{DefP1:1} is not an option, since then the join-operation results in a simplicial complex (in any sense) only if one of the join-factors is the empty complex.

The typical morphisms in the classical category {\bf ${\mathcal K}$}
of simplicial complexes with vertices in {\bf W} are the simplicial
maps as defined in \cite{30} p.~109, implying in particular
that;\label{DefMorphism}
$$
\left\{\begin{array}{lll}
{\bf Mor}_{\mathcal K}(\emptyset,_{\!}{\Sigma})\!=\!
\{\emptyset_o\}\!=\!\{0_{\emptyset^{\!},{\Sigma}}\},\\
{\bf Mor}_{\mathcal K}({\Sigma},\emptyset)=\emptyset,\
\hbox{\rm if}\ \Sigma\neq\emptyset,\\
{\bf Mor}_{\mathcal K}(\emptyset,{\emptyset})=
\{\emptyset_o\}=\{0_{\emptyset,{\emptyset}}\}=\{\hbox{\rm
id}_{\emptyset}\},
\end{array}\right.
$$
where $0_{_{\!{\Sigma,\Sigma^{\prime}}}}$ denotes $\emptyset =$
{\it the\ empty function from $\Sigma$ to $\Sigma^{\prime}$}.

\medskip
So; $ 0_{_{\!{\Sigma,\Sigma^{\prime}}}}\!\in{\bf
Mor}_{_{\hbox{\fivesy K}}} (\Sigma,\Sigma^\prime)
\Longleftrightarrow \Sigma=\emptyset. $

\medskip
If in a category $\varphi_i\in {\bf Mor} ({ R_i,S_i}),\
i=1,2,$ we put:

\medskip
\noindent
$$\varphi_1 \sqcup \varphi_2: R_1\sqcup R_2 \rightarrow
S_1\sqcup S_2 :r\mapsto
\left\{\begin{array}{ll}
{\varphi_1(r)} \ \ \hbox{if {$r\in R_1$}} \\
{\varphi_2(r)}
\ \ \hbox{if {$r\in R_2$}} \\
\end{array}\right.\
\text{where $ \sqcup :=$ ``disjoint union''.}
$$

\begin{vardef}{} \label{DefP7:K0}
{\rm (of the objects in  ${\mathcal K}\!_{_{^{o}}}\!$)}
An $($abstract$)$ simplicial complex $\Sigma$ on a vertex set ${{\bf
V}\!_{\Sigma_{}}}$ is a collection $($empty or non-empty$)$ of
finite $($empty or non-empty$)$ subsets $\sigma$  of ${{\bf
V}\!_{\Sigma_{}}}$ satisfying;
\begin{enumerate}
\item[({\bf a})] If $v\!\in\! {{\bf V}\!_{\Sigma_{}}}$,
then $\{v\}\in \Sigma$.
\item[({\bf b})] If $\sigma\in \Sigma$ and $\tau\subset \sigma$
then $\tau\in \Sigma$.
\end{enumerate}
\end{vardef}

So, $\!\{\emptyset_o\}$ is allowed as an object in {\bf ${\mathcal
K}\!_{_{^{o}}}$}.
We will write ``concept$_o$'' or \ ``concept$_{\wp}$'' when we want to stress that a concept is to be related to our modified categories.

If $
\hbox{\rm card}(_{\!}\sigma{_{^{\!\!_{\ }}}}\!_{\!})\!=\!q\
\!\!${\eightrm+}$1$ then $\dim\sigma{_{^{\!\!_{\ }}}}\!\!:=\!q$
\label{DefP7:SimpDim}and $\sigma{_{^{\!\!_{\ }}}}$ is said to be a
$q$-{\it face}$_o$ or a
$q$-{\it simplex}$_o$ of $\Sigma{_{^{\!_o}}}$ and
$\dim\Sigma{_{^{\!_o}}}\!\!:=$sup$\{\dim (\sigma{_{^{\!\!_{\
}}}}\!)\vert \sigma{_{^{\!\!_{\ }}}}\!\!\in\!
\Sigma{_{^{\!_o}}}\}$.
Writing\ $\emptyset{_{^{\!_o}}}\!$ when using $\emptyset$ as a
simplex, we get dim$(\emptyset)\!=\!-\infty$ and
$\dim(\{\emptyset{_{^{\!_o}}}\!\})\!=\!\dim(\emptyset{_{^{\!_o}}}\!)\!=\!-1.$

\begin{varnote}{}
Any simplicial complex$_{_{^{\!o}}}\!$ $\Sigma_o\!\neq\emptyset$ in
${\mathcal K}\!_{_{^{o}}}\!$ includes $\{\emptyset\!_{_{^{o}}}\!\}$
as a subcomplex.
\end{varnote}

\goodbreak

So, a typical object$_{_{^{\!o}}}\!$ in {\bf ${\mathcal K}\!_{_{^{o}}}$}
is {\bf ${\Sigma} \sqcup \{\emptyset\!_{_{^{o}}}\!\}$} or
$\emptyset$, where ${\Sigma}\in{\mathcal K}$ and {$\chi$} is a
morphism in {\bf ${\mathcal K}\!_{_{^{o}}}$} if either:

\smallskip
$ \hbox{\rm a})\ \ \chi=\varphi \sqcup \hbox{\rm
id}_{\{\emptyset_o\}}\ \hbox{\rm for\ some}\ \varphi\in {\bf
Mor}_{_{\hbox{\fivesy K}}}(\Sigma,\Sigma^\prime)\ \ \ \hbox{\rm or}
$

\smallskip $ \hbox{\rm b})\ \ \chi= 0_{\emptyset,\Sigma_o}.$

In\ particular, ${\bf Mor}_{_{\hbox{\fivesy
K}\!_{_{^{o}}}}}(\Sigma_o,\{\emptyset_o\})=\emptyset $ if and only if
$ \Sigma_o\neq\{\emptyset_o\},\emptyset. $

\bigskip
\underbar{A functor {${\mathcal{E}}$}:\ {\bf ${\mathcal{K}}_o$}$\longrightarrow${\bf ${\mathcal{K}}$}:}

\smallskip
\noindent Set ${\mathcal{E}}(\Sigma_o)=
 \Sigma_o \setminus \{\emptyset_o\}\in Obj({\mathcal K})$
and given a morphism  $\chi: {\Sigma_o}
\rightarrow {\Sigma_o^\prime}$ we define;

\medskip
${{\mathcal{E}}}(\chi)=\left\{\begin{array}{ll}\medskip
\varphi  & \text{if $\chi$ fulfills {\rm a)} above \ and}
\\
0_{\emptyset,{{\mathcal{E}}}(\Sigma_o)} & \text{if $\chi$ fulfills
{\rm b)} above.}
\end{array}\right.
$

\bigskip
\underbar{A functor {${\mathcal{E}}_o$}: ${\mathcal
K}\longrightarrow{\mathcal K}_o:$}

\smallskip
\noindent Set\ ${{\mathcal{E}}}_o(\Sigma)=\Sigma \sqcup
\{\emptyset_o\}\in{Obj}({\mathcal K}_o)$
and given
$\varphi:{\Sigma} \rightarrow {\Sigma^\prime}$, we put\
$\chi:=\varphi \sqcup \hbox{\rm id}_{\{\emptyset_o\}}$, which gives;
${\mathcal{E}}{\mathcal{E}}_o=$id$_{\hbox{\fivesy K}}$;
{\bf imE$_o\!\!= Obj({\mathcal K}_o)\setminus \{\emptyset_o\}$} and
${{\mathcal{E}}}_o{{\mathcal{E}}}=\hbox{\rm id}_{\hbox{\fivesy K}}$ except for
                ${{\mathcal{E}}}_o{{\mathcal{E}}}(\emptyset)=\{\emptyset \}.$

\medskip
Similarly, let {\bf ${\mathcal C}$} be the category of topological spaces and continuous maps. Consider the category {\bf ${\mathcal
D}\!_{_{^{\!\wp}}}\! $}\label{DefP7:CategoryDp} with objects:
$\emptyset$ together with $ X_{\wp}:=X{\bf +}\{\wp\}$ for all $X\in {Obj}({\mathcal C})$, i.e., the set $X_{\wp}:= X\sqcup \{\wp\}$
with the weak topology, $\tau\!\!_{_{X}\!_{_{^{\!\wp}}}\!\!}$, with
respect to $_{^{\!}}{ X}_{^{\!}}$ and $\{\wp\}$.

\medskip
\noindent
$f\!_{{{\wp}}}\!\! \in {\bf Mor}\!_{_{^{\hbox{\fivesy
D}\!\!_{_{^{\wp}}}\!\! }}}({ X}\!_{_{^{\wp}}}\! ,{
Y}\!\!\!_{_{^{\wp}}}\! )\ \text{if;}$
$$f\!_{{{\wp}}} =
\left\{\begin{array}{ll}
\hbox{\hbox{\rm a)}}\ \ f{\bf +} \hbox{\rm id}\!_{_{\{\!\wp\!\}}}\
(:=f\sqcup\hbox{\rm id}\!_{_{\{\!\wp\!\}}})\ \hbox{\rm with}\ f\in
{\bf Mor}\!_{_{\hbox{\fivesy C}}}\!(X,Y)\
\hbox{\rm and}\\
f\ \hbox{\rm is}\ \underline{{\hbox{\rm on}\ \!X}}\ \hbox{\rm to}\ {
Y},\ \hbox{\rm i.e}\ \hbox{\rm the\
domain\ of}\ f\hbox{\rm \ is\ the}\ {\rm whole\ of}\ \!X\!\\
\hbox{\rm and}\
{ X}\!_{_{^{\wp}}}\!\! ={ X}{\bf +} \{\wp\},\ {
Y}\!\!\!_{_{^{\wp}}}\!\! ={ Y} {\bf +} \{\wp\}\ \hbox{\rm \ {\bf
or}}
\\
\hbox{\hbox{\rm b)}}\ \ 0_{{{{\emptyset,Y\!\!_{_{^{\!\wp}}}\!\! }}}}
\ \ (=\emptyset = \hbox{\rm the\ empty\ function\ from}\
\emptyset\ \hbox{\rm to}\ { Y}\!\!_{_{^{\!\wp}}}\! ).\\
\end{array}\right.
$$

\noindent
There are functors:
\vskip-0.5cm

$${\mathcal
F}\!\!_{_{^{\!\wp}}}\!:\!
\left\{\!\!\!\begin{array}{ll}\medskip
 {\mathcal C} \longrightarrow
{\mathcal D}\!_{_{^{\!\wp}}}\hfill \\
{ X}\mapsto\!{ X}{\bf +} \{\wp\}
\end{array}\!\!,\right.\ \text{and}\
\ {\mathcal F}\!:\!
\left\{\!\!\!\begin{array}{ll}\medskip
{\mathcal D}\!_{_{^{\!\wp}}}\longrightarrow  {\mathcal C}\hfill \\
\!{X}{\bf +} \{\wp\} \mapsto{X}\medskip\\
\emptyset\mapsto\emptyset \hfill
\end{array}\right.
\text{resembling}\ {{\mathcal{E}}}_{o},\ \text{resp.}\  {{\mathcal{E}}}.
$$

\begin{varnote}
The ``${\mathcal F}\!\!_{_{^{\wp}}}$-lift topologies'',
$$
\tau_{_{\!\!{X_{_{^{\!{o}}}}}}}\!\!\!:=\! {\mathcal
F}\!\!_{_{^{\wp}}}\!(\tau_{_{\!\!{X}}})_{\!} \cup_{\!}
\{\emptyset_{}\} \!=\! \{ {\mathcal O}\!_o\!\!=\! {\mathcal
O}_{\!}\sqcup_{\!}\{\wp\}\ \!\mid\ \! {\mathcal
O}\!\in\!\tau_{_{\!\!{X}}}
\}_{\!}\cup_{\!} \{\emptyset_{}\}
$$
\vskip-0.3cm{\rm and}\vskip-0.3cm
$$
\tau_{_{\!\!{X_{_{^{\!{o}}}}}}}\!\!\!:=\! \tau_{_{\!\!{X}}}\cup
\{\!{X}_{_{^{\!{o\!}}}}\!\}\!= \!\{{\mathcal O}_{_{^{\!{o\!}}}} \mid
{\mathcal O}_{_{^{\!{o\!}}}}\!=
\!\!{X}_{_{^{\!{o\!}}}}\setminus({\mathcal N}\sqcup\{\wp\})\ \!;\
{\mathcal N}\ \hbox{\rm closed\ in}\ {X}
\}\cup\nobreak\{\!{X_{_{^{\!{o\!}}}}}\!\}\ \
$$
would also \hbox{\rm give} $_{\!}{\mathcal D}\!\!_{_{^{\wp}}}\!$
due to the domain restriction in ${\rm a})$, disallowing ``partial'' functions.
This makes $_{\!}{\mathcal D}\!\!_{_{^{\wp}}}\!$
a link between the two constructions of partial maps treated in \cite{2}
pp.~184-6.

No {extra} morphisms$ _{_{^{\!\wp}}}\! $\nobreak\ has
been allowed into ${\mathcal D}\!_{_{^{\!\wp}}}$ $({\mathcal K}_{\!o
^{\!}} )$ in the sense that the morphisms$_{_{^{\!\wp}}}\!$ are all
targets under ${\mathcal F}\!\!_{_{^{\!\wp}}}$
$
\!\!\
({{\mathcal{E}}}_{_{^{\!o}}}\! )$
except $0_{{{{\emptyset,Y\!\!_{_{^{\!\wp}}}\!\!   }}}} $ defined
in b), re-establishing $\emptyset$ as the unique initial
object.
\end{varnote}

The underlying principle for our definitions is that a concept in
${\mathcal C}\ ({\mathcal K})$ is carried over to ${\mathcal
D}_{\!\wp}\ ({\mathcal K}_o)$ by ${\mathcal F}\!\!_{\wp}\ ({{\mathcal{E}}}_o)$ with addition of definitions of the concept$_{\!o}$ for cases
which are not proper images under $ {\mathcal F}\!\!_{\wp}\ ({{\mathcal{E}}}_o)$. The definitions of the product/join operations
``$\times\!_o\!$'', ``$\ast\!\!_{{^o}}\!$'', ``$\mhatast$'' in
page~\pageref{DefP11}
as well as ``$/\!{^{^{_{o}}}}$'', ``$+\!\!_{_{^{o}}}$'' and ``$\setminus\!_{o}{_{\!}}$'' below certainly comply with this principle.

\begin{vardef} \label{DefP8:SetDiv0}
$X_{\wp1}{/^{o}}{X_{\wp2}}:=
\left\{\begin{array}{ll}
\!\hskip-0.1cm\emptyset &\text{if}\ X_{\wp1}\!= \emptyset\\
\!\hskip-0.1cm{\mathcal F}\!_\wp\bigl( {{\mathcal F}(X_{\wp1})}/
 {{{\mathcal F}}(X_{\wp2})}\bigr)
 &\text{if}\ X\!\!_{_{\wp2}}\!\!\neq \!\emptyset \ {\rm in}\
{{\mathcal D}}\!_{\!\wp} \ {\rm for}\
X\!_{\wp2}\!\subset\!X\!_{\!\wp1}.
\end{array}\right.
$

\medskip\noindent
``$/$'' is the classical ``quotient'' except that
${{{\mathcal F}}( X_{\wp1})}/{\emptyset}
:={{\mathcal F}( X_{\wp1})},$
cf. \cite{2}
p.~102.
\end{vardef}

\noindent
So, for example:
$(X_{\wp1}{/^{o}}X_{\wp2},X_{\wp2}{/^{o}}{X_{\wp2}}):=
\left\{\begin{array}{lll}
\!\hskip-0.1cm(X_{\wp1},\emptyset) & \text{if}\ X_{\wp2}=\emptyset, \\
\!\hskip-0.1cm(X_{\wp1},\{\wp\}) & \text{if}\ X_{\wp2}=\{\wp\},\\
\!\hskip-0.1cm(X_{\wp1}{/^{o}}{X_{\wp2}},\bullet) & \text{if}\ X_{\wp2}\not=\emptyset,\{\wp\}.\\
\end{array}\right.
$

\begin{vardef} \label{DefP8:SetSum0}
Let ``$+$'' denote the classical ``topological sum'', then;

\medskip
$X_{\!\wp1}+\!\!_{_{^{o}}}X_{\!\wp2}:=
\left\{\begin{array}{lll}
\!X_{\!\wp1} & \text{if}\
X_{\!\wp2}=\emptyset\\
\!X_{\!\wp2} & \text{if}\
X_{\!\wp1}=\emptyset\\
\!{{\mathcal F}\!_{\wp}}({{{\mathcal F}}(X_{\!\wp1})+ {{\mathcal
F}}(X_{\!\wp2})}) & \text{if}\
X_{\!\wp1}\neq \emptyset \neq X_{\!\wp2}
\\
\end{array}\right.
$
\end{vardef}

\smallskip
To avoid dropping out of the category ${{\mathcal D}\!_{\wp}}$ we introduce a setminus ``$\setminus{_{_{^{\!o\!}}}}$'' in ${{\mathcal D}\!_{\wp}}$, giving Proposition~\ref{PropP15:1} p.~\pageref{PropP15:1} its compact form.
Obs, that ``$\setminus{_{_{^{\!o\!}}}}$'' might be non-associative if $\{{\wp}\}$ is involved.
(``$\smallsetminus$'':= classical ``setminus''.)

\begin{vardef} \label{DefP8:PointSetMinus0}\noindent
$\!{X}_{\!\wp}\!\setminus\!{_{_{^{\!o\!}}}}{X}_{\!\wp}^{^\prime}\!:=
\left\{\begin{array}{ll}
\hskip-0.1cm\emptyset & \text{if}\ { X}_{\!\wp}\!\!=\emptyset,~{
X}_{\!\wp}\mmysubsetneqq { X}_{\!\wp}^\prime~
\text{or}~{X}_{\!\wp}^\prime=\{{\wp}\}\\
\hskip-0.1cm{\mathcal F}\!_{\wp}\bigl({{\mathcal F}}({X_{\!\wp}})
\smallsetminus
{{\mathcal F}}({X}_{\!\wp}^\prime)\bigr) & \text{else}.
\end{array}\right.
$
\end{vardef}

\medskip
\subsection{Realization of a simplicial
complex$_o$}\label{UnderavdP8:RealizationOfSimpCompl}
\label{DefP8:RealizationOfSimpCompl}

A ``point'' in a classical topological space is simply a subspace containing a single element. In a simplicial complex$_o$ a ``point'' is a subcomplex$_o$ containing only one vertex. The simplicial complexes$_o$ $\emptyset$, $\{\emptyset_o\}$ are both pointless and non-final. Any functorial definition of the {\it realization of a simplicial  complex}$_o$ forces their realizations $|\emptyset|$ and $|\{\emptyset_o\}|$ resp. to be separate pointless objects, which implies that no realization functor could target any \underbar{classical} Euclidian space of any dimension, since each of these spaces only possess one pointless subspace - namely the empty space $\emptyset$.

\goodbreak
So, to be able to create something similar to the classical realization functor one need to (externally) add a non-final object $\{\wp\}=|\{\emptyset_{_{^{\!\!o}}}\!\}|$ into the classical category of topological spaces as {\it join-unit} and (-1)-dimensional standard simplex.
If $X\!_{_{^{\wp}}}\!\! \neq\emptyset,\ \{\wp\}$, then
($X\!_{_{^{\wp}}}\!,\tau\!\!_{_{^{X_{{^{\!\wp}}}\!\!}}}$) is a
non-connected  space.
We therefore \underbar{define} ${ X}\!_{_{^{\wp}}}\!\!\neq
\emptyset,\ \{\wp\}$ to have a certain topological
``property$\!_{_{^{\wp}}}\!\!$'' if $ {\mathcal F}({
X}\!_{_{^{\wp}}})$ has the ``property'', e.g., ${ X}
\!_{_{^{\wp}}}\! $ is connected$\!_{_{^{\wp}}}$ \underbar{\rm iff}   $X$ is connected.

We will use Spanier's definition of the
\htmladdnormallink{{\it function space realization}
}{http://books.google.com/books?id=h-wc3TnZMCcC&printsec=frontcover&dq=\%22covariant+functor\%22&hl=en&source=gbs\_summary\_r&cad=0\#PPA110,M1}
${|\Sigma_{_{^{o}}}|}$\label{FunctionSpaceRealization} as given
in \cite{30}
p. 110, unaltered, except for the `` $_{_{^{o}}}$''s and the
underlined $\underline{\underbar{addition}}$ below, where $\wp:=\alpha_0$ and $\alpha_0(v)\equiv 0\ \text{for all}\ v\in {\bf
V}\!\!_{_{^{\Sigma}}}$\label{DefP8:alpha0}:
\begin{quote}
We now define a
covariant functor from the category of simplicial
complexes$_{_{^{\!o}}}$ and simplicial maps$_{_{^{\!o}}}$ to
the
category of topological spaces$_{_{^{\!o}}}$ and continuous
maps$_{_{^{\!o}}}$. Given a nonempty simplicial
complex$_{_{^{_{\!}o}}}$ ${\Sigma_{_{^{\!o}}} }$, let ${
|\Sigma_{_{^{\!o}}}|}$ be the set of all functions $\alpha$ from the
set of vertices of ${\Sigma_{_{^{\!o}}}}$ to ${\bf I}:=[0,1]\
\hbox{\rm such\ that;}$
\begin{enumerate}
\item[(a)] For any $\alpha$, $\{v\in {\bf
V}\!\!_{_{^{\Sigma_o }}} \mid
         \alpha(v)\neq 0\}$
          is a simplex$_o$ of ${\Sigma_o }$.\\ (in particular,
          $\alpha(v)\neq 0$ for only a
          finite set of vertices).
          \medskip
\item[(b)]
For any $\alpha\ \underline{\underline{\neq\alpha_0}}$,
$\sum\!\!_{_{v\in V\!\!{_{_{\Sigma_o}}}}}\hskip-0.2cm \alpha(v)=1.$\ \
\end{enumerate}
\nobreak \noindent If  ${\Sigma_o }=\emptyset$ , we define
         ${ |\Sigma_o |}=\emptyset$.\phantom{\hskip0.5cm\ space}
\end{quote}

The {\it barycentric} {\it coordinates}$_o$ define a metric
$$
d(\alpha,\beta)=\sqrt{\sum_{v\in{V}\!_{\Sigma_o}}
[\alpha({v})-\beta({v})]^{2}}
$$
on $|\Sigma_o|$ providing it with the {\it metric topology}\ \ and thereafter denoting it $|\Sigma_o|_d$.

We will equip $|\Sigma_o|$ with another, more commonly used topology and for this purpose we define the {\it closed geometric simplex$_o$}
$|\sigma_o|$\label{DefP8:Closedsigma0} of $\sigma_o\in{\Sigma_o }$
to be:
$$|\sigma_o|:=\{\alpha\in|\Sigma_o|\ \mid\ { [\alpha(\hbox{\rm v})\neq0]
\Longrightarrow [v\in \sigma_o] \}.}$$

\begin{vardef} \label{DefP8:4}
For ${\Sigma_{_{^{\!o}}} }\ne\emptyset$,
$|\Sigma_{_{^{\!o}}}|$ is topologized by setting
${|\Sigma_{_{^{\!o\!}}}|}\!:=|{{\mathcal{E}}}(\Sigma_{_{^{\!o\!}}})|+\{\alpha_0\}$. This is equivalent to giving
${|\Sigma_{_{^{\!o}}}|}$ the weak topology with respect to the
$|\sigma\!_{_{^{\!o\!}}}|$'s naturally imbedded in $\mathbb{R}^{n}{\bf +}\{\wp\}$, and we define $\Sigma_{_{^{\!o}}}$
to be connected if $|\Sigma_{_{^{\!o\!}}}|$ is, i.e., if ${\mathcal
F}(\vert\Sigma_{_{^{\!o\!}}}\vert)\simeq\vert{{\mathcal{E}}}
(\Sigma_{_{^{\!o\!}}})\vert\ \hbox{\rm is}.$
\end{vardef}

\begin{varprop} \label{PropP8}
$$\!\hskip-0.3cm|{_{\!}}\Sigma_{{{\!o\!}}}|\ \hbox{\it is always homotopy
equivalent to}\
 |{_{\!}}\Sigma_{{{\!o\!}}}|_{_{^{\!d}}}\ \
\hbox{\rm(\cite{11}
pp.~115, 226.)}.$$
$$\hskip0.9cm|{_{\!}}\Sigma_{{{\!o\!}}}|\ \hbox{\it is
homeomorphic to}\
 |{_{\!}}\Sigma_{{{\!o\!}}}|_{_{^{\!d}}}\ \underbar{iff}\ \
\Sigma_o\ \text{is locally finite}\ \text{\rm(\cite{30}
p.~119 Th.\ 8.)}.\hskip0.1cm \square
\quad
$$
\end{varprop}

\begin{varth}\label{Theorem:CountableRealization}\ {\rm(\cite{30} p.~120)}
                      If  ${\Sigma_o }$ has a realization in
                    $\mathbb{R}^{n}{\bf +}\{\wp\}$,
                    then  ${\Sigma_o }$ is countable and locally finite,
                    and $\dim{\Sigma_o }\leq n$. Conversely,
                    if ${\Sigma_o }$ is countable and locally finite,
                    and $\dim{\Sigma_o }\leq n$,
                    then  ${\Sigma_o }$ has  a realization as a closed
                    subset in  $\mathbb{R}^{2n+1}{\bf +}\{\wp\}$.
                    \hfill$\square$
\end{varth}

\subsection{Topological properties of realizations}\label{UnderavdP9:TopPropertiesOfRealizations}

\cite{6}
p.~352 ff.,
\cite{9}
p.~54 and \cite{18}
p.~171
emphasize the importance of triangulable spaces, i.e., spaces
homeomorphic to
the realization$_{_{\!^o}}\!$
(p.~\pageref{DefP8:RealizationOfSimpCompl}) of a simplicial
complex$_{_{\!^o}}\!$.
Non-triangulable topological (i.e., $0$-differentiable) manifolds
have been constructed, cf. \cite{28},
but since all continuously $n$-differentiable ($n>0$)
manifolds are triangulable, \label{triangulable}
cf. \cite{25}
p.~103 Th.\ 10.6, it is clear that essentially all spaces in practical use, in particular those within mathematical physics, are triangulable and therefore CW-complexes$_{_{\!^o}}\!$.

A topological space has the
\htmladdnormallink{{\it homotopy type}}
{http://en.wikipedia.org/wiki/Homotopy}
of the {\it realization}
(p.~\pageref{DefP8:RealizationOfSimpCompl})
{\it of a simplicial complex} (or, of a {\it polytope} for short)
if and only if (\underbar{\rm iff})
it has the type of an
\htmladdnormallink{{\it absolute neighborhood}}
{http://planetmath.org/encyclopedia/AbsoluteRetract.html}
\htmladdnormallink{{\it retract}}
{http://planetmath.org/encyclopedia/AbsoluteRetract.html}
(ANR),
which it has \underbar{\rm iff} it
has the type of a \htmladdnormallink{{\it CW-complex}}
{http://en.wikipedia.org/wiki/CW_complex},
\rm cf.\ \cite{11}
p.~226\ Th.\ 5.2.1.
Our CW-complexes$_{_{\!^o\!}}\!$ will have a $(-1)$-cell $\{\wp\}$,
as the ``relative CW-complexes'' defined in \cite{11}
p.$\ _{\!}$26, but their topology is such that they belong to the
category ${\mathcal D}_{_{\!\!^\wp}}$,
as defined in p.~\pageref{DefP7:CategoryDp}.

CW-complexes are
\htmladdnormallink{{\it compactly generated}}
{http://en.wikipedia.org/wiki/Compactly_generated_space},\
\htmladdnormallink{{\it perfectly normal spaces}}
{http://planetmath.org/encyclopedia/PerfectlyNormal.html},
cf. \cite{11},
pp.~22, 112, 242,
that are locally
\htmladdnormallink{{\it contractible}}
{http://en.wikipedia.org/wiki/Contractible},
in a strong sense and (hereditarily)
\htmladdnormallink{{\it paracompact}}
{http://en.wikipedia.org/wiki/Paracompact},
cf. \cite{11},
pp.~28-9\ Th.~1.3.2 \raise1pt\hbox{\eightrm+}\
Th.~1.3.5\ (Ex.~1, p.~33).\\
$\{\wp\}$ is called an {\it ``ideal cell''} in \cite{5},
p.~122.

\smallskip
\begin{varnorem}{\bf Notations.} \label{NotationsP9}
We have used
$\tau\!\!_{_{X}}$:=the topology of $X$.
We will also use; {\bf PID}:= Principal Ideal Domain and ``an {\bf
R-PID}-module {\bf G}'':= {\bf G} is a module over a {\bf PID} {\bf
R}.
${\bf L\mhatH S} :=$ Long $\mhatH$omology Sequence, {\bf
M-$\!$Vs} $\!:=\!$ Mayer-$\!$Vietoris sequence.
``$\simeq$''denotes ``homeomorphism'' or ``chain isomorphism''.

Let, here in Chapter~\ref{Ch3:AugmHomology},
$\Delta =\{{\Delta}^p,\partial\}$
be the classical singular chain complex.
\end{varnorem}

\begin{varnote}
In Definition\ 1 p.~\pageref{DefP1:1} we had $\{\wp\}=\{\emptyset_o\}$ and in the realization definition we had $\{\wp\}=\{\alpha_0\}$. $\emptyset$ is unquestionably the universal initial object. Our $\{\wp\}$ could be regarded as a local initial object and $\wp$ as a member identifier, since every object with $\wp$ as its unique ``$-1$-dimensional'' element, has therefore a ``natural'' membership tag, that fully complies with
\htmladdnormallink{Grothendieck's}{http://en.wikipedia.org/wiki/Grothendieck}
assumption that, for every set $\mathcal{S}$, there is a
\htmladdnormallink{``universe''}{http://en.wikipedia.org/wiki/Grothendieck_universe}
$U$ with $\mathcal{S}\in U$,
cf. \cite{21} p.~24.
Note also that $d(\alpha_0,\alpha)\equiv1$ for any $\alpha$, i.e., $\alpha_0$ is an isolated point.
\end{varnote}

\subsection{Simplicial augmental homology theory}\label{SubSecP9:2.2}

$\mhatH$ denotes the simplicial as well as the singular augmental (co)homology functor.
We begin with ten (10) lines that constitute a firm formal foundation for the simplicial augmental \underbar{pair}-(co)homology functor.

\smallskip
Two vertex-orderings are equivalent if an even permutation makes them equal.
Choose {\it ordered} $q$-simplices to generate the simplicial
\htmladdnormallink{{\it chain \nobreak{complex}}}
{http://en.wikipedia.org/wiki/Chain_complex}
$\mathbf{C}^{o}(\Sigma_o;{\bf
G})\}:=\{C^{o}\!\!_{_{\!\!^q}}(\Sigma_o;{\bf G})\}_{q}$, where the
coefficient module ${\bf G}$ is a unital
($\leftrightarrow\!1\!_{^{_{\mathbf{A}}}}\!\circ\!g=g$) module
over any commutative ring {\bf A} with unit.
We use {\it ordered simplices} instead of
\htmladdnormallink{{\it oriented}}
{http://en.wikipedia.org/wiki/Simplex}
since the
former is more natural in relation to our definition
p.~\pageref{DefP18:OrderedSimpCartProd},
of the {\it ordered simplicial cartesian product}.

Now, with {\bf0} as the additive unitelement,

\smallskip
{\bf-}${\bf C}^{o}(\emptyset;{\bf G})\equiv {\bf0}$
in all dimensions.

{\bf-}${\bf C}^{o}(\{\emptyset_o\};\!{\bf G})\equiv {\bf0}$
in all dimensions except for
$ C^{o}\!\!\!\!_{_{-1}}\!(\{\emptyset_o\};\!{\bf G})\ \! \cong\
\!{\bf G}.$

{\bf-}${\bf C}^{o}(\Sigma_o;{\bf G})\equiv \widetilde{{\bf C}}({{\mathcal{E}}}(\Sigma_o);{\bf G})\equiv$ the classical
``$\{\emptyset\!_{_{^{o\!}}}\!\}$''-augmented chain.

\medskip
Totaly analogous to the classical introduction of
\htmladdnormallink{{\it simplicial homology}}
{http://en.wikipedia.org/wiki/Simplicial_homology}
and by just hanging on to the
``$\{\emptyset\!_{_{^{o\!}}}\!\}$-augmented chains'', also when defining relative chains$_o$, we get the {\it relative
simplicial augmental homology functor for}
${\mathcal K}_{_{^{o}}}\!$-{\it pairs}, denoted $\mhatH_{\ast}$.
The graded quotient
$\mathbf{C}^{o}(\Sigma_{1o},\Sigma_{2o};\mathbf{G}):=
\mathbf{C}^{o}(\Sigma_{1o};\mathbf{G})/\mathbf{G}^{o}(\Sigma_{2o};\mathbf{G}):=
\{C^{o}_{q}(\Sigma_{1o};\mathbf{G})/C^{o}_{q}(\Sigma_{2o};\mathbf{G})\}_{q}$
induces the ``relative chains for the simplicial pair $(\Sigma_{1o},\Sigma_{2o})$''.
$$\mhatH_{i}(\Sigma_{o1},\Sigma_{o2};{\bf G})=
\left\{\begin{array}{llll}
\hbox{\rm H}_i({{\mathcal{E}}}(\Sigma_{o1}),{{\mathcal{E}}}(\Sigma_{o2});{\bf G}) &
\text{if}\ \Sigma_{o2}\neq\emptyset\\
\widetilde{{\rm H}}_i({{\mathcal{E}}}(\Sigma_{o1}),{\bf G}) & \text{if}\
\Sigma_{o1}\neq\{\emptyset_o\}, \emptyset\ \text{and}\ \Sigma_{o2}=\ \emptyset\\
\left\{\begin{array}{ll}
\cong {\bf G}  & \text{if}\ i=-1\\
={\bf0} & \text{if}\ i\neq -1
\end{array}\right.
& \text{when}\ \Sigma_{o1}\!=\!\{\emptyset_o\}\ \text{and}\
\Sigma_{o2}=\ \emptyset\\
{\bf0} & \text{for all $i$ when}\ \Sigma_{o1}=\Sigma_{o2}=
\emptyset.
\end{array}\right.
$$
${\rm H}_\ast$\ ($\widetilde{{\rm H}}_\ast)$ denotes the
classical (reduced) homology functor.

\subsection{Singular augmental homology theory}\label{SubSecP10:2.3}

The $|\sigma\!_{_{^{\!\wp\!}}}|$, defined in p.~\pageref{DefP8:Closedsigma0},
imbedded in $\mathbb{R}^{n}{\bf +}\{\wp\}$ generates a
satisfying set of\nobreak\ ``standard simplices$ \!_{_{^{\wp}}}\! $'' and ``singular simplices$\!_{_{^{\wp}}}$''. This implies, in
particular, that the ``$p$-standard simplices$\!_{_{^{\wp}}}$'',
denoted $\Delta\!^{{\wp}p}$, are defined by
$\Delta\!^{{\wp}p}\!:=\Delta\!^p{\bf +}\ \!\{\wp\}$ where
$\Delta\!^p$ denotes the usual $p$-dimensional standard simplex and
${\bf +}$ is the topological sum, i.e.
$\Delta\!^{{\wp}p}:=\Delta\!^p\sqcup\{\wp\}$ with the weak topology
with respect to $\Delta\!^p$ and $\{\wp\}$.
Now, and most important:
$\Delta\!^{^{_{\wp(\!-1\!)}}}\!\!:=\{\wp\}$.

Let $T^p$ denote an arbitrary classical singular $p$-simplex
($p\geq 0$). The ``singular $p$-simplex$\!_{_{^{\wp}}}\! $'',
denoted $\sigma\!\lower3pt\hbox{$^{^{{\ \!\!\wp}\ \!\!\!p}}$}$,
now stands for a function of the following kind, where the bar, $|$, in $ \sigma\!\lower3pt\hbox{$^{^{{\ \!\!\wp}\ \!\!\!p}}$}
\!\!\!_{_{|\!\Delta\!^p}}$ indicates ``restriction'':
$$
\sigma\!\lower3pt\hbox{$^{^{{\ \!\!\wp}\ \!\!\!p}}$}:
\Delta\!^{{\wp}p}=\Delta\!^p{\bf +}\{\wp\}\longrightarrow X{\bf +}
\{\wp\} \ \ \text{where}\ \ \sigma\!\lower3pt\hbox{$^{^{{\
\!\!\wp}\ \!\!\!p}}$}(\wp)=\wp\ \ \text{and}
$$
\centerline{$ \sigma\!\lower3pt\hbox{$^{^{{\ \!\!\wp}\ \!\!\!p}}$}
\!\!\!_{_{|\!\Delta\!^p}}\! =T^p \ \text{ for some classical
$p$-dimensional singular simplex}\ T^p\
\ \text{for all}\ p\geq 0.$}

In particular;
$$
\sigma\!\lower3pt\hbox{$^{^{{\ \!\!\wp}(\!-1\!)}}$}
:\{\wp\}\longrightarrow X\!_{_{^{\wp}}}\!= X{\bf +} \{\wp\} :\wp
\mapsto \wp.$$

\smallskip
The boundary function $\partial_{_{{^{^{\!}}}^\wp}}\!$
\label{DefBoundaryFunc} is defined by $
\partial\!_{_{{^{^{\!}}}^{\wp_{_{\!}}p}}}\!
(\sigma\!\lower1pt\hbox{$^{^{{_{_{\!}}}_{\wp_{_{\!}}p}}}$}\!)$
$:={\mathcal F}\!\!_{_{^{\!\wp}}}\!(\partial\!_p(T^p)) \ \ \hbox{\rm
if}\ p>0$ \noindent where $\partial\!_p$ is the classical singular
boundary function, and $\partial
_{\wp0}(\sigma\!\lower3pt\hbox{$^{^{{\ \!\!\wp}\ \!\!\!0}}$})
\!\equiv\! \sigma\!\lower3pt\hbox{$^{^{{\ \!\!\wp}(\!-1\!)}}$}$ for
every singular 0-simplex$\!_{_{^{\wp}}}\!$\
$\sigma\!\lower3pt\hbox{$^{^{{\ \!\!\wp}\ \!\!\!0}}$}\!\!.$
Let ${\Delta\!\!^{^{\wp}}}\!=
\!\{{\Delta}_{\!}\!^{^{_{\wp_{\!}p}}}\!\!,\partial_{\wp}\}$ denote
the singular augmental chain complex$_{_{\!}\wp}.$ Observe that;
$$|\Sigma \!_{_{^{o}}}\!| \neq\emptyset\Longrightarrow\! |\Sigma
\!_{_{^{o}}}\!|= {\mathcal F}\!\!_{\alpha_0}(|{{\mathcal{E}}} (\Sigma
\!_{_{^{o}}}\!)|)\in {\mathcal D}\!_{\alpha_0}.$$

\medskip
It is now obvious that e.g.; $\mhatH_{i}(|\Sigma_{o1}|,|\Sigma_{o2}|;{\bf G})
=\mhatH_{i}(\Sigma_{o1},\Sigma_{o2};{\bf G})$.
By the strong analogy to classical
\htmladdnormallink{{\it singular homology}}
{http://en.wikipedia.org/wiki/Singular_homology},
we omit the proof of the next lemma.

\begin{varlemma} \label{LemmaP10}
{\rm (Equivalently for co{\hatH}omology -- with the subindex $i$ as an index.)}
$$\mhatH_{i}(X_{\wp1},X_{\wp2};{\bf G})=
\left\{\begin{array}{llll}
{\rm H}_i({\mathcal F}(X_{\wp1}), {\mathcal F}(X_{\wp2});{\bf
G}) & \text{if}\ X_{\wp2}\neq \emptyset\\
\widetilde{\rm H}_i({\mathcal F}(X_{\wp1});{\bf G})
& \text{if}\ X_{\wp1}\neq \{\wp\},\ \emptyset\ \text{and}\
X_{\wp2}=\emptyset\\
\left\{\begin{array}{ll}
\cong{\bf G} & \text{if}\ i=-1 \\ ={\bf0} & \text{if}\ i\neq-1\\
\end{array}\right.
&
\text{when}\ X_{\wp1}= \{\wp\}\ \text{and}\ X_{\wp2}=\emptyset\\
{\bf0} &  \text{for all $i$ when}\ X_{\wp1}=X_{\wp2}=\emptyset. \quad
\end{array}\hskip-0.1cm\right.
$$
\end{varlemma}

\begin{varnote}
\begin{enumerate}\medskip
\item[i)]
\label{NoteP10:i}
${\Delta}_{\!}\!^{^{_{_{\!}i}}}\!(X_{1},X_{2};{\bf G})
\!\cong\!
{\Delta}_{\!}\!^{^{_{\wp_{\!}i\!}}}\!({\mathcal
F}\!\!_{\wp}(X_{1}),\ \!\!{\mathcal F}\!\!_{\wp}(X_{2});{\bf G})$
always.
So,
$$\mhatH_{-1}\!({\mathcal
F}\!\!_{\wp}\!(X_{1\!}), {\mathcal F}\!\!_{\wp}\!(X_{2});{\bf G})
\!\equiv\!{\bf0},\
\text{for any pair $(X_{1},X_{2})$ of classical spaces.}
$$

\medskip
\item[ii)]
\label{NoteP10:ii}
\hskip2cm$
{\Delta\!\!^{^{\wp}}}(X_{\wp1},X_{\wp2}\!) \simeq {\Delta}({\mathcal
F}(X_{\!\wp1\!}), {\mathcal F}(X_{\!\wp2}))$
\vskip0.3cm
except if, and only if,
$X_{\!\wp1\!}\neq_{\!}X_{\!\wp2\!} = \emptyset$, when
non-isomorphisms occur only in $\deg=-1$, i.e.,
$$ {\Delta}{\!\!}^{^{_{\wp_{\!}\hbox{\fivebf($_{\!}$-1$_{\!}$)}}}}\!(
X_{\!\wp1\!},\emptyset)\cong {\mathbb{Z}}
\ncong
{\bf0} \cong
{\Delta}_{\!}\!^{^{_{_{\!}\hbox{\fivebf($_{\!}$-1$_{\!}$)}}}}\!
({\mathcal F}(X_{\!\wp1\!}),\emptyset)$$
which, if $X_{\!\wp1\!}\neq\!\emptyset$, gives,
$$ {\mhatH}_{_{^{_{\!}0}}}\!(X_{\wp1},\emptyset) \oplus{\mathbb{Z}}
\cong
\hbox{\rm H}_{_{^{_{\!}0}}}\!({\mathcal F}(X_{\wp1}),\emptyset).$$
Remember that:
\hbox{${\bf Mor}\!\!\!\!\!\!\!_{_{_{\hbox{\fivesy
D}_{\!_{^{\wp}}}\!\!\hbox{{\fiverm
-pairs}}}}}\!\!((X\!_{_{^{\wp}}},X^\prime_{\!\!_{^{\wp}}}),
(X\!_{_{^{\wp}}},\emptyset))\!\neq\!\emptyset$
if, and only if,
$
X^\prime_{\!\!_{^{\wp}}}\!\!=\emptyset.{\!}$}
\vskip1pt

\medskip
\item[iii)]
\label{NoteP10:iii}
${\bf C}^{o}\!(\Sigma_{_{^{\!o1}}},\Sigma_{_{^{\!o2}}};{\bf G})
\approx {
\Delta\!\!^{^{\wp}}}\!(|\Sigma_{_{^{\!o1}}}|,|\Sigma_{_{^{\!o2}}}|;{\bf
G})$ connects the simplicial and singular functor$_{_{^{\!o}}}\!,$
where ``$\approx$'' stands for ``chain equivalence''.

\medskip
\item[iv)]
\label{NoteP10:iv}
$\mhatH_{0}(X_{_{\!}\wp}\raise1pt\hbox{\eightrm+}_{_{^{\!\!\!o}}}
Y_{\!\wp},\{\wp\};{\bf G}) = \mhatH_{0}(X_{_{\!}\wp},\{\wp\};{\bf
G})\oplus
\mhatH_{0}(Y_{\!\wp},\{\wp\};{\bf G}) $ \ but
\begin{eqnarray}
\mhatH_{0}(X_{_{\!}\wp}\raise1pt\hbox{\eightrm+}_{_{^{\!\!\!o}}}
Y_{\!\wp},\emptyset;{\bf G})
=\ \mhatH_{0}(X_{_{\!}\wp},\emptyset;{\bf G})\oplus\
\mhatH_{0}(Y_{\!\wp},\{\wp\};{\bf G})\ =\hskip2cm\phantom{h}&
\nonumber\\ \nonumber
=\ \mhatH_{0}(X_{_{\!}\wp},\{\wp\};{\bf G})\oplus\
\mhatH_{0}(Y_{\!\wp},\emptyset;{\bf G}).
\end{eqnarray}
\end{enumerate}
\end{varnote}

\begin{vardef} \label{DefP10}
The $p$:th {\it singular augmental homology group of}
$X_{\wp}$ with respect to ${\bf G}$ is
$\mhatH_{\!p}(X\!_{\wp};{\bf
G})\!:=\mhatH_{\!p}(X\!_{\wp},\emptyset;{\bf
G}).$
\smallskip
 The {\it coefficient group}$_{o}$
$:={\mhatH}_{\!-1}(\{\wp\},\emptyset;{\bf G})\!\simeq\!{\bf G}.$

\medskip
Using ${\mathcal F}_{\wp}\ ({{\mathcal{E}}}_{o}),$ we ``lift'' the foundational concepts of {\it homotopy}, {\it excision} and {\it point} in
${\mathcal C}\ ({\mathcal K})$ into ${\mathcal D}\!_{\wp}$-concepts
$({\mathcal K}_{o}$-concepts$)$ {\it homotopy}$_{\wp}$,
{\it excision}$_{\wp}$ and {\it point}$_{\wp}$ $(\bullet)$, respectively.
\end{vardef}

\smallskip
So, $f_{o},g_{o}\in{\mathcal D}\!_{\wp}$ are {\it homotopic}$_{\wp}$ if and only if $f_{o}=g_{o}=
0_{\emptyset,Y_o}$ or there are homotopic maps
$f_{1},g_{1}\in {\mathcal C}$ such that $f_{o}=f_{1}+\hbox{\rm
Id}_{\{\wp\}},
g_{o}=g_{1}+\hbox{\rm Id}_{\{\wp\}}.$

\smallskip
An inclusion$_o$
$$(i_o,i_{o{A_o}}): (X_o\setminus_{o}
\!U_o,A_o\setminus_o U_o) \longrightarrow (X_o,A_o),\
U_o\neq\{\wp\},$$
is an {\it excision}$_\wp$ if and only if there is an excision
$$(i,i_{A}):(X\smallsetminus U,A\smallsetminus U)
\longrightarrow (X,A)$$
such that $i_{o}=i+ \hbox{\rm Id}_{\{\wp\}}$ and
$i_{o{A_o}}=i_{_{A}}+ \hbox{\rm Id}_{\{\wp\}}.$

\smallskip
$\{{\bf P},\wp\}\in {\mathcal D}_{\wp}$ is a {\it
point}$_{\wp}$ \underbar{iff} $\{{\bf
P}\}+\{\wp\}={\mathcal F}_{\wp}(\{{\bf P}\})$ and $\{{\bf
P}\}\in{\mathcal C}$ is a {\it point}.\break
So, $\{\wp\}$ is \underbar{\bf not} a point$_{\wp}$.

\medskip
\begin{varnorem}{\underbar{\bf Conclusions}:} \label{RemarkP11:Conclusion}
$(\mhatH,{ {\partial}}\!_{\wp\!})$, abbreviated $\mhatH$, is a
{\it homology theory on the $h$-category of pairs from} ${{\mathcal
D}}\!_{\wp}$ $({\mathcal K}_{o})$,
c.f. \cite{9}
p.~117, i.e. $\mhatH$ fulfills the $h$-category analogues,
given in {\cite{9} \S\S8-9\ pp.~114-118}, of the classical \htmladdnormallink{``{\it Eilenberg-Steenrod axioms}''}
{http://en.wikipedia.org/wiki/Eilenberg-Steenrod_axioms}
from  \cite{9} \S3\ pp.~10-13.
The necessary verifications are either equivalent to the classical or completely trivial. E.g. the {\it dimension} {\it axiom} is fulfilled since $\{\wp\}$ is not a point$_{\wp}$ and the ``excision''-concept doesn't apply to a construction like $(X_o\setminus_{o}\{\wp\},A_o\setminus_{o}\{\wp\})$.

Since the exactness of the relative Mayer-Vietoris sequence of a
proper triad follows from the axioms, cf.
\cite{9}
p.~43, we will use this without further motivation, paying proper attention to Note iv above.

$ { \widetilde{\hbox{\rm H}}(X) = \mhatH({\mathcal
F}_{\wp}(X),\emptyset) } $ explains all the ad-hoc reasoning
surrounding the classical reduced homology functor $\widetilde{\hbox{\rm H}}$.
\end{varnorem}

\section{Augmental homology modules for products and joins}\label{SecP11:3}
\subsection{Background to the product- and joindefinitions} \label{SubSecP11:3.1}

What is usually called simply
\htmladdnormallink{`{\it\underbar{the}\ product topology}''}{http://en.wikipedia.org/wiki/Product_topology}
with respect to products including infinitely many factors is
actually the {\it Tychonoff product topology} which became the
dominant product topology when the {\it Tychonoff product theorem}
for compact spaces was introduced.
J.L. Kelley (1950) (J. Lo\'s and C. Ryll-Nardzewski (1954)) has
proven the
\htmladdnormallink{{\it Tychonoff product theorem} }{http://en.wikipedia.org/wiki/Tychonoff's_theorem}
(Tychonoff product theorem for
\noindent\htmladdnormallink{{\it Hausdorff spaces}}{http://en.wikipedia.org/wiki/Hausdorff_space})
to be equivalent to the
\htmladdnormallink{{\it axiom of choice} }{http://en.wikipedia.org/wiki/Axiom_of_choice}
(the
\htmladdnormallink{{\it Boolean} }{http://en.wikipedia.org/wiki/Boolean_prime_ideal_theorem}
\htmladdnormallink{{\it prime ideal theorem,} }{http://en.wikipedia.org/wiki/Boolean_prime_ideal_theorem}
which in its equivalent dual formulation is known as the
\htmladdnormallink{{\it ultra filter} }{http://en.wikipedia.org/wiki/Ultrafilter_lemma}
\htmladdnormallink{{\it lemma} }{http://en.wikipedia.org/wiki/Ultrafilter_lemma}).
Compare \cite{Brown:TenProductTopologies}.

Our ``realization'' of simplicial complexes  p.~\pageref{DefP8:alpha0} turns out to be a contravariant functor.

As was stated in p.~\pageref{FunctorialityOfJoin}:
The simplicial as well as the topological join-operation are (modulo realization, equivalent) cases of
\htmladdnormallink{colimits,}{http://en.wikipedia.org/wiki/Limit_(category_theory)}
 as being simplicial resp. topological attachments.
Restricting to topological $k$-spaces,
\cite{11}
p.~157ff remains true also with our new realization functor, i.e. {\it geometric realization preserves finite limits and all colimits.}

For a background on the join-definitions, see  \cite{Ehlers&Porter} and \cite{Fritsch&Golasinski} and
for the topological join-definitions, see the introductions in \cite{23} and \cite{36}, the definitions of which are quoted below.

Since $\emptyset$ is a simplicial zero-element with respect to simplicial join  and $\{\emptyset\!_o\}$ is the join-unit, it is absolutely impossible to functorially realize the \underbar{(augmented abstract)} simplicial complexes into the \underbar{classical} topological spaces including the \underbar{classical} Euclidian spaces.

That there should be a join-unit also in the category of
\underbar{classical} topological spaces and continuous maps
has been taken for granted. The join-unit in classical literature has been introduce with a jargong like: ``We use $\emptyset$ as join-unit'', or for example as in \cite{7} p.~135: ``...(it is a convention that $\emptyset\ast{Y}=Y,\ X\ast\emptyset=X$)''. The ``new'' topological space $\{\wp\}$ now provide general topology with a natural join-unit.

\subsection{Classical definitions of topological product and join}
\label{SubSec4.2:ClassicProdJoin}

\begin{definition}
{\rm(E.H.~Spanier
\cite{30}
p.~4)}
\label{TopProd:Spanier}
The {\it topological product} of an indexed collection of topological spaces $\{X\}_{j\in J}$ is the cartesian product $\times X_{j}$ given the topology induced by the projection maps
$p:{\times} X_{j}\longrightarrow X_{j}$ for $j\in J$.
\end{definition}

\begin{definition}
{\rm(
\cite{30}
p.~234)}
\label{TopPairProd:Spanier}
Given topological pairs $(X,A)$ and $(Y,B)$, we define their {\it product} $(X,A)\times(Y,B)$ to be the pair $X\times Y,(X\times B)\cup (A\times Y)$.
\end{definition}

\begin{definition}
{(G.W.~Whitehead
\cite{36}
p.~56)}
\label{TopJoin:Whitehead}
Let $X$ and $Y$ be spaces, which we assume to be disjoint from each other and from $X\times Y\times I$, where $I$ is the unit interval $\{t\mid 0\le t\le 1\}$.
Let $W=X\cup (X\times Y\times I)\cup Y$; we topologize $W$ by defining a subset to be open if and only if its intersection with each of the spaces $X$, $(X\times Y\times I)$ and $Y$ is open.
The {\it join of $X$ and $Y$} is the identification space $X\ast Y$ obtained from $W$ by identifying each $x\in X$ with all of the points $(x,y,0)$ and each $y\in Y$ with all of the points of $(x,y,1)$. The identification map sends $X$ and $Y$ homeomorphically into $X\!\ast\!Y$; Hence we may consider $X$ and $Y$ as subspaces of $X\ast Y$. Let $(t-1)x\oplus ty$ be the image of\nobreak$\ (x,y,t)\ \text{in}\ X\ast Y$.

The join operation is easily seen to be commutative (up to a natural homeomorphism). The join of $X$ with the empty set $\emptyset$ is $X$.
\end{definition}

\begin{definition}
{\rm(J.W.~Milnor,
\cite{23} p.~430.)}
\label{TopJoin:Milnor}
The join $A{_{_{^{1}}}}\!\circ \dots\circ\! A{_{_{^{n}}}}$ of $n$
topological spaces $A{_{_{^{1}}}}\!,\dots,\!A{_{_{^{n}}}}$  can be defined as follows. A point of the \hbox{\spaceskip1.5pt join is specified by}

\noindent
(1) $n$ real numbers $t{_{_{^{1}}}}\!,\cdots, t{_{_{^{n}}}}$ satisfying
$t{_{i}}\ge 0,$
$t{_{_{^{1}}}}\! + \cdots + t{_{_{^{n}}}}\!=1,$ and

\noindent
(2) a point $a{_{i}}\!\in A{_{i}}$ for each $i$ such that $t{_{i}}\not= 0.$
Such a point in $A{_{_{^{1}}}}\!\circ\cdots\circ A{_{_{^{n}}}}$ will be denoted by the symbol
$t{_{_{^{1}}}}\!a{_{_{^{1}}}}\!\oplus\cdots\oplus t{_{_{^{n}}}}\!a{_{_{^{n}}}}\!$
where the element $a{_{i}}$ may be chosen arbitrarily or omitted whenever the corresponding $t{_{_{^{1}}}}$ vanishes.

\noindent
By the {\sl strong topology} in $A{_{_{^{1}}}}\!\circ\cdots\circ A{_{_{^{n}}}}$
we mean the strongest topology such that the coordinate functions
$
t{_{i}}\!:A{_{_{^{1}}}}\!\circ\cdots\circ A{_{_{^{n}}}}\longrightarrow
[0,1]
$ and
$
a{_{i}}\!:t{_{i}}\!\!{^{^{_{-1}}}}\!\!(0,1]
\longrightarrow A{_{i}}
$
are continuous. Thus a sub-basis for the open sets is given by the sets of the following two types
\ (1) the set of all
$t{_{_{^{1}}}}\!a{_{_{^{1}}}}\!\oplus\cdots\oplus t{_{_{^{n}}}}\!a{_{_{^{n}}}}\!$
such that $\alpha<t{_{i}}\!<\beta,$

\noindent
(2) the set of all
$t{_{_{^{1}}}}\!a{_{_{^{1}}}}\!\oplus\cdots\oplus t{_{_{^{n}}}}\!a{_{_{^{n}}}}\!$
such that $t{_{i}}\!\not=0$ and $a{_{i}}\!\in U,$ where $U$ is
an arbitrarily open subset of $A{_{i}}\!.$

\noindent
The join of infinitely many topological spaces in the strong topology can be defined in exactly the same manner, with the restriction that all but a finite number of the $t{_{i}}$ should vanish. It is clear from the definition that the formation of finite or infinite joins in the strong topology is an associative, commutative operation.

\noindent
The strong topology is not the same as the more conventional
{\sl weak topology}, in which
$A{_{_{^{1}}}}\!\circ\cdots\circ A{_{_{^{n}}}}$
is considered as an identification space of the product of
$A{_{_{^{1}}}}\!\times\cdots\times A{_{_{^{n}}}}$
with an $(n-1)$-simplex.
\end{definition}

Let $A{_{_{^{1}}}}\!\mhatast\cdots\mhatast A{_{_{^{n}}}}$ denote $A{_{_{^{1}}}}\!\circ\cdots\circ A{_{_{^{n}}}}$ above equipped with the strong topology.

\subsection{Augmental versions of products and joins} \label{SubSecP11:3.2}
$$
\textrm{Let }\ \mytinynabla
\ \textrm{denote the above classical}\!:
\left\{\begin{array}{lll}
\textrm{topological product}\ \times\ {\rm from}\ \textrm{\cite{30}
p. 4}\ {\rm or}\\
\textrm{topological join}\ \ast\ {\rm from}\ \textrm{
\cite{36}
p.~56}\ {\rm or}\\
\textrm{topological join}\ \mhatast\ {\rm from}\ \textrm{
\cite{23}~p.~430.}
\end{array}
\right.
$$
Definition~2 p.~\pageref{DefP2:2} implies that the join-unit in the category of \underbar{classical} (abstract) simplicial complexes is $\emptyset$, the \underbar{classical} realization of which, $\emptyset$, is assumed to be the join-unit in the category of \underbar{classical} topological spaces, i.e.;
$$ X \mhatast\ \!{\emptyset} = X \ast \emptyset=X=
\emptyset\ast X=\emptyset\ \! \mhatast{X}. $$
\begin{vardef} \label{DefP11}
$$X_{\wp1}\mytinynabla\!\!\!_{{o}} X_{\wp2}:=
\left\{\begin{array}{ll}
\emptyset & if X_{\wp1}=\emptyset\ \hbox{\rm or}\ X_{\wp2}=\emptyset
\\
{{\mathcal F}_{\wp}}({{\mathcal
F}(X_{\wp1})\mytinynabla{\mathcal F}(X_{\wp2})}) & \hbox{\rm if} X_{\wp1}\neq \emptyset
\neq X_{\wp2}.\\
\end{array}\right.
$$
\end{vardef}

From now on we will delete the \raise0.6pt\hbox{$\wp$/\eightrm o}\ $\!_{_{\!}}$-indices. So, e.g. $\times, \ast, \mhatast$ now means $\times\!_{{^o}}, \ast\!_{{^o}}, \mhatast\!_{{^o}}$ respectively, while ``$X$ connected'' means ``${\mathcal F}(X)$ connected''.

\smallskip\noindent
{\bf Equivalent Join Definition.} \label{EquivJoinDef}
Put $$\emptyset{\myPsqcup1}X=X{\myPsqcup1}
\emptyset\!:=\!\emptyset$$
and
$$\{\wp\}{\myPsqcup1}X=\{\wp\},\ X{\myPsqcup1}
\{\wp\}=X,\ \text{if}\ X\ne\emptyset.$$
For $X,Y\ne\emptyset\ \text{or}\ \wp$,
let $X {\myPsqcup1} Y $ be ``$X \times Y \times(0,1]$
{\it pasted} to $X$'' along the function
$$ \varphi\!{_{_{^{1}}}}: X \times
Y \!\times\{1\} \longrightarrow X ;\ (x{_{_{^{}}}}\!,y ,1) \mapsto
x{_{_{^{}}}}, $$
i.e. the quotient set of $\big(X \times Y
\times(0,1]\big) \sqcup X ,$ under the equivalence relation
$(x,y,1)\sim x$ and let
$p_1:\big(X \times Y  \times(0,1]\big)\sqcup X
\longrightarrow X {\myPsqcup1}Y $ be the quotient function.

For $X,Y=\emptyset\ \text{or}\ \{\wp\}$ put $ X {\myPsqcup0}
\!Y := Y {\myPsqcup1}X  $ and else
$X {\myPsqcup0}\!Y :=X \times Y \times[0,1)$ {\it pasted} to the set $Y $ along  the function $$ \varphi_2: X \times Y \!\times\{0\}
\longrightarrow Y ;\ (x,y ,0) \mapsto y, $$
and let
$$ p_2:\big(X \!\times Y
\times[0,1)\big)\sqcup Y  \longrightarrow X {\myPsqcup0}\ Y  $$
be the quotient\nobreak\ function. \textrm{Put}\
$$X \!\circ\!Y  :=\big(X  {\myPsqcup1}\ Y \big) \cup
\big(X\myPsqcup0\ Y \big). $$

To any $ (x,y,t)\in X \times Y \times[0,1] $
there corresponds precisely one point $(x,y ,t)\in X {\myPsqcup1}Y
\cap X{\myPsqcup0}Y , $, if $ 0<t<1$ and
the ``equivalence class containing $x$'' if $t=1$
$(y\ \text{if}\ t=0)$, which is denoted $(x,1)\
((y,0)).$ \label{EquivJoinDef(x,1)}
This allows one to introduce ``coordinate functions''
\begin{eqnarray}
\!\xi:
X {\circ}Y \longrightarrow &[0,1],\nonumber\\
\eta_1:X{\myPsqcup1}\
\!Y  \longrightarrow &X,\nonumber\\
\eta_2: X{\myPsqcup0}Y  \longrightarrow  &Y
\nonumber\end{eqnarray} \hbox{\rm extendable\ to}
$X{\circ}Y$,
by setting
$$\eta_1(y,0):\equiv x_0 \in  X, $$
resp.
$$ \eta_2(x,1):\equiv y_0 \in  Y.  $$
Combining $p_1$ and $p_2$ we also have a projection $p$;
$$ p\!:X \sqcup \big(X \times
Y \times[0,1]\big) {\sqcup}Y \rightarrow X {\circ}Y .
$$

Let $X\mhatast Y $ denote $ X {\circ}Y$ equipped with the
smallest topology making $\xi, \eta_1,\eta_2 $
continuous and let $X {{\ast}}Y\ \text{be}\ X {\circ}Y  $ with the quotient topology with respect to $p$, i.e. the largest topology
making $p$ continuous
$(\Rightarrow \tau\!\! _{_{\tiny^{X\mtopast\!Y}}} \hskip-0.0cm
\raise1.5pt\hbox{{\eightsy {\char"1A}}}\
\tau\!\!_{_{\tiny X\ast Y\! }}).$

\begin{varnorem}{\bf Pair-definitions.} \label{DefP11:PairDef}
Let
$$({X}_{_{^{\!1}}},{X}_{_{^{\!2}}})
\mytinynabla
\!\!_{\ominus}\ \!
({Y}\!_{_{^{\!1}}},{Y}\!_{_{^{\!2}}}) \!:= \!({X}_{_{^{\!1}}}
\mytinynabla
\!\!\!_o \ \!{Y}\!_{_{^{\!1}}},({X}_{_{^{\!1}}}
\mytinynabla
\!\!\!_o \ \!{Y}\!_{_{^{\!2}}}) \ominus
({X}_{_{^{\!2}}}
\mytinynabla
\!\!\!_o \ \! {Y}\!_{_{^{\!1}}})),$$
where
${\ominus}$\ stand for ``$\cup$'' or ``$\cap$''
and
if either $X_{_{^{\!2}}}\!$ or $Y\!_{_{^{\!2}}}\ \!$ is not closed
\hbox{\rm(}open\hbox{\rm)},
$$ ({X}_{_{^{\!1}}}\!\!\ast\!_{_{^{\!\ }}}\!\! {Y}\!_{_{^{\!1}}},
({X}_{_{^{\!1}}}\!\!\ast\!_{_{^{\!\ }}}\!\!
{Y}\!_{_{^{\!2}}})\ominus ({X}_{_{^{\!2}}}\!\!\ast\!_{_{^{\!\
}}}\!\! {Y}\!_{_{^{\!1}}})):=
({X}_{_{^{\!1}}}\!\!\ast\!_{_{^{\!\ }}}\! {Y}\!_{_{^{\!1}}},
({X}_{_{^{\!1}}} \mcircast\ \! {Y}\!_{_{^{\!2}}})\ominus
({X}_{_{^{\!2}}} \mcircast\ \! {Y}\!_{_{^{\!1}}})),$$
where $\mcircast$ indicates that the subspace topology is to be used, i.e. the pair
$$({X}_{_{^{\!1}}}\!\!\ast\!_{_{^{\!\ }}}\! {Y}\!_{_{^{\!1}}},
({X}_{_{^{\!1}}}\!\circ{Y}\!_{_{^{\!2}}})\ominus
({X}_{_{^{\!2}}}\!\circ{Y}\!_{_{^{\!1}}}))$$
with the subspace
topology in the 2:nd component. (We will use similar pair-definitions for simplicial complexes
with, ``$\times\!$'' $($``$\ast_{\!}$''$)$ from
\cite{9}
p.~67 Def. 8.8
($\ \!$\cite{30}
 p.~109 Ex.\ 7) - of course, without any topologizing considerations.)
\end{varnorem}

\begin{varnote}\label{NoteP11}
\begin{enumerate}
\item[{\bf i.}]
With $t\in[0,1]$, the ``{Equivalent Join Definition}'' \label{JoinNote1}
above, suggest the notation
$(X\ast Y)^{_{^{t>0}}}=X\myPsqcup1 Y$ resp. $(X\ast Y)^{_{^{t<1}}}=X\myPsqcup0 Y$ and $(X\ast Y)^{_{^{t=1}}}=X$ resp. $(X\ast Y)^{_{^{t=0}}}=Y$.
Now, for any $s,t\in(0,1)$ we see that
$(X\ast Y)^{_{^{t\geq s}}}\!$ is homeomorphic to
the
\htmladdnormallink{{\it mapping}}
{http://en.wikipedia.org/wiki/Mapping_cylinder}
\htmladdnormallink{{\it cylinder}}
{http://en.wikipedia.org/wiki/Mapping_cylinder}
$C(q_1)$
with respect to the coordinate map
${q}_1:{X}\times{Y} \rightarrow  X $
and that
$(X\ast Y)^{_{^{t\leq s}}}$ is homeomorphic to
the mapping cylinder $C(q_2)$ with respect to the coordinate map ${q}_2:{X}\times\!{Y} \rightarrow  Y. $
\smallskip
\item[{\bf ii.}]
${X}\!_{_{^{2}}}\mhatast{Y}\!_{_{^{2}}}$
is a subspace of $X_{1} \mhatast Y_{1}$ by
\cite{2}
5.7.3 p.~163.
${X}\!_{{2}}{\ast}${\it Y}$\!_{{2}}$
is a subspace of$\!$ {\it X}$_{1}{\ast}${\it Y}$\!_{1}$ if
$X_{2},\ Y_{2}$ are both closed (open),
cf. \cite{7}
Th.~2.1(1) p.~122.
\smallskip
\item[{\bf iii.}] $({X}_{_{^{\!1}}}\!\circ{Y}\!_{_{^{\!2}}})\ \!{\cap}\ \!
({X}_{_{^{\!2}}}\!\circ{Y}\!_{_{^{\!1}}})\!=\!
{X}_{_{^{\!2}}}\!\circ{Y}\!_{_{^{\!2}}}$
and with
$\mytinynabla\!:=\mytinynabla\!\!_{\cup}$
for short in pair operations;
$({X}_{_{^{\!1}}},\{\wp\})\times
({Y}\!_{_{^{\!1}}},{Y}\!_{_{^{\!2}}})=
({X}_{_{^{\!1}}},\emptyset)\times
({Y}\!_{_{^{\!1}}},{Y}\!_{_{^{\!1}}})\ \hbox{\rm if}\
{Y}\!_{_{^{\!2}}}\neq\emptyset$.
\smallskip
\item[{\bf iv.}]
$\mhatast$ and $\ast$ are both commutative but, while
$\mhatast$
is associative by \cite{2}
p.~161, $\ast$ is not in general, cf.
p.~\pageref{P18:nonassociativity}.
\smallskip
\item[{\bf v.}]$\!$
``$\times\!\!_{_{^\cap}}\!$'' is (still, cf. {\cite{6}
p.~15},) the
\htmladdnormallink{{\it categorical product}}
{http://en.wikipedia.org/wiki/Product_(category_theory)}
on \underbar{pairs} from ${\mathcal
D}\!\!_{{\wp}}$.
\end{enumerate}
\end{varnote}

\subsection{Augmental homology for products and joins} \label{SubSecP12:3.2}

The following theorem and comments comes from \cite{30}, page 235 and it is the {\it classical Künneth formula for products}.

\centerline{------$\ast$------}
\begin{varth}[10] {\rm[The {\it \underbar{classical} relative Künneth formula for products} {(from \cite{30}
\htmladdnormallink{{ p.~235}
}{http://books.google.com/books?id=h-wc3TnZMCcC&pg=PA235&lpg=PA235&dq=\%22excisive+couple\%22&source=web&ots=LeRef7z0bO&sig=qcjk46rlISNAQPz_PRLEpFVXdas&hl=en&sa=X&oi=book_result&resnum=4&ct=result}).}]}\\
 {\it If $\{ X\times B,A\times Y\}$ is
an excisive couple in $X\times Y$ and $G$ and $G^\prime$ are modules
over a principal ideal domain $R$ such that ${\rm Tor}_1^R(G,G^\prime)=0$,
there is a functorial short exact sequence;
$$0\rightarrow [H(X,A;G)\otimes H(Y,B;G^\prime)]_q\buildrel
\mu^\prime\over\longrightarrow H_q((X,A)\times (Y,B);G\otimes G^\prime)
\longrightarrow$$
\smallskip
\hfill{$\longrightarrow[{\rm Tor}_1^R(H(X,A;G),H(Y,B;G^\prime)]_{q-1}
\longrightarrow 0$\indent

\medskip
\noindent
and this sequence is split.}}\hfill$\square$
\end{varth}

\smallskip
In particular, if the right hand side
vanishes (which always happens if $R$ is a field) then the cross product $\mu^\prime$
is an isomorphism.

\centerline{------$\ast$------}

By using Lemma $\!${\eightbf+}$\!$ Note {\bf ii} p.~\pageref{NoteP10:ii}
we convert the {\it classical Künneth formula} above,
mimicking what Milnor did
at the end of his proof of Lemma\ 2.1
\cite{23}
p.~431, where $H_0(X)=\tilde{H}_0(X)\oplus \mathbb{Z}$ and $H_r(X)=\tilde{H}_r(X)$ for $r\neq0$, is used to reach the first line in Theorem~\ref{TheoremP12:1} below. Recall that $\tilde{H}_r(X)=\mhatH_r(X,\emptyset)$ and that a pair $\{X,Y\}$ is an
\htmladdnormallink{{\hbox{\it excisive couple of subsets}}
}{http://books.google.com/books?id=h-wc3TnZMCcC&pg=PA188&lpg=PA188&dq=\%22excisive+couple\%22&source=web&ots=LeRef7z0bO&sig=qcjk46rlISNAQPz_PRLEpFVXdas&hl=en&sa=X&oi=book_result&resnum=4&ct=result}
if the inclusion chain map $\Delta\!^{\wp}(X) + \Delta\!^{\wp}(Y)\subset \Delta\!^{\wp}(X\cup Y)$ induces an isomorphism in homology.
The ``new'' object $\{\wp\}$ gives additional strength to the Künneth formula w.r.t. the \underbar{classical} Künneth formula ($\equiv$4:th line in Th.~\ref{RelativeKunnethformulaforproducts}), 
but much of the classical beauty is lost -- a loss which is regained in the join version i.e. in Theorem~\ref{TheoremP14:4}
p.~\pageref{TheoremP14:4} -- but still, this is now how the Künneth formula for products looks in any augmented environment. Th.~\ref{TheoremP14:4} is, via Note~\ref{NoteP21:iii}~{\bf iii} p.~\pageref{NoteP21:iii}, at the heart of modern algebraic topology.

\begin{theorem} \label{TheoremP12:1}\label{RelativeKunnethformulaforproducts}
{\rm[The {\it relative Künneth formula for products}.]}\\
If $\{X_{1}\!\!\times\! Y\!_{2},
X\!_{_{2\!\!}}\times\!Y\!_{1}\}$
is an \htmladdnormallink{{\it excisive couple of subsets}
}{http://books.google.com/books?id=h-wc3TnZMCcC&pg=PA188&lpg=PA188&dq=\%22excisive+couple\%22&source=web&ots=LeRef7z0bO&sig=qcjk46rlISNAQPz_PRLEpFVXdas&hl=en&sa=X&oi=book_result&resnum=4&ct=result}
{\rm(Def. \cite{30} p.~188),}
${\bf q}\geq0$, {\bf R} a {\bf PID}, and assuming $\hbox{\rm
Tor}_1^{\mathbf{R}}({\bf G},{\bf G}^\prime)=0$ for {\bf R}\hbox{-}modules ${\bf G}$ and ${\bf G}^\prime$, then;

\smallskip
\noindent$\mhatH_{q}((X_{1},X_{2})\times (Y_{1}, Y_{2});{\bf G}\
\!{{{\otimes}} _{\!}{{_{_{\mathbf{R}}}}{\bf
G}^\prime}})\mringHom{R}
$
$$\noindent\phantom{I}\hskip-1.7cm
\left\{\!\!\begin{array}{ll}
\smallskip
[\mhatH_i(X_{1};{\bf G})
\!\otimes\!\!{_{_{_{{\bf R}}}}}\!
{
{\mhatH_j(Y_{1};{\bf G}^\prime)]\lower2pt\hbox{$_{q}$}\!\oplus}
}
(\mhatH_{q}(X_{1};{\bf G})
\!\otimes\!\!{_{_{_{{\bf R}}}}}\!
\!{\bf G}^\prime)
\!\oplus\!
({\bf G}
\!\otimes\!\!{_{_{_{{\bf R}}}}}\!
\mhatH_{q}(Y_{1};{\bf G}^\prime))
\!\oplus\!{\bf T}\!_1\ \ & \text{if}\ {\bf C}_1
\\
\smallskip
{[\mhatH_i(X_{1};{\bf G})\otimes_{_{\mathbf{R}}}
\mhatH_j(Y_{1},Y_{2};{\bf G}^\prime)]}_q \oplus
({\bf G}\!\otimes_{_{\mathbf{R}}}\!
\mhatH_{q}(Y_{1},Y_{2};{\bf G}^\prime))\oplus{\bf T}_2
& \text{if}\ {\bf C}_2\hskip-0.4cm
\\
\smallskip
{[\mhatH_i(X_{1},X_{2};{\bf G})\otimes_{_{\mathbf{R}}}
\mhatH_j(Y_{1};{\bf G}^\prime)]}_q \oplus
(\mhatH_{q}(X_{1},X_{2};{\bf G}) \otimes_{_{\mathbf{R}}}{\bf G}^\prime)\oplus{\bf T}_3 & \text{if}\ {\bf C}_3,\hskip-0.4cm
\\
{[\mhatH_i(X_{1},X_{2};{\bf G})\otimes_{_{\mathbf{R}}}\mhatH_j(Y_{1},Y_{2};{\bf G}^\prime)]}_q\oplus{\bf T}_4\ \
& \text{if}\ {\bf C}_4,\hskip-0.4cm
\end{array}\hskip-0.4cm \right.\hskip-0.4cm
\eqno{{\bf(1)\hskip-0.0cm}}
\label{EqP12:1}\!\!
$$
\end{theorem}
where the torsion terms, i.e. the ${\bf T}\!_{_{^{}}}$-terms, split as those ahead of them,\nobreak\ e.g.,
$${\bf T}\!_{_{^{1}}}\!=\![\hbox{\rm Tor}_{1}^{\bf R}\bigl(
\mhatH_i(X_{1};{\bf G}), \mhatH_j(Y_{1}; {\bf
G}^\prime)\big)]_{_{q-1}}
\oplus\hbox{\rm Tor}_{1}^{\bf R}
\bigl(\mhatH\hskip-0.3cm\lower3pt\hbox{$_{q-1}$}\!(X_{1};{\bf G}), {\bf
G}^\prime\bigr)\oplus\hbox{\rm Tor}_{1}^{\bf R}\bigl({\bf
G},\mhatH\hskip-0.3cm\lower3pt\hbox{$_{q-1}$}\!(Y_{1};{\bf G}^\prime)\bigr),
$$
and
$$
{\bf T}\!_{_{^{4}}}=[\hbox{\rm Tor}_{1}^{\bf R}\hbox{\tenbf(} \mhatH_i(X_{1},X_{2};{\bf G}),\
\!\mhatH_j(Y_{1},Y_{2}; {\bf
G}^\prime)\hbox{\tenbf)}]{{\lower3pt\hbox{{\fivei q}{\fiverm-1}}}}.
$$
${\bf C}_i$, $(i=1-4)$, are ``conditions'' and should be interpreted as follows, resp.,
$$\left\{\begin{array}{ll}
{\bf C}_1:=\!``X\!_{1}\!\times\!Y\!\!_{1}\neq \emptyset,~\{\wp\}\ \text{and}\ X_{2}=\emptyset=Y_{2}\hbox{''},\\
{\bf C}_2:=\!``X\!_{1}\!\times\!Y\!\!_{1}\neq \emptyset,~\{\wp\}\ \text{and}\ X_{2}=\emptyset\neq\!Y_{2}\hbox{''},\\
{\bf C}_3:=\!``X\!_{1}\!\times\!Y\!\!_{1}\neq \emptyset,~\{\wp\}\ \text{and}\ X_{2}\neq\emptyset=Y_{2}\hbox{''},\\
{\bf C}_4:=\!``X_{1}\!\!\times\!\!Y\!_{1} = \emptyset,~\{\wp\}\ \ {\bf or}\ \ \ X_{2}\neq \emptyset\neq Y_{2}\hbox{''}.
\end{array}\right.$$
Let [\hbox{\bf...}]$_q$, as in \cite{30}
p.~235 Th.\ 10, be interpreted as \hskip0.5cm $\nobreak\
\raise1.5pt\hbox{
$\bigoplus\hbox{\tenbf...}\hskip-0.9cm
\lower3pt\hbox{$_{_{\scriptstyle i+j=q\&
i,j\geq0}}$} .$} \hfill\square
$

\begin{varlemma}\label{LemmaP12}
For a relative homeomorphism $ f\!\!:\!(X,A)\!\rightarrow \!(Y,B)\ (
$i.e.
$f\!:X\!\rightarrow\!Y$
is continuous and
$f:X\setminus A \rightarrow
Y\setminus B\ \text{is}$
a homeomorphism$)$,
let  $F\!\!:\!N\!\times{\bf I}\!\rightarrow\! N\! $ be a
strong neighborhood deformation retraction of $N$
onto $A$. If $B$ and
$f(_{\!}N )$ are closed in
$_{\!}N^\prime\!:=\!f(N\!\setminus\!
A)\cup B,$ then $B$ is a
strong neighborhood deformation retract of
$N^\prime\!$ by means of,
\vskip-0.4cm
$$\hskip-0.2cmF^\prime : N^\prime\!\!\times{\bf
I}\!\rightarrow N^\prime\ ;\
\left\{\begin{array}{ll}
{\!}\hbox{\rm({\it y},{\it t})\ $\!$\lower0pt\hbox{${{\mapsto}}$}\
$\!${\it y}} & \text{if}\
{\it y} \!\in\! {\it B}, {\it t} \!\in\! {\tenbf I} \\
{\!}\hbox{\rm({\it y},{\it t})\ $\!$\lower0pt\hbox{${{\mapsto}}$}}
{\hbox{\rm{\it f}$\circ${\it F} ({\it f}
$\!^{^{_{\hbox{\fiverm-}1}}}\! $\hbox{\ninerm(}{\ninei
y}\hbox{\ninerm)}, {\it t}) }}  &  \text{if}\ {\it y} \!\in\!
{\it f}({\it N})\setminus${\it B} ${\!}={\!} {\it
f}({\it N}\setminus{\it A}),\ {\it t}
\!\in\! {\bf I}.
\end{array}\right.
$$
\end{varlemma}
\begin{proof} {\rm(cf. p.~\pageref{PropP88i})}
$F^\prime$ is continuous by \cite{2} p.~34;~2.5.12 as being so when\nobreak\ $\hbox{\rm restricted}$ to any one of the \underbar{closed} subspaces $f(N)\times{\bf I}$, resp. ${\it
B}\times{\bf I},$ where $N^\prime\times{\bf
I}=(f(N)_{\!}\times_{\!}{\bf I})\cup({\it B}_{\!}\times_{\!}{\bf
I})$.
\end{proof}

\begin{varnote}{} \label{NoteP12}
$X$ (resp. $Y$) is a strong deformation retract of ``the mapping cylinder
with respect to its product projection'', which is homeomorphic to
$ ({X}\!\ast{Y})\!^{^{_{t\ge s}}}$ (resp. $({X}\!\ast{Y}_{_{^{\!\ }}}\!
)^{^{_{\!t\le 1-s}}}\!$), with $s,t\in[0,1]$.
Equivalently for
``$\mhatast$''-join, by the Lemma.
So,
$$^{\!}r\!\!:\!X \mhatast Y
\!\!\rightarrow\!X{\!}\ast^{\!}Y;
\left\{\begin{array}{ll}
(x,y,t)\!\mapsto\! x\ i\!f\ t\geq0.9,\\
(x,y,t)\!\mapsto\! y\ i\!f\ t\leq0.1, \\
(x,y,t)\!\mapsto\! (x,y,\frac{0.5}{0.4}(t^{\!}-{0.5})+0.5)\ \text{otherwise},
\end{array}\right.
$$
is a homotopy inverse of the identity,
i.e.,
$$\underline{\underbar{$X\ \mhatast\ Y$
and
$X\ast Y$
are homotopy equivalent.}}$$
\end{varnote}
This
\htmladdnormallink{{\it homotopy}}
{http://en.wikipedia.org/wiki/Homotopy}
is a well-known fact and it allows us to substitute ``$\ast$'' for any occurrence of ``$\mhatast$'', and vice versa, in any discussion of homology groups.

\begin{theorem}
{\rm(Analogously for $\mhatast$ by the last Note (mutatis mutandis).$)$}\\
If $
(X_{1},X_{2})\!\neq\!(\{\wp\},\emptyset)\!\neq\!(Y_{1},Y_{2})$\
and {\bf G} is an {\bf A}\hbox{-}module, then  \label{TheoremP12:2}

\bigskip \noindent
$
\medskip
\mhatH_{q}((X_{1},X_{2})\!\times\!(Y_{1},Y_{2});{\bf G})
\mringHom{A}
$\\
\noindent
$
\msmallRingHom{A}
\mhatH\!\!\!{_{_{q+1}}}\!((X\!_{_1}\!,\!X\!_{_2}\!)\ast
(\!Y\!\!_{_1}\!,\!Y\!\!_{_2}\!);{\bf G})
\oplus\
\mhatH\!_{_{q}}\!((X\!_{_1}\!,\!X\!_{_2}\!)\ast
(\!Y\!\!_{_1}\!,\!Y\!\!_{_2}\!)^{^{_{\!t\geq0.5}}}\hbox{\bf+}\
(X\!_{_1}\!,\!X\!_{_2}\!) \ast
(Y\!\!_{_1}\!,\!Y\!\!_{_2}\!)^{^{_{\!t\leq0.5}}}\!;{\bf G}){\bf=}
$
$$
\msmallRingHom{A}
{\left\{\begin{array}{ll}
\mhatH_{_{q+1}}
 (\!X\!_{_1}\!\ast Y\!\!_{_1};
{\bf G}) \oplus \mhatH_{q}(X\!_{_1};{\bf G}) \oplus
\mhatH_{q}(Y\!\!_{_1};{\bf G})\!\!\!\!\! & \text{if}\ \ \
{\bf C}_1
\\
\mhatH_{_{q+1}}((\!X\!_{_1},\emptyset)\ast
(Y\!\!_{_1},Y\!\!_{_2});{\bf G})\oplus
\mhatH_{q}(Y\!_{_1},Y\!\!_{_2};{\bf G}) & \text{if}\ \ \
{\bf C}_2
\\
\mhatH_{_{q+1}}((\!X\!_{_1},X\!_{_2})\ast
(Y\!\!_{_1},\emptyset);{\bf G})\oplus
\mhatH_{q}(X\!_{_1},X\!_{_2};{\bf G}) & \text{if}\ \ \
{\bf C}_3
\\
\mhatH_{_{q+1}}( (\!X\!_{_1},X\!_{_2})\ast
(Y\!\!_{_1},Y\!\!_{_2})  ;{\bf G}) & \text{if}\ \ \
{\bf C}_4
\indent\indent\
\\
\end{array}\right.
} \!
\eqno{\bf(2)}
\label{EqP12:2}
\nonumber
$$
where
$$\left\{\begin{array}{ll}
{\bf C}_1:=\!``X\!_{1}\!\times\!Y\!\!_{1}\neq \emptyset,~\{\wp\}\
{and}\ X_{2}\!=\!\emptyset\!=\!Y_{2}\hbox{''},\\
{\bf C}_2:=\!``X\!_{1}\!\times\!Y\!\!_{1}\neq \emptyset,~\{\wp\}\
{and}\ X_{2}\!=\!\emptyset\!\neq\!Y_{2}\hbox{''},\\
{\bf C}_3:=\!``X\!_{1}\!\times\!Y\!\!_{1}\neq \emptyset,~\{\wp\}\ {
and}\ X_{2}\!\neq\!\emptyset\!=\!Y_{2}\hbox{''},\\
{\bf C}_4:=\!``X_{1}\!\!\times\!\!Y\!_{1} = \emptyset,~\{\wp\}\
\ \ {\bf or}\ \ X_{2}\neq \emptyset\neq Y_{2}\hbox{''}.
\end{array}\right.$$
\end{theorem}

\begin{proof}
\label{DefOfJoinparts}
Splitting ${X}\!\ast\!{Y}$ at $t\!=\!0.5$ gives
$ ({X}\ast{Y})^{^{_{t\geq0.5}}}$ and $({X}\ast{Y})^{^{_{\!t\leq0.5}}}$. Now,
$({X}\ast{Y})^{^{_{t\geq0.5}}}$ is
\htmladdnormallink{{\it homeomorphic}}
{http://en.wikipedia.org/wiki/Homeomorphism}
to the ``
\htmladdnormallink{{\it mapping cylinder}}
{http://en.wikipedia.org/wiki/Mapping_cylinder}
$C(p_X)$
of the product
\htmladdnormallink{{\it projection}}
{http://en.wikipedia.org/wiki/Projection_(mathematics)}
$p_X:\ X\times Y\!\longrightarrow\! X;\ (x,y)\mapsto x$''
of which $X$ is a strong
\htmladdnormallink{{\it deformation retract.}}
{http://en.wikipedia.org/wiki/Deformation_retract}
Equivalently, $({X}\ast{Y})^{^{_{\!t\leq0.5}}}$ is homeomorphic to $C(p_Y)$ of the projection $p_Y$.

The relative {\bf M-$\!$V$\!$s} with respect to $\!$the excisive
couple of pairs
$$ \{({^{\!}}(\!X\!_{_1}\!,\!X\!_{_2}\!)_{^{\!}}\ast_{^{\!}}
(\!Y\!\!_{_1}\!,\!Y\!\!_{_2}\!){^{\!}})^{^{_{\!t\geq0.5}}}\!\!,\ \
({^{\!}}(\!X\!_{_1}\!,\!X\!_{_2}\!)_{^{\!}}\ast_{^{\!}}
(\!Y\!\!_{_1}\!,\!Y\!\!_{_2}\!){^{\!}})^{^{_{\!t\leq0.5}}}
\}
$$
splits since the inclusion of their topological sum into
$(\!X\!_{_1}\!,\!X\!_{_2}\!)_{\!}\ast_{\!}
(\!Y\!\!_{_1}\!,\!Y\!\!_{_2}\!)$
is pair null-homotopic, cf. \cite{26}
p.~141 Ex.\ 6c, and \cite{17} p.~32 Prop.~1.6.8.
Since the 1:st (2:nd) pair is acyclic if
$Y\!\!_{2}\ \!(X\!_{_2\!})\!\neq\!\emptyset$, we get
Theorem~\ref{TheoremP12:2}.

This proven
\htmladdnormallink{{\it homomorphisms}}
{http://en.wikipedia.org/wiki/Homomorphism}
remains true also when $\mhatast $ is substituted for $\ast$, by the homotopy equivalence in the last Note.
\end{proof}

Milnor finished his proof of \cite{23}
Lemma\ 2.1 p.~431 by simply comparing the r.h.s. of the ${\bf
C}_1$-case in Eq.\ 1 with that of Eq.\ 2. Since we are aiming at a
stronger result of ``natural chain equivalence'' in Theorem~\ref{TheoremP13:3} this
is not strong enough. We will therefore use the following three auxiliar results to prove our next two theorems.
We hereby avoid explicit use of ``proof by acyclic models''.

\newtheorem*{Varth}{}

\noindent
\begin{Varth}
\hskip-0.3cm{\bf5.7.4.} {\rm  (from \cite{2} p.~164.)} \label{P13:Brown5.7.4}
(${\bf E}^0\!:=\!\{e,\wp\}=\bullet$ \ denotes a point, i.e. a
$0$-disc.)\\
\indent There is a homeomorphism:
$$\nu:X \mhatast Y \mhatast {\bf E}^0
\longrightarrow
(X \mhatast {\bf E}_{_1}^0)\times
(Y \mhatast {\bf E}_{_2}^0)
$$
{\it which restricts to a homeomorphism}:
$$
X  \mhatast Y
\rightarrow
((X \mhatast {\bf E}^0_{_1})\times Y)\cup (X\times(Y \mhatast {\bf
E}^{0}_{_2})).
\qed$$
\end{Varth}

\begin{varcor}[5.7.9]  \label{CorP13:HiltonWylie}
{\rm(from \cite{17} p.~210.)}
\ \ {\rm(``$\approx$'' stands for ``chain equivalence''.)}\ \

If $\phi$:\ {\bf C} $\approx$ {\bf E} with inverse $\chi$ and
$\phi^\prime:$\ {\bf C}$^\prime\approx$ {\bf E}$^\prime$ with
inverse $\chi^\prime$, then
$$\phi\otimes\phi^\prime:{\bf C}\otimes{\bf
C}^\prime\approx {\bf E}\otimes{\bf E}^\prime\ \hbox{\it with
inverse}\ \chi\otimes\chi^\prime.\qed$$
\end{varcor}

The $\otimes$-operation in the
\htmladdnormallink{{\it monoidal category}
}{http://en.wikipedia.org/wiki/Monoidal_category}
(= tensor category) in this corollary can be substituted for any monoid-inducing operation in any category, cf. \cite{koch:frobalgand2dimtqft} Ex.~2 p.~168.

\noindent
\begin{varth}[46.2] {\rm(from \cite{26}
p.~279)}
\label{TheoremP13:Munkres}

Let $\!\ {\mathcal C}$ and ${\mathcal D}_{\!}$ be free chain complexes that vanish below a certain dimension; let $\lambda:{\mathcal C}\!\rightarrow\!{\mathcal D}$ be a chain map.
If $\lambda$ induces homology isomorphisms in all dimensions, then $\lambda$ is a chain equivalence.
\end{varth}

\begin{theorem} \label{TheoremP13:3}
$\!${\rm(The relative Eilenberg-Zilber theorem for
topological join.)}
For an excisive couple
$\{X\mhatast Y\!_{2},X_{2} \mhatast  Y\}$
from the category of ordered pairs $((X ,X_{2}\!),$ $\!(Y
,Y\!_{2}\!))\!$
of\ topological\ pairs$\!_{_{^{_{\wp}}}}{_{\!\!}},$
$$ {\bf s}{\bf (} { \Delta\!\!^{^{\wp}}} \!(\!X\!_{_{\
\!}}\!,X\!_{_2})\nobreak\otimes\nobreak{\Delta \!\!^{^{\wp}}}
\!(Y_{_{\!\!\ \!}}\!,Y\!\!_{_2}){\bf )}\ \textrm{is naturally chain
equivalent to}\ { \Delta\!\!^{^{\wp}}}\! {\bf (}(\!X\!_{_{\
\!}}\!,X\!_{_2}) \ \!\mhatast \ \! (Y\!\!_{_{\
\!}}\!,Y\!\!_{_2}){\bf )}.$$
{\rm (``{\bf s}'' stands for  suspension i.e. the suspended
chain equals the original one except that the dimension $i$ in the original chain becomes $i+1$ in the suspended chain. Th.~\ref{TheoremP13:3} is the join version of the classical Th. 9 in \cite{30}
p.~234.)}
\end{theorem}

\begin{proof}
The second isomorphism is the key one and
is induced by the pair homeomorphism in \cite{2}
5.7.4\ p.~164.
For the 2:nd last isomorphism we use \cite{17}
p.~210 Corollary\ 5.7.9 and that {\bf L{\hatH}S}-homomorphisms
are ``chain map''-induced. The last follows from the {\bf L{\hatH}S} since cones are null-homotopic.
Note that the second component in the third module is an excisive
union.
\vskip-1.3cm
\baselineskip=0.9cm
\begin{eqnarray}
{ \mhatH\!_{_{\!q\!}}(X
\ \!\!\mhatast \ \!
Y_{^{_{^{\!}}}})\ \
\mringHombb{Z}
 \ \ \mhatH\!\!\!_{_{\!q\hbox{\fivebf+}1\!}}\!(X\!_{}
\ \!\mhatast \ \!
Y\!_{}
\ \!\mhatast \ \!
\{\!{\bf v}_{_{\!}}\!,\wp\!\}, X\!_{}
\ \!\mhatast \ \!
Y\!_{})
}
\ \
\mringHombb{Z}
\nonumber\hskip2.5cm\\
{ \ \! \mhatH\!\!\!_{_{\!q\hbox{\fivebf+}1\!}}
        ((\!X\!_{}
\ \!\!\mhatast \ \!\! \{\!{\bf u}_{_{\!}},_{\!}\wp_{\!}\}\!)
\!\times\!(Y\!_{} \ \!\mhatast \ \!\!\{\!{\bf v}_{_{\!}}\!,\wp\}\!),
\hbox{\tenbf(}(\!X\!_{} \ \!\mhatast \ \!
\{\!{\bf u}_{_{\!}},_{\!}\wp_{\!}\}) \!\times\! Y\!_{}\hbox{\tenbf)}
\cup\! \hbox{\tenbf(}\!X\!_{} \!\times\! (Y\!_{} \ \!\mhatast \ \!
\{{\bf v}_{_{\!}}\!,\wp\})\hbox{\tenbf)} ) { =\ \! } }\hskip1.0cm
\nonumber\\
{ =\ } \mhatH\!\!\!_{_{\!q+1\!}}
        ((\!X\!_{}
\ \!_{\!}\mhatast \ \!_{\!} \{{\bf u}_{_{\!}},\wp\}, X) \!\times\!
(Y\!_{} \ \!_{\!}\mhatast \ \!_{\!} \{{\bf v}_{_{\!}}\!,\wp\}, Y) )\
\mringHombb{Z}
\nonumber\hskip2.7cm
\nonumber
\end{eqnarray}
\hskip1.0cm\mringHombb{Z}
\bigg[\lower0.38cm\vbox{\hsize8.5cm\scriptsize
\noindent
Motivation: The underlying chains on the l.h.s. and r.h.s. are, by
Note {\bf ii} p.~\pageref{NoteP10:ii}
isomorphic to their classical counterparts to which we apply the classical Eilenberg-Zilber Theorem.}\bigg]\
\mringHombb{Z}\hskip1.5cm
\normalbaselines
\begin{eqnarray}
\mringHombb{Z}
\ \hbox{\rm H}\!\!_{_{\!q+1\!\!}}
        ({\Delta\!\!^{^{\wp}}}(\!X\!_{}
\ \!\mhatast \ \! \{{\bf u}_{_{\!}},\wp\}, X)
\!\otimes\!\!\!_{_{_{\mathbb{Z}}}}\! {
\Delta\!\!^{^{\wp}}}(Y\!_{} \ \!\mhatast \ \! \{{\bf
v}_{_{\!}}\!,\wp\}, Y))\
\mringHombb{Z}\hskip1.5cm\nonumber\\
\hskip2.5cm\mringHombb{Z}
\ \hbox{\rm H}\!\!_{_{\!q+1\!\!}}
        ({\bf s}{\Delta\!\!^{^{\wp}}}\!
{\bf (} \!X\!_{}{\bf )} \!\otimes\!\!\!_{_{_{\mathbb{ Z}}}}\!
{\bf s}{ \Delta\!\!^{^{\wp}}}\! {\bf (}Y\!_{}{\bf )})
\mringHombb{Z}
\ \hbox{\rm H}_{_{\!q\!\!}}
        ({\bf s}[{\Delta\!\!^{^{\wp}}}\!
{\bf (} \!X\!_{}{\bf )} \!\otimes\!\!\!_{_{_{\mathbb{ Z}}}}\! {
\Delta\!\!^{^{\wp}}}\! {\bf (}Y\!_{}{\bf )}]).
\nonumber
\end{eqnarray}
\normalbaselines

\smallskip
So,
${\Delta\!\!^{^{\wp}}}\!(X
\ \!\!\mhatast \ \!
Y_{^{_{^{\!}}}})$
is naturally chain equivalent to
${\bf s}[{\Delta\!\!^{^{\wp}}}\! {\bf (} \!X\!_{}{\bf )}
\otimes\!\!_{_{_{\mathbb{Z}}}}\! { \Delta\!\!^{^{\wp}}}\!
{\bf (}Y\!_{}{\bf )}]$
by \cite{26}
p.~279 Th.\ 46.2 quoted above,
proving the non-relative Eilenberg-Zilber Theorem for\nobreak\
joins.

(The
\htmladdnormallink{{\it boundary map}}
{http://en.wikipedia.org/wiki/Chain_complex}
for the $\otimes$-complex
is given in \cite{30}
p.~228 together with the complete classical $\times$-proof in \cite{30}
p.~232 Theorem\ 6.)
\hfill$\triangleright$

\medskip
Substituting
1) ``$\times$'', 2) ``$\Delta$'', 3) ``Theorem 6'' in the original $\times$-proof for pairs given in
\cite{30}
p.\ 234, \label{RelSingHom}
with 1) ``$\ \!\mhatast \ \!$'',
2) ``${\Delta\!\!^{^{\wp}}}\!$ resp. $\!{\bf s}
{\Delta\!\!^{^{\wp}}}$'',
3) ``Theorem~\ref{TheoremP13:3},
1:st part'' respectively,
will do
since; (\cite{ZariskiAndSamuel} Th.~35 p.~184 gives the algebraic motivation.)
\vskip-1.5cm
\baselineskip=1.0cm
\begin{eqnarray}
{\bf s}\big({
\Delta\!\!^{^{\wp}}}(X\!_{_1})\otimes
{\Delta\!\!^{^{\wp}}}(Y\!_{_1})\big)/ \hbox{\twelvbf{\char"28}}{\bf
s}\big({ \Delta\!\!^{^{\wp}}}(X\!_{_1})\otimes
{\Delta\!\!^{^{\wp}}}(Y\!_{_2})\big)+ {\bf s}\big({
\Delta\!\!^{^{\wp}}}(X\!_{_2})\otimes
{\Delta\!\!^{^{\wp}}}(Y\!_{_1})\big) \hbox{\twelvbf{\char"29}}=\nonumber\\
=\! {\bf s}\big\{{ \Delta\!\!^{^{\wp}}}\!(\!X\!_{_1}\!)\otimes
{\Delta\!\!^{^{\wp}}}\!(Y\!\!_{_1}\!)/
\hbox{\twelvbf{\char"28}}\big({
\Delta\!\!^{^{\wp}}}\!(\!X\!_{_1}\!)\otimes
{\Delta\!\!^{^{\wp}}}\!(Y\!\!_{_2}\!)\big)+
\big({\Delta\!\!^{^{\wp}}}\!(\!X\!_{_2}\!)\otimes
{\Delta\!\!^{^{\wp}}}\!(Y\!\!_{_1}\!)\big)
\hbox{\twelvbf{\char"29}}\!\big\}
=\hskip0.5cm\nonumber
\end{eqnarray}
\normalbaselines

\medskip
\hfill
$
=
{\bf s}\big\{\hbox{\twelvbf{\char"28}}{
\Delta\!\!^{^{\wp}}}\!(\!X\!_{_1}\!)/ {
\Delta\!\!^{^{\wp}}}\!(\!X\!_{_2}\!)
\hbox{\twelvbf{\char"29}}\otimes \hbox{\twelvbf{\char"28}}{
\Delta\!\!^{^{\wp}}}\!(Y\!\!_{_1}\!)/ {
\Delta\!\!^{^{\wp}}}\!(Y\!\!_{_2}\!)\hbox{\twelvbf{\char"29}}\!\big\}.
\hfill
$
\hfill
\end{proof}
\begin{varcor}[4]  \label{Cor4P231:Spanier}
{\rm(from \cite{30} p.~231.)}
Given torsion-free chain complexes ${\bf C}$ and ${\bf C}^\prime$ and modules ${\bf G}$ and ${\bf G}^\prime$
\ such that\ \
$\hbox{\rm Tor}_1({\bf G},{\bf
G}^{\prime})={\bf0}$,
there is a functorial short exact sequence\\
\noindent\phantom{I}\hskip-0.2cm${\bf0}\rightarrow \![H({\bf C};{\bf G})\otimes H({\bf C}^\prime;{\bf G}^\prime)]_{_q}\!
{{\!\rightarrow\hskip-0.35cm}^{\mu^\prime}}
H_q({\bf C}\otimes {\bf C}^\prime;{\bf G}\otimes{\bf G}^\prime) \!\rightarrow\!
[\hbox{\rm Tor}_1(H({\bf C};{\bf G}),H({\bf C}^\prime;{\bf G}^\prime))]_{\!_{q-1}}\hskip-0.2cm\rightarrow\!{\bf0}$\\
and this sequence is split exact.  \hfill$\square$
\end{varcor}

Corollary~4
\cite{30}
p.~231 now gives Theorem~\ref{TheoremP14:4},
since\\
$\phantom{|}\hskip0.2cm
\mhatH\lower0pt\hbox{$\!{_{_{\star}}}$}\!(\circ) \!\cong\! {\bf
s} \mhatH\lower0.4pt\hbox{$\!\!\!{_{_{\star+1}}}$}\!(\circ)
\!\cong\! \mhatH\lower0.4pt\hbox{$\!\!\!{_{_{\star+1}}}$}\!({\bf
s}(\circ))
\ \ \ \text{\rm and}\ \ \
{\Delta\!\!^{^{\wp}}}\! {\bf (}(\!X\!_{_{\
\!}}\!,X\!{_{\!_2}}_{\!})\ast (Y\!\!_{_{\
\!}}\!,Y\!\!{_{\!_2}}_{\!}){\bf )}
\approx
{\Delta\!\!^{^{\wp}}}\! {\bf (}(\!X\!_{_{\ \!}}\!, X\!{_{_2}}_{\!})
\ \!\mhatast \ \! (Y\!\!_{_{\ \!}}\!,Y\!\!{_{_2}}_{\!}){\bf )}
$\\
by Th.~\ref{TheoremP12:2}
and \cite{26}
p.~279 Th.\ 46.2, quoted above.

\medskip
The couple
$\{X_{1}\ast Y_{2}, X_{2}\ast Y_{1}\}$
is excisive \underbar{iff}
$\{X_{1}\times Y_{2}, X_{2}\times Y_{1}\}$
is excisive, which is seen through the \hbox{\bf M-$\!$Vs}-stuffed 9-Lemma
p.~\pageref{3x3-lemma}ff
and Theorem~\ref{TheoremP12:2}
(line~4).

\goodbreak
\begin{theorem} \label{TheoremP14:4}
{\rm (The Relative K\"{u}nneth Formula for
Topological Joins; cf. \cite{30}
p.~235.)}

If
$\{X_{1}
\ \!\!\mhatast \ \! Y\!_{2}\!, X_{2}
\ \!\!\mhatast \ \! Y\!_{1}\}$
is an excisive couple in
$X_{1}\!\mhatast Y\!_{1}\!,$
{\bf R}\ a\ {\bf PID},\ {\bf G} and {\bf G}$^\prime$\
{\bf R}-modules  and
$
\hbox{\rm Tor}\!_1^{_{^{\mathbf{R}}}}\!({\bf G},{\bf
G}^{\prime})={\bf0}.
$,
Then the functorial sequences below are
$($non-naturally$)$\nobreak\ split\nobreak\ exact;
\begin{eqnarray}
{\bf0}\longrightarrow {\rlap{$_{_{_{i+j=q}}}$} {\
\raise2pt\hbox{$\bigoplus$}}}\ \ [\mhatH_{_{i}}
           (X\!_{_1},X\!_{_2};{\bf G})
\otimes\!\!_{_{\mathbf{R}}}
        \mhatH_{_{j}}
(Y\!_{1},Y\!_{2};{\bf G}^{\prime})]
\longrightarrow\hskip4cm
\nonumber
\\
\longrightarrow
\  \mhatH_{_{q+1}}
        ((\!X\!_{_1},X\!_{_2})\
\mhatast \ (Y\!_{_1},Y\!_{2}); {\bf
G}\otimes\!\!_{_{\mathbf{R}}}{\bf G}^{\prime})
\longrightarrow
\hskip3cm{{\hbox{\rm({\bf 3}\label{EqP14:3})}}}\!\!
\nonumber
\\
\hskip3cm\longrightarrow
{{\rlap{$\!\!\!\!_{_{_{i+j=q-1}}}$}}{\raise2pt\hbox{$\bigoplus$}}}\!\!
\ \ \hbox{\rm Tor}_1^{\mathbf{R}}\bigl(\mhatH_{_{i}}
           (X\!_{_1},X\!_{2};{\bf G}),
 \mhatH_{_{j}}
      (Y\!\!_{_1},Y\!_{2};
{\bf G}^{\prime})\bigr)
\longrightarrow\ {\bf0}
\hskip1.5cm\square
\nonumber
\end{eqnarray}
\end{theorem}
Analogously for the $\ast$-join.

\cite{30}
p.~247 Th.\ 11 gives the co$\mhatH$omology-analog of Theorem~\ref{TheoremP14:4}.

\medskip
$(\!X_{1},X_{2}\!)=(\{\wp\},\emptyset)$ in
Theorem~\ref{TheoremP14:4},
immediately gives our next theorem.

\begin{theorem}\label{TheoremP14:5}
{\rm[The Universal Coefficient Theorem for (co){\hatH}omology]}

$$ \mhatH_{_{i\!}} (Y\!\!_{_1},Y\!\!_{_2};{\bf G})\
\mringHom{R}
\ \big[\hbox{\rm \cite{30}
p.~214}\big]\
\mringHom{R}
\ \mhatH_{_{i\!}} (Y\!\!_{_1},Y\!\!_{_2};{\bf
R}\otimes\!_{\!_{\mathbf{R}}}\!{\bf G})\
\mringHom{R}
$$

$$
\mringHom{R}
\
\hbox{\twelvbf{\char"28}}\mhatH_{_{i\!}}(Y\!\!_{_1},Y\!\!_{_2};{\bf
R}) \otimes\!_{\!_{\mathbf{R}}}\!{\bf
G}\hbox{\twelvbf{\char"29}}\oplus {\rm Tor}^{\!\hbox{\fivebf
R}}_1\! \big(\mhatH_{_{{i\!-\!1}}}\! (Y\!\!_{_1},Y\!\!_{_2};{\bf
R}), {\bf G}\big),
$$

\noindent for {any }{\bf R}-{\bf PID} {module}\ {\bf G}.

\medskip
If\ all\ $\mhatH_{{\ast}}\!
        (Y\!_{_1},Y\!_{_2};{\bf R})$
are of finite type or {\bf G} is finitely generated, then;

\bigskip
\noindent \centerline{$ \mhatH^{{i\!}}
(Y\!\!_{_1},Y\!\!_{_2};{\bf G})
\mringHom{R}
\mhatH^{{i\!}} (Y\!\!_{_1},Y\!\!_{_2};{\bf
R}\otimes\!_{\!_{\mathbf{R}}}\!{\bf G})
\mringHom{R}
\hbox{\twelvbf{\char"28}}\mhatH^{{i\!}}(Y\!\!_{_1},Y\!\!_{_2};{\bf
R}) \!\otimes\!\!_{\!_{\mathbf{R}}}\!\!{\bf
G}\hbox{\twelvbf{\char"29}}\!\oplus\! \hbox{\rm Tor}^{\!\hbox{\fiverm R}}_1\! \big(\mhatH^{{{i{\raise1pt\hbox{\fivebf{+}}1}}}}\!
(Y\!\!_{_1},Y\!\!_{_2};{\bf R}), {\bf G}\big).
\ \square$
}
\end{theorem}

\noindent
{\bf N.B.} There can be no {\it \underbar{classical} \underbar{relative} K\"{u}nneth Formula for Topological Joins}, cf. p. \pageref{Milnor'sReducedFormula}.

\subsection{Local homology groups for products and joins}\label{SubSecP15:3.3}

Proposition\ \ref{PropP15:1} below motivates in itself the introduction of a topological
$(-1)$-object, which imposed the definition
p.~\pageref{DefP8:PointSetMinus0}
of a ``setminus'', ``$\setminus\!{_{_{^{o\!}}}}$'', in ${{\mathcal
D}\!_{\wp}}$, revealing the shortcomings of the classical boundary
definitions with respect to manifolds,
cf. p.~\pageref{DefP2:BoundaryOfSimplManif} and
p.~\pageref{DefP17:BoundaryOfHomologyManif}.
Somewhat specialized, Proposition\ \ref{PropP15:1} below is found in
\cite{12}
p.~162
and partially also in \cite{27}
p.~116 Lemma 3.3.
``${ X}\setminus x$'' usually stands for ``${ X}\setminus\{x\}_{\!}$''
and we will write $x$ for $\{x,\wp\}$ as a notational convention.
Recall the definition of $\alpha_{_{^{\!0}}}\!$, p.~\pageref{DefP8:alpha0},
and that
\hbox{$\dim{\rm Lk}\!\lower1.2pt\hbox{$_{_{\!\Sigma_{_{\!{\
}}}}}$}\!\!\!\sigma\!=\!{\dim\!{\Sigma
\!_{\!}-\!
\hbox{\eightbf\#} \sigma\!}},$} where $\#$ stands for cardinality.
\noindent
The {\it contrastar} of\ $\sigma\!\in\!\Sigma\!=\hbox{\rm
cost}_{_{\!{\Sigma}}}\!\sigma\!:=\! \{\tau\!\in\!\Sigma|\
\tau\!\not\supseteq\sigma\}.$
So, $\hbox{\rm
cost}_{_{{\!\Sigma}}}\!\emptyset\!_{o}\!\!=\!\emptyset\!$
\ and\
${\hbox{\rm cost}_{_{{\!\Sigma}}}\!\sigma\!=\!\Sigma\
\underline{\hbox{\rm iff}}\ \sigma\!\not\in\!\Sigma}.$
In the proof we use the following simplicial complexes,
${\bar{\sigma}} :=\{\tau \mid \tau
\hbox{\scriptsize$\subseteqq$}
\sigma\}
\ \text{and}\
{\dot{\sigma}}\!:=\{\tau \mid \tau
\hbox{\scriptsize$\subsetneqq$}
\sigma\}$.
So, ${ \bar{\emptyset}}\!_{o}\!=\{\emptyset\!_{o}\!\}$ and
${ \dot{\emptyset}}\!_{o}\!=\emptyset.$
${\rm Int} (\sigma) \!:=\! \{\alpha\in|\Sigma|\mid[{\bf
v}\in\nobreak\sigma ] \Longleftrightarrow [\alpha({\bf v})\neq0]
\}.$

\begin{proposition} \label{PropP15:1}
Let ${\bf G}$ be any module over a commutative ring {\bf A}$_{\!}$
with unit. With $\alpha\ \!\!\in\ \!\!\hbox{\rm Int}_{\!}\sigma$
and $\ {\!}\alpha_{\!}=_{\!} \alpha_{_{^{\!0}}}\!$
{\underbar{\hbox{\rm iff}}}\ $\sigma\!=\!\emptyset\!_{_{^{o}}}{\!}$
the following module isomorphisms are all induced by chain
equivalences
$(${\rm cf. \cite{26}
p.~279 Th.\ 46.2 quoted here in p.~\pageref{TheoremP13:Munkres}}$)$
$$\hskip-0.3cm\mhatH
\vbox{\moveleft0.0cm\hbox{\lower0.0pt\hbox{$_{{{_{\!}{i
\lower1.0pt\hbox{-} \hbox{\fivebf\#}_{\!}\sigma}}}}$}}}
\vbox{\moveleft0.05cm\hbox{$({\rm Lk}_{_{\!\Sigma}}\sigma;{\bf
G})$}}
\mringHom{A}\
\mhatH_{i}(\Sigma,\hbox{\rm cost}_{_{\!\Sigma}}\sigma;{\bf G})
\mringHom{A}\
\mhatH_{i}(\vert\Sigma\vert,\vert\hbox{\rm
cost}_{_{\!\Sigma}}\sigma\vert;{\bf G})
\mringHom{A}\
\mhatH_{i}(|\Sigma|,|\Sigma|\setminus\!{_{_{^{\!o\!}}}}\alpha;{\bf
G}),
$$
$$\hskip-0.3cm\mhatH
^{_{\!}i\raise0.5pt\hbox{\fiverm-}\hbox{\fivebf\#}_{\!}\sigma}
\!_{\!}({\rm Lk}_{_{\!\Sigma}}\sigma;{\bf G})
\mringHom{A}\ \!
\mhatH^{i}(\Sigma,\hbox{\rm cost}_{_{\!\Sigma}}\sigma;{\bf G})
\mringHom{A}\ \!
\mhatH^{i}(\vert\Sigma\vert,\vert\hbox{\rm
cost}_{_{\!\Sigma}}\sigma\vert;{\bf G})
\mringHom{A}\ \!
\mhatH^{i}(|\Sigma|,|\Sigma|\setminus\!{_{_{^{\!o\!}}}}\alpha;{\bf
G}).
$$
\end{proposition}

\begin{proof}
(Cf. $\!$definitions p.~\pageref{SecP34:Appendix}ff.) The
$``\setminus_o\!\!$''-definition p.~\pageref{DefP8:PointSetMinus0}
and
\cite{26}
Th.\ $\!$46.2 p.~279 $\!$\raise1pt\hbox{\eightbf+} $\!$pp.~194-199 Lemma 35.1-35.2 \raise1pt\hbox{\eightrm+}$\!$
Lemma\ 63.1 p.~374 gives the two last isomorphisms since
$ \vert\hbox{\rm cost}_{\Sigma}\sigma\vert $
is a deformation retract of
$ \vert\Sigma\vert\!\setminus\!_{_{o}}\alpha, $
$\hbox{\rm while}$ already on the chain level we have
$$C^{o}_{_{\star}}(\Sigma,\hbox{\rm
cost}_{\Sigma}\sigma) \ \!=\!\
C^{o}_{_{\star}}(\overline{\hbox{\rm {st}}}
_{_{\!\Sigma}}(\sigma),\ \! { \dot{\sigma}}\ast{\rm
Lk}_{_{\!\Sigma}}\sigma)
\ \!=\ \! C^{o}_{_{\star}}(\overline{{{\sigma}}}
\ast{\rm Lk}_{_{\!\Sigma}}\!\sigma,\ \!
{\dot{\sigma}}\ast{\rm
Lk}_{_{\!\Sigma}}\!\sigma) \ \!\simeq\
C^{o}_{_{\!\!\star\lower1pt\hbox{-}\#\sigma}}\! ({\rm
Lk}_{_{\!\Sigma}}\sigma).$$
\end{proof}

\begin{varlemma} \label{lemmaP15}
If $x\ (y)$ is a closed point in $X\ (Y)$, then both
$\!\{\!X\!_{_{}}\!\!\times\! (Y\!\!_{_{\
}}\!\!\setminus\!\!{_{_{^{o\!}}}} y), (X\!_{_{\
}}\!\!\setminus\!\!{_{_{^{o\!}}}} x)\!\times\!Y\}$
\ and\
$\{\!X\!_{_{}}\!\ast (Y\!\!_{_{\ }}\!\!\setminus\!\!{_{_{^{o\!}}}}
y), (X\!_{_{\ }}\!\!\setminus\!\!{_{_{^{o\!}}}}  x)\ast Y\!\!_{_{\
}}\!\}$
are excisive pairs.
\end{varlemma}

\begin{proof}
\cite{30}
p.~188 Th.\ 3, since $X\!_{_{\ }}\!\!\times\! (Y\!\!_{_{\
}}\!\setminus\!\!{_{_{^{o\!}}}}  y\!_{_{0}}\!)$
 $\big(X\!_{_{\ }}\!\!\ast\!
(Y\!\!_{_{\ }}\!\setminus\!\!{_{_{^{o\!}}}}  y\!_{_{0}}\!)\big)$
is open in $(X\!_{_{\ }}\!\!\times\! (Y\!\!_{_{\
}}\!\setminus\!\!{_{_{^{o\!}}}} y\!_{_{0}}\!)) \cup ((X\!\!_{_{\
}}\!\setminus\!\!{_{_{^{o\!}}}}  x\!_{_{0}}\!)\times\! Y\!_{_{\
}}\!\!) $
$\big(\!\text{resp.}(X\!_{_{\ }}\!\!\ast (Y\!\!_{_{\
}}\!\setminus\!\!{_{_{^{o\!}}}} y\!_{_{0}}\!)) \cup ((X\!\!_{_{\
}}\!\setminus\!\!{_{_{^{o\!}}}}  x\!_{_{0}}\!)\ast Y\!_{_{\
}}\!\!)\!\big) $, which proves the excisivity.
\end{proof}

\begin{theorem} \label{TheoremP15:6}
If $(t\!{_{_{^{1}}}}\!,\widetilde{x_{{\!}}\ast_{{\!}}y},
t\!{_{_{^{2}}}}\! )\!:=\!\{\!(x,y,t)\ \!{\mid}\
\!0_{{\!}}<_{{\!}}t\!{_{_{^{1}}}}\!\!\le \!t\!\le\!
t\!{_{_{^{2}}}}\!\!<\!1\}$ for closed points
$x\!_{\!}\in \!\!X$ and $y\!\in \!Y$ we have

\medskip
\noindent {\bf i.}\ \ \ $ \mhatH\!\!\!_{_{\!q+1\!}}({X} \
\!\mhatast \ \! {Y}\!,{X} \ \!\mhatast \ \! {Y}
\setminus\!{_{_{^{o\!}}}} (x ,y ,t);{\bf G})
\mringHom{A}\ \!
\mhatH\!\!\!_{_{\!q+1\!}}({X} \ \!\mhatast \ \! {Y}\!,{X} \
\!\mhatast \ \! {Y}\setminus\!{_{_{^{o}}}}\
\!(t\!{_{_{^{1}}}},\widetilde{x \ast y},t\!{_{_{^{2}}}});{\bf G}) \
\mringHom{A}\ \!
$
$$
\mringHom{A}\ \!
\mhatH_{_{\!q\!}}({X}\!\times\!{Y}\!,{X}\!\times\!{Y}
\setminus\!{_{_{^{o}}}}\! (x,y);{\bf G})
\mringHom{A}\ \!
\Big[\phantom{}^{\text{\scriptsize\rm Motivation: A\ }}_{\text{\scriptsize\rm simple calculation.}}\Big]
\mringHom{A}\!
$$
$$
{\mhatH_{_{\!q\!}}((X \!,X \!\setminus\!\!{_{_{^{o\!}}}} x ) \!\times\! (Y \!,Y\!\setminus\!\!{_{_{^{o\!}}}} y );{\bf G}) }
\mringHom{A} \
\Big[\phantom{}^{{\scriptsize\rm
Motivation:\ Th.~{\ref{TheoremP12:2}\!
}}}_{{\scriptsize\rm
\ \ p.~\pageref{TheoremP12:2}\
line\ four.}} \Big]
\mringHom{A} \
\mhatH\!\!\!_{_{\!q+1\!}}(
        (X \!,
X \!\setminus\!\!{_{_{^{o\!}}}}   x )
 \
\!\mhatast \ \!
(Y \!,Y \!\setminus\!\!{_{_{^{o\!}}}} y );{\bf G}).
$$

\vskip0.2cm
\noindent {\bf ii.}\ \
$ \hskip1.5cm \mhatH\!\!\!_{_{\!q+1\!}}
        (\!X\!_{}
\ \!\mhatast \ \! Y\!_{}, X\!_{} \ \!\mhatast \ \! Y\!_{}
\setminus\!{_{_{^{o}}}}(y ,{0});{\bf G})\
\mringHom{A}\ \!
\mhatH\!\!\!_{_{\!q+1\!}} ((X \!,\emptyset)\ \!\mhatast \ \! (Y
\!,Y \!\setminus\!\!{_{_{^{o}}}} y \!\ ) ;{\bf G})
$\\
and equivalently for the $(x ,{1})$-points.

\smallskip
$_{\!}$All isomorphisms are induced by
\htmladdnormallink{chain $($homotopy$)$ equivalences}
{http://en.wikipedia.org/wiki/Homotopy_category_of_chain_complexes}
{\rm(cf. \cite{26} p.~279 Th.\ 46.2 quoted here in p.~\pageref{TheoremP13:Munkres})}.

Analogously for co$\mhatH$omology and for $``\ast$'' substituted\ for $``\mhatast$''.
\end{theorem}

\begin{proof}
{\bf i.}
\indent $\left\{\begin{array}{ll}
&
 A:=X {\myPsqcup1}Y \setminus\!{_{_{^{\!o}}}} \{(x_0,y_0,t)\mid\ t\!{_{_{^{1}}}}\le t<1 \}
\\
&
B:=X {\myPsqcup0}\ Y \!\setminus\!{_{_{^{\!o}}}}\ \!\!
\{(x_0,y_0,t)\mid\  0<t\le
t_2\}
\\
\end{array}\right.
\indent
\Longrightarrow $

\medskip
$
\hskip3.0cm
\Longrightarrow\!
\indent
\left\{\begin{array}{ll}
&
\!\!\!\!\!\! A\cup B =\!X\!_{} \ \!\mhatast \ \!
Y\!_{}\setminus\!{_{_{^{\!o}}}}\ \!(t\!{_{_{^{1}}}},\widetilde{x
\ast y}, t\!{_{_{^{2}}}} )
\\
& \!\!\!\!\!\! A\cap B=\!X \!\times\!Y \!\times\!(0,1)
\!\setminus\!{_{_{^{\!o}}}}\ \!\! \{x\!{_{_{^{0}}}}\}\!\times\!
\{y\!{_{_{^{0}}}}\!\}\!\times\!(0,1)\hbox{\tenbf,}
\\
\end{array}\right.
$

\medskip
\noindent with\ \
$ x\!{_{_{^{0}}}}\times y\!{_{_{^{0}}}}\times(0,1)\!:=
\!\{x\!{_{_{^{0}}}}\!\}\times \{y\!{_{_{^{0}}}}\!\}\times \{t\
\!\vert\ \! t\in(0,1) \}\cup\{\wp\}
\ \hbox{\rm and}\
(x\!{_{_{^{0}}}},y\!{_{_{^{0}}}},t\!{_{_{^{0}}}})\!:=\!
\{(x\!{_{_{^{0}}}},y\!{_{_{^{0}}}},t\!{_{_{^{0}}}}),\wp\}.
$

\medskip\noindent
Now, using the null-homotopy in the relative {\bf M-$_{\!}$Vs} with
respect to $ \{ ( X {\myPsqcup1}Y,A)$,
$(X {\myPsqcup0}Y,B)_{\!}\}$
and the resulting splitting of it, see Note~{\bf i} p.~\pageref{NoteP11}, and the involved pair deformation
retractions as in the proof of Th.~\ref{TheoremP12:2}
we get

\medskip
\noindent $
{\mhatH\!\!\!_{_{\!q+1\!}}
        (\!X\!_{}
\ \!\mhatast \ \!
Y\!_{}, X\!_{}
\ \!\mhatast \ \!
Y\!_{} \ \!\setminus\!{_{_{^{\!o}}}}\
\!\!(t\!{_{_{^{1}}}},\widetilde{x\!{_{_{^{0}}}}\! \ast
y\!{_{_{^{0}}}}\!}, t\!{_{_{^{2}}}} ) } ) \ \!
\mringHom{A}
\ \mhatH_{_{\!q\!}}
        (X \!\times\!Y \!\times\!(0,1),
X \!\times\!Y \!\times\!(0,1) \ \!\setminus{_{_{^{\!o}}}}\ \!
\{x\!{_{_{^{0}}}}\!\}\!\times\!
\{y\!{_{_{^{0}}}}\!\}\!\times\!(0,1)\ ) \ \!
\mringHom{A}
$
$$
\mringHom{A}
\ \Big[ {{
^{\textrm{\scriptsize Motivation: The underlying pair on the
r.h.s.}}_{
{{\textrm{\scriptsize is a pair deformation retract of that
on the l.h.s.}}}}}}
\Big]\
\mringHom{A}
$$
$$
\mringHom{A}
\ \mhatH_{_{\!q\!}}
        (X \!\times\!Y \!\times\!
\{t\!{_{_{^{0}}}}\!,\wp\}, X \!\times\!Y
\!\times\!\{t\!{_{_{^{0}}}}\!,\wp\} \ \!\setminus\!{_{_{^{o}}}}\
\!\! (x\!{_{_{^{0}}}}\!,\! y\!{_{_{^{0}}}}\!,\!t\!{_{_{^{0}}}})\ )\
\mringHom{A}
\ {\mhatH_{_{\!q\!}}
        (X \!\times\!Y \!,
X \!\times\!Y \!\setminus\!{_{_{^{\!o}}}}\! (x\!{_{_{^{0}}}}\!,\!
y\!{_{_{^{0}}}}\!)\ {\!}\!) }=
$$
$$
=
\mhatH_{_{\!q\!}}
        (X \!\times\!Y \!,
(X \!\times\!(Y \!\setminus\!\!{_{_{^{o}}}} y\!{_{_{^{0}}}}\!)) \cup
((X \!\setminus\!\!{_{_{^{o}}}}x\!{_{_{^{0}}}}\!) \!\times\!Y )\! \
) \ \!=
{ \mhatH_{_{\!q\!}}(
        (X \!,
X \!\setminus\!\!{_{_{^{o}}}}x\!{_{_{^{0}}}}\!) \!\times\! (\!Y \!,Y
\!\setminus\!\!{_{_{^{o}}}} y\!{_{_{^{0}}}}\!) \ ). }
\eqno{\hbox{$\triangleright$}}
$$
$$
\leqno{\bf ii.}
\ \
\left\{\begin{array}{ll}
&
\!\!\!\!\!\!\! A:=X {\myPsqcup1}Y \hbox{\tenbf,}
\\
&
\!\!\!\!\!\!\!B:=X {\myPsqcup0}Y \!\setminus\!{_{_{^{\!o}}}}\ \!\!
X \!\times\!\{y_0\}\times[0,1)
\\
\end{array}\right.
\Longrightarrow
\left\{\begin{array}{ll}
&
\!\!\!\!\!\!\! A\cup B =\!X\!_{} \ \!\mhatast \ \! Y\!_{}\
\!\setminus\!{_{_{^{\!o}}}}\ \!(y\!{_{_{^{0}}}},0)
\\
&
\!\!\!\!\!\!\! A\cap B= X \!\times\! (Y \!\setminus\!{_{_{^{\!o}}}}\
\!\! y\!{_{_{^{0}}}})\!\times\!(0,1)
\\
\end{array}\right.
$$
where\\ $ (x\!{_{_{^{0}}}}, y\!{_{_{^{0}}}},t)\in X
\!\times\!\{y\!{_{_{^{0}}}}\!\}\!\times\![0,1) $
is independent of $x\!{_{_{^{0}}}}$ and $
(x\!{_{_{^{0}}}},y\!{_{_{^{0}}}},t\!{_{_{^{0}}}}\!)\!:=
\{(x\!{_{_{^{0}}}},y\!{_{_{^{0}}}},t\!{_{_{^{0}}}}\!), \wp\}. $

\smallskip
Now use Th.~\ref{TheoremP12:1}
p.~\pageref{TheoremP12:1}
line\ 2 and that the r.h.s. is a pair deformation retract of the
l.h.s.;
$$
(X \!\times\!Y \!\times\!(0,1) , X \!\times\! (Y \
\!\setminus{_{_{^{\!o}}}}\ \! y\!{_{_{^{0}}}}\!)\!\times\!(0,1)) \
\lower4pt\hbox{$\widetilde{\phantom{{..}}}$} \
(X \!\times\!Y  , X _{\!}\times_{\!} (Y \ \!\setminus{_{_{^{\!o}}}}\
\! y\!{_{_{^{0}}}})) =
(X \!,\emptyset) \!\times\! (Y \!,Y \!\ \!\setminus\!{_{_{^{o}}}}\
\!y\!{_{_{^{0}}}}\!).
$$\end{proof}

\begin{vardef}
Put
${\rm Hip}\!\!{{{\lower3.5pt\hbox{\fivebf G}}}}\!X \!=\!
\{\hbox{\rm{The homologically$_{_{\!\hbox{\fivebf G}}}$ instable points in $X$}}\}
:= \{x\!\in\!X\mid
\mhatH_{i}(X,X\setminus_{\!_{^o}}x;\hbox{\tenbf
G})=0\ \text{for all}\ i\ \!\in {\mathbb{Z}}\}.$
\end{vardef}
So, for products and joins we get
$$
{\rm Hip}\!{{{\lower3.5pt\hbox{\fivebf
G}}}}{\!}(X_{_{^{_{\!}1}}}\!\
\mytinynabla
\!X_{_{^{_{\!}2\!}}})
\supset
(X_{_{^{_{\!}1}}}\!\mytinynabla\!{\rm Hip}\!{{{\lower3.5pt\hbox{\fivebf
G}}}}{\!}X_{_{^{_{\!}2\!}}})\cup (({\rm
Hip}\!{{{\lower3.5pt\hbox{\fivebf
G}}}}{\!}X_{_{^{_{\!}1}}})\mytinynabla\! X_{_{^{_{\!}2\!}}}).$$

\subsection{Homology manifolds under product and
join}$\!$\label{SubSecP16:3.4}

\begin{vardef} \label{DefP16}
$\emptyset$ is said to be a weak {homology$_{_{\mathbf{G}}}\!$-$(-\infty)$-manifold}.
Any nonempty
\htmladdnormallink{${\bf T}_1$-{\it space}
}{http://planetmath.org/encyclopedia/PerfectlyNormal.html}
($\Leftrightarrow\!$ all points are closed)
$X\!_{\!}\in\!{{\mathcal D}}_{\!\!\wp}$ is a  weak {\it
homology}$_{_{\mathbf{G}}}\!$ $n$-manifold \hbox{\rm (}{\it
n}-whm$_{_{^{\!\hbox{\fivebf G}}}\!}$\hbox{\rm)} if, for some {\bf
A}-module\ ${\mathcal R}$;

\medskip
\noindent\phantom{I}\hskip-0.3cm
{
\vbox{\hsize=0.9 \hsize
{
{{
\noindent
$
\mhatH_{i}(X,X\setminus\!_{o}{_{\!}}x;{\bf G})
\ \!=\ \!
{0\ \ {\rm if}\ \hbox{$i\!\ne\! n$} \ {\rm for\ all}\
\hbox{$\wp\!\ne\! x\!\in\! X$} , \hskip2.7cm {\bf (4.i)} }
$

\medskip\noindent
${
\mhatH_{_{^{n}}}(X,X\setminus\!_{o}{_{\!}}x;{\bf G})\
\!\cong
{\bf G}\oplus {\mathcal R}\ {\rm for\ some}\ \hbox{$ \wp\!\ne\!
x\!\in\! X$}\ {\rm if}\ X\ne\{\wp\}.\
\hskip1.0cm \hbox{\bf (4.ii$^\prime$)}
}\!$
}}}
}
\hskip0.7cm {\raise0.35cm\hbox{({\tenbf 4$^\prime$})}}

\medskip\noindent
An {\it n}-whm$_{_{^{\!\hbox{\fivebf G}}}\!}$ $X$ is {\it
joinable}
({\it n}-jwhm$_{_{^{\!\hbox{\fivebf G}}}\!}$)
if {\bf (4.i)} holds also for $x=\wp$.
}

\noindent
An {\it n}-jwhm$_{_{^{\!\hbox{\fivebf G}}}\!}$ $X$
is a {\it weak} {\it homology} {\it n}-{\it
sphere}$_{_{^{\!\hbox{\fivebf G} }}\!}$ $\!(${\it n}-whsp$)$
if
$\mhatH\!\!\!\!_{_{^{n-1}}}\!(X\!\setminus_{o}{_{\!}}x;{\bf
G})\!=\!0\ \forall\ _{\!}x\!_{\!}\in\!\! X.$
\end{vardef}

\begin{vardef}\label{DefP16:Technical}
\  $X$ is acyclic$\!_{{_{\bf G}}}\!\!$ if
$\mhatH_{i}(X,\emptyset;{\bf G})=0$ for all $i \in {\mathbb{Z}}$.
So, $\{\wp\}\ (=|\{\emptyset_o{_{\!}}\}|)$ isn't acyclic$\!_{{_{\bf
G}}}$.

An arbitrary $X$ is {\it ordinary}$\!_{_{\bf G}}\!$ if (${\bf 4}^{\prime}$) implies that $\mhatH_{i}(X\!\setminus_{o}x;{\bf G})=0,$
 $\ \forall\ \!i\ge n\,$ {and} $\forall x\in X.$
 \label{DefP16:ordinary}

$X$ \hbox{\rm is} {\it locally} {\it weakly} {\it
direct}$_{\mathbf{G}}$ if
$\mhatH_{i}(X,X\setminus\!_{o}x;{\bf G})
\cong
{\bf G}\oplus {\mathcal Q}$
{for some {\it i}}, {some {\bf A}-module}\ ${\mathcal Q}$\
and some $\wp\!\ne\! x\!\in\! X$. ({This definition was introduced to avoid troublesome tensor annihilations in case of general $\mathbf{G}$.})
\end{vardef}

\begin{varnote}{1.} \label{NoteP17:1}
{\rm Lk}$_{_{\!\Sigma}}\sigma=\{\emptyset\!_{_{^{o}}}\}$ for any maxidimensional simplex $\sigma$ in any any simplicial complex $\Sigma$. So, for such a $\sigma$ and if $\alpha\!\in\! \hbox{\rm Int}\sigma$ we get,
$\mhatH \hskip-0.3cm\raise-0.55pt\hbox{${{_{_{_{\dim\!{\Sigma
}}}}}\!}$}
(|\Sigma|,|\Sigma|\setminus\!_{o}\!\alpha;{\bf
G})\cong{\bf G}$ by Proposition~\ref{PropP15:1}
p.~\pageref{PropP15:1}
and
Lemma p.~\pageref{LemmaP10}.
Quasi-manifolds are ordinary by Note~1
p.~\pageref{NoteP29:1}.
\end{varnote}

\begin{theorem} \label{CorP16:ToTh6}
For compact triangulable spaces $X_{1},X_{2}$
$(X_{i}\neq \emptyset, \{\wp\});$

\medskip
\noindent {\bf i.}\label{CorP16:iiToTh6} \
$
X_{1}\times X_{2} $
$(n_{1}+n_{2})\text{\rm -whm}_{\mathbf{k}}$
$\Longleftrightarrow
X_{1},X_{2} $ both \text{\rm whm}$_{\mathbf{k}}.
$

\smallskip
\noindent
{\bf ii.}\label{CorP16:iiiToTh6}
If $n_{_{^1}}\!\! +n_{_{^2}}\! > n _{_{^i}}\ i=1,2, $
then
 $$[X_{1}\!\times X_{2}
(n_{1}\!+n_{2}\!)\hbox{\rm -jwhm}_{\mathbf{k}}]\!\Longleftrightarrow [X_{1},X_{2} \text{both}\ n_{_{^i}}\text{\rm -jwhm}_{\mathbf{k}}\text{and acyclic}_{\mathbf{k}}].
$$
{\rm (Since by Eq.\ {\bf 2} p.~\pageref{EqP12:2};}
{$[\mhatH_{_i}(X_{_{1}}\!\!\times\!X_{_{2}};{\mathbf{k}})=0 \
\hbox{\rm for}\ i\not=n_{1}+n_{2}]
\Longleftrightarrow [X_{_{1}}, X_{_{2}}\ \!\hbox{\rm both}$\ ${\it
acyclic}_{\mathbf{k}}]
\!\Longleftrightarrow\![X_{_{1}}\!\times\! X_{_{2}}\ \hbox{\rm
acyclic}_{\mathbf{k}}].$
{\rm So $X_{_{1}}\!\times\! X_{_{2}}$ is never
a \text{\rm whsp}$_{\mathbf{k}})$.}
}

\smallskip
\noindent {\bf iii.}\label{CorP16:iToTh6} \ \
$
X_{1}\!\!\ast\!X_{2}\
(n_{1}\!+n_{2}\!+1)\text{\rm-whm}_{\mathbf{k}}
$
$ \Longleftrightarrow
X_{1},X_{2}\! $ {both} n$_{_{{\!i}}}\!$\text{\rm-jwhm}$_{\mathbf{k}}
$
$
\Longleftrightarrow X_{1}\!\ast\!X_{2}\!$
\text{\rm jwhm}$_{\mathbf{k}}.
$

\smallskip
\noindent {\bf iv.}\label{CorP16:ivToTh6}
$X_{1},X_{2}$\ are {both} \text{\rm  whsp}$_{_{\mathbf{k}}\!}$
\underbar{iff} $X_{1}\!\!\ast\! X_{2}\!$ is a \text{\rm  whsp}.
\end{theorem}

\begin{proof}
\noindent [Originating from \cite{10}.]
({\bf i-iii}) Use Th.~\ref{TheoremP12:1},
\ref{TheoremP14:4}-\ref{TheoremP15:6}
and the weak directness$_{\mathbf{k}}$ to
transpose non-zeros from one side to the other, using
Th.~\ref{TheoremP15:6}{\bf.ii} only for joins, i.e., in particular, with $\varepsilon =
0\ \hbox{\rm or}\ 1$ depending on whether
$\mytinynabla =\times\ \hbox{\rm or}\ \ast$ respectively, use

\smallskip
\noindent
{\bf(}$\mhatH
\vbox{\moveleft0.0cm\hbox{\lower0.0pt\hbox{$_{{{_{\!}{p
\lower0.0pt\hbox{\fivebf+} \varepsilon}}}}$}}}
(X\!_{_{^{1}}}\!\mytinynabla\!X\!_{_{^{2}}},
X\!_{_{^{1}}}\!\mytinynabla\!X\!_{_{^{2}}} \!\setminus\!\!_{o}
\widetilde{(x\!_{_{^{1}}},x\!_{_{^{2}}}\!)}
;\mathbf{k})
\mringHom{k}${\bf)}
$
\mhatH _{_{^{p\!\!}}}
((X\!_{_{^{1}}},X\!_{_{^{1}}}\!\setminus\!{_{\!}}_{o}x\!_{_{^{1}}})\times
(X\!_{_{^{2}}},X\!_{_{^{2}}}\!\setminus\!{_{\!}}_{o}x\!_{_{^{2}}});\mathbf{k})
\mringHom{k}
$
$
[\mbox{Eq.}~1\ p.~\pageref{EqP12:1}]
\mringHom{k}
$
\vskip-0.23cm
\vskip-0.23cm
$$
\hskip2cm\mringHom{k}
{\rlap{$\!_{_{_{i+j=p\!\!}}}$}{\raise2pt\hbox{\ $\bigoplus$}}}\ \
[\mhatH_{i}
({X\!_{_{^{1}}},X\!_{_{^{1}}}\!\setminus\!{_{\!}}_{o}x\!_{_{^{1}}}};{\mathbf{k}})
\!\otimes\!_{\mathbf{k}}\!
        \mhatH_{j}
({X\!_{_{^{2}}},X\!_{_{^{2}}}\!\setminus\!{_{\!}}_{o}x\!_{_{^{2}}}\!};\mathbf{k})].
\hskip3cm\triangleright
$$

\medskip\noindent
{\bf iv.} Use, by the Five Lemma, the chain equivalence
of the second component in the first and the last item of
Theorem~\ref{TheoremP15:6}{\bf.i} and the {\bf M-$_{\!}$Vs} with respect to
${\!}\{(X\ast (Y\setminus\!{_{_{^{o\!}}}} y\!_{_{0}}\!)),((X\setminus\!{_{_{^{o\!}}}} x\!_{_{0}}\!)\ast Y)\}$.
\end{proof}

\begin{varnote}{2.} With $\varepsilon =
0\ \hbox{\rm or}\ 1$ depending on whether
$\mytinynabla =\times\ \hbox{\rm or}\ \ast$ respectively, Th.~\ref{TheoremP15:6} will handle any case encountered in practice, since e.g., for any module $\mathbf{G}$:

\smallskip
\noindent
$
\hbox{\tenbf(}\mhatH
\vbox{\moveleft0.0cm\hbox{\lower0.0pt\hbox{$_{{{_{\!}{p
\lower0.0pt\hbox{\fivebf+} \varepsilon}}}}$}}}
(X\!_{_{^{1}}}\!\mytinynabla\!X\!_{_{^{2}}},
X\!_{_{^{1}}}\!\mytinynabla\!X\!_{_{^{2}}} \!\setminus\!\!_{o}
\widetilde{(x\!_{_{^{1}}},x\!_{_{^{2}}}\!)}
;\mathbf{G})
\mringHombb{Z})$
$
\mhatH _{_{^{p\!\!}}}
((X\!_{_{^{1}}},X\!_{_{^{1}}}\!\setminus\!{_{\!}}_{o}x\!_{_{^{1}}})\times
(X\!_{_{^{2}}},X\!_{_{^{2}}}\!\setminus\!{_{\!}}_{o}x\!_{_{^{2}}});\mathbf{G})
\mringHombb{Z}\ \!
$\\
\centerline{$\mringHombb{Z}
\Big[\lower3.5pt\vbox{{\hbox{\scriptsize Lemma\ p.~\pageref{LemmaP10}}\vskip-0.1cm{\hbox{\scriptsize + Eq.~1\ p.~\pageref{EqP12:1}}}}}\Big]\ \!
\mringHombb{Z}\ \!
$}
\vskip-0.23cm
\vskip-0.23cm
$$
\mringHombb{Z}
{\rlap{$\!_{_{_{i+j=p\!\!}}}$}{\raise2pt\hbox{\ $\bigoplus$}}}\ \
[\mhatH_{i}
({X\!_{_{^{1}}},X\!_{_{^{1}}}\!\setminus\!{_{\!}}_{o}x\!_{_{^{1}}}};{\mathbb{Z}})
\!\otimes\!_{_{\mathbb{Z}}}\!
        \mhatH_{j}
({X\!_{_{^{2}}},X\!_{_{^{2}}}\!\setminus\!{_{\!}}_{o}x\!_{_{^{2}}}\!};\mathbf{G})] \ \ \oplus\hskip6cm
$$
$$ \hskip3cm
\oplus\ \ {\rlap{$\!\!\!\!_{_{_{i+j=p-\!1\!\!}}} $}
{\raise2pt\hbox{$\bigoplus$}}}\ \
\hbox{\rm Tor}_1^{^{_{\mathbb{Z}}}}\bigl(
\mhatH_{i}
({X\!_{_{^{1}}},X\!_{_{^{1}}}\!\setminus\!{_{\!}}_{o}x\!_{_{^{1}}}}\!;{\mathbb{Z}}),
 \mhatH_{j}
({X\!_{_{^{2}}},X\!_{_{^{2}}}\!\setminus\!{_{\!}}_{o}x\!_{_{^{2}}}}\!;
\mathbf{G})\big).
$$
\end{varnote}
\begin{varnote}{3.}
If all the involved homology groups
are of finite type with $\mathbb{Z}$ as coefficient module, then the Universal Coefficient Theorem implies that the above field-versions of Th.~\ref{TheoremP17:7.1} and Th.~\ref{CorP16:ToTh6} are true also with $\mathbb{Z}$ as coefficient module, cp. \cite{Ho} Proposition~2.4 p.~181-182.
\end{varnote}

\begin{vardef} \label{DefP17}
$\emptyset$ is said to be a {homology$_{_{\mathbf{G}}}\!$-$(-\infty)$-manifold} and
$X\!=\!{\bullet}{\bullet}$ is a {\it homology}$_{_{\mathbf{G}}}\!$
$0$-{\it manifold}.
Any other (than $\emptyset$ or ${\bullet}{\bullet}$) connected, locally compact Hausdorff space
$X\!_{\!}\in\!{{\mathcal D}}_{\!\!\wp}$ is a ({\it singular})
homology$\!_{_{\mathbf{G}}}\!$
$n$-{\it manifold} \hbox{\rm (}$n$-hm$_{_{^{\!\hbox{\fivebf G}
}}\!}$\hbox{\tenbf)} if

\smallskip
\noindent\phantom{I}\hskip-0.3cm
{
\vbox{\hsize=0.92 \hsize
{
{{
\noindent
$\mhatH_{i}(X,X\setminus{_{\!}}_{o}x;{\bf G})
\ \!=\ \!$
0\ \
if\ $i\!\ne\! n\ \forall$
$\wp\!\ne\! x\!\in\! X,
\hskip4.7cm\hbox{({\bf 4.i})}
$
\vskip0.2cm \noindent
${
\mhatH_{_{^{n}}}(X,X\setminus\!{_{\!}}_{o}x;{\bf G})
\!\cong\!
0\ {or}\ {\bf G}\ \forall \ \wp\!\ne\! x\!\in\! X \hbox{\rm and
$={\bf G}$ for some}\ x\!\in\!X.
\hskip1.0cm\hbox{({\bf 4.ii})}
}$
}}}}}
\ {\raise0.45cm\hbox{\hsize1cm\hskip0.3cm({\bf 4})}}

\medskip
The {\it boundary} is defined to be: \label{DefP17:BoundaryOfHomologyManif}
\ \ $\!{\rm Bd}_{_{\mathbf{G}}}\!X\!\!:=\! \{x\!\in\!X\mid
\mhatH_{n}(X,X\setminus\!{_{\!}}_{o}x;{\bf G})=0\}.
$
{If} {\rm Bd}$_{_{\mathbf{G}}}X\!\!\ne\!\emptyset$,
$X$ is said to be a {\it homology}$_{_{\mathbf{G}}}\!$ $n$-{\it manifold with boundary}.

A compact $n$-manifold ${\mathcal S}$ is {\it
orientable}$_{_{{^{^{\!\!\hbox{\fivebf G}}}}}}\!$ {if}
$\mhatH\!_{_{{n}}}\!({\mathcal S},{\rm Bd}{\mathcal S};{\bf G})\
\widetilde{\hbox{\rm=}}\ {\bf G}.$
An $n$-manifold is {\it orientable}$_{_{{^{^{\!\!\hbox{\fivebf G}}}}}}\!\!$ if all its compact $n$-submanifolds are orientable;  otherwise it is {\it non}-{\it orientable}$_{_{{^{^{\!\hbox{\fivebf G}}}}}}\!._{\!}$
Orientability of $\emptyset$ is not defined.

An $n$-hm$\!_{_{^{\!\hbox{\fivebf G} }}\!\!}$ $X$ is {\it joinable}
if $(4)$ holds also for $x\!=\wp.$

An {\it n}-hm$_{_{^{\!\hbox{\fivebf G} }}\!}$
$X{\!}\!\not={\!}\emptyset$
is a {\it homology$_{_{\mathbf{G}}}{\!}{\!}$ n-sphere}
$(n$-hsp$_{_{^{\!\hbox{\fivebf G} }}\!})$ if,
for all $x{\!}\in{\!} X,$
$\mhatH_{i}(X,X{\!}\setminus\!_{o}{{\!}}x;
{\bf G}){\!}={\!}{\bf G}$ if $i=n$ and $0$ otherwise.
So, a triangulable $n$-hsp$_{_{^{\!\hbox{\fivebf G} }}\!}$ is a compact space.
\end{vardef}

When $\mytinynabla$ in Theorem~\ref{TheoremP17:7}, is interpreted throughout as $\times$, the symbol ``hm$_{\mathbf{k}}\!$'' on the r.h.s. of Th.~\ref{TheoremP17:7}.1 temporarily excludes
$\emptyset,\{\emptyset_o\}$ and $\bullet\bullet$, and we assume
$\epsilon:=0$.
When $\mytinynabla$ is interpreted throughout as $\ast$, put $\epsilon:=1$, and let the
symbol ``hm$_{\mathbf{k}}\!$'' on the right hand side of
Th.~\ref{TheoremP17:7}.1 be limited to ``any \underbar{compact
joinable} homology$_{\mathbf{k}}$
$n_{i}$-manifold''.

\begin{theorem}\label{TheoremP17:7}
For compact triangulable spaces $X_{1}, X_{2}$ and any field {\bf k};

\smallskip\noindent
{\bf 1.}\label{TheoremP17:7.1}
$X_{1}\diamond X_{2}$ is a {\it homology}$_{_{\mathbf{k}}}$
$(n_{1}+n_{2}+\epsilon)$
$\hbox{-manifold}
\Longleftrightarrow X_{i}$
is a $n_{_{i}}$-hm$_{_{\mathbf{k}}},\ \text{for}\ i=1,2.$

\smallskip\noindent
\noindent {\bf 2.}\label{TheoremP17:7.2} \
${\rm Bd}_{_{\mathbf{k}}} (\bullet\times X)= \bullet\times({\rm Bd}_{_{\mathbf{k}}} X).\  \text{Else}\ {\rm Bd}_{_{\mathbf{k}}} (X_{1}\diamond X_{2})= (({\rm Bd}_{_{\mathbf{k}}}
X_{1})\diamond X_{2})\cup
(X_{1}\diamond ({\rm Bd}_{_{\mathbf{k}}} X_{2})).$

\smallskip\noindent
\noindent {\bf 3.}\label{TheoremP17:7.3} \
$X_{1}\diamond X_{2}\text{ is\ orientable}_{\mathbf{k}}
\Longleftrightarrow\text{ both } X_{1},X_{2}$ \text{are orientable}$_{\mathbf{k}}.$
\end{theorem}

\begin{proof}
Th.\ \ref{TheoremP17:7} is trivially true for
$\!X_{i}\!\times_{\!}\bullet$ and
$X_{i}\!_{\!}\ast_{\!}\{\wp\}.$
Otherwise, exactly as in the proof for Theorem~\ref{CorP16:ToTh6} above, adding for {\rm
\ref{TheoremP17:7}.1}\ that for Hausdorff-like spaces $($:= all
compact subsets are locally compact$)$, in particular for Hausdorff
spaces, $X_{1}\!\ast X_{2}\!$ is locally compact
(Hausdorff) if and only if both
$X_{1}\!,X_{2}\!$ are compact (Hausdorff), cf. \cite{4}
p.~224.
\end{proof}

\begin{varnote}{4.} \label{NoteP17:2}
\noindent$\!
(X_{_{^{_{\!}1}}}\!\mytinynabla\!X_{_{^{_{\!}2\!}}}, {\rm Bd}\!{{{\lower3.5pt\hbox{\fivebf
G}}}}(X_{_{^{_{\!}1}}}\!\mytinynabla\!X_{_{^{_{\!}2\!}}}))
=\big[{{\lower3pt\hbox{\scriptsize \ref{TheoremP17:7}.2\ above
+}}\atop{\raise3pt\hbox{\scriptsize pair-def.\
p.~\pageref{DefP11:PairDef}}}}\big]
$
$=
(X_{_{^{_{\!}1\!}}},{\rm Bd}\!{{{\lower3.5pt\hbox{\fivebf
G}}}}{\!}_{\!}X_{_{^{_{\!}1\!}}})
\mytinynabla\!
(X_{_{^{_{\!}2\!}}},{\rm Bd}\!{{{\lower3.5pt\hbox{\fivebf
G}}}}{\!}_{\!}X_{_{^{_{\!}2\!}}}).
$
\end{varnote}

\noindent{\bf Test case.} Using join and quotient to calculate the $\mhatH$omology groups of the $\ominus$-figure.
($\mhatH_{i}\!(\ominus;{\bf G})\!=\!0$ if $i\!\ne\! 1$.)
Prop.~{\bf iii} p.~\pageref{PropP88iii} gives the second equality, the join-definition the two next, while Eq.~3 p.~\pageref{EqP14:3} and Lemma p.~\pageref{LemmaP10} take care of the rest in the following calculation:

\normalbaselines \noindent $\hskip-0.1cm
\mhatH_{_{1}}\!(\ominus;{\bf G})
= $ $ \mhatH_{_{1}}\!( {{| {{ \bullet\hskip-0.09cm
\raise1.4pt\hbox{$^{{\nearrow^{\hskip-0.15cm\lower1.45pt\hbox{$\bullet$}}}}$}
\hskip-0.35cm{{\mmyline}}\raise1.7pt\hbox{{$_{{\!\!_{\!}\hbox{$\bullet$}}}$}}
\hskip-0.48cm
\lower2.1pt\hbox{${{_\searrow}}$}\!_{_{_{_{\hskip-0.1cm\hbox{$\bullet$}}}}}
}} |}} \lower2pt\hbox{/} _{^{^{|\raise1.5pt\hbox{$_{\ \!
\!{^{^{\hbox{\bf .}}}\hskip-0.13cm{\raise0.65pt\hbox{\bf .}}
\hskip-0.1cm_{_{\hbox{\bf .}}}}}$}|}}} ;{\bf G}) =
\mhatH_{_{1}}\!( {{ \bullet\hskip-0.09cm
\raise1.4pt\hbox{$^{{\nearrow^{\hskip-0.15cm\lower1.45pt\hbox{$\bullet$}}}}$}
\hskip-0.35cm{{\mmyline}}\raise1.7pt\hbox{{$_{{\!\!_{\!}\hbox{$\bullet$}}}$}}
\hskip-0.48cm
\lower2.1pt\hbox{${{_\searrow}}$}\!_{_{_{_{\hskip-0.12cm\hbox{$\bullet$}}}}}
}}\ ,\ \bullet\!\bullet\!\bullet\ ;{\bf G}) $ $ =
\mhatH_{_{1}}\!((\bullet)\ \!{\ast}\ (\bullet\!\bullet\!\bullet)\
,\ \bullet\!\bullet\!\bullet\ ;{\bf G}) = [{\rm Def.}\ p.~\pageref{DefP11:PairDef}]\hskip-0.15cm$\\
$
= \mhatH_{_{1}}\!((\bullet,\{\emptyset_o\})\
\!{\ast}_{_{\!^\cup}}\ \!({\bullet\!\bullet\!\bullet},\emptyset) ;{\bf
G})= $ [Eq.~3 p.~\pageref{EqP14:3}]
$ =\! {\rlap{$_{_{_{i+j=0}}}$} {\ \raise2pt\hbox{$\bigoplus$}}}\
\! \mhatH_{_{i}}(\bullet,\{\emptyset_o\};{\bf R})\otimes
\mhatH_j(\bullet\!\bullet\!\bullet,\emptyset;{\bf G}) =$

\noindent $ = \mhatH_{_{0}}(\bullet,\{\emptyset_o\};{\bf
R})\otimes
\mhatH_{_{\!0}}\!(\bullet\!\bullet\!\bullet,\emptyset;{\bf G}) = $
[Lemma p.~\pageref{LemmaP10}] $ ={\bf R}\otimes({\bf G}\oplus{\bf G}) $ $ ={\bf
G}\oplus{\bf G}. $

\part{Merging Combinatorics, Topology \& Commutative Algebra} \label{PartIII}%

\section{Relating General Topology to Combinatorics} \label{SecP18:I}

\subsection{Realizations for simplicial products and joins} \label{SubSecP18:I:1}

\ \ \ The $k$-{\it ifikation} $k({\it X})$ or ${\it X}_k$ of ${\it X}$ is ${\it X}$
with its topology enlarged to the weak topology with respect to all continuous maps with compact Hausdorff domain.
Generally, this topology differs from the
\htmladdnormallink{{\it compactly}}
{http://en.wikipedia.org/wiki/Compactly_generated_space}
\htmladdnormallink{{\it generated}}
{http://en.wikipedia.org/wiki/Compactly_generated_space}
enlargement
${\it X}_{C}$, but they coincide for Hausdorff spaces.
For any two simplicial complexes $\Sigma_1$ and $\Sigma_2$, $|\Sigma_1|{\times}|\Sigma_2|$,
$|\Sigma_1|{\ast}|\Sigma_2|$ and $|\Sigma_1|{\mhatast}|\Sigma_2|$ are all Hausdorff.
Put
$${X}\bar{\times}\ \!Y\!:=k({X}{{{\times}}}\ \!Y).$$
For CW-complexes, this is a proper topology-enlargement only if
none of the two underlying complexes is locally finite and at
least one is uncountable.
Let ${X}\bar{\ast}Y$ be the quotient space with respect to
$p:(X\bar\times Y)\times{\bf I}
\rightarrow X\circ Y$
from p.~\pageref{SubSecP11:3.1}.
Related subtleties are examined in \cite{Fritsch&Golasinski} and \cite{20}.
Simplicial ${\times}$ and ${\ast}$ both ``commute'' with realization by turning into $\bar{\times},\bar{\ast}$, respectively.
Unlike $\ast$, Def. p.~\pageref{SubSecP11:3.1},
$\bar{\ast}$ is associative for arbitrary topological spaces, cf.
\cite{33}~\S3.\label{P18:nonassociativity}

\begin{definition} \label{DefP18:OrderedSimpCartProd}
{\rm(cf. \cite{9}
Def.\ 8.8 p. 67.)}
Given ordered simplicial complexes $\Delta{\!^{^{_{\prime}}}}\!$ and
$\Delta{\!^{^{_{\prime\prime}}}}\!$ i.e. the vertex sets ${\bf
V}\!\!_{{^{\Delta^{^{_{\!{{\!}}\prime}}}}}}\!\!$ and ${\bf
V}\!\!_{{^{\Delta^{^{_{\!\!\prime{_{\!}}\prime}}}}}}\!\!$ are
partially ordered so that each simplex becomes linearely ordered
resp.
{\it The Ordered Simplicial Cartesian Product}
$\Delta{\!^{^{_{\prime}}}}\!
{\raise1.0pt\hbox{\ninesy{\char"02}}}
\Delta{\!^{^{_{\prime\prime}}}}\!\!$ of
${\Delta{\!^{^{_{\prime}}}}}\!$ and
$\Delta{\!^{^{_{\prime\prime}}}}\!\!$ $($\hbox{\rm triangulates}
$|_{\!}\Delta{\!^{^{_{\prime}}}}{\!}|
\bar{\raise0.5pt\hbox{\ninesy{\char"02}}}|_{\!}\Delta{\!^{^{_{\prime\prime}}}}\!|$
\hbox{\rm and}$)$ is defined as
$$\!{\bf V}_{{^{{\Delta^{^{_{\!\prime}}}}
\!\times \Delta^{^{_{\!\!\prime{_{\!}}\prime}}} }}}:=
\{_{\!}(v\!\raise1pt\hbox{$_{i}$}\!^{_{\!}\prime},
v\!\raise1pt\hbox{$_{_{^{j}}}$}\!^{_{\!}\prime\prime})_{\!}\}
\!\!=\!\!{\bf
V}\!\!_{{^{\Delta^{^{_{\!\prime}}}}}}\!\!{\raise1.0pt\hbox{\ninesy{\char"02}}}
{\bf V}\!\!_{{^{\Delta^{^{_{\!\!\prime{_{\!}}\prime}}}}}}\! .$$
Put
$w\!_{_{^{i,j}}}\!_{\!}{_{\!}}:=\!(v\!\raise1pt\hbox{$_{i}$}\!^{_{\!}\prime},
v\!\raise1pt\hbox{$_{_{^{j}}}$}\!^{_{\!}\prime\prime}).$
Simplices in
$\Delta{\!^{^{_{\prime}}}}\!{\raise1.0pt\hbox{\ninesy{\char"02}}}\Delta{\!^{^{_{\prime\prime}}}}\!\!$
are sets $\{ w_{i_0,j_0},$ $w_{i_1,j_1},...,
w_{i_k,j_k}\},$
with
$w\!_{_{^{{i_{\!s\!}},{j_{\!s}}}}}\!\!\!
\neq\!w_{{i_{s+1}},{j_{s+1}}}$ and
$v_{i_0}^\prime \le v_{i_1}^\prime\le... \le
v_{i_k}^\prime\
(v_{j_0}^{\prime\prime}\le v_{j_1}^{\prime\prime}
\le...\le\!\!v_{j_k}^{\prime\prime})$ where\ $v_{i_0}^\prime,
v_{i_1}^\prime,..., v_{i_k}^\prime$
$( \text{resp.}\ v_{j_0}^{\prime\prime},
v_{j_1}^{\prime\prime},...,v_{j_k}^{\prime\prime})$ is a sequence of
vertices, with repetitions possible,
constituting a simplex in $\Delta^\prime$
(\text{resp.} $\Delta^{\prime\prime})$.
\end{definition}

\begin{lemma} \label{LemmaP18}
{\rm(cf. \cite{9}
p.~68.)}
If $p_{{{i}}}\!:\Sigma_{_{^{\!1}}}\!\!\times\!\Sigma_{_{^{\!2}}}
\!\!\rightarrow\! \Sigma_{{{i}}} $ is the simplicial projection, then
$$\eta:=(|p_{_{^{\!1}}}\!|,|p_{_{^{\!2}}}\!|),|\Sigma_1\!\times\!
\Sigma_2|{\rlap{$\hookrightarrow$}{\ \rightarrow}}|\Sigma_1|\bar{\times}|\Sigma_2|$$
triangulates $|\Sigma_1|\bar{\times}|\Sigma_2|$.

If $L_1$ and $L_2$ are subcomplexes of $\Sigma_1$ and $\Sigma_2$,
then $\eta$ carries $|L_1\times L_2|$ onto $|L_1|\bar{\times}|L_2|$.

This triangulation has the property that, for each vertex $B$ of \
$\Sigma_2$, say, the correspondence $x\rightarrow(x,B)\ {\sl is\ a}$
simplicial map of $\Sigma_1$ into $\Sigma_1\times \Sigma_2$.

Similarly for joins, with
$
 \eta \!:|\Sigma_{_{^{\!1}}}\!\ast
\Sigma_{_{^{\!2}}}\!|\
{\rlap{$\hookrightarrow$}{\ \rightarrow}}\ |\Sigma_{_{^{\!1}}}\!|\
\!\bar{\ast}\ \!|\Sigma_{_{^{\!2}}}\!|
$
defined in the proof.
\end{lemma}

\begin{proof}
{\bf($\times$)}. The simplicial projections $
p_{{{i}}}\!:\Sigma_{_{^{\!1}}}\!\times\Sigma_{_{^{\!2}}}
\!\!\rightarrow\! \Sigma_{{{i}}} $ gives realized continuous maps
$|p_{i}|,\ \text{for}\ i\!=\!1,2.$ The map $\eta:=
\!(|p_{_{^{\!1}}}\!|,|p_{_{^{\!2}}}\!|)\!:\!
|\Sigma_{_{^{\!1}}}\!\!\times\!\Sigma_{_{^{\!2}}}\!|
\!\!\rightarrow\!
|\Sigma_{_{^{\!1}}}\!|\bar\times|\Sigma_{_{^{\!2}}}\!| $ is
bijective and continuous, cf. \cite{2}
2.5.6 p.~32 \raise1pt\hbox{\eightrm+} Ex.~12,
14 p. 106-7.
The topology $\tau\!\!{_{_{^{|\!\Sigma_{_{^{\!1\!\!}}}\times\Sigma_{_{^{\!2\!\!}}}|}}}}\!\!$
$(\text{resp.}\ \tau\!\!{_{_{^{|\!\Sigma_{_{^{\!1\!\!}}}|
\!\bar\times\!|\!\Sigma_{_{^{\!2\!\!}}}|}}}}\!)\!$ is the weak
topology with respect to $\!$the compact subspaces
$\{|\Gamma\!_{_{^{\!1}}}\!\times\Gamma\!_{_{^{\!2}}}\!|\}
_{_{^{\!\Gamma\!_{{^{i}}}\!\subset\Sigma_{{^{_{\!}i}}}}}}\!\! $
(resp.\ $
\{\cong\
\!|\Gamma\!_{_{^{\!1}}}\!|\times|\Gamma\!_{_{^{\!2}}}\!|\}
_{_{^{\!\Gamma\!_{{^{i}}}\!\subset\Sigma_{{^{_{\!}i}}}}}}\!),\!\hbox{\rm\
cf.}$
\noindent
\cite{11}
p.~246 Prop.~A.2.1.
Since,
$$(|p_{_{^{\!1}}}\!|,|p_{_{^{\!2}}}\!|)(
|\Gamma\!_{_{^{\!1}}}\!\times\Gamma\!_{_{^{\!2}}}\!| \cap A) \!=\!
\nobreak
(|\Gamma\!_{_{^{\!1}}}\!|^{^{_{\!}}}\times^{^{_{\!}}}|\Gamma\!_{_{^{\!2}}}\!|)
\cap(|p_{_{^{\!1}}}\!|,|p_{_{^{\!2}}}\!|)(A),$$
the inverse $(|p_{_{^{\!1}}}\!|,|p_{_{^{\!2}}}\!|)^{-1}$ is continuous, i.e.,  $(|p_{_{^{\!1}}}\!|,|p_{_{^{\!2}}}\!|)$ is a homeomorphism.
$\hfill \triangleright$

\noindent
({$\ast$}). This follows from
\cite{34}
p.~99, using the map $\eta$ below.
With
$\Sigma_{_{^{\!1}}}\!\ast \Sigma_{_{^{\!2}}}\!:=
\{\sigma\!{_{_{^1 \!\!}}}\cup \sigma\!{_{_{^2\!\!}}}\ \vert\
\sigma\!{_{_{i \!\!}}}\in \Sigma_{_{^{\!i}}}\!,\ \text{for}\ i=1,2\}$,
and if
$
\sigma=
\{v_{_{^{\!1}}}^{^{_{_{\prime}}}}\! \
\!\hbox{\tenbf,}\!\dots\!\hbox{\tenbf,}\ \!
 v_{_{^{\!\hbox{\fiverm q}}}}^{^{_{_{\prime}}}}\! \ \hbox{\tenbf,}\
v^{^{_{_{\prime\prime}}}}\!\!\!\!\!\raise2pt\hbox{$
_{_{^{\!\raise0pt\hbox{\fiverm q{\fivebf+}1} }}}^{^{_{_{\ }}}}\! $}
\ \!\hbox{\tenbf,}\!\dots\!\hbox{\tenbf,}\ \! v^{^{_{_{\prime\prime}}}}\!\!\!\!\!\raise2pt\hbox{$
_{_{^{\!\raise0pt\hbox{\fiverm q{\fivebf+}r} }}}^{^{_{_{\ }}}}\!
$}\}
\!\in\! V_{\Sigma_1\ast\Sigma_2}$,
then put
$$\{t_{_{^{\!1}}}^{^{_{_{\ }}}}\!v_{_{^{\!1}}}^{^{_{_{\prime}}}}\! \
\!\hbox{\tenbf,}\!\dots\!\hbox{\tenbf,}
t_{_{^{\!\hbox{\fiverm q}}}}^{^{_{_{\ }}}}\! v_{_{^{\!\hbox{\fiverm q}}}}^{^{_{_{\prime}}}}\! \ \hbox{\tenbf,}\
t\raise2pt\hbox{$ _{_{^{\!\raise0pt\hbox{\fiverm q{\fivebf+}1}
}}}^{^{_{_{\ }}}}\! $} v^{^{_{_{\prime\prime}}}}\!\!\!\!\!\raise2pt\hbox{$
_{_{^{\!\raise0pt\hbox{\fiverm q{\fivebf+}1} }}}^{^{_{_{\ }}}}\! $}
\ \!\hbox{\tenbf,}\!\dots\!\hbox{\tenbf,}\ \! t\raise2pt\hbox{$
_{_{^{\!\raise0pt\hbox{\fiverm q{\fivebf+}r} }}}^{^{_{_{\ }}}}\! $}
v^{^{_{_{\prime\prime}}}}\!\!\!\!\!\raise2pt\hbox{$
_{_{^{\!\raise0pt\hbox{\fiverm q{\fivebf+}r} }}}^{^{_{_{\ }}}}\!
$}\}:=
\alpha_{\sigma}:{\bf V}_{\Sigma}\rightarrow [0,1];\
\alpha_{\sigma}\!(v)\!=
\left\{\begin{array}{l}
{t_{\hbox{\sevenrm v}}}
\ \text{if}\ v\in \sigma,\\
0\ \text{otherwise}.
\end{array}\right.
$$
Set $ t^{_{^{{\prime}}}}\!:=\sum_{_{_{_{\hskip-0.4cm{1\le i \le
q}}}}}\hskip-0.2cm t_i$\ \ \  and
\hskip0.1cm
$\!\ t^{_{^{{\prime\prime}}}}\!:=1-t^{_{^{{\prime}}}}\!=\ \
\sum_{_{_{_{\hskip-0.5cm{ q+1\le j \le q+r}}}}}\hskip-0.7cm t_j.
$\\
Now:
$$
 \eta \!:|\Sigma_{_{^{\!1}}}\!\ast
\Sigma_{_{^{\!2}}}\!|\
{\rlap{$\hookrightarrow$}{\ \rightarrow}}\ |\Sigma_{_{^{\!1}}}\!|\
\!\bar{\ast}\ \!|\Sigma_{_{^{\!2}}}\!|\ ;\
\{t_{_{^{\!1}}}^{^{_{_{\ }}}}\!v_{_{^{\!1}}}^{^{_{_{\prime}}}}\! \
\!\hbox{\tenbf,}\dots\hbox{\tenbf,}\ \! t_{_{^{\!\hbox{\fiverm
q}}}}^{^{_{_{\ }}}}\! v_{_{^{\!\hbox{\fiverm
q}}}}^{^{_{_{\prime}}}}\! \ \hbox{\tenbf,}\ t\raise2pt\hbox{$
_{_{^{\!\raise0pt\hbox{\fiverm q{\fivebf+}1} }}}^{^{_{_{\ }}}}\! $}
v^{^{_{_{\prime\prime}}}}\!\!\!\!\!\raise2pt\hbox{$
_{_{^{\!\raise0pt\hbox{\fiverm q{\fivebf+}1} }}}^{^{_{_{\ }}}}\! $}
\ \!\hbox{\tenbf,}\dots\hbox{\tenbf,}\ \! t\raise2pt\hbox{$
_{_{^{\!\raise0pt\hbox{\fiverm q{\fivebf+}r} }}}^{^{_{_{\ }}}}\! $}
v^{^{_{_{\prime\prime}}}}\!\!\!\!\!\raise2pt\hbox{$
_{_{^{\!\raise0pt\hbox{\fiverm q{\fivebf+}r} }}}^{^{_{_{\ }}}}\!
$}\} \ \mapsto
$$
$$\!\!\hskip0.0cm \mapsto \hbox{\tenbf(} t^{_{^{{\prime}}}}\!\{ {{ t\raise2pt\hbox{$
_{_{^{\raise0pt\hbox{\fiverm 1} }}}^{^{_{_{\ }}}}\! $}} \over { \
t^{{{{\prime}}}} }} v^{^{_{_{\prime}}}}\!\!\!\raise2pt\hbox{$
_{_{^{\!\raise0pt\hbox{\fiverm 1} }}}^{^{_{_{\ }}}}\! $} \
\!\hbox{\tenbf,}\!\dots\!\hbox{\tenbf,}\ \! {{ t\raise2pt\hbox{$
_{_{^{\raise0pt\hbox{\fiverm q} }}}^{^{_{_{\ }}}}\! $} } \over { \
t^{{{{\prime}}}} }} v^{^{_{_{\prime}}}}\!\!\!\raise2pt\hbox{$
_{_{^{\!\raise0pt\hbox{\fiverm q} }}}^{^{_{_{\ }}}}\! $}\} \
\hbox{\tenbf,}\ t^{_{^{{\prime\prime}}}}\!\{ {{ t\raise2pt\hbox{$
_{_{^{\!\raise0pt\hbox{\fiverm q{\fivebf+}1} }}}^{^{_{_{\ }}}}\! $}
} \over { \ t^{_{^{{\prime\prime}}}} }} v^{^{_{_{\prime\prime}}}}\!\!\!\!\!\raise2pt\hbox{$
_{_{^{\!\raise0pt\hbox{\fiverm q{\fivebf+}1} }}}^{^{_{_{\ }}}}\! $}
\ \!\hbox{\tenbf,}\!\dots\!\hbox{\tenbf,}\ \! {{ t\raise2pt\hbox{$
_{_{^{\!\raise0pt\hbox{\fiverm q{\fivebf+}r} }}}^{^{_{_{\ }}}}\! $}
} \over { \ t^{_{^{{\prime\prime}}}} }} v^{^{_{_{\prime\prime}}}}\!\!\!\!\!\raise2pt\hbox{$
_{_{^{\!\raise0pt\hbox{\fiverm q{\fivebf+}r} }}}^{^{_{_{\ }}}}\!
$}\}\hbox{\tenbf)},\ \hbox{where}\
(tx,(1-t)y)\!:=\widetilde{(x,y,t)}.
$$

(\cite{36}
(3.3) p.~59 is useful, when deling moore explicitly with the compact subspaces of $|\Sigma_{1}\!|\
\!{\ast}\ \!|\Sigma_{2}\!|$.)
\end{proof}
\normalbaselines

\cite{11} pp.\ 303-4 describes a category theoretical version of Milnor's original realization procedure for simplicial sets in \cite{24}, see \S~\ref{CatDefRealization} here in  p.~\pageref{CatDefRealization}.
{N.B.} Any complex with a nonempty vertex set in the category of augmented simplicial complexes possesses exactly two ``pointless'' subcomplexes $\emptyset$ and $\{{\emptyset}_o\}$. Any functorial realization of the latter, the join-unit, must fulfill $\mhatbar\{{\hat\emptyset}_o\}\mhatbar=|\{{\emptyset}_o\}|=\{\wp\}\in{\mathcal
D}\!_{_{^{\!\wp}}}$ defined in p.~\pageref{DefP7:CategoryDp}. Plainly, if the source is ``augmented'' so must also the target be.
The equality of the outmost spaces in \ref{product} resp. \ref{join} below also follows from Lemma~\ref{LemmaP18} above.
When we regard an ordered simplicial complex $\Sigma$ as an
(augmented) {\it simplicial set} it will be denoted by
$\mhatSigma$,
cf. p.~\pageref{UnderavdP4:RelCombToLogics}ff.
Let
$\mhatSigma_{_{^{\!1}}}\!\hat\times\mhatSigma_{_{^{\!2}}}\!$
be the (augmented) semi-simplicial product of \
$\mhatSigma_{_{^{\!1}}}\!$
and
$\mhatSigma_{_{^{\!2}}}\!$,
while
{$\mhatbar
\mhatSigma
\mhatbar$
} is the Milnor realization of $\mhatSigma$ from \cite{24}.
\cite{11}
p.~160 Prop.~4.3.15 + p.~165
Ex.\ 1\ $\!\!$+2 gives the homeomorphisms
\begin{equation} \label{product}
\vert\Sigma_{1}\vert
{{\bar\times}}
\vert\Sigma_{2}\vert
\simeq
\hat{\vert}
\hat\Sigma_{1}\hat{\vert}\bar\times\hat{\vert}\hat\Sigma_{2}
\hat{\vert}
\simeq
\hat{\vert}
\hat\Sigma_{1}\hat\times\hat\Sigma_{2}
\hat{\vert}
\simeq
\vert
\Sigma_{1}
{{\times}}\
\Sigma_{2}\vert.
\end{equation}
Similarly:
\begin{equation} \label{join}
\vert\Sigma_{1}\vert
{\bar\ast}
\vert\Sigma_{2}\vert
\simeq
\hat{\vert}
\hat\Sigma_{1}\hat{\vert}\bar\ast
\hat{\vert}\hat\Sigma_{2}\hat{\vert}
\simeq
\hat{\vert}
\hat\Sigma_{1}\hat\ast\hat\Sigma_{2}
\hat{\vert}
\simeq
\vert
\Sigma_{1}
{{\ast}}
\Sigma_{2}\vert,
\end{equation}
with the join of simplicial sets, $\hat\ast$,
as well as the motivations taken from \cite{Fritsch&Golasinski}.

The Milnor realization $\mhatbar\Xi\mhatbar$
of \underbar{any} simplicial set $\Xi$
is {triangulable} by \cite{11}
p.~209 Cor.\ 4.6.12. E.g., the augmental singular complex
$\Delta\!^{\wp}\!(X)$\label{DefP18:RealizationOfSimpSet}
with respect to any topological space $X$, is a simplicial set
and, cf.~\cite{24}
p.~362 Th.\ 4,
the map
$j:\mhatbar
\Delta\!^{\wp}\!(X)
\mhatbar\rightarrow X
$
is a weak  homotopy equivalence i.e. induces isomorphisms in
homotopy groups, and $j$ is a true homotopy equivalence if $X$ is of
CW-type,
cf. \cite{11}
pp.~76-77,\ 170,\ {189ff},\ 221-2.
Also the identity map $k(X)\hookrightarrow\hskip-0.3cm\rightarrow X$ is at least a weak  homotopy equivalence.


\subsection{Local homology for simplicial products and joins} \label{SubSecP19:I:2}
\normalbaselines

Lemma~\ref{LemmaP18} p.~\pageref{LemmaP18} and
\cite{7}
Th.~12.4 p.\ 89, also implies that
\begin{equation} \label{LocalTop}
\eta\!:\!(|\Sigma_{1}\hbox{\lower1pt\hbox{\mynabla{7}}}
\Sigma_{2}\!|,
|\Sigma_{1}\hbox{\lower1pt\hbox{\mynabla{7}}}
\Sigma_{2}\!|\setminus\!_{o}{_{\!}}
{(\!\widetilde{\alpha_{1}\!,\alpha_{2}\!})}\})\
\rlap{$\longrightarrow$}{\raise4pt\hbox{$\ \simeq$}}\ \ \
(|\Sigma_{1}\!|\ _{\!}{{\mybarnabla{7}}} |\Sigma_{2}\!|,
|\Sigma_{1}\!|\ _{\!}{{\mybarnabla{7}}} |\Sigma_{2}\!|
\setminus\!_{o}{_{\!}} (\alpha_{1}\!,\alpha_{2}\!))
\end{equation}
is a homeomorphism if
$\eta(\!\widetilde{\alpha_{1}\!,\alpha_{2}\!})\!=\!
(\!\alpha_{1}\!,\alpha_{2}\!)$.

We conclude, cf. \cite{34}
p.~99, that the following four topological spaces
\begin{equation} \label{RealJoin}
k(|{\Sigma}|\ast|\Delta|)=
|{\Sigma}|\bar\ast|\Delta|=
|\Sigma\ast\Delta|=
k(|{\Sigma}|\mhatast |\Delta|),
\end{equation}
are homeomorphic, since the underlying spaces are all Hausdorff and their respective compact subspaces are homeomorphic, i.e., they are homeomorphic compactly generated Hausdorff spaces.

\medskip
Moreover, if
${\Delta^{^{_{\!{{\!}}\prime}}}}\!\!\subset\!{\Sigma}^{^{_{\prime}}}\!$,
${\Delta^{^{_{\!{{\!}}\prime\prime}}}}\!\!\subset\!{\Sigma}^{^{_{\prime\prime}}}$
then,
${\vert{\Delta^{^{_{\!{{\!}}\prime}}}}_{\!}\vert} \bar{\ast}\
\!\!\vert\Delta{\!^{^{_{\prime\prime}}}}\!\vert $
is a subspace of
${\vert{\Sigma}^{^{_{\prime}}}\!\vert} \bar{\ast}\
\!\!\vert{\Sigma}^{^{_{\prime\prime}}}\!\vert$ and
$\ \ \dim (\Sigma\times \Delta)=\dim \Sigma+ \dim \Delta
\ \ \textrm{while}\ \
\dim(\Sigma\ast \Delta)=\dim \Sigma+ \dim \Delta +1.$

\medskip
With $\alpha_i\!\in \hbox{\rm Int}\sigma\!_{i}\!\subset\!|\Sigma_{i}|$ and
$
\left\{\begin{array}{l}
{(\widetilde{\alpha_{_{^{\!1}}},\alpha_{_{^{\!2}}}})}
:=\eta^{{_{-1}}}\!(\alpha_{1}\!,\alpha_{_{\!2\!}})
\!\in\hbox{\rm Int}\sigma
\subset|\Sigma_{_{^{\!1}}}\times\Sigma_{_{^{\!2}}}|
\label{DefP19:csigma}\\
c_{\sigma}\!:=\!\dim\sigma_{1}\!+
\dim\sigma_{2}\!-\dim\sigma,
\end{array}\right.
$\\
we conclude that,
$$ c_{\sigma}\!\ge\!0\ \text{\rm and}\ [c_{\sigma}\!\!=\!0\ \
\ \textrm{if and only if}\ \ \ \sigma\ \ \textrm{is a maximal
simplex in}\ \
\bar\sigma_{1}\times\bar\sigma_{2}\!\subset\!\Sigma_{1}\times\Sigma_{2}].\
$$

\begin{corollary}
\label{CorP19:ToTh6}
{\rm (to \S\S \ref{SubSecP12:3.2}-\ref{SubSecP15:3.3})}
Let ${\bf G}$, ${\bf G}^{\prime}$ be arbitrary modules over a {\bf PID} {\bf R} such that
$\hbox{\rm Tor}_1^{\mathbf{R}}({\bf G},{\bf G}^\prime)=0$,
then, for any
$\underline{\emptyset_o\ne\sigma}\in\Sigma_{1}\times\Sigma_{2}$ with
$\eta\big(\hbox{\rm Int}(\sigma_{})\big)\!\subset \hbox{\rm
Int}(\sigma_{1})\times\hbox{\rm Int}(\sigma_{2})$,

$$ {\underline{\underline{\mhatH_{i+c_{\!\sigma}+1}
                  ({\rm Lk}_{_{\Sigma_{1}\times\Sigma_{2}}}\sigma; {\bf
G}\otimes_{_{\bf R}}{\bf G^{\prime}})}}}
\mringHom{R}
\hskip10.5cm
$$
\vskip-0.2cm
\begin{eqnarray}
\mringHom{R}
{\rlap{$_{_{_{{{{{p+q=i}\atop{p,q\ge -1}}}}}}}$} {\ \
\raise2pt\hbox{$\ \bigoplus$}}}\
[\mhatH_{p}
           ({\rm Lk}_{\Sigma_{1}}\sigma_{1};{\bf G})
\otimes_{_{\bf R}}
        \mhatH_{q}
({\rm Lk}_{\Sigma_{2}}\sigma_{2});{\bf G^{\prime}})]
\oplus\hskip2cm
\\
\phantom{II}\hskip2cm{\rlap{$\oplus$}{\ {\rlap{$_{_{_{_{_{{{p+q=i-1}\atop{p,q\ge
-1}}}}}}}$} {\ \ \raise2pt\hbox{$\ \bigoplus$}}}}}
{\rm Tor}_1^{\bf R} \bigl(\mhatH_{p}
           ({\rm Lk}_{\Sigma_{1}}\sigma_{1};{\bf G}),\
 \mhatH_{q}
      ({\rm Lk}_{\Sigma_{2}}\sigma_{2};{\bf G^{\prime}})\bigr)
      \mringHom{R}
     \nonumber
\end{eqnarray}
$$\hskip5cm\mringHom{R}\
{\underline{\underline{\mhatH_{i+1}
        ({\rm Lk}_{_{\Sigma_{1}\ast\Sigma_{2}}}(\sigma_1\cup\sigma_2); {\bf
G}\otimes_{_{\bf R}}{\bf G^{\prime}})}}}.$$

\nobreak
\smallskip
\indent
So, if $\emptyset\!_{_{^{o}}}\!\ne\!\sigma$ and $c_{\sigma}\!=\!0$
then
$$\mhatH_{i} ({\rm Lk}\!\!\!\!\!\!
 \lower1.1pt\hbox{$_{_{\Sigma_{1}\!\times\Sigma_{2}}}$}
\!\!\!\! \sigma;{\bf G})
\ \mringHom{R}\
\mhatH_{i} ({\rm Lk}\!\!\!\!\!\!
 \lower1.1pt\hbox{$_{_{\Sigma_{1}\!\ast\Sigma_{2}}}$}\!\!(\sigma_{1}
\cup
\sigma_{2})
                ;{\bf G})$$
                and, since $\hbox{\rm Tor}_1^{\mathbf{R}}({\bf G},{\bf G}^\prime)=0$,
\begin{eqnarray}
\noindent
\mhatH_{_{0}}
 ({\rm Lk}\!\!\!\!\!\!
 \lower1.1pt\hbox{$_{_{\Sigma_{1}\!\times\Sigma_{2}}}$}
\!\!\!\! \sigma; {\bf G}\otimes{\bf G}^{\prime})
\ \mringHom{R}\hskip8cm
\\
\ \mringHom{R}\
{\mhatH_{_{0}}({\rm Lk}\!\!
 \lower1.1pt\hbox{${_{_{\Sigma_{1}}}}$}\!\!\!\sigma_{1}; {\bf G})
\otimes\mhatH\!\!_{_{-\!1}}\! ({\rm Lk}\!\!
 \lower1.1pt\hbox{${_{_{\Sigma_{2}}}}$}\!\!\!\sigma_{2};
{\bf G}^{\prime}) \oplus \mhatH\!\!_{_{-\!1}}\! ({\rm Lk}\!\!
\lower1.1pt\hbox{${_{_{\Sigma_{1}}}}$}\!\!\!\sigma_{1};
 {\bf G})}
\otimes
        \mhatH_{_{0}}({\rm Lk}\!\!
 \lower1.1pt\hbox{${_{_{\Sigma_{2}}}}$}\!\!\!\sigma_{2};
       {\bf G}^{\prime}).
       \hskip1cm
       \nonumber
       \end{eqnarray}
\end{corollary}

\begin{proof}
Note that
$\sigma\!\ne\emptyset_o\!\Rightarrow\!\sigma_{_{\!\!j}}\!\!\ne\emptyset_o,\
j\!=\!1,2.$
The isomorphisms of the underlined modules are, by
Proposition~\ref{PropP15:1}\ p.~\pageref{PropP15:1},
just a simplicial  version of
Theorem~\ref{TheoremP15:6}.{\bf i} p.~\pageref{TheoremP15:6},
and holds even without the {\bf
PID}-assumption.
Prop.~\ref{PropP34:1}\ p.~\pageref{PropP34:1},
and
Theorem~\ref{TheoremP14:4}\ p.~\pageref{TheoremP14:4},
gives the second isomorphism even for
$\sigma_{1}\!\!=\!\emptyset\!_{_{^{o}}}\ \hbox{\rm and}\ \!\!/\
\!\!\hbox{\rm or}\ \sigma_{2}\!\!=\!\emptyset\!_{_{^{o\!}}}.$
\end{proof}
The above module homomorphisms concerns only simplicial homology, so, it should be possible to prove them purely in terms of simplicial homology. This is, however, a rather cumbersome task, mainly due to the fact that $\Sigma_1\!\times\!\Sigma_2$ is not a subcomplex of $\Sigma_1\!\ast\Sigma_2$.\\ A somewhat more explicit proof is provided in p.~\pageref{LinkBevisS88}.

\medskip
Lemmas~\ref{LemmaP19:1} and \ref{LemmaP20:2} below are related
to the defining properties for
quasi-manifolds, resp. pseudomanifolds, cf. p.~\pageref{SubSecP27:III:1}.
We write out all seven items mainly just to be able to see the
details once and for all.

\smallskip
In the $\ast$-case, lemma~\ref{LemmaP19:1} below is trivially
true if any $\Sigma_{i}\!\!=\!\{\emptyset_{\!o}\}$.
``codim\lower3.5pt\hbox{$^{\hbox{$\sigma\ge 2$}}$}''
($\Rightarrow$
$\dim{\rm Lk}_{_{^{\!\Sigma}}}\!\sigma\!\geq\!1$)
means that a maximal simplex, say $\tau$, containing  $\lower3.5pt\hbox{$^{\hbox{$\sigma$}}$}$
fulfills $\dim\tau\!\ge\dim\sigma+\!2$.

\begin{lemma}
\label{LemmaP19:1}
Read $\mytinynabla$ below
as ``$\times\!$'' or throughout as ``$\ast$''.
When $\mytinynabla$ is $\ast$-substituted, $\Sigma_{{i}}$ is assumed to be connected \underbar{or}
$0$-dimensional.
${\bf G}_{{1}},{\bf G}_{{2}}$ are ${\bf A}$-modules such that $\hbox{\rm Tor}_1^{\mathbf{A}}({\bf G}_{1},{\bf G}_{2})=0$.
Now, if $\dim\Sigma_{i}\geq0$ and $v_i\!:=\dim \sigma_{i}\ (i=\ ,1,2)$
then {\bf D}$_1$-{\bf D} are all \hbox{equivalent};\\
$($
${\rm Int} (\sigma) \!:=\! \{\alpha\in|\Sigma|\mid[{\bf
v}\in\nobreak\sigma ] \Longleftrightarrow [\alpha({\bf v})\neq0]
\}$
and
${\bar{\sigma}} :=\{\tau \mid \tau
\hbox{\scriptsize$\subseteqq$}
\sigma\}$.
$)$

\smallskip
\noindent{\bf D}$_1)$\ $\mhatH_{0} ({\rm
Lk}_{(\Sigma_{1}{\lower1pt\hbox{\mytinynabla}}\Sigma_{2})}\sigma; {\bf G}_1\otimes{\bf
G}_2)=0$
                       for $\emptyset_o \neq \sigma\in
                   \Sigma_{1}\mytinynabla\Sigma_{2}$,
                    whenever codim$\sigma\ge 2.$

\smallskip
\noindent{\bf D$_1^{\prime})$}\
            $\mhatH_{0}
                      ({\rm Lk}_{\Sigma_{i}}\sigma_{i};{\bf G}_i)=0$
                           for $\emptyset_o \neq \sigma_{i}\in \Sigma_{i}$,
                     whenever codim$\sigma_{i}\ge2\ (i=1,2)$.

\smallskip
\noindent{\bf D$_2)$} $\mhatH_{v+1}
{(\Sigma_{1}\!\mytinynabla\!
\Sigma_{2}}, \hbox{\rm
cost}_{_{\!\Sigma_{1}\!\!\nabla\!\Sigma_{2}}}\!\!\sigma; {\bf
G}_1\otimes{\bf G}_2)=0$
                     for $\emptyset_o \neq \sigma\in
                   \Sigma_{1}\mytinynabla\!\Sigma_{2}$,
                   if
                   codim$\sigma\ge 2.$

\smallskip
\noindent{\bf D$_2^{\prime})$} $\mhatH\hskip-0.2cm\lower0.1cm\hbox{$_{v_i+1}$}\!
                    (\Sigma_{i},\hbox{\rm
                     cost}\lower2pt\hbox{\sixrm {\char"06}}_{i}\!\sigma\!_{i}\!;{\bf G}_i)=0$
                           for $\emptyset_o \neq \sigma_{i}\in \Sigma_{i}$,
                     whenever codim$\sigma_{i}\ge2\ (i=1,2)$.

\smallskip
\noindent{\bf D$_3)$}\ \ $\!\mhatH\hskip-0.2cm\lower0.1cm\hbox{$_{v+1}$}\!
                       (|\Sigma_{1}\mytinynabla\!\Sigma_{2}|,
                         |\Sigma_{1}\mytinynabla\!\Sigma_{2}|
                         \setminus_{_{\!^o}}\!\alpha;
                        {\bf G}_1\!\otimes{\bf G}_2)=0$
                     for all $\alpha_0\ne\alpha_{\!}\in\hbox{\rm Int}(\sigma)$
                           if codim$\ \!\sigma\!_{\!}\ge2.$

\smallskip
\noindent{\bf D$_3^{\prime})$}\
                  $\mhatH_{v_i+1}
                      (|\Sigma_{i}|,|\Sigma_{i}|\setminus_{o}{_{\!}}
                      \alpha_i;{\bf G}_i)=0$
                     for $\alpha_0\ne \alpha_i\in \hbox{\rm Int}(\sigma_{i})$,
                            if codim$\sigma_{i}\ge2\ (i=1,2)$.

\smallskip
\noindent {\bf D$)$} $\mhatH_{v+1}
                   (\ \!|\Sigma_{1}|\bar{\mytinynabla\!}|\Sigma_{2}|,
                    |\Sigma_{1}|\bar{\mytinynabla\!}|\Sigma_{2}|
                    \setminus_{o}{_{\!}}
                    (\alpha_1,\alpha_2);
                   {\bf G}_{1}\!\otimes{\bf
                   G}_{2})\!\!=\!0\!$
                    for all $\alpha_0\ne{(\widetilde{\alpha_1,\alpha_2})}\in{\rm Int}(\sigma_{ })$\\
            $
             \subset|\Sigma_{1\!}
            \mytinynabla\!\Sigma_{2}|$
if codim$\sigma\ge2$, where
$\eta\!:|\Sigma_{1}\!\mytinynabla\!\Sigma_{2}|\ \!
\rlap{$\rightarrow$}{\raise4pt\hbox{$\simeq$}}\
|\Sigma_{1}|{{\bar{\mytinynabla\!}}}|\Sigma_{2}|$ and
${\eta(\widetilde{\alpha_1,\alpha_2})}\!=\!(\alpha_1,\alpha_2)$
{\rm($k$-ifikations never effect the homology modules).}
\end{lemma}

\begin{proof}
By the homogenity of the interior of
$\vert\bar\sigma\!_{_{^{1}}}\!\times\bar\sigma\!_{_{^{2}}}\!\vert$
we only need to deal with simplices $\sigma$ fulfilling
$c_{_{^{\!\sigma}}}\!\!\!=\!0$.
By
Prop.~1 p.~\pageref{PropP15:1}
and Eq.~\ref{LocalTop} above,
all non-primed, resp. primed, items are equivalent among
themselves.
{\bf D$\!_1$} $\Leftrightarrow$ {\bf D$\!_1^{\ \! \prime}$} by
Corollary~\ref{CorP19:ToTh6} above, since $\hbox{\rm Tor}_1^{\mathbf{A}}({\bf G}_{1},{\bf G}_{2})=0$.\\
For joins of finite complexes, this is done explicitly in \cite{12} p.~172.
\end{proof}

\vfill\break
\normalbaselines


\subsection{Connectedness for simplicial complexes} \label{SubSec:Connectedness}
$$\phantom{.}$$
\vskip-0.8cm
Any $\Sigma$ is representable as $\Sigma=
\bigcup\!\!\!\!\!\!\!\!\!\!_{_{_{_{_{\sigma ^m \in \Sigma}}}}}\!\!
\overline {\sigma ^m},$ where $\overline {\sigma ^m}$ denotes the
simplicial complex generated by the maximal\ simplex $\sigma^m$.

\begin{definition}
\label{DefP20:1}
Two maximal faces $\sigma,\tau\!\in\!\Sigma$ are {\it strongly
connected} if they can be connected by a finite sequence $\sigma=
{{\delta\!_{_{0}}}},..,{\delta\!_{_{i}}},..,{{\delta\!_{_{^{q}}}}}\!
= \tau $ of maximal faces with
$\#({\delta\!_{_{i}}}\!\cap{\delta\!_{_{i+1}}}\!)=
\hbox{max}\!\!\!\!\!\!\!\!\!\!\!\! _{_{_{{{0\le j\le q}}}}}
\#{\delta\!_{_{^{j}}}}\!-_{\!}1$ for consecutives.
Strong connectedness imposes an equivalence relation among the
maximal faces, the equivalence classes of which defines the {\it
maximal strongly connected components} of $\Sigma,$ {\rm cf.\
\cite{2}
p.\ 419ff.}
A complex $\Sigma$ is said to be {\it strongly connected} if each pair of of
its  maximal simplices are strongly connected.
A {\it submaximal face}  has exactly one vertex less then some maximal
face  containing it.
\end{definition}

Note that strongly connected complexes are pure, i.e.,
${_{^{\!\!}}}\sigma {_{^{\!\!\ }}} \!\!\in\!\Sigma{_{^{\!\ }}}\
\!\!\hbox{\rm maximal}
\Rightarrow\dim\!\sigma\nobreak=\nobreak\dim\!\Sigma.$

\begin{lemma}
\label{LemmaP20:2}
{\rm(This lemma concerns the defining properties for pseudomanifolds.
\cite{10}
p.\ 81 gives a proof, valid for finite-dimensional complexes.)}

\smallskip
\noindent
{\bf A}$)$ If $d_i\!:=\dim\Sigma_{i}\geq0$
then $\Sigma_{1} \!\times\! \Sigma_{2}$
is pure
$\Longleftrightarrow$ $\Sigma_{1}$ and $\Sigma_{2}$  are
both pure.

\smallskip
\noindent{\bf B}$)$ If $\dim\sigma_{i} ^m\!\geq1$ for each maximal
simplex $\sigma_{i}^m \!\in\! \Sigma_{i}$ then

\smallskip
Any submaximal face in $\Sigma_{1}
 \!\times\!
\Sigma_{2}\!$ lies in at most (exactly) two maximal faces
$\Longleftrightarrow $
Any submaximal face in $\Sigma_{i}\!$ lies in at most
(exactly) two maximal faces of $\Sigma_{i}$, where $i=1,2.$

\smallskip
\noindent{\bf C}$)$ For $d_i>0$;
$\Sigma_{1}\!\times\! \Sigma_{2}$ strongly connected
$\Longleftrightarrow$ $\Sigma_{1},\Sigma_{2}$ both
strongly connected.
\end{lemma}

\begin{note} \label{NoteP20:1}
Lemma~\ref{LemmaP20:2} is true also for $\ast$ with exactly the same reading but now with no other restriction than that $\Sigma_{_
{\!i}}\!\ne\emptyset$ and this includes, in particular, heading {\bf B}.
\end{note}

\begin{definition}\label{DefP20:2}
$ \Delta\smallsetminus \Delta\!^{^{\!_{o}}}
\!:=\{\delta\in\Delta\ \mid\ \delta\!\not\in\!\Delta\!^{^{\!_{o}}}
\} $ is {\it connected as a poset} $(${\rm partially ordered set}$)$
with respect to $\!$simplex inclusion if, for every pair
$\sigma,\tau\!\in\!\Delta\smallsetminus \Delta\!^{^{\!_{o}}}$, there is a chain
$\sigma\!=\!\sigma\!_{_{^{0}}}\!,\sigma\!_{_{^{1}}}\!,
...,\sigma\!_{_{^{k}}}\!\!=\!\tau$, where
$\sigma\!_{i}\!\!\in\! \Delta\smallsetminus \Delta\!^{^{\!_{o}}} $ and
$\sigma\!_{i}\!\!\subseteq\!\sigma\!_{_{^{\!i+1}}}\!\!$ or
$\sigma\!_{i}\!\!\supseteq\!\sigma\!_{_{^{\!i+1\!}}}$.
\end{definition}

\begin{note} \label{NoteP20:2}
$\!$(cf. \cite{12}
p.\ 162.)
Let $\Delta \mmysupsetneqq \Delta\!^{^{\!_{o}}} \mmysupsetneqq
\{\emptyset_{_{^{\!o}}}\!\}$.
$\Delta_{\!}\smallsetminus \Delta\!^{^{\!_{o}}}\!$
             is connected as a poset \underbar{iff}
             $|\Delta^{}|\ \!{{{\setminus}\!_{_{^{o}}}}}|\Delta\!^{^{\!_{o}}}|$
             is pathwise connected.
             When $\Delta\!^{^{\!_{o}}}=\{\emptyset_{_{^{\!o}}}\!\},$
             then the notion of connectedness as a
             poset\nobreak\ is\nobreak\ equivalent to
             \hbox{the usual one for $\Delta^{}$.}
             Let $\Delta\!^{^{\!_{(\!p\!)}}}\!\!:=
             \{\sigma\in\Delta\ \!\vert\ \!\#\sigma\le p+1\}$; then
             $|\Delta^{}|$ is connected \underbar{iff}
             $|\Delta\!^{^{\!_{(\!1\!)}}}\!|$
             is connected.
\end{note}

\begin{lemma}
\label{LemmaP20:3}
{\rm(\cite{12}
p.\ 163)}
                    { $\Delta\smallsetminus \Delta\!^{^{\!_{o}}}$
                     is connected as a poset
                     \underbar{\rm iff}
                    for each pair of maximal $\hbox{\rm simplices}\
              \sigma,\tau\in\Delta\smallsetminus \Delta\!^{^{\!_{o}}}$,
                there is a chain,
                $\sigma\!=\!\sigma\!_{_{^{0}}}\!\supseteq\!\sigma\!_{_{^{1}}}\!
                 \subseteq\!\sigma\!_{_{^{2}}}\!\supseteq\!...\!
                    \subseteq\!\sigma\!_{_{^{2m}}}\!\!=\!\tau\!$  in
                $\Delta\smallsetminus \Delta\!^{^{\!_{o}}}$,
               where the $\sigma_{_{{\!\!2i}}}\!$s are maximal faces and
                $\sigma\!_{_{^{2i}}}\!{\raise0.5pt\hbox{\eightmsbm \char"72}}
                 \sigma\!_{_{^{2i+1}}}\!$ and
               $\sigma\!_{_{^{2i+2}}}\!
{\raise0.5pt\hbox{\eightmsbm \char"72}}
\sigma\!_{_{^{2i+1}}}\!$
                    are situated in different components of
                    ${\rm Lk}\!_{_{^{\Delta}}}\!\!\sigma\!_{_{^{2i+1}}}
                   (i\!=\!0,1,..., m\!-\!1).$}
\quad
\end{lemma}

\begin{lemma}
\label{LemmaP20:4}
{\rm (A. Bj\"orner 1995.)}
{\rm (A direct consequence of Lemma~\ref{LemmaP20:3}.)}
{Let $\Sigma$ be a finite-dimensional simplicial
           complex, and assume that {\rm Lk}$_{_{\!\Sigma}}\sigma$ is
connected for all $\sigma\!\in\!\Sigma$, inkluding\ \
$\emptyset_{_{^{\!o}}}\!\!\in\Sigma$, such that $\dim${\rm
Lk}$_{_{^{\!\Sigma}}}\!\sigma\!\geq\!1.$ Then $\Sigma$ is pure and
strongly connected.} \quad
\end{lemma}

\goodbreak
\section{Relating Combinatorics to Commutative Algebra} \label{SecP21:II}

\subsection{Definition of Stanley-Reisner rings}$\!$ \label{SubSecP21:II:1}

\medskip
\centerline{---$\ast\ast\ast$--- What we are aiming at.
---$\ast\ast\ast$---}
\noindent
\S{\bf \ref{SubSecP22:II:2}.}
Through the Stanley-Reisner Functor below, attributes like
Buchsbaum, Cohen-Macaulay and Gorenstein on (gra\-ded) rings and
algebras (13H10 in MSC2000), became relevant also within Combinatorics as the classical
definition (Def.~2 p.~\pageref{DefP2:2}) of simplicial complexes was
altered to Def.~1 p.~\pageref{DefP1:1}.\nolinebreak\ Now this
extends to General/Algebraic Topology, cf.
p.~\pageref{SubSecP22:II:2}.
Theorem~\ref{TheoremP23:8} p.~\pageref{TheoremP23:8}
indicates that simplicial homology manifolds can be inductively generated.

 \noindent
\S{\bf \ref{SubSecP28:III:2}.}
Corollary~\ref{CorP30:1} p.~\pageref{CorP30:1},
tells us that there is no $n$-manifold ($\ne\bullet$) with an $(n-2)$-dimensional boundary.
The examples on p.~\pageref{ExampleP32}
implicitly raise the arithmetic-geometrical question: Which
boundary-dimensions are accessible for quasi-$n$-manifolds with respect to different coefficient-modules?

\noindent
\S{\bf \ref{SubSecP32:III:4}.}
Corollary~\ref{CorP33:2} p.~\pageref{CorP33:2},
confirms Bredon's conjecture in \cite{1}
p.\ 384, that homology manifolds with ${\mathbb{Z}}$ as coefficient-module are locally orientable, and so, these manifolds have either an empty, a $(-1)$-dimen\-sional or an $(n-1)$-dimensional boundary.

\centerline{---$\ast\ast\ast$--- \indent ---$\ast\ast\ast$---}

\begin{definition} \label{DefP21}

A subset $s\!\!\subset\!\!{\bf W}\!\!\supset\!\! V_\Delta$ is said
to be a {\it non}-{\it simplex} (with respect to\ {\bf W}) of a simplicial complex $\Delta$, denoted
$s~{\propto}~\Delta,\
\text{if}\ s\!\not\in\!\Delta$ but
 $\dot s\!=\!{(\bar s)}^{(\dim{s})-1}\!\!\!\!\subset\!\Delta$
(i.e.,\ the $(\dim s-1)$-dimensional skeleton of $\bar s$, consisting of all proper subsets of $s$, is a subcomplex of $\Delta)$. For a simplex $\delta = \{ v_{i_1},\dots, v_{i_k}\}$ we define $m_{\delta}$ to be the square-free monic monomial
 $m_{\delta} :=1_{_{\bf A}}\!\!\cdot v
 _{i_1}\!\!\cdot\!\dots\!
\cdot v_{i_k}\in {\bf A}[{\bf W}]$ where $ {\bf A}[{\bf
W}]$ is the graded polynomial algebra on the variable set ${\bf W}$
over the commutative ring ${\bf A}$ with unit $1\!_{_{\bf A}}.$
\hbox{\rm So},\ $m_{\emptyset\!_o}\!\!=1\!_{_{\bf A}}.
$

Let ${\bf A} \langle  \Delta \rangle  := {\bf A}[{\bf W}]/{{\bf
I}_{\Delta}}$ where {\bf I}$_{\Delta}$
is the ideal generated by $\{m_\delta  \mid \delta
{\propto}\Delta \}$.
${\bf A} \langle  \Delta \rangle $ is the {\it face ring} or
{\it Stanley-Reisner $($St-Re$)$ ring} of $\Delta$ over ${\bf A}$.
Frequently ${\bf A}={\bf k}$, a field, e.g. a prime field, the real numbers $\mathbb{R}$ or the rational numbers $\mathbb{Q}$.
\end{definition}

{\it Stanley-Reisner} (St-Re) {\it ring theory} is a basic tool
within combinatorics, where it supports the use of commutative
algebra.
The generating set for the ideal ${\bf I}_{\Delta}$ above is minimal, while in the traditional definition ${\bf I}_{\Delta}$ is generated by all squarefree monomials
$p_\delta=v_{i_1}\!\!\cdot\!\dots\!\cdot v_{i_k}\in {\bf A}[{\bf W}]$ such that $\{ v_{i_1},\dots, v_{i_k}\}\not\in\Delta$, i.e. ${\bf I}_{_{^{\!\Delta}}}\!\!=
\hbox{\tenbf(}\{p_\delta \mid \delta \!\not\in\!\Delta
\}\hbox{\tenbf)}$.
This latter definition of $p_\delta$ in
${\bf I}_{_{^{\!\Delta}}}\!\!=
\hbox{\tenbf(}\{p_\delta \mid \delta \!\not\in\!\Delta
\}\hbox{\tenbf)}$
leaves the questions
``$\!{\bf A} \langle  \{\emptyset_{_{\!^o}}\!\} \rangle =\
\!?$'' and ``$\!{\bf A}\langle \
\!\emptyset\ \!\rangle =\ \!?$'' open,
while Note~\ref{NoteP21}.{\bf ii} below, shows that:
${{\bf A}} \langle  \{\emptyset_{_{\!^o}}\!\}
\rangle \cong {\bf A}\!\ne\! {\bf0}$
=``{\it The trivial ring}''$= {\bf A}\langle \ \!\emptyset\
\!\rangle $.

\begin{note}\label{NoteP21}
\!\hskip-0.1cm{\bf i.}
Let $\overline{{\bf P}}$ be the set of all finite subsets of a set
${\bf P}$, like the $\bar{s}$ above. Now,
${\bf A} \langle  \ \! \overline{{\bf P}}\ \! \rangle ={\bf
A[{\bf P}]}.$
$ \overline{{\bf P}}$ is known as the ``{\it the full simplicial complex on ${\bf P}$}''. The set of natural numbers {\bf N} gives
$\overline{{\bf N}}$, named the ``{\it infinite simplex}''.

Another useful concept is the {\it costar}, defined by
${\rm cost}_{_{\!{\Sigma}}}\!\sigma\!:=\! \{\tau\!\in\!\Sigma\mid\
\tau\!\not\supseteq\sigma\}$.
We have
$\Delta=\mbfcupcap{\bigcap}{s\propto\Delta}
\hbox{\rm cost}\!\!\!_ {_{_{\overline{\bf W}}}}\!\!{{{s}}\!_{_{^{\ }}}}\!$, for any universe {\bf W} containing the vertex set $V_\Delta$.

The cone points in a simplicial complex $\Delta$ are those vertices \  ${\bf v}\!\in\!V_\Delta$ that are contained in every maximal simplex of $\Delta$. So, $\delta\!\!_{_{\Delta}}\!:=\!\{{\bf v}\in\! V_\Delta\!\mid {\bf v}\ \text{cone point}\}\!\in\!\Delta.$
The set of cone points in $\Delta$ is also characterized as the intersection of all maximal simplices in $\Delta$, and thus resembling the
\htmladdnormallink{{\it radical}}
{http://en.wikipedia.org/wiki/Jacobson_radical}
in algebra.

Many moore basic simplicial constructions are collected in  pp.~\pageref{SecP34:Appendix}ff.

\smallskip
\noindent {\bf ii.}\label{NoteP21:ii}
${\bf A} \langle  \Delta \rangle
\ \!\cong\ \!
\lower0.0pt\hbox{\sixrm${ {{\lower0.0pt\hbox{${\hbox{\eightbf
A}}$}}[{\lower0pt\hbox{\eighti V}}\!\!_{_{\Delta}}] \over
(\{m_{\delta}
\in{\bf A}
[{\lower0pt\hbox{\eighti V}}\!\!_{_{\Delta}}]
\mid\ {\delta}\ \!{\propto}\ \!\!{{\Delta}} \})}}$} \raise3pt\hbox{\ \ },$ if $\Delta\neq \emptyset$,
${\{{\emptyset}\!_o\!\}}.$
           So, the choice of the universe {\bf W}
           isn't all that critical.

If  $\Delta=\emptyset$, then
           the set of non-simplices equals $\{\emptyset\!_o\}\!$,
since
           $\emptyset_o\!\not\in\emptyset$, and
           ${\overline{\emptyset}\!_o\!\!}
           ^{_{((\dim\emptyset_{\!o})-1)}}\hskip-0.1cm  =\!
           {\{\emptyset_o\}}^{_{(-2)}}\!\!=\emptyset
           \subset\emptyset$, implying ${\bf A}\langle \ \!\emptyset\ \!\rangle \!=\!{\bf0}\!=
           $``{\it The trivial ring}'' since $m_{\emptyset\!_{o}}\!\!=\!1_{_{{\bf A}}}$.

Since $\emptyset \in \Delta$ for
           every simplicial complex $\Delta\neq \emptyset$,
           ${\{{\bf v}\}}$ is a non-simplex of
           $\Delta$ for every
           ${\bf v}\!\in\!{\bf W}\setminus {\it V}\!_{\Delta}$,
           i.e.,
           $$[{\bf v}\not\in{\it V}\!_{\Delta}\neq\emptyset
           {{{{]}{\Longleftrightarrow}}}{[}}\{{\bf v}\}
            {\propto}\Delta
           \neq \emptyset].$$
           So
${\bf A} \langle  \{\emptyset_{_{\!^o}}\!\} \rangle = {\bf A}$
since $\{\delta\mid{\delta}
{\propto}\Delta\}= {\bf W}.$

\medskip
\noindent {\bf iii.}\label{NoteP21:iii}
${\bf k} {\langle }
\Delta_{_{^{_{\!}1\!}}}\ast\Delta_{_{^{_{\!}2\!}}}
{\rangle \!}
$
$\ \!\cong\ \!$
${\bf k} {\langle }
\Delta_{_{^{_{\!}1\!}}} {\rangle }
\ \!{\otimes}\ \!
{\bf k} {\langle } \Delta_{_{^{_{\!}2\!}}} {\rangle }
$
(R. Fröberg, 1988)
$
\ \hbox{\rm with}\
{\bf I}
{\lower2.0pt\hbox{\sevenrm {\char"01}$_{_{1}}\!${\sevensy
{\char"03}}\sevenrm {\char"01}$_{_{2}}$
}}
\!=
\hbox{\rm (} \{m_{\delta}\ \!\mid\ \!\![
{\delta}\ {\propto}\
\Delta_{_{^{_{\!}1\!}}}\!\!\lor
{\delta}\ {\propto}\
\Delta_{_{^{2}}}\!]
\land[\delta\!\notin\!{\Delta_{_{^{_{\!}1\!}}}\!
\!\ast\!{\Delta_{_{^{2}}}\!}}]\}{\bf )} $
by \cite{10}
Ex. 1\ p.\ $\!$70.\ \
(Here the parenthesis $``(\cdot)":=$``the ideal generated by $\cdot$")

\medskip
\noindent {\bf iv.}\label{NoteP21:iv}
If $\Delta_i\neq\emptyset\ \text{for}\ i=1,2$, it is well known that

{\bf a.}
${\bf I}_{{{\Delta}_1}\cup{{\Delta}_2}}
{={\bf I}\!_{{\Delta}_1}\!\cap{\bf I}\!_{{\Delta}_2}
=
(\{m=\hbox{\rm Lcm}(m_{\delta_1},m_{\delta_2})\ \!\mid\ \!
\delta_i{{\propto}\Delta_i}\ \text{for}\ i=1,2
\})};$

\indent {\bf b.} ${\bf I}
_{{{\Delta}_1} \cap{{\Delta}_2}}=
{\bf I}\!_{{{\Delta}_1}}
\!+
{\bf I}\!_{{\Delta}_2}
\!= (\{m_{{\delta}}\ \mid\ \delta{{\propto}\ \!\Delta_1}\lor \delta{{\propto}\ \!\Delta_2}\} )$ in {\bf A}[{\bf W}].\\
${\bf I}_{\Delta_1} \!\!\cap {\bf I}_{\Delta_2}$
and
${\bf I}_{\Delta_1} \!\!+ {\bf I}_{\Delta_2}$
are generated by a set (no restrictions on its cardinality) of
square-free monomials, if both ${\bf I}{_{\Delta_1}}$ and ${\bf I}{_{\Delta_2}}$ are square-free.
These squarefree monomially generated ideals form a distributive
sublattice $({\mathcal J}\ \!^{\circ};\cap,+,{\bf A[W]}_+)$, of the
ordinary lattice structure on the set of ideals of the polynomial
ring {\bf A[W]}, with a counterpart, with reversed lattice order,
called the {\it squarefree monomial rings {with unit}}, denoted ({\bf
A}$^{\!\circ}${\bf[W]}$;\cap,+,{\bf I})$. {\bf A[W]}$_+$ is the
unique homogeneous maximal ideal and zero element. We can use the
ordinary subset structure to define a distributive lattice structure on ${\Sigma}_{\mathbf{W}}^\circ:=$ ``The set of {non-empty}
simplicial complexes over {\bf W}'', with $\{\emptyset_o\}$ as zero
element and denoted (${\Sigma}_{\mathbf{W}}^\circ;{\cup},\cap,\{\emptyset_o\})$.
The {\it Weyman/Fr\"oberg/Schwartau Construction} eliminates the
above square-free demand, cf. \cite{32}
p. 107.
\end{note}

\begin{proposition} \label{PropP21}
$\!{\!}$The Stanley-Reisner ring assignment functor defines a
mo\-nomorphism on distributive lattices
$({\Sigma}_{\mathbf{ W}}^\circ,{\cup},\cap,\{\emptyset_o\})\longrightarrow({\bf I}^{\!\circ}{\bf[W]},\cap,+,{\bf I})$, which is an
isomorphism for finite {\bf W}.
\quad
\end{proposition}

\subsection{Buchsbaum complexes are weak manifolds} \label{SubSecP22:II:2}

Combinatorialists call a finite simplicial complex a {\it Buchsbaum} ({Bbm}$_{\!_{\mathbf{k}}}\!$) ({\it Cohen-Macaulay}
$({\rm CM}_{\!_{\mathbf{k}}})$, {\it Gorenstein} $({\rm Gor}_{\!_{\mathbf{k}}})$) complex if its Stanley-Reisner ring is a Buchsbaum (Cohen-Macaulay, Gorenstein) ring.
The following three propositions, found in \cite{31} together with proof-references, can be used as "definitions" when restricted to finite simplicial complexes. The proof of Proposition~\ref{PropDef:Gor}.iv is given by A. Björner as an appendix in \cite{10}.

\begin{proposition} \label{PropDef:Bbm}
{\rm(Schenzel)  (\cite{31} Th. 8.1)} Let $\Sigma$ be a finite simplicial complex and let  $\mathbf{k}$ be a field. Then the following are equivalent:\\
{\rm({\bf i})} $\Sigma$ is Buchsbaum over $\mathbf{k}.$\\
{\rm({\bf ii})} $\Sigma$ is pure, and $\mathbf{k}\langle\Sigma\rangle_{\mathfrak{p}}$ is Cohen-Macaulay for all prime ideals ${\mathfrak{p}}$ different from the unique homogeneous maximal ideal i.e., the irrelevant ideal.\\
{\rm({\bf iii})} $\mbox{ For all } \sigma\in \Sigma,\ \sigma\ne \emptyset \mbox{ and } i< \dim(${\rm Lk}$_{\Sigma}\sigma),$
$ {\hat{\mathbf{H}}}_{i}({\rm Lk}_{\Sigma}\sigma;\mathbf{k})=0.$\\
{\rm({\bf iv})} $\mbox{ For all } {\alpha}\in|\Sigma|,\ \alpha\ne \alpha_0 \mbox{ and } i<\dim\Sigma,{\hat{\mathbf{H}}}_{i}(|\Sigma|,|\Sigma|\setminus_o\alpha ;\mathbf{k})=0.$\hfill$\square$
\end{proposition}

\goodbreak
\begin{proposition} \label{PropDef:CM}
{\rm(\cite{31} Prop. 4.3, Cor. 4.2)} Let $\Sigma$ be a finite simplicial complex and let  $\mathbf{k}$ be a field. Then the following are equivalent:\\
{\rm({\bf i})} $\Sigma$ is Cohen-Macaulay over $\mathbf{k}.$\\
{\rm({\bf ii}) (Reisner)} $\mbox{ For all } \sigma\in \Sigma \mbox{ and } i< \dim(${\rm Lk}$_{\Sigma}\sigma), {\hat{\mathbf{H}}}_{i}({\rm Lk}_{\Sigma}\sigma;\mathbf{k})=0.$\\
{\rm({\bf ii}) (Munkres)} $\mbox{ For all } {\alpha}\in|\Sigma| \mbox{ and } i<\dim\Sigma,
{\hat{\mathbf{H}}}_{i}(|\Sigma|,|\Sigma|\setminus_o\alpha;\mathbf{k})=0.$
\hfill$\square$
\end{proposition}

\begin{proposition} \label{PropDef:Gor}
{\rm(\cite{31} Th. 5.1)} Let $\Sigma$ be a finite simplicial complex, $\mathbf{k}\mbox{ a field or } \mathbb{Z} \mbox{ and }\Gamma:={\rm core}\Sigma:= {\rm Lk}_{_{^{\!\Sigma}}}\delta_{_{^{\!\Sigma}}}.$  Then the following are equivalent:\\
{\rm({\bf i})} $\Sigma$ is Gorenstein over $\mathbf{k}.$\\
{\rm({\bf ii})} $\mbox{ For all }  {\sigma}\in\Gamma,\ {\hat{\mathbf{H}}}_{i}(\mathrm{Lk}_\Gamma\sigma;\mathbf{k})\cong
\begin{cases}
\mathbf{k}\ \ \mbox{if}\ i=\dim({\rm Lk}_{\Gamma}\sigma),\\
\mathbf{0}\ \ \mbox{if}\ i<\dim({\rm Lk}_{\Gamma}\sigma)
\end{cases}$\\
{\rm({\bf iii})} $\mbox{ For all }  \alpha\in|\Gamma|,\ {\hat{\mathbf{H}}}_{i}(|\Gamma|,|\Gamma|\setminus_o \alpha;\mathbf{k})=\begin{cases}
\mathbf{k}\ \ \mbox{if}\ i=\dim\Gamma,\\
\mathbf{0}\ \ \mbox{if}\ i<\dim\Gamma.
\end{cases}$\\
{\rm({\bf iv})} {\rm(A. Björner)} $\Sigma$ is Gorenstein $_\mathbf{k}  \Longleftrightarrow\Sigma$ is C-M $_\mathbf{k}, \mbox{ and } \Gamma$ is an orientable pseudomanifold without boundary.\\
{\rm({\bf v})} Either {\rm({\bf1})}  $\Sigma=\{\emptyset\},\ \bullet,\ \bullet\bullet$, or
{\rm({\bf2})} $\Sigma$ is Cohen-Macaulay  over  $\mathbf{k}$, $\dim\Sigma\geq 1,$ and the link of every $(\dim(\Sigma)-2)$-face is either a circle or a line with two or three vertices and  $ \tilde \chi (\Gamma)=(-1)^{\dim\Gamma}$.
\hfill$\square$
\end{proposition}

Through Prop.~\ref{PropP15:1} p.~\pageref{PropP15:1},
the following definitions for arbitrary modules and
topological spaces are consistent with the original.

\begin{definition} \label{DefP22}

\noindent
{$X$ is ``{Bbm}$_{_{\mathbf{G}}}$'' (``{CM}$_{_{\mathbf{G}}}$'', 2-$``${CM}$_{_{\mathbf{G}}}$'') if} $\!X$ is an
$n$-{whm}$_{_{\mathbf{G}}}\!$ ($n$-{jwhm}$_{_{\mathbf{G}}}$,
$n$-whsp$_{_{{\bf G}\!}}).$

A simplicial complex $\Sigma$ is defined to be
``$\ \!${\rm Bbm}$_{_{\mathbf{G}}}$'',
``{CM}$_{_{\mathbf{G}}}$''
resp. 2-$``${CM}$_{_{\mathbf{G}}}$'' if $|\Sigma|$ is.
In particular,
$\Delta$ is $2$-``{\rm CM}$\!_{_{\bf G}}$'' \underbar{iff}\ \
$\Delta$ is ``{\rm CM}$\!_{_{\bf G}}$''
and\ \
$ \mhatH\!\!\!\!_{_{n\!_{_{}}\!-\!1}}\! (\hbox{\rm
cost}\!_{_{\Delta}}\!{{\delta\!_{_{\ }}}}\!;{\bf G})=0,\ \! \forall\
\delta\!\in\!\Delta, $
\hbox{cf. \cite{31} Prop.~3.7 p.~94.}

(The {\it n} in
``$n$-{\rm Bbm}$_{_{\mathbf{G}}}$'' (``$n$-{\rm CM}$_{_{\mathbf{G}}}$'' resp. 2-``{$n$-CM}$_{_{\mathbf{G}}}$'') is deleted since any interior point $\alpha$ of a
realization of a maxi-dimensional simplex gives
$ \mhatH\!\!\!\!\!\raise0.5pt\hbox{{$_{_{_{\dim\Delta}}}$}}
\!\!(|\Delta|,|\Delta|\setminus_{_{\!o}\!}\alpha;{\bf G}) ={\bf G}.
$
Recall that $\sigma\in\Sigma$ is said to be {\it maxi-dimensional} if $\dim\sigma =
\dim\Sigma$. The $\vert\{\emptyset\!_{o}\}\vert$-extended {``$\setminus_{_{\!o}\!}$}'' is defined in p.~\pageref{DefP8:PointSetMinus0}.)
\end{definition}

So we are simply renaming
$n$-{whm}$_{_{\mathbf{G}}},$
$n$-{jwhm}$_{_{\mathbf{G}}}$\nobreak\ and
$n$-whsp${_{_{\mathbf{G}}}}_{\!}$
to
$_{\!}$``{Bbm}$_{_{\mathbf{G}}}$'',
$\!``${CM}$_{_{\mathbf{G}}}$'' resp.
2-$``${CM}$_{_{\mathbf{G}}}$'',
where the quotation marks indicate that we are not limited to compact spaces nor to just ${\mathbb{Z}}$ or {\bf k} as coefficient modules.

\begin{varnorem} \label{RemarkP22:N.B.}
{\bf N.B.} The definition p.~\pageref{DefP18:RealizationOfSimpSet}
of
$\mhatbar
\Delta\!\!^{^\wp}\!({\it X})
\mhatbar$
provides each topological space $_{\!}X_{\!}$ with a
Stanley-Reisner ring, which, with respect to
``{Bbm}$_{_{\mathbf{G}}}$''-,
$\!``${CM}$_{_{\mathbf{G}}}$''-
and
2-$\!``${CM}$_{_{\mathbf{G}}}$''-ness, is triangulation invariant. \end{varnorem}

\begin{proposition}
\label{PropP22:1}
{The following conditions are equivalent}:

{\bf a.} $\Delta$ is $\hbox{\rm ``Bbm}\!_{_{\bf G}}"$,

\smallskip
{\bf b.} {\rm(Schenzel, \cite{14}
p.\ 96)} $\Delta$ is pure and ${\rm
Lk}\!_{_{{\Delta}}}\!\delta$ is
``CM$\!_{_{\bf G}}"\ \text{for all}\
\emptyset\!_{o}\ne{\delta}\in{\Delta}$,

\smallskip
{\bf c.} {\rm(Reisner, \cite{3} {5.3.16.(b)}
p.\ 229)} $\Delta$ is pure and $ {\rm
Lk}\!_{_{{\Delta}}}\!\!{{{\bf v}\!_{_{ }}}}\ \hbox{\it is}\
\hbox{\rm ``CM}\!_{_{\bf G}}\!\!"\ \text{for all}\ {\bf v}\!\in\!{\it V}_{_{\!\!\!\Delta}}. $
\end{proposition}

\begin{proof}$\!$
Use Prop.\ \ref{PropP15:1} p.~\pageref{PropP15:1}.
and Lemma~\ref{LemmaP20:4} p.~\pageref{LemmaP20:4}.
and then use Eq. {\bf I} p.~\pageref{EqP34:I}.
\end{proof}

\begin{varex}\label{Ex:HomologyVsSt-Re}
When limited to compact polytopes  and a field {\bf k} as
coefficient module, we add, from \cite{31}
p.\ 73,
the following Buchsbaum-equivalence using {\it local\ cohomology};

\smallskip
{\bf d.} {\rm(Schenzel)} \  $\Delta$ {\it is Buchsbaum}
\underbar{\it iff}\ \
$\dim{~
_{{\lower2.1pt\hbox{\tenbf--}\hskip-0.25cm}} \hbox{\rm {H}} }
^{^{_{i}}}\!\!
\raise0.8pt\hbox{$_{_{\!{\mathbf{k}
\hbox{\fivebf[}\Delta\hbox{\fivebf]}_{{\bf+}\!\!}}} }$}
({\bf k} \langle  \Delta \rangle )
\le \infty
\ if\
0\le i < \dim{\bf k} \langle  \Delta \rangle ),
$
in which case
${
_{{\lower2.1pt\hbox{\tenbf--}\hskip-0.25cm}} \hbox{\rm {H}} }
^{^{_{i}}}\!\!
\raise0.8pt\hbox{$_{_{\!{\mathbf{k}
\hbox{\fivebf[}\Delta\hbox{\fivebf]}_{{\bf+}\!\!}}} }$}
({\bf k} \langle  \Delta \rangle )
\cong
{\mhatH}_{_{i\!-\!1}}\!(\Delta;{\bf k})$;
for proof
cf.\ \cite{32}
p.\ 144.
Here, ``$\dim$'' is the
\htmladdnormallink{{\it Krull }}
{http://en.wikipedia.org/wiki/Krull_dimension}
\htmladdnormallink{{\it dimension},}
{http://en.wikipedia.org/wiki/Krull_dimension}
 which for
Stanley-Reisner rings is simply ``1 + the simplicial dimension''.

\smallskip
For
$\Gamma\!\!_{_{^{1}}}, \Gamma\!\!_{_{^{2}}}$ finite and {CM}$
_{\bf k}$,
a {\it Künneth formula} for ring theoretical local cohomology follows; (``$\cdot_+$'' indicates the unique homogeneous maximal ideal i.e., the irrelevant ideal of ``$\cdot$''.)

\smallskip\noindent
{\bf 1.}\indent
$
{\underbar{H}}^{^{_{q}}}
\vbox{\moveleft0.7cm\hbox{
\lower4.5pt\hbox{$_{_{\!(
{\bf k}\langle\Gamma\!_{_{^{1}}}\!\rangle
\otimes\hskip0.05cm
{\bf k}\langle\Gamma\!_{_{^{2}}}\!\rangle
)_{{\bf+}\!\!} }}$}
}}
$
$\hskip-1.0cm
({\bf k}{\langle\Gamma\!\!_{_{^{1}}}\rangle}
\otimes
{\bf k}{\langle\Gamma\!\!_{_{^{2}}}\rangle})
\cong $
\vskip-0.15cm
$$
\ \!\cong\ \ \!\!
\Big[ {^{
\hbox{\small Motivation: {\bf k}$\langle\Gamma\!\!_{_{^{1}}}\rangle\otimes${\bf k}$\langle\Gamma\!\!_{_{^{2}}}\rangle
\cong
{\bf k}\langle\Gamma\!\!_{_{^{1}}}\!\ast \Gamma\!\!_{_{^{2}}}\rangle$}\
\hbox{\small\rm by Note~\ref{NoteP21}{\sevenbf .iii}\ p.~\pageref{NoteP21:iii}, and
} }
_{ \hbox{\small\rm
then, use Th.~\ref{TheoremP14:4} Eq.~3 p.~\pageref{TheoremP14:4}
 and  Theorem~\ref{CorP16:iToTh6}{\sevenbf .i} p.~\pageref{CorP16:iToTh6}.
} }} \Big]
\cong
$$
$$\hskip8cm
\cong
{\rlap{$\!\!_{_{_{{{i+j=q}}}}}$} {\
\raise2pt\hbox{$\bigoplus$}}}\hskip0.2cm
{\underbar{H}}^{^{_{i}}}\!
\vbox{\moveleft0.4cm\hbox{
\lower4.5pt\hbox{$_{_{\!
\hbox{\sixbf k}\langle\Gamma\!_{_{^{1}}}\!\rangle
_{{\bf+}\!\!} }}$}
}}
\hskip-0.5cm
( {\bf k}
{
\raise1pt\hbox{$\langle$}  \Gamma\!\!_{_{^{1}}}\! \raise1pt\hbox{$\rangle$}
} )
\otimes
{\underbar{H}}^{^{_{j}}}\!
\vbox{\moveleft0.4cm\hbox{
\lower4.5pt\hbox{$_{_{\!
\hbox{\sixbf k}\langle\Gamma\!_{_{^{2}}}\!\rangle
_{{\bf+}\!\!} }}$}
}}
\hskip-0.5cm
( {\bf k}
{
\raise1pt\hbox{$\langle$}  \Gamma\!\!_{_{^{2}}}\! \raise1pt\hbox{$\rangle$}
} ).
$$

\noindent
{\bf 2.}
Put:
$$ \beta_{_{\mathbf{G}}}\!(X):= \hbox{\rm inf}\{j\ \!\vert\ \!\exists
x;\ \! x\!\in\! X \land \ \! \mhatH_{j}\!
(X,X\setminus_{_{\!o}}\!\! x;{\bf G})\ne0\}.$$
For a finite $\Delta$, $\beta_{_{{\bf k}}}\!(\Delta)$ is related to
the concepts ``depth of the ring
$ {\bf k}\langle  \Delta \rangle $''
and ``C-M-ness of ${\bf k} \langle  \Delta \rangle  $'' through
$\beta_{_{{\bf k}}}\!(\Delta)=\hbox{\rm depth} ( {\bf k}\langle
\Delta \rangle  )-\!1,$
in \cite{3}
Ex.\ 5.1.23 p.\ $\!$214 and
\cite{31}
p.\ $\!$142 Ex.\ $\!$34. See also \cite{Mi} and \cite{27}.
\end{varex}

\begin{proposition}
\label{PropP22:2}

{\bf a)}
$[{\Delta}\ \hbox{\rm is ``CM}\!_{_{\bf G}}\!\!\hbox{''}]
\!\Leftrightarrow\! [\mhatH_{_{i}}({\Delta},\!\hbox{\rm
cost}\!_{_{\Delta}}\!{{\delta\!_{_{ }}}}; {\bf G})=0\ \forall\
\!\delta\!\in\!\Delta \ \hbox{\rm and} i\leq\!n\!-\!1].$

\medskip
{\bf b)}~{\rm(cf. \cite{31}
p.\ 94)}
$[{\Delta}\ \hbox{\rm is 2-``CM}_{_{\mathbf{G}}}\!\!\!\hbox{''}]
\!\Leftrightarrow\! [\mhatH_{_{i}}(\hbox{\rm
cost}\!_{_{\Delta}}\!{{\delta\!_{_{ }}}}; {\bf G})\!=\!0~\text{for all}~ \delta\!\in\overline{\hbox{\teni$\!$ V}}_{_{\!\!\!\Delta}} \ \text{and}\ \!i\!\leq\!n\!-\!1].$
\end{proposition}

\begin{proof}
Use the {\bf L{\hatH}S} (long homology sequence) with respect to
$\!({\Delta}, \hbox{\rm cost}\!_{_{\Delta}}\!{{\delta\!_{_{
}}}})$, Prop.~\ref{PropP15:1} p.~\pageref{PropP15:1}
and the fact that $\hbox{\rm
cost}\!_{_{\Delta}}\!{{\emptyset\!_{_{o}}}}\!\!=\!\emptyset$
resp.
$\hbox{\rm cost}\!_{_{\Delta}}\!{{\delta}}\!\!=\!\Delta$ if
$\delta\not\in\Delta.$
\end{proof}

\noindent Put\label{DefP23:SimplSkeleton}
$${ \Delta\!^{^{_{(_{\!}p_{\!})}}}\!\!:= \{\delta\in\Delta\ \!\vert\
\!\#\delta\le p+1\}},\quad
{\Delta\!^{^{_{p\!}}}\!\!:=
\Delta\!\!^{^{_{(_{\!}p_{\!})}}}\!\!\setminus
\Delta\!\!^{^{_{(_{\!}p\!-\!1_{\!})}}} }\!\!\!,\quad
{\Delta\!^{^{_{\prime\!}}}\!:=
\Delta\!^{^{_{(\!n{_{_{^{_{\!}}}}}\!-\!1_{\!})}}}}\!\!\!,\quad
n\!:=\!\dim\Delta.$$
So, in particular
$ {
\Delta\!^{^{_{(\!n{_{_{^{_{}}}}}\!\!_{\!})}}}\!\!=\!\Delta
}.$

\begin{theorem}
\label{TheoremP23:8}

$${\bf a)}\
\Delta\ \text{is}\ ``\text{CM}_{\mathbf{G}}\!\!\text{''}\ \Longleftrightarrow\ \
\Delta\!^{^{_\prime}}\ \text{is}\  2-\!``\text{CM}_{_{^{\mathbf{G}}}}\!\!\!\text{''}\ \text{and}\
\mhatH\hskip-0.25cm_{_{n-1}}\hskip-0.0cm(\Delta,\hbox{\rm
cost}\!_{_{^{\!\Delta}}}\!\!\delta ;{\bf G})=0\
\text{for all}\ \delta\!\in\!\Delta.$$
$${\bf b)}\
\Delta\ \text{is}\ ``\text{CM}_{\mathbf{G}}\!\!\text{''}
\Longleftrightarrow\ \ \Delta\!^{^{_\prime}}\ \text{is}\
2-\!``\text{CM}_{_{^{\mathbf{G}}}}\!\!\!\text{''}\ \text{and}\
\mhatH\hskip-0.25cm_{_{n-1}}\hskip-0.0cm(\hbox{\rm
cost}\!_{_{^{\!\Delta}}}\!\!\delta ;{\bf G})=0
\text{for all}\ \delta\!\in\overline{\hbox{\teni
V}}_{_{\!\!\!\Delta}}.\ \ \ \ $$
\end{theorem}

\begin{proof}
Proposition~\ref{PropP22:2} above
together with the fact that adding or deleting $n$-simplices does
not effect homology groups of degree $\leq
n-2$. See Proposition~\ref{PropP35:3e}.{\tenbf e} p.~\pageref{PropP35:3e}.
\end{proof}

Theorem~\ref{TheoremP23:8}
is partially deducible, using commutative algebra, from \cite{15}
pp.\ 358-360.

Our next corollary was originally, for ${\bf G}={\bf k}$, ring
theoretically proven by T. Hibi. We will essentially keep Hibi's formulation, though using:
$$
{\Delta\!^{^{\!_{o}}}\!} \!:=\!
\Delta\smallsetminus
\{\tau\!\in\!\Delta\ \!\vert\ \!\tau\!\supset\!\delta_{i} \
\hbox{\rm for\ some}\ i\!\in\!{\bf I}\}
=
\mbfcupcap{\bigcap}{\!_{i\in{\bf I}}}
{\rm cost}\!_{_{\Delta}}\!\!{{\delta}\!_{i}}\!.
$$

\begin{corollary}\label{CorP23}
{\rm(\cite{14}
Corollary\ p.\ 95-96)}
Let $\Delta$ be a pure simplicial complex of dimension $n$ and
$\{{\delta}\!_{i}\!\}\!_{_{^{i\in\lower0.5pt\hbox{\eightbf
I}}}}$, a finite set of faces in $\Delta$ satisfying
$\delta_{i}\cup\delta_{j}\!\notin\!\Delta$ for all $i\!\ne\!j$.
Set,
$\Delta\!^{^{\!_{o}}}\!:=
\mbfcupcap{\bigcap}{_{\!i\in{\bf I}}} \hbox{\rm cost} \!_{_{\Delta}}\!\!{{\delta}_{i}}. $

${\bf a)}$  If $\Delta$ is $\hbox{\rm ``CM}\!_{_{\bf G}}\!\!\hbox{''}\!$ and
$\dim\Delta\!^{^{_{\!o}}}\!<\!n$, then
$\dim\Delta\!^{^{_{\!o}}}\!=\!n\!-\!1$ and $\Delta\!^{^{_{\!o}}}$ is
$\hbox{\rm ``CM}\!_{_{\bf G}}\!\!\hbox{''}\!$.

\medskip
${\bf b)}$ If ${\bar{\hbox{\rm
st}}}\!_{_{\Delta}}\!{{\delta\!_{_{i}}}}$ is $\hbox{\rm
``CM}\!_{_{\bf G}}\!\!$'' for all ${i\!\in\!{\bf I}}$ and
$\Delta\!^{^{_{\!o}}}$ is $\hbox{\rm ``CM}\!_{_{\bf G}}\!\!\hbox{''}\!$ of
dimension $n$, then $\Delta$ is
$\hbox{\rm ``CM}\!_{_{\bf G}}\!\!\hbox{''}\!$.
\end{corollary}

\normalbaselines

\begin{proof}
\noindent{\bf a}) \noindent$
[\delta_{i}\raise0.5pt\hbox{\scriptsize$\cup$}
\delta_{j}\!\!\notin\!\Delta$ for all
$i\!\ne\!j\!\in\!\lower0.2pt\hbox{\bf I}] \!\Longleftrightarrow\!
[\hbox{\rm cost}\!_{_{\Delta}}\!\!{{\delta}\!_{i}}
\cup
\big(
\mbfcupcap{\bigcap}{\!{j\not=i}}
\hbox{\rm
cost}\!_{_{\Delta}}\!\!{{\delta}\!_{_{^{j}}}}\!\big)\!=\!\Delta]\!
\Longrightarrow
\langle
{\hbox{\rm cost}\!\!_{_{\Delta^{^{{\!{_{\!\prime}}}}}}}
\!\!{{\delta}\!_{i}}}  \cup
\big(
\mbfcupcap{\bigcap}{\!{j<i}}
\hbox{\rm
cost}\!\!_{_{\Delta^{^{{\!{_{\!\prime}}}}}}}
\!\!{{\delta}\!_{_{^{j}}}}\! \big)
\!=\! \big(\hbox{\rm cost}\!_{_{\Delta}}\!\!{{\delta}\!_{i}}
\cup \big(
\mbfcupcap{\bigcap}{\!{j<i}}
\hbox{\rm
cost}\!_{_{\Delta}}\!\!{{\delta}\!_{_{^{j}}}}\big) \big)
\cap
\Delta\!^{^{_{\!\hbox{$_{\prime}$}}}}
\!=\!\Delta\!^{^{_{\!\hbox{$_{\prime}$}}}} \rangle
\Longrightarrow
 \dim{\Delta\!^{^{\!_{o}}}}\!\!=\!n\!-\!1\Longrightarrow
$
$
{\Delta\!^{^{\!_{o}}}} \!\!=\!
{\Delta\!^{^{\!_{o}}}}\! \cap
\Delta\!^{^{_{\!\hbox{$_{\prime}$}}}} \!=\!
\big(
\mbfcupcap{\bigcap}{_{\!i\in{\bf I}}}
\hbox{\rm
cost}\!_{_{\Delta}}\!\!{{\delta}\!_{i}}\!\big)
\cap
\Delta\!^{^{_{\!\hbox{$_{\prime}$}}}}
=
\mbfcupcap{\bigcap}{_{\!i\in{\bf I}}}
{\hbox{\rm cost}\!\!_{_{\Delta^{^{{\!{_{\!\prime}}}}}}}
\!\!{{\delta}\!_{i}}}\!. $

By Th.~\ref{TheoremP23:8}
we know that $\Delta\!^{^{_{\!\hbox{$_{\prime}$}}}}$ is 2-$\!``{\rm
CM}\!_{_{\bf G}}\!\!\hbox{''}\!$, implying that $ {\hbox{\rm
cost}\!\!_{_{\Delta^{^{{\!{_{\!\prime}}}}}\!}}
\!{{\delta}}}\!_{i}$ is $``{\rm CM}\!_{_{\bf
G}}\!\!$'' for all $i\!\in\lower0.2pt\hbox{\bf I}. $

Induction using the
\hbox{\bf M-$\!$Vs} with respect to $ (\hbox{\rm
cost}\!_{_{\!\Delta^{^{{\!{_{\!\prime}}}}}}}
\!\!{{\delta}\!_{i}},
\mbfcupcap{\bigcap}{\!_{j<i}}
\hbox{\rm
cost}\!\!_{_{\Delta^{^{{\!{_{\!\prime}}}}}}}
\!\!{{\delta}\!_{_{^{j}}}}\!) $ gives $
\mhatH_{i}\!(\Delta\!^{^{_{\!o}}};{\bf G})\!=\!0$ for all $i\!<\!n\!_{_{\Delta}}\!\!\!-\!\!1. $

For links, use Prop.~\ref{PropP34:2}$a$+$b$ p.~\pageref{PropP34:2}.
E.g.,
$$ {{\rm Lk}\!\!\raise0.5pt\hbox{$_{_{{
\Delta\!\!^{^{_{\circ}}}}}}$}\!\!\delta}
=
\!{{\rm Lk}\!\!\!\!\!
\raise-0.5pt\hbox{$_{_{{ \Delta\!\!^{^{_{\circ}}}\!\!\cap
{\Delta^{^{{\!{_{\prime}}}}}}}}}$}\!\!\!\!\!\delta}
=
 {{\rm
Lk}\hskip-0.4cm\lower0.5pt\hbox{$ {_{_{\hskip-0.3cm(
\cap
\hskip-0.15cm_{_{_{{{i\in\lower0.2pt\hbox{\fivebf I}}}}}}
\hskip-0.2cm\hbox{\sevenrm
cost}\!_{_{\!\Delta}}\!\!{{\!\delta}\!_{i}}\!) }}} $}
\raise0pt\hbox{$_{_{{ \!\cap
{\Delta^{^{_{\!{\prime}}}}}}}}$}\hskip-0.4cm\delta} \ \
=
[{\rm Prop.~\ref{PropP34:2}.a\ p.~\pageref{PropP34:2}}]=
\mbfcupcap{\bigcap}{_{\!i\in{\bf I}}}\
{{\rm Lk}\!\!\!\!\!\raise0.75pt\hbox{$_{_{_{ \hbox{\sevenrm
cost}\!\!_{_{\!\Delta^{^{{\!\!{_{\prime}}}}}}}\!{{\!\delta}\!_{i}}
}}}$}\!\!\!\!\!\delta}$$
where $ {\hbox{\rm
cost}\!\!\!_{_{\Delta^{^{{\!{_{\prime}}}}}\!}}
\!{{\delta}}}\!_{i}\ \hbox{\rm and\ so}\ {{\rm
Lk}\!\!\!\!\!\raise0.75pt\hbox{$_{_{_{ \hbox{\sevenrm
cost}\!\!_{_{\!\Delta^{^{{\!\!{_{\prime}}}}}}}\!{{\!\delta}\!_{i}}
}}}$}\!\!\!\!\!\delta}\ \ \hbox{\rm is}\ ``{\rm CM}\!_{_{\bf
G}}\!\!\hbox{''}\!,\ \text{for all}\ {i\!\in\lower0.2pt\hbox{\bf I}}.
\hfill\triangleright$

\noindent{\bf b}) $ \Delta\!=\!{\Delta\!^{^{\!_{o}}}}
\mbfcupcap{\bigcup}{_{\!j\in{\bf I}}}
\overline{\hbox{\rm {st}}}\!_{_{\Delta}}\!\delta\!_{_{^{j}}}$
and
$$\overline{\hbox{\rm {st}}}\!_{_{\Delta}}\!\delta \cap
{\Delta\!^{^{\!_{o}}}} = \Big[\!\!\phantom{o}^{
\hbox{\scriptsize Motivation:  Eq. {\bf II}+{\bf III}}} _{\hbox{\scriptsize p.~\pageref{EqP34:II}ff plus that}\
\hbox{\small$\delta_{i}\cup\delta_{j}\!\notin\Delta$}}\Big] =
\overline{\hbox{\rm {st}}}\!_{_{\Delta}}\!\delta \cap {\hbox{\rm
cost}\!_{_{\Delta}}\!\!{{\delta}}} = \dot{\delta}\ast {{\rm
Lk}\!_{_{\Delta}}\!\!{{\delta}}}. $$
So, by Th.~\ref{CorP16:iToTh6}.i, p.~\pageref{CorP16:iToTh6}
$ \ \!\overline{\hbox{\rm {st}}}\!_{_{\Delta}}\!\delta \cap
{\Delta\!^{^{\!_{o}}}}\ \hbox{\rm is} $
\noindent $\hbox{\rm ``CM}\!_{_{\bf G}}\!\!\hbox{''}$ \underbar{iff} $\
\!\overline{\hbox{\rm {st}}}\!_{_{\Delta}}\!\delta$ is.
Induction, using the \hbox{\bf M-$\!$Vs}$\ \!$ with respect to
$(\overline{\hbox{\rm {st}}}\!_{_{\Delta}}\!\delta\!_{i},
{\Delta\!^{^{\!_{o}}}}
\mbfcupcap{\bigcup}{_{j<i}}
\overline{\hbox{\rm {st}}}\!_{_{\Delta}}\!\delta\!_{_{^{j}}}),$
gives $ \mhatH_{i}\!(\Delta\!^{^{_{\ }}}\!;{\bf G})\!=\!0$ for all $i\!<\!n\!_{_{\Delta}}.$
$\hbox{\rm End\ as\ in\ {\bf a}.}
$
\end{proof}

Now, we will show that our generalized definition of 2-``CM''-ness p.~\pageref{DefP22} is consistent with K. Baclawski's original definition in \cite{Baclawski} p.~295.

\begin{lemma}
\label{LemmaP23:1}
$
\Delta\ \hbox{\rm ``CM}_{_{\!\!\mathbf{G}}}\!\!\!\hbox{''}\Rightarrow
\left\{\hskip-0.2cm\begin{array}{ll}
{\bf a})\ \mhatH\!_{_{i}}(\!{\rm
cost}\!_{_{\Delta}}\!{{\delta\!_{_{ }}}}; {\bf G})=0\
\text{for all}\ \delta\ \in\!\Delta & \text{if}\ i\leq n-2,\\
{\bf b})\ \mhatH\!_{_{i}}( \!\hbox{\rm cost}\hskip-0.55cm
_{_{{\rm cost}\!\!_{_{\Delta}}\!\!\!{{{\hbox{\sixi
{\char"0E}}}\!_{_{\mathrm{2}}} } }}}
\!{{\delta_{1}}}\!; {\bf G})=0\ \text{for all}\
\delta_{1},\delta_{2}\in\!\Delta\!\!\! & \text{if}\ i\leq n-3. \\
\end{array}\right.
$
\end{lemma}

\medskip
\noindent{\bf Proof.} {\bf a}) Use
Proposition~\ref{PropP15:1} p.~\pageref{PropP15:1}, the
definition of ``CM''$\!$-ness and the {\bf L{\hatH}S} with respect
to\ $\!$ $({\Delta},\hbox{\rm cost}\!_{_{\Delta}}\!{{\delta\!_{_{
}}}})\!$,\ {which\ reads;}

\smallskip
$\lower4.5pt\hbox{${
\hskip-0.15cm
{\buildrel{\beta}_{_{^{1\ast}}} \over \rightarrow}\
{
\mhatH\!_{_{n\!\!}} ({\Delta}
,{\rm cost}\!_{_{\!\Delta\!}}\!{{\delta\!_{_{
}}}};\mathbf{G})}
\ {\buildrel{\delta}_{_{^{1\ast}}} \over \rightarrow}\
{ \mhatH\!\!\!\!_{_{n\!-\!1\!\!}} ({\rm cost}\!_{_{\!\Delta\!}}\!{{\delta\!_{_{ }}}};\mathbf{G})}
{\buildrel{\alpha}_{_{^{1\ast}}} \over \rightarrow}\ \ \!
{\mhatH\!\!\!\!_{_{n\!-\!1\!\!}} (\Delta
;\mathbf{G})}
\ \!{\buildrel{\beta}_{_{^{1\ast}}} \over \rightarrow} \ \!
\mhatH\!\!\!\!_{_{n\!-\!1\!\!}} (\Delta
,{\rm cost}\!_{_{\!\Delta\!}}\!{{\delta\!_{_{
}}}};\mathbf{G}) \ \!
{\buildrel{\delta}_{_{^{1\ast}}} \over \rightarrow}
{ \mhatH\!\!\!\!_{_{n\!-\!2\!\!}} ({\rm cost}_{_{\!\Delta\!}}\!{{\delta\!_{_{ }}}};\mathbf{G}) } \
{\buildrel{\alpha}_{_{^{1\ast}}} \over \rightarrow}\hskip-0.35cm\dots
}$}$

\hfill$\triangleright$

\medskip
\noindent$\!${\bf b}) Apply Prop.~\ref{PropP35:3} p.~\pageref{PropP35:3}.{\bf a}+{\bf b} to the
{\bf M-$\!$V$\!$s} with respect to
$\! (\!{{\rm cost}\!_{_{\Delta}}{{\delta_{1}}}}\!,
{{\rm cost}\!_{_{\Delta}}{{\delta_{2}}}}\!),$ then use {\bf a}), i.e.,

\smallskip
$
\lower0.5pt\hbox{${
\cdot\cdot=\cdot { \buildrel \beta_{_1\ast} \over \rightarrow }\
{ \mhatH_{_{\!n}}\! ({\rm cost}\!_{_{\Delta}}\!(\!{{\delta_{1}}}\cup
{{\delta_{2}}}\!);{\bf G}) } \ {\buildrel
\delta_{_1\ast} \over \rightarrow }\
\mhatH\!\!\!\!_{_{n\!-\!1}}\! (\!{\rm cost}\!\!\!\!\!
_{_{\hbox{\sixrm cost}\!\!_{_{\Delta}}\!\!{{\delta_{2} } }}}\!
\!\!{{\delta_{1}}};{\bf G}) \ {\buildrel
\alpha_{_1\ast} \over \rightarrow }\
\mhatH\!\!\!\!_{_{n\!-\!1}}\! ({\rm
cost}\!_{_{\Delta}}\!{{\delta_{1}}};{\bf G})
\oplus \mhatH\!\!\!\!_{_{n\!-\!1}}\! ({\rm
cost}\!_{_{\Delta}}\!{{\delta_{2}}};{\bf G}) \
{\buildrel \beta_{_1\ast} \over \rightarrow}}$}$

$\lower0.5pt\hbox{$
{
{\buildrel \beta_{_1\ast} \over \rightarrow }\ {
\mhatH\!\!\!\!_{_{n\!-\!1}}\! ({\rm cost}\!_{_{\Delta}}\!(\!{{\delta_{1}}}\cup{{\delta_{2}}}\!);
\mathbf{G}) } \ {\buildrel \delta_{_1\ast} \over
\rightarrow }\ \mhatH\!\!\!\!_{_{n\!-\!2}}\!
(\!{\rm cost}\!\!\!\!\! _{_{\hbox{\sixrm
cost}\!\!_{_{\Delta}}\!\!{{\delta_{2} } }}}\!
\!\!{{\delta_{1}}}\!;\mathbf{G}) \ {\buildrel
\alpha_{_1\ast} \over \rightarrow }\
\mhatH\!\!\!\!_{_{n\!-\!2}}\! ({\rm cost}\!_{_{\Delta}}\!{{\delta_{1}}};\mathbf{G}) \
\oplus\ \mhatH\!\!\!\!_{_{n\!-\!2}}\! ({\rm cost}\!_{_{\Delta}}\!{{\delta_{2}}};\mathbf{G}) \
{\buildrel \beta_{_1\ast} \over \rightarrow}\hskip-0.1cm\dots }
$}\square$

\medskip
\noindent{\bf Observation.} \label{ObservationP23}
To turn the implication in Lemma~\ref{LemmaP23:1} into an equivalence we just have
to add
$\mhatH\!_{_{i}}\! (\Delta;\!{\bf G})\!=\!0\ \hbox{\rm for}\
i\leq n\!-\!1 $, giving us the equivalence in;

\medskip
\noindent $\Delta\ \!\hbox{\rm ``CM}\!_{_{\bf G}}\!\!\!\hbox{''}
\Longleftrightarrow
\left\{\begin{array}{ll}
\!\!{\bf i)}\ \ \ \mhatH_{_{i}}\! (\Delta;{\bf G})\!=\!0 &
\text{for}\ i\!\leq\!n\!-\!1 \\
\!\!{\bf ii)}\ \ \mhatH_{_{i}}(\!\hbox{\rm
cost}\!_{_{\Delta}}\!{{\delta\!_{_{ }}}}; {\bf G})=0\
\text{for all}\ \delta\!\!\in\!\Delta & \text{for}\ i\!\leq\!n\!-\!2 \\
\!\!{\bf iii)}\ \mhatH_{_{i}}( \!\hbox{\rm cost}\!\!\!\!\!\!\!\!
_{\lower0.8pt\hbox{$ {_{\hbox{\sevenrm
cost}\!\!_{_{\Delta}}\!\!{{\delta_{2} } }}} $}}\!
\!{{\delta_{1}}}\!;{\bf G})\! =\!0\ \text{for all}\
\!\delta_{1}\!,\delta_{2}\!\!\in\!\Delta &
\text{for}\ i\!\leq\!n\!-\!3 \\
\end{array}\right.
$
\hfill$({\bf 1})$

\smallskip\noindent
The property $\mhatH\!\!\!_{_{^{n-1}}}\!(\hbox{\rm
cost}\!_{_{^{\!\Delta}}}\!\!\delta ;{\bf G})\!=0$ for all
$\delta\!\in\! \Delta$
allows one more step in the proof of Lemma~\ref{LemmaP23:1}.b, i.e.,

\medskip
\noindent [$\Delta$ is
$\hbox{\rm ``CM}\!_{_{\bf G}}\!\!\!\hbox{''}$
\indent and \indent
$\mhatH\!\!\!_{_{^{n-1}}}\!(\hbox{\rm
cost}_{_{^{\!\!\Delta}}}\!\!\delta ;{\bf G})\!=0\ \text{for all}\
\!\delta_{\!\!}\in_{\!\!}\Delta$]
$\Longleftrightarrow$
$$\hskip0.0cm\Longleftrightarrow\left\{\begin{array}{ll}
\!\!{\bf i)}\ \ \ \!\mhatH_{_{i}}\!
(\Delta;{\bf G})\!=\!0 & \text{for}\ i\!\leq\!n\!-\!1,\\
\!\!\!{\bf ii)}\ \ \mhatH_{_{i}}\!(\!\hbox{\rm
cost}\!_{_{\Delta}}\!{{\delta\!_{_{ }}}}; {\bf G})=0\ \text{for all}\ \delta\in\Delta & \text{for}\ i\!\leq\!n\!-\!1,\\
\!\!\!{\bf iii)}\ \mhatH_{_{i}}\!( \!\hbox{\rm
cost}\!\!\!\!\!\!\! _{\lower0.8pt\hbox{$ {_{\hbox{\sevenrm
cost}\!\!_{_{\Delta}}\!\!{{\delta_{2} } }}} $}}\!
\!{{\delta_{1}}}\!;{\bf G})\! =\!0\ \text{for all}\
\delta_{1}\!,\delta_{2}\!\in\!\Delta &
\text{for}\ i\!\leq\!n\!-\!2.
\end{array}\right.
\hskip2.2cm({\bf 2}\label{EQP24:2})$$

\begin{varnote} \label{NoteP24}
{\bf a}) Item {\bf i}) in Eq.\ 2 follows from {\bf ii}) and {\bf iii}) by the {\bf M-$\!$V$\!$s} above.

\noindent
{\bf b}) The l.h.s. in Eq.\ 2 is, by definition, equivalent to
$\Delta$ being 2-``CM''.
\end{varnote}

\smallskip
Since,
$$[{\Delta}\ \!\!\hbox{\rm ``CM}\!_{_{\bf G}}\!\!\hbox{''}]
\!\Longleftrightarrow\! [\mhatH_{_{i}}({\Delta},\!\hbox{\rm
cost}\!_{_{\Delta}}\!{{\delta\!_{_{ }}}};{\bf G})=0\ \text{for all}\ \!\delta\!\in\!\Delta \ \hbox{\rm and}\
\!i\!\leq\!n\!-\!1],
$$
it is obvious that item {\bf iii}) in Eq.\ 1 above, by the {\bf L{\hatH}S} with respect
to
$\!({\Delta},\!\hbox{\rm cost}\!_{_{\Delta}}\!{{\delta\!_{_{ }}}})$,
is totally superfluous as far as the equivalence is concerned but nevertheless it becomes quite useful when substituting $\hbox{\rm
cost}\!_{_{\Delta}}\!\!{{\delta\!_{_{}}}}$
for every occurrence of $\Delta$ and, using that $ \hbox{\rm
cost}\!\!\!\!\!\!\! _{\lower0.8pt\hbox{$ {_{\hbox{\sixrm
cost}\!\!_{_{\Delta}}\!\!{{\delta\!_{_{\ }}}} }} $}}
\!\!{{\delta_{1}}} = \hbox{\rm cost}\!
_{_{{\Delta}}}\!\!{{\delta\!_{_{\ }}}} $ for $ {{\delta_{1}
}}\!\!\notin \!\hbox{\rm cost}\! _{_{{\Delta}}}\!\!{{\delta\!_{_{\
}}}}\!\! $, we get:

\centerline{${\rm cost}\!_{_{\Delta}}\!{{\delta\!_{_{}}}} \
\!\hbox{\rm ``CM}\!_{_{\bf G}}\!\!\!\hbox{''}\ \text{for all}\
\!\delta\!\!\in\!\Delta\!\Longleftrightarrow
\left\{\begin{array}{ll}
\!\hskip-0.2cm{\bf i})\ \ \ \mhatH_{i}\! (\hbox{\rm
cost}\!_{_{\Delta}}\!{{\delta\!_{_{}}}} ;{\bf G})\!=\!0 &
\!\!\!i\!\leq\!n\!_{_{\delta}}\!-\!1\ \ \text{for all}\
\!\delta\!\!\in\!\Delta  \\
\!\hskip-0.2cm{\bf ii})\ \mhatH_{i}(\!\hbox{\rm cost}\hskip-0.5cm_{\lower0.8pt\hbox{$ {_{\hbox{\sixrm
cost}\!\!_{_{\Delta}}\!\!{{\delta\!_{_{\ }}}} }} $}}
\!\!{{\delta_{1}}}; {\bf G})=0&
\!\!\!i\leq\!n\!_{_{\delta}}\!-\!2\ \text{for all}\
\delta,\delta_{1}\!\!\in\! {\Delta}.\hskip-0.3cm
\end{array}\right.
$\hskip0.2cm{({\bf 3})}\label{EqP24:3}}

\medskip
\begin{varnorem}{\bf Remarks\ 1)} \label{RemarkP24:1}
From the {\bf L{\hatH}S} with respect to $\!({\Delta},\hbox{\rm
cost}\!_{_{\Delta}}\!{{\delta_{_{
}}}}),$ Theorem~\ref{CorP16:iiiToTh6}.ii p.~\pageref{CorP16:iiiToTh6}
and Note\ 1 p.~\pageref{NoteP29:1}
we conclude that:
If ${{\Delta\!}}$ is $2$-$\!$``CM$\!_{_{\bf G}}\!\!\hbox{''}$
then Note\ 1 p.~\pageref{NoteP17:1}
plus Proposition~\ref{PropP15:1} p.~\pageref{PropP15:1}
implies that $\mhatH_{n_{\Delta}}({\Delta};{\bf G})\not={\bf
0}.$
\end{varnorem}

\begin{varnorem}{\bf2}) \label{RemarkP24:2}
It is always true that
$$n\!_{_{\Delta}}\!\!-\!1\ \!\le
n_{_{\!\raise1pt\hbox{${_{{{\tau}}}}$}} } \!\le
n\!_{_{\raise1pt\hbox{${_{{{\sigma}}}}$}} } \!\le n\!_{_{\Delta}}\
\ \ \textrm{if}\ \tau\!\subset\!\sigma,$$
where
$n\!_{_{\sigma}}\!\!\!:=\! \dim\hbox{\rm
cost}_{_{\!\Delta}}\!\sigma \ \ \hbox{\rm and}\ \
n\!_{_{\Delta}}\!:=\nolinebreak\dim{{\Delta}}.$

Any $\Delta$ is representable as
$\Delta=
\mbfcupcap{\bigcup}{\!\!_{\sigma ^m \!\in\! \Delta}}
\overline {\sigma ^m}$, where $\overline {\sigma ^m}$ denotes the simplicial complex generated by the maximal\ simplex $\sigma^m\in \Delta$.

The {\it cone points} in a simplicial complex $\Delta$ are those vertices \  ${\bf v}\!\in\!V_\Delta$ that are contained in every maximal simplex of $\Delta$. Set $\delta\!\!_{_{\Delta}}\!:=\!\{{\bf v}\in\! V_\Delta\!\mid {\bf v}\ \text{cone point}\}\!\in\!\Delta.$
So, $\delta\!\!_{_{\Delta}}=
\mbfcupcap{\bigcap}{\!\!_{\sigma ^m \!\in\! \Delta}}
\overline{\sigma ^m}$.

Let the subindex $m$ in ${\sigma\!_m}$ denote that the simplex is maxi-dimensional, i.e.,
$\dim{{{{\sigma\!_m}}}}\!\!=\!n\!_{_{\Delta}}$.
\quad
Set $\bullet\ \!\bullet\!\!\mthickline\!\bullet:=
\{{\emptyset}\!_{_{^{o}}}\!,\{{\bf v}\!_{_{^{1}}}\!\},\{{\bf
v}\!_{_{^{2}}}\!\}, \{{\bf v}\!_{_{^{3}}}\!\},\{{\bf
v}\!_{_{^{2}}},{\bf v}\!_{_{^{3}}}\!\} \}$.
Now,

\medskip\noindent
$
\left\{\begin{array}{ll}
 {
\!\hskip-0.2cm[n\!_{_{{\bf v}}}\!\!\!=\!\!n\!_{_{\Delta}}
\hbox{\rm and}\ \hbox{\rm cost}_{_{\!\!\Delta}}\!\!{\bf v} \
\hbox{\rm pure}\ \text{for all}\ \!{\bf v}\!_{_{}}\! \in\!\hbox{V}\!\!_{_{\Delta}}] \!\Longleftrightarrow\! [n\!_{_{{\sigma}}}\!\!\!
=\!\!n\!_{_{\Delta}}
\hbox{\rm and}\
\hbox{\rm cost}\!_{_{\Delta}}\!{{\sigma\!_{_{ }}}}\
\hbox{\rm pure}\ \forall\ \!\emptyset\!_{_{^{o}}}\!\!\!\ne\!\!
\sigma\!\in\!\Delta] \!\Rightarrow
}\\
\noindent
\!\hskip-0.2cm\Longrightarrow [\ ^{_{\!}}\!\! \hbox{\rm
cost}_{_{\!\!\Delta}}\!\!\sigma \
(=
\bigcup\hskip-0.35cm\lower3pt\hbox{$_{_{{{{\bf v\in\sigma}}}}}$}
\hbox{\rm
cost}\!_{_{\Delta}}\!\!{{\bf v}\!_{_{^{\ }}}}) \
\hbox{\rm pure}\ \text{for all}\ \!\emptyset\!_{_{^{o}}} \! \ne
\sigma\!_{_{}}\in\Delta \not=\bullet\
\!\bullet\!\!\!\mthickline\!\bullet ] \Longrightarrow [{\Delta} \
\hbox{\rm is\ pure}] \Longrightarrow
\\
\noindent
{
\!\hskip-0.2cm\Longrightarrow\langle
[{\bf v}\!\in\!\hbox{\rm V}\!_{_{\!\Delta}}
}
\ \text{is a cone point}
]\
{ \Longleftrightarrow
[n\!_{_{{\bf v}}}\!\! =\!n\!_{_{\Delta}}\!\!-\!1]\ \!
\Longleftrightarrow [{\bf
v}\!\in\!{{\sigma\!_m}}\!\!\!\in\!\Delta\ \hbox{\rm if}\
\dim{{{{\sigma\!_m}}}}\!\!=\!n\!_{_{\Delta}}] \rangle }.
\end{array}\right.
$

\medskip
Since,
$
n\!_{_{\varphi}}:=\dim\hbox{\rm cost}_{_{\!\Delta}}\!\varphi\!
=n\!_{_{\Delta}}\!\!\!-\!1 \!\Longleftrightarrow
\emptyset\!_{_{^{o}}} \!\ne \varphi \!\subset\!
{{\sigma\!_m}}\!\!\!\in\!\Delta\ \text{for all}\
\!{\sigma\!_m}\!\!\in\Delta$ with
${\dim{{{{\sigma\!_m}}}}\!\! = n\!_{_{\Delta}}{_{\!}},}$ we
conclude that; if $\Delta$ is pure, then $\varphi$ contains only cone points.

\medskip
So,
$$ [\Delta\ \hbox{\rm pure\ and}\ n\!_{_{{\sigma}}}\!
\!=\!\dim{{\!\Delta}}\ \text{for all}\ \sigma\in\Delta] \Longleftrightarrow
[\Delta\ \hbox{\rm pure}\ \text{\rm and has no cone points}].$$

We also note that
$$[n\!_{\rm v}\! =n\!_{_{\Delta}}\!\!-\!1\ \!
\forall\ {\bf v}\!\in{V}_{_{\!\!\!\Delta}}] \Longleftrightarrow
[{\it V}_{_{\!\!\!\Delta}}\ \text{is finite and}\
\Delta=\overline{\!{\it V}}_{_{\!\!\!\Delta}}],$$
where  $\overline{\!{\it V}}_{_{\!\!\!\Delta}}:=\{\sigma\subset{V}_{_{\!\!\!\Delta}} \mid\ \! \sigma\ \text{finite}\}$, i.e., the {\it full
complex} with respect to\ ${\it V}_{_{\!\!\!\Delta}}$.
\end{varnorem}

\goodbreak
\begin{theorem} \label{TheoremP24:9}
Each one of the following two double conditions are equivalent to
``$\Delta$ is $2$-$\!``{\rm CM}\!_{_{\bf G}}\!\!\hbox{''}\hbox{''};$

\medskip\noindent
{\bf a)}\indent
$\left\{\begin{array}{ll}
{\bf i.}\ \hbox{\rm cost}\!_{_{\Delta}}\!{{\delta\!_{_{}}}}\ is
\ ``{\rm CM}\!_{_{\bf G}}\!\!\hbox{''}\!,\ \!\forall\
\delta\!_{_{}}\in\!\Delta\\
{\bf ii.}\ \ \!n\!_{_{\delta}} \!:=\! \dim\hbox{\rm
cost}\!_{_{\Delta}}\!{{\delta\!_{_{}}}}=\dim\Delta
=:n\!_{_{\Delta}}\ \!\forall\ \emptyset\!_{_{^{o}}} \!\! \ne\!
\delta\!_{_{}}\in\!\Delta.
\end{array}\right. $

\medskip\noindent
\noindent{\bf b)}\indent
$ \left\{\begin{array}{ll}
{\bf i.}\ \hbox{\rm cost}\!_{_{\Delta}}\!{{{\bf v}\!_{_{}}}}\ {\sl
is}\ ``{\rm CM}\!_{_{\bf G}}\!\!\hbox{''}\!,\ \ \forall\ {\bf
v}\!_{_{}}\in\!\hbox{\rm
V}\!_{_{\Delta}}\\
{\bf ii.} \ \{\bullet\
\!\bullet\!\!\mthickline\!\bullet\}\ne{\Delta}\ \hbox{\it has\ no\
cone\ points}.
\end{array}\right. $
\end{theorem}

\begin{proof}
With no dimension collapse in Eq.~3 above, this equation is equivalent to Eq.~2.
\end{proof}

\begin{vardef}\label{DefP24}
$ \Delta\setminus[\{{\bf v}\!_{_{^{1}}},\dots,{\bf
v}\!_{_{^{p}}}\}]:= \{ \delta\in\Delta\vert\ \! \delta\cap\{{\bf
v}\!_{_{^{1}}},\dots,{\bf v}\!_{_{^{p}}}\}=\emptyset \}. $
\quad$(\Delta\!\!\setminus\![\{{\bf v}\}] \!=\! \hbox{\rm
cost}\!_{_{\Delta}}\!\!{{\bf v}} $.)
\end{vardef}

Permutations and partitions within $\{{\bf v}\!_{_{^{1}}},\dots,{\bf
v}\!_{_{^{p}}}\}$ does not effect the result:

\begin{lemma}
\label{LemmaP24:2}
$ \Delta\!\!\setminus\![\{{\bf v}\!_{_{^{1}}}\!,\dots,\!{\bf
v}\!\!_{_{{p+q}}}\}] \!=\! (\Delta\!\!\setminus\! [\{{\bf
v}^{{_{\prime}}}\!\!\!_{_{^{1}}}\!,\dots,\!{\bf
v}^{{_{\prime}}}\!\!\!_{_{^{p}}}\}]) \!\!\setminus\! [\{{\bf
v}^{{_{\prime\prime}}}\!\!\!\!_{_{^{1}}}\!,\dots,\! {\bf
v}^{{_{\prime\prime}}}\!\!\!\!_{_{^{q}}}\}]$\ \ {and}

$$\Delta\!\setminus[\{{\bf v}\!_{_{^{1}}}\!,\dots,\!{\bf
v}\!_{_{^{p}}}\}]
=
\bigcap\hskip-0.45cm\lower4pt\hbox{$_{_{{{{i=1,p}}}}}$}
\hbox{\rm cost}\!_{_{\Delta}}\!\!{{{\bf v}}\!_{i}}\!
=
\hbox{\rm cost}\hskip-0.3cm_{_{{\rm cost
{\bf v}}\!_{_{^{2}}}\hskip-0.3cm_{_{ {^{{{\bf .}}}\!_{\!}{{\bf
.}}_{_{\!{\bf .}}}} }}\!}}\hskip0.0cm{{{\bf v}}\!_{_{^{1}}}}
\hskip-0.7cm{\lower10pt\hbox{\tiny\fiverm cost$\!\!_{_\Delta}\hskip-0.05cm{\bf v}\!_p$}}
\eqno{\square}$$
\end{lemma}

\begin{varnorem}{\bf Definition.} \label{DefP24:Alternative}
(From \cite{Baclawski} p.~295.)
For $k\in{\bf N}$,
$\Delta$ is $k$-``{\rm CM}$\!_{_{\bf G}}\!\!\hbox{''}\!$ if for every subset
$T\subset\hbox{\rm V}_{_{\!\!\!\Delta}}$ such that $\#T=k-1$, we have:
$$\left\{\begin{array}{ll}
{\bf i.}\ \Delta\setminus[T]\ is\ \hbox{``\rm CM}\!_{_{\bf
G}}\!\!\hbox{''}\!,\\
{\bf ii.}\ \
\dim\Delta\setminus[T]=\dim\Delta\ \quad(=:\!n\!_{_{\Delta}}\!\!=:\!n).
\end{array}\right.
$$
\end{varnorem}

\smallskip\hskip1.5cm
\framebox{\vbox{\hsize11.8cm
Changing ``$\#T\!=k-1$'' to ``$\#T\!<k$'' does not
alter the extension of the last definition
(Iterate in Th.~\ref{TheoremP24:9}{\bf b} mutatis mutandis.).
So, for ${\bf G}\!=\!{\bf k}$, a field, our definition of $k$-``{\rm
CM}$\!_{_{\bf G}}\!\!\hbox{''}$ is equivalent to Kenneth Baclawski's
original definition in \cite{Baclawski} p.~295.
}}

\subsection{{Stanley-Reisner rings for simplicial products are Segre products}}  \label{SubSecP25:II:3}
\begin{vardef} \label{DefP25}
The {\it Segre product} of the graded ${\bf A}$-algebras $R_1$ and $R_2$, denoted
$R=\sigma\!_{_{\bf A}}(R_1,R_2)$ or $R=\sigma(R_1,R_2)$, is defined through;\\
\centerline{$[R {\rlap{{\rlap{$]_{_p}$}{\ \ $:=$}}}{\ \ \ \ \ [}}
R_1]_{_p}{ \otimes_{_{\bf A}}}[R_2]_{_p}\!,\  \forall p\!\in\!{\bf
N}.$}
\end{vardef}

\begin{varex}  \label{ExampleP25}
{\bf 1})
The trivial Segre product, $R_1{\bf 0}R_2$, is equipped with the
trivial product,
i.e., every {product} of elements, both of which lacks ring term,
equals 0.

\smallskip\noindent
{\bf 2})
The ``canonical'' $\!$Segre product, $R_1\underline\otimes R_2$, is
equipped with a product induced by extending $(\!$ linearly and
distributively$)$ the
{\it componentwise multiplication}
on simple homogeneous \hbox{elements: If $_{\!}m_1^\prime
\!_{\!}\otimes_{\!} m_1^{\prime\prime}$}
$\in \bigl[R_1\underline\otimes R_2\bigr]_\alpha$ and $m_2^\prime
\otimes m_2^{\prime\prime} \in \bigl[R_1\underline\otimes
R_2\bigr]_\beta$ then
$$(m_1^\prime \otimes m_1^{\prime\prime}) (m_2^\prime \otimes
m_2^{\prime\prime}):= m_1^\prime m_2^\prime \otimes
m_1^{\prime\prime} m_2^{\prime\prime}\in
 \bigl[R_1\underline\otimes R_2\bigr]_{\alpha+\beta}.$$
{\bf 3})
The ``canonical'' generator-order sensitive Segre product,
$R_1{\bar\otimes} R_2$, of two graded ${{\bf k}}$-standard algebras
$R_1$ and $R_2$ presupposes the existence of a uniquely defined
partially ordered minimal set of generators for $R_1$ $(R_2)$ in
$[R_1]_{_1}\ ([R_2]_{_1})$ and is equipped with a product induced by
extending $($distributively and linearly$)$ the following operation
defined on simple homogeneous elements, each of which now are
presumed to be written, in product form, as an increasing chain of
the specified linearly ordered generators:

If $(m_{11} \otimes m_{21}) \in \bigl[R_1{\bar\otimes}
R_2\bigr]_\alpha$ and $(m_{12} \otimes m_{22}) \in \bigl[R_1{
\bar\otimes} R_2\bigr]_\beta$ then
$$(m_{11} \otimes m_{21}) (m_{12} \otimes m_{22}):= ((m_{11} m_{12}
\otimes m_{21} m_{22})\in \bigl[R_1{\bar\otimes}
R_2\bigr]_{\alpha+\beta}$$
if by ``pairwise'' permutations, $(m_{11} m_{12}, m_{21} m_{22})$ can
be made into a chain in the product ordering, and 0 otherwise. Here,
$(x,y)$ is a pair in $(m_{11} m_{12}, m_{21} m_{22})$ if $x$ occupy
the same position as $y$ counting from left to right in $m_{11}
m_{12}$ and $m_{21} m_{22}$ respectively.
\end{varex}

\begin{varnote}{1.} \label{NoteP25:1}
(\cite{32}
{\rm p.}\ 39-40$)$ Every Segre product of $R_1$ and $R_2$ is
module-isomorphic by definition and
{so}, they all have the same Hilbert series.
$\!$The {\it Hilbert series} of a graded ${\bf k}$-algebra
$R\!=\!\bigoplus\!\!_{_{^{i\ge 0}}}\!R_i\!$ is
$${\bf Hilb}_{R}(t)\!:=\!\sum_{i\ge 0}(\hbox{\rm H}(R,i))t^i\!:=\!
\sum_{i\ge 0}(\hbox{\rm dim}_{\mathbf{k}}R_i)t^i,$$
where $\dim$ is {\it Krull} {\it dimension} \rm and H stands for the
{\it The Hilbert function}.

\smallskip
\noindent \noindent {\it {2}.}
If $R_1$, $R_2$ are graded algebras finitely
generated (over ${\bf k}$) by $x_1,\dots,x_n\in[R_1]_{_1},
y_1,\dots,y_m\\ \in[R_2]_{_1}$, resp., then $R_1{\underline\otimes}
R_2$ and $R_1{\bar\otimes} R_2$ are generated by $(x_1\otimes
y_1),\dots,(x_n\otimes y_m)$, and
$$\dim R_{_1}{\bar\otimes} R_{_2}\!=\!\dim R_{_1}
{\underline\otimes}
R_{_2} \!=\dim R_{_1} \!+ \dim R_{_2}\!-\!1.$$

\smallskip\noindent
{\it {3}.}
The generator-order sensitive case covers all cases above.
In the theory of Hodge Algebras and in particularly in its
specialization to {\it Algebras with Straightening Laws}  (ASLs), the generator-order is the
\noindent
main issue, cf. \cite{3} \S7.1 p.~289ff.
and \cite{16} p.~123ff.
\end{varnote}

\cite{10}
p.\ 72 Lemma gives a reduced (Gr\"obner) basis, $C^\prime\cup D$,
for ``{\bf I}'' in $${\bf k} {\langle}
\Delta_{_{^{_{\!^{\!}}1\!}}}\!\times\Delta_{_{^{_{\!}2\!}}}
{\rangle }\
\cong\raise1pt\hbox{ ${\bf
k}[V\!\!\!_{_{^{\Delta_1}}}\!\!\!\times\!
V\!\!\!_{_{^{\Delta_2}}}\!]$} / \lower1pt\hbox{${\bf I}$}$$
\mbox{with}
$$C^\prime\!\!:=\!\{ w_{_{^{\!\lambda,\mu}}}w\!_{_{^{\nu,\xi}}}|\
\lambda\!<\!\nu\land \mu\!>\!\xi\},
w_{\!_{^{\lambda,\mu}}}\!\!:= \!(v_{\!_{^{\lambda}}}\!,
v\!_{_{^{\mu}}}\!), v_{\!_{^{\lambda}}}\!\!\in\!
V\!\!\!_{_{^{\Delta_{^{\!}1}}}}\!,
v\!_{_{^{\mu}}}\!\!\in\!V\!\!\!_{_{^{\Delta_2}}}\!\!$$ where the
subindices reflect the assumed linear ordering on the factor
simplices and with $\overline{p_i}$ as the projection down onto the $i$:th factor;

\medskip
\goodbreak
\noindent
$D:={ \Bigl\{\hbox{\eightbf w}\ \!=\
\!w_{_{^{{\!\lambda_{_{^{\!1}}}}\!\!,{{\mu_{_{^{\!1}}}}}}}} \!\!\!
\cdot\dots\cdot
w_{_{^{{\!\lambda_{_{^{\!k}}}}\!\!,{{\mu_{_{^{\!k}}}}}}}} \mid}$

\noindent
${\Bigl[\bigl[\bigl[\{\overline{p_1}(\hbox{\eightbf w})\} \propto {\Delta_1}\bigr]\
\!\land\ \! \bigl[\{\overline{p_2}(\hbox{\eightbf w})\}
\in{\Delta_2} \bigr] \ \!\land\ \! {{\lambda_1<\dots<\lambda_k}
\brack {\mu_1\le\dots\le \mu_k}}\bigr]\lor
\bigl[\bigl[ \{\overline{p_1}(\hbox{\eightbf w})\} \in{\Delta_1}
\bigr]\ \!\land
}$

\medskip
\noindent
$
\land
\ \! \bigl[\{\overline{p_2}(\hbox{\eightbf w})\} \!\propto\! {\Delta_2}\bigr]
\!\land\!{{\lambda_1\le\dots\le \lambda_k} \brack {\mu_1<\dots<
\mu_k}}\bigr] \!\lor\!
\bigl[\bigl[\{\overline{p_1}(\hbox{\eightbf w})\} \!\propto\! {\Delta_1}\bigr]
\!\land\! \bigl[\{\overline{p_2}(\hbox{\eightbf w})\} \!\propto\! {\Delta_2}\bigr] \land
{{\lambda_1<\dots<\lambda_k}\brack{\mu_1<\dots<\mu_k}}\bigr]\Bigr]\Bigr\}.
$

\medskip
Now,
$
C^\prime\cup D=\{m_\delta\mid\ \! \delta\
{{{\propto}}}\ \!
{\Delta_{1}\!
\!\times\!{\Delta_{2}}\}}$\label{DefP25:ProduktS-R}
and the identification
$ (v_{\lambda}, v_{\mu}) \leftrightarrow
v_{\lambda}\otimes v_{\mu}$
gives, see \cite{10}
Theorem 1\ p.\ $\!$71, the following graded ${{\bf k}}$-algebra
isomorphism of degree zero;
${\bf k} \langle
\Delta_{1}\times\Delta_{2}
\rangle
\cong
{\bf k} \langle
\Delta_{1} \rangle
{\bar\otimes}
{\bf k} \langle
\Delta_{2}\rangle$,
which, in the Hodge Algebra terminology, is {\it the} \underbar{\it discrete} {\it algebra with the same data} as
${\bf k} \langle
\Delta_{1}\rangle
{\underline\otimes}
{\bf k} \langle
\Delta_{2}\rangle$
cf.\ \cite{3}\
\S\ 7.1.\
If the discrete algebra is ``C-M'' or Gorenstein (defined below), so
is the original by
\cite{3}
Corollary\ 7.1.6.
Any finitely generated graded {\bf k}-algebra has a Hodge
Algebra structure, see \hbox{\cite{16}
p.\ 145.}

\subsection{Gorenstein complexes without cone points
are homology spheres} \label{SubSecP26:II:4}

\begin{definition}\label{DefP26:1}
The {\it cone points} in $\Delta$ are those vertices \  ${\bf v}\!\in\!V_\Delta$ that are contained in every maximal simplex $\sigma^m$. Set $\delta\!\!_{_{\Delta}}\!:=\!\{{\bf v}\in\! V_\Delta\!\mid {\bf v}\ \text{cone point}\}\!=
\bigcap\hskip-0.4cm_{_{_{_{_{\sigma^m \!\in\! \Delta}}}}}\!\!
\overline {\sigma ^m}$.
\end{definition}

\begin{definition}\label{DefP26:2}
\text{\rm (from \cite{31}
Th.~5.1 p.~65)}
Let $\Sigma$ be an arbitrary $($finite$)$ complex and put;
                   $\Gamma:=\hbox{\rm core}\Sigma
                   :=
                  \{\sigma\in\Sigma\ \mid\ \sigma\ \hbox{contains\ no\
                  cone\ points}\}.
                  $
\end{definition}

Then;
$\emptyset$ is not Gorenstein$_{\!_{\bf G}}\!$, and
$\emptyset\not=\Sigma$ is Gorenstein$_{\!_{\bf G}}\!$ (Gor$_{\!_{\bf
G}}\!$) if \

\medskip
$
\mhatH_{i}(|\Gamma|,|\Gamma|\setminus_{_{\!o}\!} \alpha;
              {\bf G})=$
$
\left\{\begin{array}{ll}
0 & if\ i\ne \dim\Gamma
\\
{\bf G} & if\ i= \dim\Gamma \\
\end{array}\right.
\ \ \forall\ \!\alpha\in |\Gamma|.
$

\smallskip
\begin{varnote} \label{NoteP26}
$\Sigma$ Gorenstein$_{\!_{\bf G}}\!\!$ $\Longleftrightarrow \Sigma$
finite and $|\Gamma|$ is a
homology$_{_{^{\!\hbox{\fivebf G}}}}\!$
sphere as defined in p.~\pageref{DefP17}.
Now;
${\bf v}$ is a cone point \underbar{iff} $\
\overline{{\hbox{\rm{st}}}}_{_{^{\!\Sigma\!}}}\!{\bf v}\!=\!\Sigma$
and so,
$$\overline{{\hbox{\rm{st}}}}_{_{^{\!\Sigma\!}}}\!
{{\delta}}_{_{^{\!\Sigma}}}\!
\!=\!
\Sigma= ({\rm core}\Sigma)\ast {\bar{\delta}}_{_{^{\!\Sigma}}}
\ \text{and}\
{\rm core}\Sigma
= {\rm Lk}_{_{^{\!\Sigma}}} \delta_{_{^{\!\Sigma}}}\! :=
\{\tau\in \Sigma\mid [\delta_{_{^{\!\Sigma}}}\cap \tau
=\emptyset]\land [\delta_{_{^{\!\Sigma}}}\cup \tau\in \Sigma]\}.$$
\end{varnote}

\begin{proposition}
\label{PropP26:1}
{\rm (\cite{10} p.\ 77.)} If \ ${\mathbf{G}}$ is a field ${\mathbf{k}}$ or the integers ${\mathbb{Z}}$ then,
$$\hskip2.0cm\Sigma\!_{_{^{1}}}
 \!\ast\Sigma \!_{_{^{2}}}\ \text{Gorenstein}\!_{_{^{\mathbf{G}}}}
\ \ \Longleftrightarrow\ \ \Sigma\!_{_{^{1}}},\Sigma\!_{_{^{2}}}\
\text{both Gorenstein}\!_{_{^{\mathbf{G}}}}.
\hskip1.5cm\square
$$

\end{proposition}

Gorensteinness is, unlike
``{Bbm}$_{_{\mathbf{G}}}\!\!\!\hbox{''}$-,
$\!``${CM}$_{_{\mathbf{G}}}\!\!\hbox{''}\!$-
and
2-$``${CM}$_{_{\mathbf{G}}}\!\!\hbox{''}$-ness,
triangulation sensitive and in particular, the Gorensteinness for
products is sensitive to the partial orders, assumed in the
$\hbox{\rm definition,\ given\ to}$ the vertex sets of the factors;
for $\{m_\delta\mid\ \!
\delta{\ {\propto}\
}\ \!
\Delta_{_{^{_{\!}1\!}}}\!
\!\times\!{\Delta_{_{^{_{\!}2\!}}}\!}\},
$
see p.~\pageref{DefP25:ProduktS-R}.
In \cite{10}
p.\ 80, this product is represented in the form of matrices, one for each pair $(\delta_{_{^{\!1}}},\delta_{_{^{\!2}}})$ of maximal
simplices $\delta\!_{i}\!\in\!\Delta_{i},\ i=1,2.$
It is then easily seen that a cone point must occupy the upper left
corner in each matrix or the lower right corner in each matrix. So a product ($\dim\Delta_{_{i}}\ge1$) can never have more than two cone
points. For Gorensteinness to be preserved under product the factors must have at least one cone point to preserve even $``{\rm
CM}_{\!_{\bf G}}\!\!\hbox{''}$-ness, by Theorem~\ref{CorP16:iiiToTh6}.ii
p.~\pageref{CorP16:iiiToTh6}.

Bd(core$(\Delta_{_{1}}\!\!\times\! \Delta_{_{2}}))=\emptyset$
forces each $\Delta_{_{i}}\!$ to have as many cone points as
$\Delta_{_{1}}\!\!\times\! \Delta_{_{2}}$.

\begin{proposition}
\label{PropP26:2}
{\rm (For proof see \cite{10}
p.\ 83ff.)} Let ${\mathbf{G}}$ be a field ${\mathbf{k}}$ or the integers ${\mathbb{Z}}$ and let $\Delta_{_{1}}\!,\ \!\Delta_{_{2}}$ be two arbitrary finite simplicial complexes
with $\dim\!\Delta_{_{i}}\!\geq\!1\ (i\!=\!1,2)$ and
a linear order defined on their vertex sets $V_{\Delta_1}, V_{\Delta_2}$, respectively, then
$$\Delta_{_{1}}\!\!\times\!\Delta_{_{2}}
{\rm Gor}_{\!_{\bf G}}
\Longleftrightarrow\!\Delta_{_{1}},\ \!\Delta_{_{2}}
{\it both\ Gor}_{_{\mathbf{G}}}$$
{\bf and}\quad either condition {\bf I} {\bf or} {\bf
II} holds:

\smallskip
\noindent{\bf (I)} $\Delta_{_{1}},\Delta_{_{2}}$
has one cone point each, either both minimal {\bf or} both\
maximal.

\smallskip\noindent
{\bf (II)}\ $\Delta_{_{1}},\Delta_{_{2}}$\
has two cone points each,
one minimal {\bf and} the other maximal.\ $\square$
\end{proposition}

\begin{varex}
\label{ExampleP26}
Gorensteinness is character sensitive!

Let $\Gamma:=\hbox{\rm core}\Sigma\!=\!\Sigma$ be
a 3-dimensional
$\hbox{\rm Gor}\!_{_{^{\mathbf{k}}}}\!$ complex, where ${{\bf
k}}$ is the prime field ${\mathbb{Z}}_{_{^{\!\hbox{\fivebf p}}}}\!$ of
characteristic {\it p}.
This implies, in particular, that $\Gamma$ is a
homology$\!_{_{^{{\mathbb{Z}}}}}$ $3$-manifold. Put
$\mhatH_{i}:=\mhatH_{i}(\Gamma;{\mathbb{Z}})$, then
$\mhatH_{0}\!=0$,
$\mhatH_{3}={\mathbb{Z}}$ and $\mhatH_{2}$
has no torsion by
Lemma~\ref{LemmaP28:1.i}.{\bf i} p.~\pageref{LemmaP28:1.i}.
Poincare' duality and \cite{30}
p.~244 Corollary~4 gives
$
\mhatH_{_{^{1}}}\!\!=\mhatH_{_{^{2}}}\!\oplus {\bf
T}\mhatH_1$,
where {\bf T}$\circ :=$
the torsion-submodule of $\circ$.
So $\Sigma=\Gamma$ with a pure torsion
$\mhatH_{1}={\mathbb{Z}}_{\bf p}$, say, gives us an
example  of a presumptive character sensitive Gorenstein complex.
Examples of such orientable compact combinatorial manifolds without
boundary are given by the 3-dimensional projective space ${{\mathbb{RP}}^3} $ and the lens space ${\bf L}(n,k)$ where
$\mhatH_1\!({\mathbb{RP}}^3;{\mathbb{Z}})= {\mathbb{Z}}_{\bf 2}$
and
$ \mhatH_1\!({\bf L}(n,k);{\mathbb{Z}})= {\mathbb{Z}}_{\bf n}$.

So, ${\mathbb{RP}}^3\ast\bullet$(${\bf L}(n,k)\ast\bullet)$ is
$\hbox{\rm Gor}\!{\lower2pt\hbox{\fivebf k}}$ for $\hbox{\rm char}{\bf k}\ne2$ (${\rm char}{\bf k}\ne n$)  while it is not even
Buchsbaum for $\hbox{\rm char}{\bf k}=2$ (resp. ${\rm char}{\bf k}=n$), cf. \cite{26}
p. 231-243 for details on ${\mathbb{RP}}^3$ and ${\bf L}(n,k)$ and  \cite{31}
Prop.\ 5.1 p.\ 65 or \cite{10}
p.\ 75 for Gorenstein equivalences.
A Gorenstein$_{\bf k}$ $\Delta$ is not in
general shellable, since if it were, $\Delta$ would have been
CM$_{_{^{\mathbb{Z}}}}$
but ${\bf L}(n,k)$ and ${\mathbb{RP}}^3$ is not.
Indeed, in 1958 M.E. Rudin publish a paper \cite{29}; \hbox{\it An unshellable triangulation of a tetrahedron}.

Other examples are given by Jeff Weeks' computer program ``SnapPea'' hosted at\\
\htmladdnormallink{http://geometrygames.org/SnapPea/}{http://geometrygames.org/SnapPea/}.
E.g.
$_{\!}\mhatH_{_{1}}\!({\Sigma}\!\!\!\!{{{\lower3.5pt\hbox{\fivebf
fig8}}}}\!(5,1);{\mathbb{Z}})\!=\!\nobreak {\mathbb{Z}}_{\mathbf{5}}$
for the old tutorial example of SnapPea, here denoted
$ {\Sigma}\!\!\!\!{{{\lower3.5pt\hbox{\fivebf fig8}}}}\!(5,1),$
i.e. the Dehn surgery filling with respect to $\!$a figure eight
complement with diffeomorphism kernel generated by (5,1).
$\!{\Sigma}\!\!\!\!{{{\lower3.5pt\hbox{\fivebf fig8}}}}\!(5,1)
\ast\bullet\ \hbox{\rm is}$
Gor$\!_{\mathbf{k}}$ if char{\bf k}$\ne5$ but not even {\rm
Bbm}$_{\mathbf{k}}$ if char${\bf k}=\!5$.
\end{varex}

\section{{Simplicial manifolds}} \label{section:SecP27:III} \label{SecP27:III}

\subsection{{{Definitions}}} \label{subsection:SubSecP27:III:1}\label{SubSecP27:III:1}

\normalbaselines
We will make extensive use of Proposition~\ref{PropP15:1} p.~\pageref{PropP15:1}
without explicit notification.
A topological space will be called a {\it $n$-pseudomanifold} or a
{\it quasi-$n$-manifold} if it can be triangulated into a simplicial
complex that is a {\it $n$-pseudomanifold} resp. a {\it
quasi-$n$-manifold}.

\begin{vardef}{\bf $\!\!$1.} \label{DefP27:1}
An $n$-dimensional pseudomanifold is a locally finite $n$-complex
${\Sigma}$ such that;\hfill\break
($\alpha$) ${\Sigma}$ is pure,\ \
i.e. the maximal simplices in ${\Sigma}$ are all $n$-dimensional.\\
{\bf($\beta$)} Every $(n-1)$-simplex of ${\Sigma}$ is the face of at
most two $n$-simplices.\\
{\bf($\gamma$)} If $s$ and $s'$ are $n$-simplices in ${\Sigma}$,
there is a finite sequence $s=s_0,s_1,\ldots s_m=s'$ of
$n$-simplices in ${\Sigma}$ such that $s_i\cap s_{i+1}$ is an
$(n-1)$-simplex for $0\le i<m$.

The {\it boundary}, {\rm Bd}${\Sigma}$, of an $n$-dimensional
pseudomanifold ${\Sigma}$, is the subcomplex generated by those
$(n-1)$-simplices which are faces of exactly one $n$-simplex in
${\Sigma}$.
\end{vardef}

\begin{vardef}{\bf $\!$2.}\label{DefP27:2}
$\Sigma$ is a {\it quasi-$n$-manifold} if it is an $n$-dimensional,
locally finite complex fulfilling;\\
\hbox{\tenbf($\alpha$)} {\bf $\Sigma$} is pure.
\indent\rm($\alpha$ is redundant since it follows from
$\gamma$ by Lemma~\ref{LemmaP20:4} p.~\pageref{LemmaP20:4}.
)\\
\hbox{\tenbf($\beta$)} Every $(n-1)$-simplex of {\bf $\Sigma$}
                           is the face of at most two
                            $n$-simplices.\\
\hbox{\tenbf($\gamma$)} {\rm Lk}$_{_{\Sigma}}\!\sigma$ is connected
i.e. $\mhatH _{0}({\rm Lk}_{_{\Sigma}}\!\sigma;{\mathbf{G}})=0$ for all
$\sigma\!\in\!\Sigma$, s.a. dim$\sigma\!<\!n-1~(\equiv\dim{\rm
Lk}_{_{^{\!\Sigma}}}\!\sigma\!\geq\!1$).

\smallskip
The {\it boundary with respect to {\bf G}}, denoted ${\rm
Bd}_{_{\mathbf{G}}}\!\Sigma,$
of a quasi-$n$-manifold $\Sigma$, is the set of simplices {${\rm
Bd}_{_{\mathbf{G}}}\Sigma \!:=\!\{\sigma\!\in\! \Sigma\ \mid\
\mhatH_{n}(\Sigma,\hbox{\rm
cost}_{_{\Sigma}}\sigma;{\mathbf{G}})=0\}$}, where {\bf G} is a unital module  over a commutative ring {\bf A}.
(According to $\beta$ in Definition 1 and 2, there are no other $0$-manifolds than $\bullet$ and $\bullet_{\!}\bullet$, due to the presens of the $(-1)$-dimensional simplex $\emptyset_o$.)
\end{vardef}

\begin{varnote}{$\!$1.} \label{NoteP27:1}
$\!\Sigma$ is (locally) finite $\!\Longleftrightarrow\!$
$\vert\Sigma\vert$ is (locally) compact.
If $X=|\Sigma|$ is a homo\-logy$_{_{\!\!\hbox{\fivebf R}}}\!$
$n$-manifold ($n$-hm$_{_{\!\!\hbox{\fivebf R}}}\!$) we will call
$\Sigma$ a $n$-hm$_{_{\!\!\hbox{\fivebf R}}}\!$. Now,
by Theorem \ref{TheoremP14:5} p.~\pageref{TheoremP14:5};
$\Sigma$ is a $n$-hm$_{_{\!\!\hbox{\fivebf G}}}\!$
for any {\bf R}-$_{\!}{\bf PID}$\ {\rm module}\
{\bf G}.

\smallskip
\rm \cite{30}\ p.\ 207-8 {\bf +} p.\ 277-8 treats the case
$_{\!}{\bf R}={\bf G}={\mathbb{Z}}:=\{\hbox{\rm the integers}\}_{\!}$,
A simplicial complex $\Sigma$ is called a hm$_{_{\!\!\hbox{\fivebf
G}}}$ if $\ \!|\Sigma|$ is. The $n$ in $n$-hm$_{_{\!\!\hbox{\fivebf
G}}}$ is deleted, since now $n\!=\dim\Sigma.$
From a purely technical point of view we really don't need the
``locally finiteness''-assumption, as is seen from
Theorem~\ref{CorP16:ToTh6} p.~\pageref{CorP16:ToTh6}.
\end{varnote}
\nobreak
\begin{vardef}{\bf $\!$3.} \label{DefP27:3}
Let ``manifold'' stand for pseudo-, quasi- or homology manifold.
\noindent A compact $n$-manifold, ${\mathcal S}$, is {\it
orientable}$_{_{{^{^{\!\!\hbox{\fivebf G}}}}}}\!\!\ { if}\ \!$
$\mhatH\!_{_{{n}}}\!({\mathcal S},{\rm Bd}{\mathcal S};{\mathbf{G}})\tilde ={\mathbf{G}}.$
An $n$-manifold is {\it orientable}$_{_{{^{^{\!\!\hbox{\fivebf
G}}}}}}\!\!$ if all its compact $n$-submanifolds are orientable
$\!-\!$ else, {\it non-orientable}$_{_{{^{^{\!\hbox{\fivebf
G}}}}}}$.
Orientability is left undefined for $\emptyset$.
\end{vardef}
\goodbreak

\begin{varnote}{$\!$2.} \label{NoteP27:2}
For a classical homology$_{_{\!\!{\mathbb{Z}}}}\!$ manifold
$X\!\neq\!\bullet$
{s.a. {\bf Bd}$_{_{{\mathbb{Z}}}}\!X=\emptyset$};
{\bf Bd}$_{_{{\mathbb{Z}}}}{\mathcal F}_{\wp}(X) =\emptyset$ if $X$ is compact and orientable and \hbox{{\bf Bd}$_{_{{\mathbb{Z}}}}{\mathcal F}_{\wp}(X) \!=\!$ $\{\wp\}$} else.
$\bullet$ is the only compact orientable manifold with boundary
$=\{\wp_{\!}\}\ \!(\{_{\!}\emptyset\!_{_{^{o}}}\})$.
\end{varnote}

\begin{vardef}{\bf $\!$4.} \label{DefP27:4}
$\{{\bf B}_{{\fivebf G},j}^{\Sigma}
\}_{j\in{\bf I}}$
is {\it the set of strongly connected boundary components of}
$\Sigma$ if
$\{{\bf B}_{{\mathbf{G}},j}^{\Sigma}\}_{j\in{\bf I}}$
{is the maximal strongly connected components of} ${\rm
Bd}_{\mathbf{G}}\Sigma$ from Definition~\ref{DefP20:1}
p.~\pageref{DefP20:1}
$(\Rightarrow{\bf B}_{{{\mathbf{G}},j}}^{\Sigma}$ is pure
and if $\sigma$ is a \underbar{maximal} simplex in ${\bf
B}_{j}$, then;
${\rm Lk} _{{\rm Bd}_{\mathbf{G}}\Sigma}\sigma={\rm Lk}_{{\bf B}_j} \sigma=\{\emptyset_o\})$.
\end{vardef}

\begin{varnote}{$\!$3.} \label{NoteP27:3}
$\emptyset, \{\emptyset_o\},$ and $0$-dimensional complexes
with either one, ${\bullet},$ or two, ${\bullet}{\bullet},$ vertices
are the only manifolds in dimensions $\le0$, and the
$|$1-manifolds$|$ are finite/infinite $1$-circles and (half)lines,
while [${\Sigma}$ is a quasi-2-manifold]$\ \Longleftrightarrow
[\Sigma$ is a homology$_{\mathbb{Z}}$
$2$-manifold].$\ $ Def. ${\bf 1}.\gamma$ is
paraphrased by ``$\Sigma$ is {\it strongly
connected}'' and $\bullet\bullet$-complexes, though strongly
connected, are the only non-connected manifolds.
Note also that
${\bf S}^{-1} :=\{\emptyset_o\}$ is the boundary
of the $0$-ball, $\bullet$, the double of which is the $0$-sphere, $\bullet\bullet$. Both the $(-1)$-sphere
$\{\emptyset_o\}$ and the $0$-sphere $\bullet\bullet$
has, as preferred, empty boundary.
\end{varnote}

\subsection{Auxiliaries}  \label{subsection:SubSecP28:III:2} \label{SubSecP28:III:2}

\begin{lemma}
\label{LemmaP28:1}
For a finite n-pseudomanifold $\Sigma$;

\noindent
{\bf i.}\label{LemmaP28:1.i}
{\rm $($cf. \cite{30}
p.\ 206 Ex.\ {\bf E}2.$)$}

\smallskip
$\mhatH_{n}\!(\Sigma,{\rm Bd}\Sigma;{\mathbb{Z}}) \!=\!{\mathbb{Z}}$
and $\mhatH_{n-1}\!(\Sigma,{\rm Bd}\Sigma;{\mathbb{Z}})$ has no torsion, or

$\mhatH_{n}\!(\Sigma,{\rm Bd}\Sigma;{\mathbb{Z}}) \!=\!0$ and
the torsion submodule of
$\mhatH_{n-1}\!(\Sigma,{\rm Bd}\Sigma;{\mathbb{Z}})$ is
isomorphic to ${\mathbb{Z}}_2.$

\smallskip\noindent
{\bf ii.}\label{LemmaP28:1.ii}
$\Sigma$ orientable$_{_{\mathbf{G}}}$ $\Longleftrightarrow\!$ $\Sigma$
orientable$_{_{{\mathbb{Z}}}}$ or ${\rm Tor}_1^{^{_{\mathbb{Z}}}}\bigl({\mathbb{Z}}_2,{\mathbf{G}}
\bigr)={\mathbf{G}}.$

So, any complex $\Sigma$ is orientable with respect to ${\mathbb{Z}}_{2}.$

\smallskip
\noindent
{\bf iii.}\label{LemmaP28:1.iii}
$\mhatH^{{n}}\!(\Sigma,{\rm Bd}\Sigma_{_{}};_{\!}{\mathbb{Z}})
\!=\!{\mathbb{Z}}\ ({\mathbb{Z}}_{2}\!)$ when $\Sigma$ is $({ non}\hbox{-}\!){
orientable}.$
\end{lemma}

\begin{proof}
\noindent
By conditions $\alpha$, $\beta$ in Def. $\!$1 p.~\pageref{DefP27:1},
a possible relative n-cycle in ${C}_{{n}}^{{o}}\!(\Sigma,{\rm
Bd}\Sigma_{_{}};{\mathbb{Z}}_{\bf m})$ must include
all oriented n-simplices all of which with coefficients of one and the same value.
When the boundary function is applied to such a possible relative n-cycle the result is an $(n-1)$-chain that includes all oriented $(n-1)$-simplices, not supported by the boundary, all of which with coefficients 0 or
${\rlap{\raise2pt\hbox{\fiverm+}}{_{\!\!}\lower2pt\hbox{\ -}}}2c\in{\mathbb{Z}}_m.$

So, $\mhatH_{{n}}\!(\Sigma,{\rm Bd}\Sigma_{_{}};{{\mathbb{Z}}_2}) \!=\!{\mathbb{Z}}_2$ --- always, and
$$\mhatH_{{n}}\!(\Sigma,{\rm Bd}\Sigma;{\mathbb{Z}}) \!={\mathbb{Z}}\ (0)
\ \Longleftrightarrow\
\mhatH_{{n}}\!(\Sigma,{\rm Bd}\Sigma ;{\mathbb{Z}}_m) \!=\!{\mathbb{Z}}_m\
(0) {\rm if}\ {\bf m}\ne2.$$
The Universal Coefficient Theorem (=
Th.\ref{TheoremP14:5} p.~\pageref{TheoremP14:5}) now gives;

\medskip
$
{\mathbb{Z}}_2=
\mhatH_{_{n\!}} (\Sigma,{\rm Bd}\Sigma_{_{}};{\mathbb{Z}}_2)
\ \! {{{_{\mathbb{Z}}}}\atop{{\raise2pt\hbox{$\cong$}} }} \ \!
\mhatH_{{n\!}} (\Sigma,{\rm Bd}\Sigma; {\mathbb{Z}}\otimes\!_{_{\mathbb{Z}}}\!{\mathbb{Z}}_2)
\ \! {{{_{\mathbb{Z}}}}\atop{{\raise2pt\hbox{$\cong$}} }} \ \!
$

\smallskip
\hfill
{$ \ \! {{{_{\mathbb{Z}}}}\atop{{\raise2pt\hbox{$\cong$}} }} \ \!
\big[\mhatH_{_{n\!}}(\Sigma,{\rm Bd}\Sigma_{_{}};{\mathbb{Z}})
\otimes\!_{_{\mathbb{Z}}}\!{\mathbb{Z}}_{_{^{2}}}\big]
\oplus
\hbox{\rm Tor}^{\mathbb{Z}}_1\!
\big(\mhatH_{_{{n\!-\!1}}}\! (\Sigma,{\rm Bd}\Sigma_{_{}};{\mathbb{Z}}), {{\mathbb{Z}}}_{_{^{2}}})
$}
\indent and,

\medskip
$
\mhatH_{_{n\!}} (\Sigma,{\rm Bd}\Sigma_{_{}};{\mathbb{Z}}_{_{^{\mathbf{m}}}})
\ \! {{{_{\mathbb{Z}}}}\atop{{\raise2pt\hbox{$\cong$}} }} \ \!
$
{$ \mhatH_{_{n\!}} (\Sigma,{\rm Bd}\Sigma_{_{}}; {\mathbb{Z}}\otimes\!_{_{\mathbb{Z}}}\!{\mathbb{Z}}_{_{^{\mathbf{ m}}}})
\ \! {{{_{\mathbb{Z}}}}\atop{{\raise2pt\hbox{$\cong$}} }} \ \!
$}

\smallskip
\hfill
{$ \ \! {{{_{\mathbb{Z}}}}\atop{{\raise2pt\hbox{$\cong$}} }} \ \!
\big[\mhatH_{_{n\!}}(\Sigma,{\rm Bd}\Sigma_{_{}};{\mathbb{Z}})
\otimes\!_{_{\mathbb{Z}}}\!{\mathbb{Z}}_{_{^{\mathbf{m}}}}\big]
\oplus
\hbox{\rm Tor}^{\mathbb{Z}}_1\!
\big(\mhatH_{_{{n\!-\!1}}}\! (\Sigma,{\rm Bd}\Sigma_{_{}};{\mathbb{Z}}), {{\mathbb{Z}}}_{_{^{\mathbf{m}}}}), \indent\indent
$} \noindent

\smallskip
\noindent where the last homology module in each formula,
by \cite{30}
p.\ 225 Cor.\ 11, can be substituted by its torsion submodule. Since
$\Sigma$ is finite, $\mhatH_{_{{n\!-\!1}}}\! (\Sigma,{\rm
Bd}\Sigma_{_{}};{\mathbb{Z}}) ={\bf C}_{_{^{_{\!}1}\!}}\oplus {\bf
C}_{_{^{_{\!}2}\!}}\oplus...\oplus {\bf C}_{_{^{\hbox{\rm s}}\!}} $
by The { Structure Theorem} for Finitely Generated Modules over {\bf
PID}s,
cf. \cite{30}
p.\ 9.
Now, a simple check, using \cite{30}
p.\ 221 Example\ 4, gives {\bf i}, which gives {\bf iii}
by \cite{30}
p.\ 244 Corollary.
Theorem~\ref{TheoremP14:5} p.~\pageref{TheoremP14:5}
and {\bf i} implies {\bf ii}.
\end{proof}

Proposition~\ref{PropP22:1} p.~\pageref{PropP22:1}
together with Proposition~\ref{PropP15:1} p.~\pageref{PropP15:1}
gives the next Lemma.

\begin{lemma}
\label{LemmaP28:2}
{\bf i.}\label{LemmaP28:2.i}
$\!\!\Sigma$ is a $n$-hm$_{_{^{\!\hbox{\fivebf G}}}}\!$
$_{\!}{{_{\Longrightarrow}}\atop{^{\
\!\not\!\!\!\Longleftarrow}}}\!$  $\Sigma$ is a quasi-$n$-manifold$\
_{\!}{{_{\Longrightarrow}}\atop{^{\
\!\not\!\!\!\Longleftarrow}}}\Sigma$ is an $n$-pseudomanifold.

\medskip
\noindent {\bf ii.}\label{LemmaP28:2.ii}
$\!\Sigma$ is a $n$-hm$_{_{^{\!\hbox{\fivebf G}}}}\!$ \underbar{iff}
it is a $\!$``Bbm$_{_{\!^{\mathbf{ G}\!\!\!}}}\hbox{''}\!$
pseudomanifold
$ {and}\
\mhatH\!\!\!\!\!\!\lower1pt\hbox{$_{_{^{\!{n\lower1pt\hbox{-}\hbox{\fivebf\#}\sigma\!}}}}$}
\!\!\!\!\raise1pt\hbox{\eightbf(}{\rm
Lk}\!_{_{_{{{\Sigma}}}}}\!\!\!\sigma;\!{\mathbf{G}}\raise1pt\hbox{\eightbf)}\!=_{\!}0\ _{\!} {or}\ \! {\mathbf{G}}\
\!\forall\ \!\sigma{\!}\not=\!\emptyset_{_{^{\!o}}}\!.
\hfill\square_{\!}$
\end{lemma}

Def.~\ref{DefP20:1} p.~\pageref{DefP20:1}.
makes perfect sense even for non-simplicial posets like $\Gamma
\smallsetminus
\Delta$ (Def.~\ref{DefP20:2} p.~\pageref{DefP20:2})
which allow us to say that
$\Gamma{\raise1.5pt\hbox{\eightmsbm\char"72}}\Delta$
is or is not {\it strongly connected} ({\it as a poset}) depending on whether
$\Gamma{\raise1.5pt\hbox{\eightmsbm\char"72}}\Delta\!^{^{\!_{\ }}}$
fulfills Def.~\ref{DefP20:1} p.~\pageref{DefP20:1}
or not. Now, for quasi-manifolds $\Gamma{\raise1.5pt\hbox{\eightmsbm\char"72}}\Delta\!^{^{\!_{\ }}}$
{\it connected as a poset}
(Def.~\ref{DefP20:1} p.~\pageref{DefP20:2})
is equivalent to
$\Gamma{\raise1.5pt\hbox{\eightmsbm\char"72}}\Delta\!^{^{\!_{\ }}}$
{\it strongly connected} which is a simple consequence of
Lemma~\ref{LemmaP20:3} p.~\pageref{LemmaP20:3}
and the definition of quasi-manifolds, cf. \cite{12}
p.\ 165 Lemma\ 4 \& 5.
$\Sigma\smallsetminus \hbox{\rm
cost}_{_{\!\Sigma}}\!\sigma$ is connected as a poset for any simplicial complex $\Sigma$ and any
$\sigma\ne\emptyset_o$. $\Sigma\smallsetminus {\rm Bd}{{\Sigma}}$ is strongly connected for any pseudomanifold $\Sigma$.
$\Delta\!^{{p}}$ is defined at the end of p.~\pageref{DefP23:SimplSkeleton}.

If $\Delta\subset\Gamma\subset\Sigma$ the {\bf L{\hatH}S} with respect to any simplicial complex $\Sigma$, with $\dim\Sigma=n$ reads;
$$
\dots0\rightarrow\mhatH_{n}(\Gamma,\Delta;{\mathbf{G}})\rightarrow
\mhatH_{n}(\Sigma,\Delta;{\mathbf{G}})\rightarrow
 \mhatH_{n}(\Sigma,\Gamma;{\mathbf{G}})\rightarrow\dots
 \eqno{({\bf1})}
$$
With $\sigma\subset\tau$ then $\Delta\!:=\hbox{\rm cost}\!_{_{\Sigma}}\sigma\subset
\Gamma\!:=\hbox{\rm
cost}\!_{_{\Sigma\!}}\tau\!$ we get {\bf a}
in the next Corollary and (1) becomes:

\smallskip\noindent
$\dots{\bf0}{
{ \buildrel \beta_{_1\ast} \over \longrightarrow }\
\hat{\mathbf{H}}\hskip-0.05cm_{{n}}\! ({\rm
cost}\!_{_{\Sigma}}\!{{\tau}},{\rm
cost}\!_{_{\Sigma}}\!{{\sigma}};{\bf G})  \ {\buildrel
\delta_{_1\ast} \over \longrightarrow }\
\hat{\mathbf{H}}\hskip-0.05cm_{{n}}\! (\Sigma,{\rm
cost}\!_{_{\Sigma}}\!{{\sigma}};{\bf G})
 \ {\buildrel
\alpha_{_1\ast} \over \longrightarrow }\
\hat{\mathbf{H}}\hskip-0.05cm_{{n}}\! (\Sigma,{\rm
cost}\!_{_{\Sigma}}\!{{\tau}};{\bf G}) \
{\buildrel \beta_{_1\ast} \over \longrightarrow}}\dots$,\\ which,  after excising everything outside $\overline{{\mathrm{st}}}\!\!
\lower1.1pt\hbox{${_{_{\Sigma}}}$}\!\sigma$ in the first homology group, becomes: \\
\smallskip\noindent
$\dots{\bf0}{
{ \buildrel \beta_{_1\ast} \over \longrightarrow }\
\hat{\mathbf{H}}_{n}\!
(\mathrm{cost}\hskip-0.2cm_{_{\hskip-0.1cm\lower1.5pt\hbox{${
\over{^{\mathrm{st}}}}$}\hskip-0.15cm{\lower0.4pt\hbox{${
_{_{_{_{\Sigma}}\!\!\sigma}}}$}}}}\hskip-0.2cm{{\tau}},{\rm
cost}\hskip-0.2cm_{_{\hskip-0.1cm\lower1.5pt\hbox{${
\over{^{\mathrm{st}}}}$}\hskip-0.15cm{\lower0.4pt\hbox{${
_{_{_{_{\Sigma}}\!\!\sigma}}}$}}}}\hskip-0.2cm{{\sigma}};{\bf G})
\ {\buildrel
\delta_{_1\ast} \over \longrightarrow }\
\hat{\mathbf{H}}\hskip-0.05cm_{{n}}\! (\Sigma,{\rm
cost}\!_{_{\Sigma}}\!{{\sigma}};{\bf G})
 \ {\buildrel
\alpha_{_1\ast} \over \longrightarrow }\
\hat{\mathbf{H}}\hskip-0.05cm_{{n}}\! (\Sigma,{\rm
cost}\!_{_{\Sigma}}\!{{\tau}};{\bf G}) \
{\buildrel \beta_{_1\ast} \over \longrightarrow}}\dots$
\hfill{({\bf2})}

Next theorem is given in {\rm \cite{12}
p.\ 166} (3.2 Hauptlemma) with the integers ${\mathbb{Z}}$ as coefficients.
\begin{theorem} \label{TheoremP28:10} {\rm(For proof see \cite{12} p.\ 166 plus Th.~\ref{TheoremP14:5}
p.~\pageref{TheoremP14:5})} Let ${\rm Tor}^{\!{\mathbb{Z}}}_1\! \big(\mhatH_{_{{i-1}}}\! (\Gamma,\Delta;{\mathbb{Z}}), {\mathbf{G}}\big)=0$.
For a finite quasi-$n$-manifold $\Sigma$ with subcomplexes
$\Delta\!\subset\Gamma \mmysubsetneqq \Sigma$ and
$\dim\Gamma\!=\!\dim\Sigma=n$,
then $\mhatH_{_{{n}}}\! (\Gamma,\Delta;{\mathbf{G}})=0$ if
$\Sigma\smallsetminus \Delta\!^{^{\!_{\ }}}$ is strongly connected, which obviously is equivalent to
$\mhatH_{n}(\Sigma,\Delta;{\mathbf{G}})\rightarrow
 \mhatH_{n}(\Sigma,\Gamma;{\mathbf{G}})$ is injective
in the relative {\bf L{\hatH}S} with respect to
$\!(\Sigma,\Gamma,\Delta)$. \hfill$\square$
\end{theorem}

With $\Gamma\!=\hbox{\rm cost}\!_{_{\Sigma}}\tau\ \hbox{\rm and}\
\Delta\!=\hbox{\rm
cost}\!_{_{\Sigma\!}}\emptyset\!_{_{^{o}}}\!\!=\emptyset$ resp.
$\hbox{\rm cost}\!\!\!\! {\lower1.0pt\hbox{${ _{_{{\!\!\hbox{\fiverm
Bd}{{{\Sigma}} }}}}\! }$}} \!\tau$ we get {\bf b}
in the next Corollary.

\begin{corollary} \label{corollary:CorP29:1} \label{CorP29:1}
{\rm (${\mathbf{G}}$ as in Th.~\ref{TheoremP28:10}.)}
If $\Sigma$ is a finite quasi-$n$-manifold and
$\sigma\mmysubsetneqq \tau\!\in\!\Sigma$, then;

\smallskip
\noindent ${\bf a.}\label{CorP29:1.a}\ \
\mhatH_{n}(\Sigma,\hbox{\rm
cost}\!\!_{_{_{^{\Sigma}}}}\!\!\sigma;{\mathbf{G}}) \longrightarrow
\mhatH_{n}(\Sigma,\hbox{\rm
cost}\!\!_{_{_{^{\Sigma}}}}\!\!\tau;{\mathbf{G}})\!$ is injective if
$\Sigma\smallsetminus
\hbox{\rm cost}\!\!_{_{_{^{\Sigma}}}}\!\!\sigma$ is strongly connected.

\medskip\noindent
{\bf b.}\label{CorP29:1.b}
$\mhatH_{n}(|\Sigma|\setminus_o\alpha;{\mathbf{G}})
\ \raise3pt\hbox{${{{\lower3pt\hbox{\fivebf A}}}}\atop{\cong}$}\
\mhatH_{n}\!(\hbox{\rm cost}{\lower1.0pt\hbox{${
_{_{{\!\!\hbox{\rm }\!\!_{_{_{\ }}}\!\!\!{{{\Sigma}} }}}}\!
}$}}\tau;{\mathbf{G}}) = \mhatH_{n}\!(\hbox{\rm
cost}{\lower1.0pt\hbox{${ _{_{{\!\!\hbox{\rm }\!\!_{_{_{\
}}}\!\!\!{{{\Sigma}} }}}}\! }$}}\tau,\hbox{\rm cost}\!\!\!\!
{\lower1.0pt\hbox{${ _{_{{\!\!\hbox{\fiverm Bd}\!\!_{_{_{\
}}}\!\!\!{{{\Sigma}} }}}}\! }$}} \!\tau;{\mathbf{G}}) =0$, if
$\alpha\in\hbox{\rm Int}(\tau)$ for any $ \tau\in \Sigma.$
\rm (Cf. the proof of Proposition~\ref{PropP15:1} p.~\pageref{PropP15:1}.)
\quad
\end{corollary}

\begin{varnote}{1.} \label{NoteP29:1}
{\bf a}\label{NoteP29:1.a}
above, implies that the boundary of a quasi-$n$-manifold is a subcomplex and
\ {\bf b}\label{NoteP29:1.b}
that simplicial quasi-manifolds are ``ordinary'', i.e. $\mhatH_{n}({\rm
cost}\!\!_{_{_{^{\Sigma}}}}\!\!\delta;{\mathbf{G}})\!=0$ as defined in p.~\pageref{DefP16:ordinary}
and so, for quasi-manifolds,
${\rm Bd}_{_{\!{\mathbf{G}}}} \!\Sigma\not=\emptyset
\Longleftrightarrow
\mhatH_{n}(\Sigma;{\mathbf{G}})=0,$ since $\mhatH_{n}(\Sigma;{\mathbf{G}}) =\mhatH_{n}(\Sigma,{\rm
cost}\!\!_{_{_{^{\Sigma}}}}\!\!\delta;{\mathbf{G}})\!=0$
\underbar{iff} $\delta\in{\rm Bd}_{_{\!{\mathbf{G}}}}
\!\Sigma$ by the {\bf LHS}.
Also recall Proposition~\ref{PropP15:1} p.~\pageref{PropP15:1}.
\end{varnote}

\begin{corollary} \label{corollary:CorP29:2} \label{CorP29:2}
\noindent
{\bf i.}\label{CorP29:2.i}
$_{\!}$If $\Sigma$ is a finite quasi-$n$-manifold with $\#{\bf
I}\!\geq\!2$ then;\\
{\hfill$\mhatH\!_{_{{{n\!\!}}}} ({\Sigma},
\mbfcupcap{\cup}{j\not=i}
{\bf B}_{_{\!\!{j\!}}} ;\!{\mathbf{G}}) \!=\!0 $
and
$ \mhatH\!_{_{{{n\!\!}}}}({\Sigma},{\bf B}_{_{\!{i}}\!};\!{\mathbf{G}}) \!=\!0 $,
{\rm (with}
${\bf B}_{_{\!\!{j\!}}}:={\bf B}_{_{\hskip-0.1cm{\mathbf{G}\!,_{\!}j}}}^{^{_{\Sigma\!\!}}}$ {\rm from
Definition~4 p.~\pageref{DefP27:4}).}\hfill}

\smallskip\noindent
{\bf ii.} \label{CorP29:2.ii}
Both $ \mhatH\!_{_{{{_{\!}n\!\!}}}}(\Sigma;_{\!}{\mathbb{Z}})$ and
$\mhatH\!_{_{{{_{\!}n\!\!}}}}(\Sigma,{\rm Bd}\Sigma;_{\!}
{\mathbb{Z}})$ equals $0$ or ${\mathbb{Z}}$.
\end{corollary}

\begin{proof}
{\bf i.}\ \
$
\mbfcupcap{\cup}{j\not=i}
{\bf B}_{j} \mmysubsetneqq \hbox{\rm
cost}_{_{_{\!{\Sigma}}}}\!\!\sigma\ \hbox{\rm for\ some}\ {\bf
B}_{_{\!{i}}}\!\hbox{-}\hbox{\rm maxidimensional}\ \sigma\!\in\!{\bf
B}_{_{\!{i}}}$ and vice versa. \hfill$\triangleright$

\medskip\noindent
{\bf ii.}
$\!\![\dim\!\tau\!=\!n] \Rightarrow\!
[[\mhatH\!_{_{{{n\!\!}}}}(\Sigma,\hbox{\rm
cost}\!_{_{\Sigma}}\!\tau;{\mathbf{G}})\!=\! {\mathbf{G}}] \land [{\rm
Bd}\Sigma\subset\hbox{\rm cost}\!_{_{\Sigma}}\!\tau]]\ \hbox{\rm
and} \ \!\Sigma\!\setminus\! {\rm Bd}\Sigma$ strongly con\-nected.
Theorem~\ref{TheoremP28:10}  and the {\bf L{\hatH}S} gives the
injections $ \mhatH\!_{_{{{_{\!}n\!\!}}}}(\Sigma;_{\!}{\mathbf{G}})\!\hookrightarrow\! \mhatH\!_{_{{{_{\!}n\!\!}}}}(\Sigma,{\rm
Bd}\Sigma;_{\!}{\mathbf{G}})\!\hookrightarrow\!{\mathbf{G}}.$
\end{proof}

\begin{varnote}{2.} \label{NoteP29:2}
$\!\Sigma_{_{^{\!{1\!}}}}\!\ast_{\!}\Sigma_{_{^{\!{2\!}}}}\!\not\subset\!\Sigma_{_{^{\!{i\!}}}}
$ locally finite
$\Longleftrightarrow$
$\!\Sigma_{_{^{\!{1\!}}}}$, $\!\Sigma_{_{^{\!{2\!}}}}\!$
both finite, and;
$\Sigma$ Gorenstein $\Longrightarrow \Sigma$ finite.
\end{varnote}

\begin{theorem} \label{theorem:TheoremP29:11} \label{TheoremP29:11}
{\bf i.a} \label{TheoremP29:11.i.a}
{\rm(cf. \cite{12}
p.\ 168, \cite{13}
p.\ 32.)} Let ${\mathbf{G}}$ denote a field ${\mathbf{k}}$ or ${\mathbb{Z}}$. $\Sigma$ is a quasi-manifold \underbar{iff}\ $\Sigma\!=\!
\bullet \bullet$ or $\Sigma$ is connected and ${\rm
Lk}\!_{_{\Sigma}}\!\sigma$ is a finite quasi-manifold for all\ \
$\emptyset\neq\sigma\in\Sigma$.

\noindent
{\bf i.b.}\label{TheoremP29:11.i.b}
$\Sigma$ is a homology$_{_{\mathbf{G}}}\!\ n $-manifold
\underbar{\hbox{\rm iff}}\ $\Sigma\!=\!\bullet{\!}\bullet\ {or}$
$\mhatH_{_{{\!{0\!}}}} ({{{{{\Sigma}}}}};\!{\mathbf{G}})\!=\mathbf{0}$\
and ${\rm Lk}\!_{_{\Sigma}}\!\sigma$ is a finite \hbox{\rm
CM}$\!_{_{\mathbf{G}}} \!$-homology$_{_{\mathbf{G}}}\!$
$(n-\#\sigma)$-manifold$\ \ \forall\ \
\emptyset\!_{_{^{o}}}\!\neq\!\sigma\!\in\!\Sigma$.

\smallskip
\noindent {\bf ii.}\label{TheoremP29:11.ii}
\ \ $\Sigma$ is a {quasi-manifold} $\Longrightarrow$ $
{\rm Bd}
\! _{_{_{\mathbf{G}}}}\!({\rm Lk}\!_{_{_{\Sigma}}}\!\sigma\!)\!
=
{\rm Lk}\!\!\!\!\! _{_{_{\mathrm{Bd}\!_{_{{\mathbf{G}}}}\!\!\!\Sigma}}}\!\!\! \sigma
\ \text{if}\ \sigma\in{\rm Bd}\! \raise0.5pt\hbox{$_{_{_{\mathbf{G}}}}$}\!\!\!\Sigma \ \text{and}\
{\rm Bd}_{_{_{\mathbf{G}}}}\!
_{\!}({\rm Lk}\!_{_{_{\Sigma}}}\!\!\sigma\!)\!
\equiv\emptyset\ \text{else}.
$

\smallskip\noindent
{\bf iii.} {$\Sigma$ is a quasi-manifold\
$\Longleftrightarrow$ ${\rm Lk}\!_{_{\Sigma}}\!\sigma$ is a
pseudomanifold} $\forall\  \sigma\in\Sigma,$  including
$\sigma=\emptyset\!_{{o}}.$
\end{theorem}

\begin{proof}
{\bf i.}
A simple check confirms all our claims for $\dim\Sigma\leq1$, cf.
Note\ 3 p.~\pageref{NoteP27:3}.

So, assume $\dim\Sigma\geq2$ and note that
\hbox{
$\!\sigma\!\in\!{\rm Bd}\! _{_{_{\mathbf{G}}}}\!\Sigma $
\ \underbar{iff}\ \
${\rm Bd}
\! \raise0.5pt\hbox{$_{_{_{\mathbf{G}}}}$}\!\!
(_{\!}{\rm Lk}
_{_{_{\Sigma}}}\sigma)\!\not=\!\emptyset$.}

\smallskip
\noindent {\bf i.a.} $\!(\!\Leftarrow\!$)
That ${\rm Lk}\!_{_{\Sigma}}\!\sigma$, with $\dim{\rm
Lk}\!_{_{\Sigma}}\!\sigma\!=\!0$, is a quasi-0-manifold implies
definition condition 2$\beta$
p.~\pageref{DefP27:2}
and since the other ``links'' are all \underbar{connected} condition
2$\gamma$ follows.\hfill $\triangleright$

\smallskip\noindent
$(\Rightarrow$) Definition condition 2$\beta$ p.~\pageref{DefP27:2}
implies that 0-dimensional links are $\bullet$ or $\bullet \bullet$
while Eq.\ {\bf I} p.~\pageref{EqP34:I}  gives the necessary
connectedness of
`links of links'${{\!}}$, cf.  Lemma~\ref{LemmaP20:4} p.~\pageref{LemmaP20:4}.
\hfill$\triangleright$

\smallskip
\noindent{\bf i.b.}
Lemma~\ref{LemmaP28:2.ii}.\ $\!${\bf ii} above plus
Proposition\ \ref{PropP22:1} p.~\pageref{PropP22:1}
and Eq.~{\bf I} p.~\pageref{EqP34:I}.
\hfill$\triangleright$

\medskip\noindent
{\bf ii.} Pureness is a local property, i.e. $\Sigma$\
{pure}\ $\Longrightarrow$ {\rm Lk}$\!_{_{\Sigma}}\!\sigma$ pure. Put
$n:=\dim{\Sigma}$.

Now;
$$\varepsilon\!\in\!\hbox{\rm Bd$\!\!$\raise0.5pt\hbox{$_{_{_{\mathbf{G}}}}$}$\!\!$ (Lk}\!_{_{\Sigma}}\!\sigma) \!\Leftrightarrow\! 0\!=\!
\mhatH\hskip-0.6cm_{_{_{\!{n\lower1pt\hbox{-}\hbox{\fivebf\#}\sigma
\lower1pt\hbox{-}\hbox{\fivebf\#}\varepsilon}}}} \hskip-0.3cm({\rm
Lk}\!\!\!\!_{_{_{\mathbf{Lk}\!_{_{\Sigma}}\!\!\sigma}}}\!\!\!\! \varepsilon;{\mathbf{G}}) \!=\!
\big[{^{\hbox{\spaceskip2.7pt
\sevenrm Eq.~{\sevenbf I} p.~\pageref{EqP34:I}
}} _{ \varepsilon\in\hbox{\sevenrm
Lk}\!_{_{\Sigma}}\!\sigma}}\hskip-0.1cm\big] \!=\!
\mhatH\hskip-0.4cm_{_{_{\!{n\lower1pt\hbox{-}
\hbox{\fivebf\#}(\sigma\cup\varepsilon)}}}} \hskip-0.4cm{\bf (}{\rm
Lk}\!\!_{_{_{{{\Sigma}}}}}\!\! (\sigma\cup\varepsilon);{\mathbf{G}}{\bf
)}\ \& \ \varepsilon\!\in\!{\rm Lk}\!_{_{\Sigma}}\!\sigma.$$
So;
$$\varepsilon\!\in\!\hbox{\rm Bd$\!\!$\raise0.5pt\hbox{$_{_{_{\mathbf{G}}}}$}$\!\!$ (Lk}\!_{_{\Sigma}}\!\sigma)
\!\Longleftrightarrow\!
[\sigma\cup\varepsilon \!\in\!{\rm Bd}\!\raise0.5pt\hbox{$_{_{_{\mathbf{G}}}}$}\!\!{{\Sigma}}\ _{\!}\hbox{\rm and}\ \varepsilon\!\in\!{\rm
Lk}\!_{_{\Sigma}}\!\sigma]
\!\Longleftrightarrow\!
$$
$$
\!\Longleftrightarrow\!
[\sigma_{\!}\cup_{\!}\varepsilon \!\in\!{\rm
Bd}\!\raise0.5pt\hbox{$_{_{_{\mathbf{G}}}}$}\!\!{{\Sigma}}\ \hbox{\rm
and}\ \sigma\cap\varepsilon_{\!}=_{\!}\emptyset]
\!\Longleftrightarrow\!
[\varepsilon\!\in\!{\rm Lk}\!\!\!\!\!
\raise0.5pt\hbox{$_{_{_{\mathbf{Bd}\Sigma}}}$}
\!\!\!\sigma].$$

\smallskip\noindent
{\bf iii.} ($\Rightarrow$)
Lemma~\ref{LemmaP28:2.i}\ $\!${\bf i} p.~\pageref{LemmaP28:2.i} and
Th.~\ref{TheoremP29:11} {\bf i.{\bf a}}.\hfill\break
{$(\Leftarrow)$} All links are connected, except $\bullet\bullet$.
\end{proof}

\begin{varnote}{3.} \label{NoteP29:3}
$\!$For a finite n-manifold $\Sigma$ and a n-submanifold $\Delta$,
put ${\mathcal U}:=\vert\Sigma\vert\setminus\vert\Delta\vert$,
implying that $\vert{\rm Bd}_{_{\!{\mathbf{G}}}}\!\Delta\vert\cup{\mathcal U}$
\noindent is the polytope of a subcomplex, $\Gamma$, of $\Sigma$
i.e. $\vert\Gamma\vert\!=\! \vert{\rm Bd}_{_{\!{\mathbf{G}}}}\!\Delta\vert\cup{\mathcal U}$, and ${\rm
Bd}_{_{\!{\mathbf{G}}}}\!\Sigma\subset\Gamma$, cf. \cite{26}
p.\ 427-429.
Consistency  of Definition~3 p.~\pageref{DefP27:3}
follows by excision in simplicial $\mhatH$omology since;
$\mhatH\!_{_{{{n\!\!}}}}(\Sigma, {\rm Bd}_{_{\!{\mathbf{G}}}} \!\Sigma;{\mathbf{G}})\!\hookrightarrow \mhatH\!_{_{{{n\!\!}}}}
(\Sigma,\Gamma;{\mathbf{G}})
\cong
\mhatH\!_{_{{{n\!\!}}}}(\Sigma\setminus{\mathcal
U},\Gamma\setminus{\mathcal U};{\mathbf{G}})
\!=
\mhatH\!_{_{{{n\!\!}}}}(\Delta,{\rm Bd}_{_{\!{\mathbf{G}}}}\!\Delta;{\mathbf{G}}).$
E.g., $\overline{\hbox{\rm
st}}{\!_{_{\Sigma}}}\!\sigma\!=\bar\sigma\ast{\rm
Lk}{\!_{_{\Sigma}}}\!\sigma$, is an orientable quasi-manifold
if\nobreak\ $\Sigma_{_{^{\!{\mathrm{q}}}}}$ {\rm is},
as is ${\rm Lk}{\!_{_{\Sigma}}}\!\sigma$ by
Theorem~\ref{TheoremP29:11}
and
Theorem~\ref{TheoremP31:12} below.
Moreover, $\delta\!\subset\!\sigma\!\Longrightarrow\!
\overline{\hbox{\rm st}}{\!_{_{\Sigma}}}\!\delta$ non-orientable if
$\overline{\hbox{\rm st}}{\!_{_{\Sigma}}}\!\sigma$ {\rm is}.
\end{varnote}

\begin{corollary} \label{corollary:CorP30:1} \label{CorP30:1}
For any quasi-$n$-manifold $\Sigma$ except infinite $1$-circles;

\smallskip
\centerline{$ \dim{\bf B}_{_{{\!{j}}}}\!\!\ge\!_{\!} n_{_{{\!{\
}}}}\!\!\!-\!2 \ \Longrightarrow\ \dim{\bf B}_{_{{\!{j}}}}\!\!=\!
_{\!}n_{_{{\!{\ }}}}\!\!\!-\!1 . $
}
\end{corollary}

\begin{proof}
Check $n \!\le\! 1$ and then assume that $n \!\ge\! 2$. If
$\dim\sigma\!=\dim{\bf B}_{_{{\!{j}}}}\!\!= n_{_{{\!{\ }}}}\!\!-2 $
and
$\sigma\!_{_{^{\!{\ }}}}\!\!\in\! {\bf B}_{_{{\!{j}}}}\!$
then;
$ {\rm Lk}\!\!\!\!\! \raise0.5pt\hbox{$_{_{_{\mathbf{Bd}\!_{_{\mathbf{G}}}\!\!\!\Sigma}}}$}\!\!\! \sigma\!_{_{^{\!{\
}}}}\!\!=\! {\rm Lk}\!\! \raise0.5pt\hbox{$_{_{_{\mathbf{B}_{_{{\!{j}}}}\!}}}$}\!\!\! \sigma\!_{_{^{\!{\ }}}} \!\!=
\!\{{\emptyset}_{_{^{\!o}}}\!_{\!}\} $.
\indent
Now, by Th.~\ref{TheoremP29:11}.{\bf ii},
$$ {\rm Lk}\!\!\!\!\! \raise0.5pt\hbox{$_{_{_{\mathbf{Bd}\!_{_{\mathbf{G}}}\!\!\!\Sigma}}}$}\!\!\! \sigma\!_{_{^{\!{\ }}}} =
{\rm Bd}_{_{\mathbf{G}}}\!({\rm Lk}\!_{_{\Sigma}}\sigma\!_{_{^{\!{ }}}})
=
\Big[{^{\hbox{\scriptsize$\rm{Lk}\!_{_{\Sigma}}\!\sigma\!_{_{^{\!{\
}}}}$\text{\scriptsize\rm is, by Th.~\ref{TheoremP29:11}\hbox{\bf .i}, a finite
quasi-}}}
_{\hbox{\scriptsize 1-manifold i.e. (a circle or) a
line.}}}\Big]
=
(\emptyset\ \hbox{\rm or})\ \bullet\bullet. $$
\hskip8cm Contradiction!
\end{proof}

Denote $\Sigma$ by $\Sigma\!\!\raise0.9pt\hbox{$_{_{_{{\rm
ps}}}}$}\!$, $\!\Sigma_{_{^{\!\hbox{\fiverm q}}}}\!$ and
$\Sigma_{_{^{{_{\!}}\hbox{\fiverm h}\!}}}$ when it is assumed to be
a $\hbox{\rm pseudo$\ \!$\hbox{-,}}$ quasi- resp. a ho\-mo\-logy manifold.
Note also that;  $\!\sigma\!\in\!{\rm Bd}\!$
\raise0.5pt\hbox{$_{_{_{\mathbf{G}}}}$}$\!\!\!\Sigma
_{_{^{\!\hbox{\fiverm q}}}}$
$\Longleftrightarrow$
$
{\rm Bd}
\! \raise0.5pt\hbox{$_{_{_{\mathbf{G}}}}$}\!\!
_{\!}({\rm
Lk}\!\!\raise0.5pt\hbox{$_{_{_{\Sigma}}}$}\!\!\sigma\!)_{_{^{\!\hbox{\fiverm
q}}}}\!\! =
{\rm Lk}\!\!\!\!\! \raise0.5pt\hbox{$_{_{_{\mathbf{Bd}\!_{_{\mathbf{G}}}\!\!\!\Sigma}}}$}\!\!\! \sigma\!\not=\emptyset$
by Theorem~\ref{TheoremP29:11}.ii above.

\begin{corollary} \label{CorP30:2}
{\bf i.} \label{CorP30:2.i}
If ${{\Sigma}}_{_{^{\!\hbox{\fiverm\ }}}} \!$ is finite
quasi-$n$-manifold with
$ {\rm Bd}\raise0.8pt\hbox{$_{_{_{\mathbf{G}}}}$}\!\!{{\Sigma}}_{_{^{\!\hbox{\fiverm q}}}} \!\!=
\mbfcupcap{\cup}{{j\in\lower0.5pt\hbox{\sixbf I}}}
\hbox{\bf B}_{_{\!\!^{j}}} \!$
and
$-1\!\le\!\dim{\bf B}_{_{{\!{i}\!}}}\! <\!\dim\Sigma\!-\!1$ for some
$i\!\in\! {\bf I}$ then,
$ \mhatH\!_{_{{{n}_{_{{\!{\ }}}}\!\!\!\!}}} ({\Sigma},{\rm Bd}
\raise0.8pt\hbox{$_{_{_{\mathbf{G}}}}$}\!\!{{\Sigma}};{\mathbf{G}})
\!=\hbox{\rm 0}.$
In particular, $\Sigma_{_{^{\!\hbox{\fiverm \ }}}}\!$
is non-orientable$_{_{\mathbf{G}}}\!$.

\smallskip
\noindent \hbox{\bf ii.}\label{CorP30:2.ii}
If
$ {\rm Bd}\raise0.8pt\hbox{$_{_{_{\mathbf{G}}}}$}\!\!{{\Sigma}}_{_{^{\!\hbox{\fiverm q}}}} \!\!=
\mbfcupcap{\cup}{{j\in\lower0.5pt\hbox{\sixbf I}}}
{\bf B}_{_{\!\!^{j}}} \!$
with
$\dim{\bf B}_{_{\!{i}}}:=\!\dim {\rm Bd}\raise0.8pt\hbox{$_{_{_{\mathbf{G}}}}$}\!\!{{\Sigma}}_{_{^{\!\hbox{\fiverm q}}}},$
then ${\bf B}_{_{\!{i\!}}}$ is a pseudomanifold.

\smallskip\noindent
\hbox{\bf iii.}\label{CorP30:2.iii}
$\mhatH\hskip-0.15cm{_{_{n\lower1pt\hbox{-}i}}}\! (\Sigma,{\rm
Bd}\Sigma\!\!_{_{\mathrm{ ps}}}\!;{\mathbf{G}}^{\!\prime})\!\!=\!
\mhatH\hskip-0.15cm{_{_{n\lower1pt\hbox{-}i}}}\! (\Sigma,{\rm
Bd}\raise0.8pt\hbox{$_{_{_{\mathbf{G}}}}$}\!\!\Sigma\!_{_{\mathrm{q}}};{\mathbf{G}}^{\!\prime})\!\!=\!
\mhatH\hskip-0.15cm{_{_{n\lower1pt\hbox{-}i}}}\! (\Sigma,{\rm
Bd}\raise0.8pt\hbox{$_{_{_{\mathbf{G}}}}$}\!\!\Sigma_{_{{\rm h}}};{\mathbf{G}}^{\!\prime}),\ i=0,1$
even if ${\mathbf{G}}\!\neq\!{\mathbf{G}}^{\!\prime}.$

\noindent
Orientability is independent of both $\mathbf{G}$ and $\mathbf{G}^{\!\prime}$, as long
as $ \hbox{\rm Tor}_1^{\mathbb{Z}}\bigl(\mathbb{Z}_{2},\!{\mathbf{G}}^{\prime} \bigr)\!\ne{\mathbf{G}}^{\prime} $
{\rm(Lemma~\ref{LemmaP28:1}\ ii p.~\pageref{LemmaP28:1})}.
Moreover,
${\rm Bd}\!\lower3.4pt\hbox{{{{\fivebf
G}}}}\!_{\!}\Sigma\!_{_{\hbox{\sixrm q}}}
\!\!=
{\rm Bd}\!\lower3.4pt\hbox{{{{\fivebf
G}}}}\!_{\!}\Sigma_{_{\hbox{\sixrm h}}}\!
$
always, while
$({\rm Bd}\Sigma\!\!_{_{\hbox{\sixrm
ps}}})^{^{{\!{n\!\lower1.0pt\hbox{-}\!\hbox{\fiverm 1}}}}}
\!\!\!\!=
({\rm Bd}\!\lower3.4pt\hbox{{{{\fivebf
G}}}}\!_{\!}\Sigma\!_{_{\hbox{\sixrm
q}}})^{^{{\!{n\!\lower1.0pt\hbox{-}\!\hbox{\fiverm 1}}}}} and
$
$({\rm Bd}\Sigma\!\!_{_{\hbox{\sixrm
ps}}})^{^{{\!{n\!\lower1.0pt\hbox{-}\hbox{\fiverm 2}}}}}
\!\!\!\!=
({\rm Bd}\!\lower3.4pt\hbox{{{{\fivebf
G}}}}\!_{\!}\Sigma\!_{_{\hbox{\sixrm
q}}})^{^{{\!{n\!\lower1.0pt\hbox{-}\hbox{\fiverm 2}}}}}
$
except for infinite $1$-circles in which case
$({\rm Bd}\Sigma\!\!_{_{\hbox{\sixrm
ps}}})^{^{{\!{n\!\lower1.0pt\hbox{-}\hbox{\fiverm 2}}}}}
\!\!=\emptyset\not=\{\emptyset_{_{\hbox{\sixrm o}}}\}=
({\rm Bd}\!\lower3.4pt\hbox{{{{\fivebf
G}}}}\!_{\!}\Sigma\!_{_{\hbox{\sixrm
q}}})^{^{{\!{n\!\lower1.0pt\hbox{-}\hbox{\fiverm 2}}}}}$.

\smallskip\noindent
\hbox{\bf iv.}\label{CorP30:2.iv}
For any finite orientable$_{_{\mathbf{G}}}\!\!$ quasi-$n$-manifold
$\Sigma_{_{^{\!\hbox{\fiverm \ }\!}}}\!$, each boundary component
${\bf B}_{_{\!{i}\!}}\!:=\!{\bf
B}\!\raise0.6pt\hbox{$_{_{\!{_{\mathbf{G},i}}}}$}\!\!\!\!^{^{_{\Sigma\!}}}\ \!\!
\nsubseteqq\!
\{\emptyset{_{_{^{\!o\!}}}}\!\}
$
in Def.~4 p.~\pageref{DefP27:4},
is an orientable $(n-1)$-pseudomanifold without boundary.

\smallskip\noindent
\ $\!$\hbox{\bf v.}\label{CorP30:2.v}\
$ \hbox{\rm Tor}_1^{^{_{\mathbb{Z}}}\!}\!({\mathbb{Z}}_{_{^{\!2}}}\!,\!{\mathbf{G}})\!=0$
$\!\Longrightarrow\!$
${\rm Bd}\raise0.6pt\hbox{$_{_{_{\mathbf{G}}}}$}\!\!\!\Sigma_{_{^{\!{{\rm q}}}}} \!= {\rm
Bd}{_{_{\mathbb{Z}}}}\Sigma_{_{^{\!{{\rm q}}}}}.$

\smallskip\noindent
{\bf vi.}\label{CorP30:2.vi}
{$_{\!} {\rm Bd}\Sigma\!\!_{_{\mathrm{ps}}} \!\!= \!{\rm
Bd}\raise0.5pt\hbox{$_{_{_{\mathbb{Z}_{^{^{\!2}}}}}}$}\!\!\!\Sigma_{_{^{{\rm q}}}} \!\!\subseteq \!{\rm
Bd}\raise0.8pt\hbox{$_{_{_{\mathbf{G}}}}$}\!\!\Sigma_{_{^{\!q}}}
\!\!\subseteq \!{\rm Bd}\raise0.5pt\hbox{$_{_{_{\mathbb{Z}}}}$}\!\!\Sigma_{_{^{{\rm q}}}}\! $ with}
{equality if
$ {\rm Bd}\!\raise0.6pt\hbox{$_{_{_{\mathbf{Z
}}}}$}\!\!{{\Sigma}}\!=\!\emptyset$
or
$ \dim\!{\bf B}\!\raise0.6pt\hbox{$_{_{\!{_{\mathbb{Z},_{\!}j}}}}$}\!\!\!\!^{^{_{\Sigma\!\!}}}\ \!\! =\!n\!-\!1\
\forall j\!\in\!{\bf I}, $
except if $\Sigma$ is infinite and
$ {\rm Bd}_{_{\!^{\mathbb{Z}_{_{^{\!2}}}}}}\!\!\Sigma_{_{^{\rm q}}}\!
\!=\!\{\emptyset_{_{^{\!o}}}\!_{\!}\} \!\neq\!\emptyset \!=\!{\rm
Bd}\Sigma\!\!_{_{\mathrm{ps}}} {{\!}}$}
{\rm (by Lemma~1\ $\!$\hbox{\bf
i}+\hbox{\bf ii} p.~\pageref{LemmaP28:1}
plus
Th.~\ref{TheoremP14:5} p.~\pageref{TheoremP14:5}
since
$\emptyset\!_{_{^{o}}}\!\in{\rm Bd}\!\lower3.4pt\hbox{{{{\fivebf
G}}}}\!_{\!}\Sigma\!\lower2.4pt\hbox{{{{\fiverm
q}}}}{\!}\neq\emptyset$
if $\Sigma\!\lower2.4pt\hbox{{{{\fiverm q}}}}_{\!} $ is infinite.}
{\rm {If $\Sigma$}$_{_{^{\mathrm{h}}}}\!\!$
{\hbox{\rm Gorenstein}}\hbox{$_{_{\!^{\mathbb{Z}_{_{^{\!2}}}}}}\!\!\!$}
{then};
${\rm Bd}_{_{\!^{\mathbb{Z}_{_{^{\!2}}}}}}\!\!\Sigma_{_{^{\rm h}}}$
=
${\rm Bd}\Sigma\!_{_{^{\mathrm{ps}}}}\!). $
}
\end{corollary}

\normalbaselines
\begin{proof}
\noindent \normalbaselines {\bf i.}
$\dim{\bf B}_{_{\!{i}}}\!\!<\dim\Sigma\!-\!2$ from Cor.\ 1 above,
gives the 2:nd equality
and Cor.\ 1.{\bf a} p.~\pageref{CorP29:1.a}
\noindent gives the injection-arrow in;\indent

$ \mhatH\!_{_{{{n}_{_{{\!{\ }}}}\!\!\!\!}}} ({\Sigma},{\rm Bd}
\raise0.8pt\hbox{$_{_{_{\mathbf{G}}}}$}\!\!{{\Sigma}};{\mathbf{G}}) \!=
\mhatH\!_{_{{{n}_{_{{\!{\ }}}}\!\!\!\!}}} ({\Sigma},{\bf
B}_{_{\!{i}}}\cup(\
\mbfcupcap{\cup}{j\not=i}
{\bf B}_{_{\!{j}}});{\mathbf{G}}) \!=
$

\smallskip
\hfill
$= \mhatH\!_{_{{{n}_{_{{\!{\ }}}}\!\!\!\!}}}
({\Sigma},
\mbfcupcap{\cup}{j\not=i}
{\bf B}_{_{\!{j}}} ;{\mathbf{G}})
 \hookrightarrow
\mhatH\!_{_{{{n}_{_{{\!{\ }}}}\!\!\!\!}}} ({\Sigma},\hbox{\rm
cost}_{_{_{\!{\Sigma}}}}\!\!\!\sigma;{\mathbf{G}})\!=\!\hbox{\bf 0}$
\ \
if\ \ $\sigma_{{\!}}\in_{{\!}}{\rm Bd}\!_{_{_{\mathbf{G}}}}\!\!\Sigma\
\!^{_{_{\setminus}}}
\mbfcupcap{\cup}{j\not=i}
{\bf
B}_{_{\!\!{j\!}}}\!^{^{\!_{\ }}}.$ \indent$\triangleright$

\smallskip
\noindent {\bf ii.} The claim is true if $\dim \Sigma\le 1$ and
assume it is true for dimensions $\le_{\!} n_{\!}-_{\!}1$.
$\alpha$ and $\gamma$ are true by definition of ${\bf
B}_{_{\!{i}}}\!$ so only $\beta$ remains.
If $\sigma\!\in\!
{\bf B}\raise0.5pt\hbox{$_{_{\!^{i}}}$}\!$
with $\dim \sigma \!=\! \dim{\bf B}_{_{\!{i}}}\!-_{\!}1$,
then
$\sigma \!\notin\!
\mbfcupcap{\cup}{j\not=i}
{\bf B}_{_{\!\!{j\!}}}$
and so,
${\rm Bd}\!_{_{_{\mathbf{G}}}}\!\! ({{\rm
Lk}\!_{_{_{\!\Sigma\!\!}}}\!\sigma})\!_{_{_{{\rm q}}}}\!
\!=\!
{{{{\rm Lk}_{_{{\!\!\!\!\!\!\!\lower1.3pt\hbox{$_{^{{\hbox{\sixrm
Bd}\!\!_{_{_{\mathbf{G}}}}\!\!\!\Sigma}}}$}}
}}\!\!\!{\sigma\!^{_{^{\mathbf{\ }\! }}}}}}}\! \!=\!
{{{{\rm Lk}_{_{{\!\!\!\lower1.3pt\hbox{$_{^{{\mathrm{B}_{^{_{\!i}}}}}}$}} }}\!\!\!\sigma}}} $,
implying, by the induction assumption, that the r.h.s.
is a $0$-pseudomanifold  i.e. $\bullet$ or $\bullet\bullet.$
\hfill$\triangleright$

\smallskip\noindent
\normalbaselines
{\bf iii.}
Use Note~1.{\tenbf a} p.~\pageref{NoteP29:1.a},
plus Corollary~\ref{CorP30:1} above and the fact that; [$
\sigma\!\!\in\!\Sigma^{^{{\!{n\lower1.0pt\hbox{-}\!\hbox{\fiverm 1}
\lower0.7pt\hbox{}}}}} \!_{\!}\cap{\rm Bd}\Sigma $ \ \
{\underbar{iff}} \ \
\hbox{${{{{\rm
Lk}\raise0.5pt\hbox{$\!\!_{_{_{{\Sigma}_{_{^{\!\hbox{\rm }}}}}
}}$}\!\!{\sigma^{_{^{\mathbf{}}}}\!_{\!}}}}}=\bullet]$.}
\hfill$\triangleright$

\smallskip\noindent
{\bf iv.}
{\bf i} and {\bf ii} implies that only orientability  and
${\rm Bd}_{\!}(_{\!}{\bf
B}\raise0.5pt\hbox{$_{_{\!^{j}}}$}\!)_{_{^{\mathrm{ ps}}}}\!=\
\! \emptyset$ remains.
${\rm Bd}_{\!}(_{\!}{\bf
B}\raise0.5pt\hbox{$_{_{\!^{j}}}$}\!)_{_{^{\mathrm{ps}}}}\!\neq\ \! \emptyset$
\ \underbar{iff}\ \
$\mhatH_{n-1}((_{\!}{\bf
B}\raise0.5pt\hbox{$_{_{\!^{j}}}$}\!)_{_{^{\mathrm{ps}}}}\!;{\mathbf{G}})= {\bf 0}
$
by Note~1 p.~\pageref{NoteP29:1},
and so,
${\rm Bd}
\! \raise0.5pt\hbox{$_{_{_{\mathbf{G}}}}$}\!\!
(_{\!}{\rm Lk}
\!\!\raise0.5pt\hbox{$_{_{_{\Sigma}}}$}\!\!\sigma\!)_{\!_{^{{\rm
q}}}}\!\not=\!\emptyset$
\ \underbar{iff}\ \
$\exists {\bf B}\raise0.5pt\hbox{$_{_{\!^{j}}}$} \!\ni\!\sigma$ by
{\bf iii}.
From Def.~4 p.~\pageref{DefP27:4}
we get,
$\dim{\bf B}_{{\!{s}}}\!= \dim{\bf B}_{{\!{t}}}\!= {n}_{_{{\!{\
}}}}\!\!-\!1,\ s,t\in{\bf I}
$
$\Longrightarrow
\dim{\bf B}_{{\!{s}}}\cap {\bf B}_{{\!{t}}}\leq n-3\ i\!f
s\!\neq\! t,$
which, since
$({\rm Bd}\Sigma\!\!_{_{\hbox{\sixrm ps}}})^{^{{\!{n-i}}}}
\!\!\!\!=
({\rm Bd}\!\lower3.4pt\hbox{{{{\fivebf
L}}}}\!_{\!}\Sigma\!_{_{\hbox{\sixrm q}}})^{^{{\!{n-i}}}}\ i=1,2,
$
gives;

\smallskip\noindent
$ \dots\hskip-0.2cm\longrightarrow\hskip-0.3cm
\underbrace{ \mhatH\!_{_{{{n}_{_{{\!{\ }}}}\!\!}}}\! ({\rm
Bd}\raise0.8pt\hbox{$_{_{_{{\bf L}}}}$}\!\!{{\Sigma}},
\mbfcupcap{\cup}{j\not=i}
{\bf
B}_{_{\!{\mathbf{L}^{\!},^{\!}j}}} ;\!{\mathbf{G}}^{\!\prime})}
_{=\ \!0\ \!\hbox{\eightrm for$\ $dimensional$\ $reasons.}}
\vbox{\moveleft0.2cm\hbox{\textbf{}$\longrightarrow$}}
\underbrace{ \mhatH\!_{_{{{n}_{_{{\!{\ }}}}\!\!}}}\! ({{\Sigma}},
\mbfcupcap{\cup}{j\not=i}
{\bf
B}_{_{\!{\mathbf{L}^{\!},^{\!}j}}} ;{\mathbf{G}}^{\!\prime})}
_{\hskip-0.2cm=\ \!0\ \! \hbox{\eightrm by Cor.~\ref{CorP29:2}
p.~\pageref{CorP29:2}
}}
\longrightarrow
\hskip-0.2cm
\underbrace{ \mhatH\!_{_{{{n}_{_{{\!{\ }}}}\!\!}}}\!({{\Sigma}},
{\rm Bd}\raise0.5pt\hbox{$_{_{_{{\bf L}}}}$}\!\!{{\Sigma}} ;{\mathbf{G}}^{\!\prime})}
_{=\hbox{{\eightbf G} \eightrm by assumption if \hbox{\eightbf
G}$^{\!\prime}$=\hbox{{\eightbf G}}}}
\hskip-0.5cm\hookrightarrow\!
$

\smallskip\noindent
$
\hookrightarrow\!
\mhatH\!\!_{_{^{ \!{n}_{_{{\!{\
}}}}\!\!\lower1pt\hbox{\tenbf-}1}}}\!( {\rm
Bd}\raise0.8pt\hbox{$_{_{_{\mathbf{G}}}}$}\!\!{{\Sigma}},
\mbfcupcap{\cup}{j\not=i}
{\bf
B}_{_{\!{\mathbf{L}^{\!},^{\!}j}}} ;{\mathbf{G}}^{\!\prime})
=\!
\big[{{\hbox{\scriptsize For dimensio-}}\atop{{\hbox{\scriptsize
nal reasons.}}}}\big]\!=
\mhatH_{\!\!n\hbox{\fivebf-}1\!}({\bf B}_{_{\!{\mathbf{L}^{\!},^{\!}i}}} ;{\mathbf{G}}).$

\smallskip
In the above truncated relative {\bf LHS} with respect to
$(\Sigma ,{\rm Bd}\raise0.8pt\hbox{$_{_{_{\mathbf{G}}}}$}\!\!{{\Sigma}},
\mbfcupcap{\cup}{j\not=i}
{\bf B}_{_{\!{\mathbf{L}^{\!},^{\!}j}}}), $
the choice of ${\bf L}$ is irrelevant by {\bf iii}.
Lemma~\ref{LemmaP28:1.ii}.{\bf ii} p.~\pageref{LemmaP28:1.ii}
gives the
orientability$_{_{^{\!\hbox{\fivebf G}}}}$.
If ${\mathbf{G}}^{\!\prime}\!=\!{\mathbb{Z}}$ in the {\bf LHS}, the injection
gives our claim.
Otherwise
$  \hbox{\rm Tor}_1^{^{_{\mathbb{Z}}}\!}\!({\mathbb{Z}}_{_{^{\!2}}}\!, \!{\mathbf{G}})\!=\!{\mathbf{G}} $
and we are done. \hfill$\triangleright$

\smallskip\noindent
{\bf v.} \label{CorP30:2.v.proof}
$
\ \ {\rm Bd}\raise0.5pt\hbox{$\!_{_{_{\mathbf{G}}}}$}\!\!\!\Sigma\!_{_{^{{\rm q}}}} \!\!\not\ni_{\!} \sigma\!\in\!
{\rm Bd}\raise0.5pt\hbox{$\!_{_{_{\!{\mathbb{Z}}}\!}}$}\!\!\Sigma\!_{_{^{{\rm q}}}}\!
 \Longleftrightarrow\ \!
\mhatH \vbox{\moveleft0.40cm\hbox{$
_{_{_{\!{n \lower1.0pt\hbox{-} \hbox{\fivebf\#}\sigma}}}}
$}}
\!\!\!\!
(({{\rm Lk}\!_{_{_{\!\Sigma\!\!}}}\!\sigma}
)\!_{_{^{{\rm q}}}}\! ;_{\!}{\mathbf{G}})
\not= 0 =
\mhatH
\vbox{\moveleft0.40cm\hbox{$
_{_{_{\!{n -\#\sigma}}}}
$}}
\!\!\!\!
(({{\rm Lk}\!_{_{_{\!\Sigma\!\!}}}\!\sigma}
)\!_{_{^{{\rm q}}}}\! ;_{\!}{\mathbb{Z}})
\Longleftrightarrow {\rm Bd}\!\raise0.5pt\hbox{$_{_{_{\mathbf{G}}}}$}\!\!\!$
$
\!({\rm Lk}\!\raise0.5pt\hbox{$_{_{_{\Sigma}}}$}
\!\!\sigma\!)\!_{_{^{{\rm q}}}}\!
=\emptyset\not=\
$
\hbox{\rm Bd$\!\!$}
\raise0.5pt\hbox{$_{_{_{\mathbb{Z}}}}$}
$\!
\!({\rm Lk}\!\raise0.5pt\hbox{$_{_{_{\Sigma}}}$}
\!\!\sigma\!_{_{^{\!{\ }}}}\!_{\!})\!_{_{^{{\rm q}}}}\!.
$

\smallskip\noindent
$({\rm Bd}
({\rm Lk}\!\raise0.5pt\hbox{$_{_{_{\Sigma}}}$}
\!\!\delta\!_{_{^{\!{\ }}}}\!_{\!})\!_{_{^{\mathrm{ ps}}}}\!
)\!^{^{_{(n-\#\delta)-i}}}
\hskip-0.2cm=
({\rm Bd} \!\!\raise0.5pt\hbox{$_{_{_{\mathbf{L}}}}$}\!
({\rm Lk}\!\raise0.5pt\hbox{$_{_{_{\Sigma}}}$}
\!\!\delta\!_{_{^{\!{\ }}}}\!_{\!})\!_{_{^{{\rm q}}}}\!
)\!^{^{_{(n-\#\delta)-i}}}\!\!,
i=1,2$,
 by
{\tenbf iii}
$\Longrightarrow$
$
n-\#\sigma-3\ge
$
$
\dim{\hbox{\rm Bd$\!\!$}
\raise0.5pt\hbox{$_{_{_{\mathbb{Z}}}}$}
\!
({\rm Lk}\!\raise0.5pt\hbox{$_{_{_{\Sigma}}}$}
\!\!\sigma\!_{_{^{\!{\ }}}}\!_{\!})\!_{_{^{{\rm q}}}}\!
}$
$\Longrightarrow$
$
0\neq \mhatH
\vbox{\moveleft0.40cm\hbox{$
_{_{_{\!{n \lower1.0pt\hbox{-} \hbox{\fivebf\#}\sigma}}}}
$}}
\!\!\!\!
(({{\rm Lk}\!_{_{_{\!\Sigma\!\!}}}\sigma}\raise1.5pt\hbox{\eightbf
{\char"29}}\!_{_{^{{\rm q}}}}\! ;_{\!}{\mathbf{G}})
$
$=\!\!\big[\!{{\scriptsize\rm Th.~\ref{TheoremP14:5}
}\atop{{\scriptsize\rm p.~\pageref{TheoremP14:5}
}}}\!\big]\!\!=\!$
$
\underbrace{ \mhatH \vbox{\moveleft0.40cm\hbox{$
_{_{_{\!{n \lower1.0pt\hbox{-} \hbox{\fivebf\#}\sigma}}}}
$}}
\!\!\!\!
(({{\rm Lk}\!_{_{_{\!\Sigma\!\!}}}\sigma}\raise1.5pt\hbox{\eightbf
{\char"29}}\!_{_{^{{\rm q}}}}\! ;_{\!}{\mathbb{Z}})
}
_{=\ \!0\ \!
\hbox{\eightrm by$\ \!$\hbox{\eightrm assumption}.}}
\!\!\! \otimes{\mathbf{G}}
\oplus  \hbox{\rm Tor}_1^{^{_{\mathbb{Z}}}\!}(
\mhatH \vbox{\moveleft0.40cm\hbox{$
_{_{_{\!{n \lower1.0pt\hbox{-}
\hbox{\fivebf\#}\sigma}\lower1.0pt\hbox{-}1}}}
$}}
\hskip-0.45cm
(({{\rm Lk}\!_{_{_{\!\Sigma\!\!}}}\sigma}\raise1.5pt\hbox{\eightbf
{\char"29}}\!_{_{^{{\rm q}}}}\! ;_{\!}{\mathbb{Z}}),{\mathbf{G}})
\!=\!
\big[{{\hbox{\scriptsize\rm For dimensio-}}\atop{{\hbox{\scriptsize\rm
nal reasons.}}}}\big]\!=\!
 \hbox{\rm Tor}_1^{^{_{\mathbb{Z}}}\!}\!
(\phantom{\ }
\mhatH \vbox{\moveleft0.40cm\hbox{$
_{_{_{\!{n \lower1.0pt\hbox{-}
\hbox{\fivebf\#}\sigma}\lower1.0pt\hbox{-}1}}}
$}}
\hskip-0.45cm
(({{\rm Lk}\!_{_{_{\!\Sigma\!\!}}}\sigma}\raise1.5pt\hbox{\eightbf
{\char"29}}\!_{_{^{{\rm q}}}}\!,
\hbox{\rm Bd$\!\!$}
\raise0.5pt\hbox{$_{_{_{\mathbb{Z}}}}$}
\!
({\rm Lk}\!\raise0.5pt\hbox{$_{_{_{\Sigma}}}$}
\!\!\sigma\!_{_{^{\!{\ }}}}\!_{\!})\!_{_{^{{\rm q}}}}\!
;_{\!}{\mathbb{Z}}),{\mathbf{G}})\!=\!0.$
\underbar{Contradiction!}

(-
$  \hbox{\rm Tor}_1^{^{_{\mathbb{Z}}}\!}\!({\mathbb{Z}}_{_{^{\!2}}}\!, \!{\mathbf{G}})\!=0$
and the torsion module of
$
\mhatH \vbox{\moveleft0.40cm\hbox{$
_{_{_{\!{n \lower1.0pt\hbox{-}
\hbox{\fivebf\#}\sigma}\lower1.0pt\hbox{-}1}}}
$}}
\!\!\!\!
(({{\rm Lk}\!_{_{_{\!\Sigma\!\!}}}\sigma}\raise1.5pt\hbox{\eightbf
{\char"29}}\!_{_{^{{\rm q}}}},
\hbox{\rm Bd}
\raise0.5pt\hbox{$_{_{_{\mathbb{Z}}}}$}
\!
\!({\rm Lk}\!\raise0.5pt\hbox{$_{_{_{\Sigma}}}$}
\!\!\sigma)\!_{_{^{{\rm q}}}} ;{\mathbb{Z}}) $ is either 0 or homomorphic
to ${\mathbb{Z}}_{_{^{\!2}}}$ by
Lemma~\ref{LemmaP28:1.i}.{\bf i} p.~\pageref{LemmaP28:1.i}
give the contradiction, since only the torsion submodules matters
in the torsion product by \cite{30}
Corollary\ 11 p.\ 225.) \hfill$\triangleright$

\medskip\noindent
{\bf vi.}
\ \ \hbox{\bf iv} \ and that ${\rm
Bd}\raise0.6pt\hbox{$_{_{_{\!\!\hbox{\fivebf
G}}}}$}\!\!\Sigma_{_{^{\!\hbox{\fiverm q}}}} \subseteq {\rm
Bd}\raise0.5pt\hbox{$_{_{_{\mathbb{Z}}}}$}\!\!\Sigma_{_{^{\!\hbox{\fiverm q}}}}$, by
Corollary~\ref{CorP29:2} p.~\pageref{CorP29:2}
plus
Theorem~\ref{TheoremP14:5} p.~\pageref{TheoremP14:5}.
\end{proof}

\goodbreak
\normalbaselines

\subsection{Products and joins of simplicial manifolds} \label{SubSecP31:III:3}
$$\phantom{.}$$
\vskip-0.8cm
Let in the next theorem, when
$\mytinynabla$,
all through, is interpreted as $\times$, the word ``manifold(s)'' in
{\bf \ref{TheoremP31:12}.1}
temporarily excludes $\emptyset,\{\emptyset_o\}$ and
$\bullet\bullet$.

When
\mytinynabla,
all through, is interpreted as $\ast$ let the word ``manifold(s)''
in {\bf \ref{TheoremP31:12}.1}
stand\nobreak\ for \underbar{finite}``pseudo-manifold(s)''
(``quasi-manifold(s)'',
cf.\ \cite{12}
4.2\ pp.\ 171\hbox{-}2).
We conclude, with respect to joins, that
Th.\ \ref{TheoremP31:12} is trivial if either
$\Sigma_{1}$ or $\Sigma_{2}$
equals $\{\emptyset_{\!o}\}$
 and that $\Sigma_{1}$ and $\Sigma_{2}$ must be
\underbar{finite}
since otherwise, their join $\Sigma_{1} \ast \Sigma_{2}$ won't be locally \underbar{finite}.
Moreover, set $\epsilon=0$ if
$\mytinynabla=\times$ and  $\epsilon=1$ if
$\mytinynabla=\ast$. Only quasi-manifold boundaries and orientability depend on the choice of coefficient module in the next theorem. Also  disposing off all torsion terms by using a field $\mathbf{k}$ as coefficient module, we are essentially back in well known  cases proven for instance as Theorem~10 in \cite{10} p.~82 (pseudomanifolds) and in \cite{12} p.~171 as ``Satz'' (joins of quasi-manifolds).

\begin{theorem} \label{theorem:TheoremP31:12} \label{TheoremP31:12}
If ${\mathbf{V}}_{_{\Sigma_{^{i\!}}}}\ne\emptyset$ then$;$

\smallskip
\noindent {\bf 1.}\
$\Sigma_{1}\mytinynabla
\Sigma_{2}\!$ is a $(n_{1}\!+n_{2}\!+\epsilon)$
$
\hbox{-manifold}\ \! \Longleftrightarrow \Sigma_{i}$ is a
$n_{_{i\!}}$-manifold.

\smallskip
\noindent {\bf 2.}
${\rm Bd} (\bullet\times \Sigma)= \bullet\times ({\rm Bd}
\Sigma)$. {Else;} ${\rm Bd} (\Sigma_{1}\!
\mytinynabla
\Sigma_{2})= (({\rm Bd} \Sigma_{1})
\mytinynabla\!
\Sigma_{2})\cup (\Sigma_{1}
\mytinynabla\!
({\rm Bd} \Sigma_{2})).\
$

\smallskip
\noindent {\bf 3.}\
If\ any\ side\ of\ \ref{TheoremP31:12}.1\ holds$;$
$\Sigma_{1}
\mytinynabla
\Sigma_{2}$ is\ orientable$_{_{\mathbf{k}}}$
$\Longleftrightarrow$
$\Sigma_{1},\Sigma_{2}$ are both orientable$_{_{\mathbf{k}}}.$
\end{theorem}

\begin{proof}
{\bf (\ref{TheoremP31:12}.1)}\
[{\rm Pseudomanifolds}.] Lemma~\ref{LemmaP20:2} p.~\pageref{LemmaP20:2}.
\hfill$\triangleright$

\medskip\hskip1cm
[{\rm Quasi-manifolds}.] Lemma~\ref{LemmaP19:1} p.~\pageref{LemmaP19:1} + Lemma~\ref{LemmaP20:2} p.~\pageref{LemmaP20:2}.
\hfill$\triangleright$

\smallskip\noindent
Proposition~\ref{PropP15:1} p.~\pageref{PropP15:1}
and Theorem~\ref{TheoremP17:7} p.~\pageref{TheoremP17:7}
gives the rest, but we still add a semi-combinatorial proof - as follows.

\smallskip
\noindent{\bf (\ref{TheoremP31:12}.2)} [{Quasi-manifolds}$;\times$]
Put $n:=\dim{\Sigma_{1}\!\times\!\Sigma_{2}}=
\dim{\Sigma_{1}\!+\!\dim\Sigma_{2}}= n_{_1}+n_{_2}$.\nobreak\
The invariance of local $\mhatH$omology within
$\hbox{\rm Int}\sigma_{1}\!\times \hbox{\rm
Int}\sigma_{2}$
implies, through Prop.\ 1 p.~\pageref{PropP15:1},
that, without loss of generality, we will only  need to study
simplices with
$c_{_{^{\!\sigma}}}\!\!=0$ (Def. p.~\pageref{DefP19:csigma}).
Put $v\!:=\!\dim\sigma\!=\!
\dim{\sigma_{_{^{_{\!}1\!\!}}}}\!+\dim\sigma_{_{^{_{\!}2\!}}}\!=:\nobreak
v_{_{^{_{\!}1\!\!}}}\!+v_{_{^{_{\!}2\!}}}.$
Recalling from Theorem~\ref{TheoremP29:11.i.a}.1.a p. \pageref{TheoremP29:11.i.a} that now, all links are finite quasi-manifolds,
``{$\sigma\!\in\!{\rm Bd}_{_{\mathbf{k}}}{\Sigma_{^{_{^{\!}1}}}\!\!\times\!\Sigma_{^{_{^{\!}2}}}}\!$
$\Longleftrightarrow\!$ $\sigma\!_{^{_{1}}}\!\in\!{\rm Bd}_{_{\mathbf{k}}}{\Sigma_{^{_{\!1}}}}$ or $\sigma_{^{_{\!2}}}\!\!\in\!{\rm
Bd}_{_{\mathbf{k}}} {\Sigma_{^{_{^{\!}2}}}}_{\!}$}''
follows from Corollary~\ref{CorP19:ToTh6} p.~\pageref{CorP19:ToTh6}
which, after deletion of a trivial torsion term through Note~\ref{NoteP17:1} p.~\pageref{NoteP17:1}, gives;
$$
\mhatH\!\!\!\!{_{_{n\lower0pt\hbox{\sevenrm-}v\lower0pt\hbox{\sevenrm-}1}}}\!
({\rm
Lk}\!\!\!\!\!\!\!_{_{_{\Sigma_{1}\!\times\!\Sigma_{2}}}}\!\!\!\!\!\sigma;
{\mathbf{k}})
\ {\lower1pt\hbox{${^{{{_{\mathbf{k}}}}\atop{\hbox{$\cong$}}}}$}}\
\mhatH\hskip-0.2cm{_{_{n\lower0pt\hbox{\sevenrm-}v\lower0pt\hbox{\sevenrm-}1}}}\!
({\rm
Lk}\hskip-0.3cm\!_{_{_{\Sigma_{1}\!\times\!\Sigma_{2}}}}\!\!\!\!\!\sigma;
{\mathbf{k}}\otimes{\mathbf{k}})
\ {\lower1pt\hbox{${^{{{_{\mathbf{k}}}}\atop{\hbox{$\cong$}}}}$}}\
\mhatH\!\!\!\!
{_{_{n\!_{_1}\!\!\lower0pt\hbox{\sevenrm-}v\!_{_1}\!\!
\lower0pt\hbox{\sevenrm-}1}}}\!
                 ({\rm Lk}\!\!_{_{_{{\Sigma_{1}}}}}\!\!\!\sigma\!_{_1};{\mathbf{k}})
\otimes\!{_{_{\mathbf{k}}}}\ \! \mhatH\!\!\!\!
{_{_{n\!_{_2}\!\!\lower0pt\hbox{\sevenrm-}v\!_{_2}\!\!
\lower0pt\hbox{\sevenrm-}1}}}\!
           ({\rm Lk}\!\!_{_{_{{\Sigma_{2}}}}}\!\!\!\sigma\!_{_2};{\mathbf{k}})
\hskip0.9cm\diamondsuit\hskip0.3cm\triangleright
$$
[{\rm Pseudomanifolds}$;\times$]
Reason as for Quasi-manifolds but restrict to sub\-maximal simplices
only, the links of which are $0$-pseudo\-manifolds.
\hfill$\triangleright$

\noindent
[{Pseudo- and Quasi-manifolds}$;\ast$]
Use Corollary~\ref{CorP19:ToTh6} p.~\pageref{CorP19:ToTh6}
as for products,
cf.\ \cite{12}
4.2\ pp.\ 171\hbox{-}2.\hfill$\triangleright$
$$\hskip-0.2cm
\noindent
\left.\begin{array}{llr}
\leftline{${\bf (\ref{TheoremP31:12}.3)}
\hskip0.5cm
(\Sigma_1\!
\mytinynabla\!
\Sigma_2, {\rm Bd}_{_{\mathbf{k}}}{\!}_{\!}(\Sigma_1\!
\mytinynabla\!
\Sigma_2))=$}\\
\noindent=
\big[{^{\hbox{\eightrm Motivation:}}
_{\hbox{\eightrm  Th.\ \ref{TheoremP31:12}.2}}}\big]
=
(\Sigma_{_{^{_{\!}1}}}\!
\mytinynabla\!
\Sigma_2,
\Sigma_1\!
\mytinynabla\!
{\rm Bd}_{_{\mathbf{k}}}{\!}_{\!}\Sigma_2
\cup
{\rm Bd}_{_{\mathbf{k}}}{\!}_{\!}\Sigma_1\!
\mytinynabla\!
\Sigma_2)
=
\Big[{^{\hbox{\scriptsize\rm Motivation: According to}}
_{\hbox{\scriptsize\rm the pair-definition p.~\pageref{DefP11:PairDef}.}}}\Big] =
\\
\hfill=(\Sigma_{_{^{1\!}}},{\rm Bd}_{_{\mathbf{k}}}{\!}\Sigma_{_{^{1\!}}})
\mytinynabla\!
(\Sigma_{_{^{_{\!}2\!}}},{\rm Bd}_{_{\mathbf{k}}}{\!}\Sigma_{_{^{2\!}}}).
\end{array}\right.
$$
Use
Eq.~1\ p.~\pageref{EqP12:1}
(resp. Eq.~3 p.~\pageref{EqP14:3})
and Note~3 p.~\pageref{NoteP29:3}.
\end{proof}

\begin{varnote}
In general, the equivalence
``{$\sigma\!\in\!{\rm Bd}_{_{\mathbf{G}}}{\Sigma_{^{_{1}}}\times\!\Sigma_{^{_{2}}}}\!$
$\Longleftrightarrow\!$ $\sigma\!_{^{_{1}}}\!\in\!{\rm Bd}_{_{\mathbf{G}}}{\Sigma_{^{_{1}}}}$ or $\sigma_{^{_{\!2}}}\in\!{\rm
Bd}_{_{\mathbf{G}}} {\Sigma_{^{_{2}}}}_{\!}$}'' for quasi-manifolds (and homology manifolds) can be investigated
through the general case of Corollary~\ref{CorP19:ToTh6} p.~\pageref{CorP19:ToTh6} which, after deletion of a trivial torsion term through Note~1 p.~\pageref{NoteP17:1}, gives;
$$
\indent
\mhatH\!\!\!\!{_{_{n\lower0pt\hbox{\sevenrm-}v\lower0pt\hbox{\sevenrm-}1}}}\!
({\rm
Lk}\!\!\!\!\!\!\!_{_{_{\Sigma_{1}\!\times\!\Sigma_{2}}}}\!\!\!\!\!\sigma;
{\mathbf{G}})
\ {\lower1pt\hbox{${^{{{_{\mathbb{Z}}}}\atop{\hbox{$\cong$}}}}$}}\
\mhatH\!\!\!\!
{_{_{n\!_{_1}\!\!\lower0pt\hbox{\sevenrm-}v\!_{_1}\!\!
\lower0pt\hbox{\sevenrm-}1}}}\!
                 ({\rm Lk}\!\!_{_{_{{\Sigma_{1}}}}}\!\!\!\sigma\!_{_1};{\mathbb{Z}})
\otimes\!\!{_{_{\mathbb{Z}}}}\ \! \mhatH\!\!\!\!
{_{_{n\!_{_2}\!\!\lower0pt\hbox{\sevenrm-}v\!_{_2}\!\!
\lower0pt\hbox{\sevenrm-}1}}}\!
           ({\rm Lk}\!\!_{_{_{{\Sigma_{2}}}}}\!\!\!\sigma\!_{_2};{\mathbf{G}})\oplus \hskip0.5cm\phantom{I}
\hskip0.5cm{\lower0.4cm\hbox{$\clubsuit$}}
$$
\vskip-0.4cm
$$
\phantom{II}\hskip1.5cm
\oplus \hbox{\rm Tor}_1^{\mathbb{Z}}\bigl( \mhatH\!\!\!\!
{_{_{n\!_{_1}\!\!\lower0pt\hbox{\sevenrm-}v\!_{_1}\!\!
\lower0pt\hbox{\sevenrm-}2}}}\!
                  ({\rm Lk}\!\!_{_{_{{\Sigma_{1}}}}}\!\!\!\sigma\!_{_1};{\mathbb{Z}})
,\mhatH\!\!\!\!
{_{_{n\!_{_2}\!\!\lower0pt\hbox{\sevenrm-}v\!_{_2}\!\!
\lower0pt\hbox{\sevenrm-}1}}}\!
           ({\rm Lk}\!\!_{_{_{{\Sigma_{2}}}}}\!\!\!\sigma\!_{_2};{\mathbf{G}})\bigr).
\phantom{III}\hskip3.0cm
$$
The coefficient module in the last theorem was restricted to a field ${\mathbf{k}}$, but it is possible to work with arbitrary coefficient modules and infinite manifolds with respect to orientation in Th.~\ref{theorem:TheoremP31:12}.3 as follows.
By Cor.~\ref{CorP30:2.iv}.iii p.~\pageref{CorP30:2.iv}
one can without loss of generality confine the study to
pseudomanifolds, and then choose the coefficient module to be,
say, a field ${\bf k}$ (char${\bf k}\neq2$).
Since any finite maxi-dimensional submanifold, i.e. a submanifold of maximal dimension, in
$\Sigma_{_{^{_{\!}1\!}}}\!\times\!\Sigma_{_{^{_{\!}2\!}}}$
\noindent ($\Sigma_{_{^{_{\!}1\!}}}\!\ast\Sigma_{_{^{_{\!}2\!}}},$
cf. \cite{36}
(3.3) p.\ $_{\!}$59)
can be embedded in the product (join) of two finite
maxi-dimensional submanifolds, and vice versa,
one might confine, without loss of generality, the attention to finite maxi-dimensional submanifolds ${\tensy S}_{1}$, ${\tensy
S}_{2}$ of $\Sigma_{_{^{_{\!}1\!}}}$ resp.
$\Sigma_{_{^{_{\!}2\!}}}$. Now, use
\hbox{Eq.~1\ p.~\pageref{EqP12:1}
{\bf(}\hbox{\rm Eq.~3\ p.~\pageref{EqP14:3}}{\bf)}
and Note~3 p.~\pageref{NoteP29:3}.}$\!$
\end{varnote}

\begin{varex}{\bf 1.}
\label{ExampleP32}\label{ExampleP32:1}
For a triangulation $\mathcal{C}$ of a two-dimensional cylinder
${\rm Bd}\lower1.3pt\hbox{$_{_{{\mathbb{Z}}}}$}\!
\mathcal{C}\!_{_{^{\!\hbox{\fiverm h}}}}\!\!=\! \hbox{\rm two\ circles}$.
\indent
${\rm Bd}\lower1.3pt\hbox{$_{_{{\mathbb{Z}}}}$}\!\bullet
=\{\emptyset_o\!\}.$ By Th.\ \ref{TheoremP31:12}
${\rm Bd}\lower1.3pt\hbox{$_{_{{\mathbb{Z}}}}$}\!
(\mathcal{C}\!\ast\bullet)_{_{{\rm q}}}
\!\!=\!
\mathcal{C}\cup(\{\hbox{\rm two\ circles}\}\ast\bullet) $.
So, ${\rm Bd}_{_{\mathbb{Z}}}\!(\mathcal{C}\ast\bullet)_{_{{\rm
q}}}$,
{\bf R}$^{^{_{3}}}$-
realizable as a pinched torus (or inwards as in p.~\pageref{BundaryCylCone}), is a 2-pseudomanifold but not a
quasi-manifold.

\item[\bf 2.]
\label{ExampleP32:2}
Cut, twist and glue along the cut to turn the cylinder $\mathcal{C}$ above
into a Möbius band ${\mathcal M}$.

``The boundary (with respect to ${\mathbb{Z}}$) of the cone of ${\mathcal M}$'' =
$${\rm Bd}\lower1.3pt\hbox{$_{_{{\mathbb{Z}}}}$}\! ({{\mathcal
M}}\!\ast\bullet)_{_{{\rm q}}}
\!\!=\!
({\mathcal M}\ast\{\emptyset_o\!\})\cup(\{\hbox{\rm a\
circle}\}\ast\bullet)
\!\!=\!
{{\mathcal M}}\cup \{\hbox{\rm a\ $2$-disk}\}$$
which is a well-known representation of the real projective plane
$\mathbb{RP}^2$,
i.e. a homology$\!_{_{\mathbb{Z}_{_{{\bf
p}}}}}$
$2$-manifold  \underbar{with} boundary
$\{\emptyset_o\!\}\neq\emptyset$ if  ${\bf p}\neq2$, cf. \cite{13}
p.~36.

$\mathbb{Z}_{\bf p}:=\!$ {\it The prime-number field
modulo {\bf p}}, i.e. of characteristic {\bf p}.

\indent
$
{\mathbb{RP}}^{_2}\#{\mathbb{S}}^{_{^{2}}}\!\!
\!=\!
{\mathbb{RP}}^{_2}_{\!}
\!=\!
{{\mathcal M}}\
\mbfcupcap{\cup}{\hskip0.05cm\rm Bd}
\{\hbox{\rm a\ $2$-disk}\}
$
confirms the obvious --- that the $n$-sphere is the unit element
with respect to the connected sum of two $n$-manifolds, cf.
\cite{26}
p.\ 38ff   + Ex.\ 3 p. \ 366.
Leaving the ${\mathbb{S}}^{_{^{2}}}\!$-hole empty turns ${\mathbb{RP}}^{_2}_{\!}$ into ${\mathcal M}$.

Let
``\mbfcupcap{\cup}{\hskip0.05cm\rm Bd}''
denote ``union through identification of boundaries''.
$${\rm Bd}\lower1.3pt\hbox{$_{_{{\mathbb{Z}}}}$}\! ({{\mathcal
M}}\!\ast\!{\mathbb{S}}^{_{^{1}}})_{_{{\rm q}}}
\!\!=\!
({\mathcal M}\ast\emptyset)\cup({\mathbb{S}}_{\!}^{_{^{1}}}\ast{\mathbb{S}}^{_{^{1}}})
\!\!=\!
{\mathbb{S}}_{\!}^{_{^{3}}}.
$$
\ \textrm{So,}\
$$ {\rm Bd}\lower1.3pt\hbox{$_{_{\mathbb{Z}}}$}\!({{\mathcal M}}_{1}{\!}\ast_{\!} {\mathcal
M}_{2}\!)_{_{{\rm q}}}\!{\!}
={\!}({\mathbb{S}}_{\!}^{_1}\!\ast {{\mathcal M}}_{2}{\!})_{_{{\rm
q}}}\
\mbfcupcap{\cup}{\hskip0.05cm\rm Bd}
(_{\!}{\mathcal M}_{1}\!\ast {\mathbb{S}}_{\!}^{_1})_{_{{\rm q}}},
\ \textrm{is a quasi-$4$-manifold.}
$$
Is
${\rm Bd}\lower1.3pt\hbox{$_{_{{\mathbb{Z}}}}$}\!({{\mathcal
M}}_{1}{\!}\ast_{\!} {\mathcal M}_{2}\!)_{_{{\rm q}}}$
orientable or is
${\rm Bd} _{_{{\mathbb{Z}}}} \!({\rm Bd}_{_{{\mathbb{Z}}}}
\!({{\mathcal
M}}_{1}{\!}\ast_{\!} {\mathcal M}_{2}\!)_{_{{\rm
q}}}\!)_{_{{\rm q}}}\!$
$\lower1.5pt\hbox{\rlap{$^{^{_{\hbox{\sevenrm?}}}}$}{$_{\!\!}$=}}
\{\emptyset_{\!o}\!\}$?

\medskip
\item[\bf 3.]
\label{ExampleP32:3}
Let ${\mathbb{RP}}^2$,
${\mathbb{RP}}^4$ be triangulations of
the 2-, resp. 4-dimensional real projective plane/space, then
${\rm Bd}\lower1.3pt\hbox{$_{_{{\mathbb{Z}}}}$}\! {\mathbb{RP}}\raise4pt\hbox{{\fiverm 2}}_{\!}\!\!_{_{{\rm h}}}\!\!=\! {\rm
Bd}\lower1.3pt\hbox{$_{_{{\mathbb{Z}}}}$}\!{\mathbb{RP}}\raise4pt\hbox{{\fiverm 4}}_{\!}\!_{_{{\rm h}}}\!\!=\!
\{\emptyset_o\!\}.$
So, by Theorem~\ref{TheoremP31:12} p.~\pageref{TheoremP31:12},
$$
\left\{\begin{array}{ll}
{\rm Bd}_{_{\mathbb{Z}_{{\!\hbox{\fivebf
p}}}}}\! ({\mathbb{RP}}\raise4pt\hbox{{\fiverm
2}}_{\!}\!\ast{\mathbb{RP}}\raise4pt\hbox{{\fiverm 4}}_{\!})_{_{{\rm
h}}} \!\!=\!{\mathbb{RP}}\raise4pt\hbox{{\fiverm 4}}\cup{\mathbb{RP}}\raise4pt\hbox{{\fiverm 2}}, &{\rm if}\ {\bf p}\!\neq\!2,\\
{\rm Bd}_{_{\mathbb{Z}_{{\!\hbox{\fivebf
p}}}}}\!
({\mathbb{RP}}\raise4pt\hbox{{\fiverm 2}}_{\!}\!\ast{\mathbb{RP}}\raise4pt\hbox{{\fiverm 4}}_{\!})_{_{{\rm q}}} \!=\!\emptyset,\
&{\rm if}\ {\bf p}\!=\!2.
\end{array}\right.
$$

Note that $\dim{\rm Bd}_{_{\mathbb{Z}_{{\!\hbox{\fivebf
p}}}}}\! ({\mathbb{RP}}\raise4pt\hbox{{\fiverm 2}}_{\!}\!\ast{\mathbb{RP}}\raise4pt\hbox{{\fiverm 4}}_{\!})_{_{{\rm h}}} \!\!=\!4$
while
$\dim({\mathbb{RP}}\raise4pt\hbox{{\fiverm 2}}_{\!}\!\ast{\mathbb{RP}}\raise4pt\hbox{{\fiverm 4}}_{\!})_{_{{\rm h}}} \!\!=\!7$, cf.
Corollary~\ref{CorP30:1} p.~\pageref{CorP30:1}.

${\rm Bd}_{_{\mathbb{Z}}}\!({\mathbb{RP}}
\raise4pt\hbox{{\fiverm 2}}_{\!}\!\ast\bullet\bullet)_{_{{\rm
q}}}\!\!=\!{\bullet\bullet}$
 \ \ and  \ \
${\rm Bd}_{_{\mathbb{Z}}}\!({\mathbb{RP}}\raise4pt\hbox{{\fiverm 2}}_{\!}\!\ast\bullet)_{_{{\rm q}}}\!
{=\!{\mathbb{RP}}\raise4pt\hbox{{\fiverm 2}}_{\!}\cup\bullet}$\ .

\medskip
\item[\bf 4.]
\label{ExampleP32:4}
(cp. \cite{22}
p.~198 Ex.~16 (Surgery)) If
${\mathbb{E}}^{_{^m}}\!\!:=$ ``the $m$-unit ball'' and $n:=p+q,\ p,q\ge0$ then,
$${\mathbb{S}}^{_{^n}}\!\! = {\rm Bd} {\mathbb{E}}^{^{_{n+1}}}\!\!\!\simeq {\rm Bd} ({\mathbb{E}}^{^{_p}}\!\!\ast {\mathbb{E}}^{^{_q}}\! ) \ \!\!\simeq
{\mathbb{E}}^{^{_{p}}}\!\!\!\ast {\mathbb{S}}^{^{_{q-1}}}\!\!
\cup {\mathbb{S}}^{^{_{p-1}}}\!\! \ast{\mathbb{E}}^{^{_{q}}}
\!\!\simeq
{\rm Bd} ({\mathbb{E}}^{^{_{p+1}}}\!\!\times {\mathbb{E}}^{^{_{q}}})
\!\!\simeq
{\mathbb{E}}^{^{_{p+1}}}\!\!\times\! {\mathbb{S}}^{^{_{q-1}}}\!\!\cup {\mathbb{S}}^{^{_{p}}}\!\!\times {\mathbb{E}}^{^{_{q}}},
$$
by Theorem~\ref{TheoremP31:12} and
Lemma p.~\pageref{LemmaP18}.
\quad
${\mathbb{S}}^{^{_{n+1}}}\simeq
{\mathbb{S}}^{^{_{p}}}\ast {\mathbb{S}}^{^{_{q}}}
\ \textrm{and}\
{\mathbb{E}}^{^{_{n+1}}}\simeq
{\mathbb{E}}^{^{_{p}}}\ast{\mathbb{S}}^{^{_{q}}}$
also hold.

See also \cite{26}
p.\ 376 for some non-intuitive manifold examples.
Also \cite{32} pp.\ 123-131 gives insights on different aspects of different kinds of simplicial manifolds.
\end{varex}

\begin{proposition} \label{proposition:PropP32}\label{PropP32}
If ${{\Sigma}}_{_{^{\!\hbox{\fiverm q}}}} \!$ is finite and
$\ -1\!\le\!\dim{\bf B}_{_{{\!{\mathbf{G},i}\!}}}\!
<\!\dim\Sigma\!-\!1$
then ${\rm Lk}_{_{\!\Sigma}}\!{\delta}$ is
non-orientable$_{_{^{\!\hbox{\fivebf G}}}}\!\!$
\hbox{for all $ {\delta}\!\in\!{\bf B}_{_{{\!{\mathbf{G},i}\!}}}.$}
{\rm(Note that
Cor.~2.{\bf i} p.~\pageref{CorP30:2.i}
is the special case ${\rm
Lk}_{_{\!\Sigma}}\!{\emptyset}_{_{^{\!o}}}\!=\Sigma$.)}
\end{proposition}
\normalbaselines

\begin{proof}
If $\sigma_{{\!}}\in_{{\!}}{\rm Bd}\!_{_{_{\mathbf{G}}}}\!\!\Sigma\
\!^{_{_{\setminus}}}
\mbfcupcap{\cup}{j\not=i}
{\bf
B}_{_{\!\!{j\!}}}\!^{^{\!_{\ }}},$
then
$\dim{\rm Bd}\!_{_{_{\mathbf{G}}}}\!\! ({{\rm
Lk}\!_{_{_{\!\Sigma\!\!}}}\!\sigma})\!_{_{_{{\rm q}}}}\!$
$
\le n-\hbox{\eightrm\#}_{\!}\sigma-3$
in $\overline{\hbox{\rm st}}{\!_{_{\Sigma}}}\!\sigma,$
cf. Note~3 p.~\pageref{NoteP29:3}.
\end{proof}

\subsection{Boundaries of simplicial homology
manifolds}  \label{SubSecP32:III:4}
$$\phantom{.}$$
\vskip-0.8cm
In this section we will work mainly with finite simplicial complexes.

The coefficient module plays, through the
St-R ring functor, a more delicate role in commutative ring theory than it does here in our $\mhatH$omology theory.
So, having a CM-complex, its St-R ring may not be CM if the coefficient module is changed into one that isn't a CM-ring.

\begin{lemma}
\label{LemmaP32:1}
{\bf i.}
$\Sigma$ is a
homology$\!_{_{\mathbf{G}}}\!$ manifold if and only if$\ \!;$
\noindent$\{\Sigma\!=\!\bullet\bullet$
\ or\
$[\Sigma$ is connected and  {\rm Lk}$\!_{_{\Sigma}}\!\!{\bf
v}$ is a finite ${\rm CM}\!_{_{\mathbf{G}}}$
pseudo manifold for all vertices\ ${\bf
v}\in\!{V}_{_{^{\!\Sigma}}}.]$
\underbar{{\bf and}}
\noindent
$\{
{\rm Bd}_{_{\mathbf{G}}}\!\Sigma_{_{^{{\rm q}}}} \!\!=\! {\rm
Bd}_{_{{\mathbb{Z}}}}\!\Sigma_{_{^{{\rm q}}}}\! $
or else;
$[[{\rm Bd}_{{_{{\mathbb{Z}}}}_{_{^{\!\hbox{\fivebf
2}}}}}\!\!\!\Sigma_{_{^{{\rm q}}}}\!\!= {\rm Bd}_{_{\mathbf{G}}}\!\Sigma_{_{^{{\rm
q}}}}\!\!=\emptyset\not=\!\{\emptyset_{_{^{\!o}}}\!\}\!= {\rm
Bd}_{_{{\mathbb{Z}}}}\!\Sigma_{_{^{{\rm q}}}}]
\ or\
[\ \!\exists\  \{\emptyset{_{_{^{\!o\!}}}}\!\}
\nsupseteqq
{\bf B}_{_{\hskip-0.1cm{\mathbb{Z}\!,_{\!}j}}}^{^{_{\Sigma\!\!}}}
\subset
\Sigma^{^{{\!{n\!\lower1.0pt\hbox{-}\hbox{\fiverm 3}}}}}\!\!
\ and\
\hbox{\rm Tor}_1^{\mathbb{Z}}
({\mathbb{Z}}_{_{^{\!2}}}\!,\!{\mathbf{G}})\!=\!{\mathbf{G}}]]
\}\}.$

\smallskip\noindent
{\bf ii.}
For a {\rm CM}$\!_{_{\mathbb{Z}}}\!$ homology$_{_{{\mathbb{Z}}}}\!$
manifold $\Delta$,
${\rm Bd}_{_{\mathbf{G}}}\!\!\Delta_{_{^{\!\hbox{\fiverm h}}}}\!\!=\!
{\rm Bd}_{_{\mathbb{Z}}}\!\Delta_{_{^{\!\hbox{\fiverm h}}}}$
for any module ${\mathbf{G}}$. So, for any homology$_{_{{\mathbb{Z}}}}\!$
manifold $\Sigma$,
${\rm Bd}_{_{\mathbf{G}}}\!\Sigma=\emptyset,\{\emptyset_{_{^{\!o}}}\!\}$
or $\dim{\bf B}_{_{{\!{_{\!}i}\!}}}=(n\!-\!1)$
for boundary components.
\end{lemma}

\begin{proof}
{\bf i.}
Prop.~\ref{PropP22:1} p.~\pageref{PropP22:1},
Lemma~\ref{LemmaP28:2.ii}.{\bf ii} p.~\pageref{LemmaP28:2.ii}
and the proof of Corollary~\ref{CorP30:2}.{\bf v} p.~\pageref{CorP30:2.v.proof}.\\
{\bf ii.}
Use Theorem~\ref{TheoremP14:5} p.~\pageref{TheoremP14:5}
and {\bf i} along with ${\mathbf{G}}_{\!}=_{\!}{\mathbb{Z}}_{_{^{\!2}}}$
in Corollary~\ref{CorP30:2.iii}.{\bf iii} p.~\pageref{CorP30:2.iii}.
\end{proof}
\vskip-10pt

\goodbreak

\begin{lemma}  \label{lemma:LemmaP32:2} \label{LemmaP32:2}
For a finite $\Sigma_{_{^{\!\hbox{\fiverm q}}}}\!$, the boundary homomorphism
$$ \delta\!_{_{n\!_{_{\Sigma}}}}\!:\!
\mhatH\!_{_{n\!_{_{\Sigma}}}}\!\! ( {\Sigma},{\rm
Bd}\raise0.8pt\hbox{$_{_{_{\mathbf{G}}}}$}\!\!{{\Sigma}}; {\mathbf{G}})
\!\rightarrow \mhatH\!\!\!\!_{_{{{n\!_{_{\Sigma}}\!\!-\!1\!\!}}}}
({\Sigma},{\bf B}_{_{\!{i}}}\cap (\
\mbfcupcap{\cup}{j\not=i}
{\bf B}_{_{\!\!{j}}}); {\mathbf{G}})
$$
in the relative {\bf
M-$\!$Vs}$_{_{\!^{o}}}$
with respect to $ \{ ({\Sigma},{\bf
B}_{_{\!{i}}}), ({\Sigma},\
\mbfcupcap{\cup}{j\not=i}
{\bf B}{_{_{\!{j}}}}{_{\!}}) \}
$
is injective if $\#{\bf I}\geq2$ in
Def.~4\nolinebreak\ p.~\pageref{DefP27:4}.

\nobreak So,
$$
\! [\mhatH\!\!\!\!_{_{{{n\!_{_{\Sigma}}\!\!-\!1\!\!}}}}
({\Sigma},\{\emptyset_o\};{\mathbf{G}})\!=\!0] \!\Longrightarrow\!
[ \mhatH\!_{_{n\!_{_{\Sigma}}}}\!\! ( {\Sigma},{\rm
Bd}\!\raise0.6pt\hbox{$_{_{_{\mathbf{G}}}}$}\!\!{{\Sigma}}; {\mathbf{G}})
\!=\!0\ \text{or}\ {\rm Bd}\!\raise0.6pt\hbox{$_{_{_{\mathbf{G}}}}$}\!\!{{\Sigma}}\ \text{is strongly connected}],
$$
as is the case e.g. if $\Sigma\!\ne\!\bullet, \bullet\bullet$ is a {\rm CM}$\!_{_{\mathbf{G}}}\!$ {quasi-n-manifold}.

\smallskip
If $\Sigma$ is a finite {\rm CM}$_{_{\mathbb{Z}}}\!$
quasi-n-manifold then ${\rm Bd}\raise0.8pt\hbox{$_{_{_{\mathbb{Z}}}}$}\!\!{{\Sigma}}= \emptyset$ or it is\ strongly\ connected and
$\dim({\rm Bd}\raise0.8pt\hbox{$_{_{_{{\mathbb{Z}}}}}$}\!\!{{\Sigma}})=
n_{_{\!\Sigma}}\!\!-1$ $
($
\hbox{\rm by Lemma~\ref{LemmaP32:1}}$)$.
\end{lemma}

\begin{proof}
Use the relative {\bf M-$\!$V$_{\!}$s}$_{_{\!^{o}}}\!$ with respect
to $\{ ({\Sigma},{\bf B}_{_{\!{i}}}), ({\Sigma},
\mbfcupcap{\cup}{j\not=i}
{\bf B}_{_{\!{j}}}) \},$ $ \dim({\bf
B}{_{_{\!{i}}}}\cap\
\mbfcupcap{\cup}{j\not=i}
{\bf B}_{_{\!{j}}})\le
n_{_{\!\Sigma}}-3$ and Corollary~\ref{CorP29:2.i}.{\bf i} p.~\pageref{CorP29:2.i}.
\end{proof}

\begin{theorem}{\bf i.}\label{TheoremP33:13}
If $\Sigma$ is a finite orientable$_{_{^{\!\!\hbox{\fivebf
G}}}}\!\!$ {\rm CM}$\!_{_{\mathbf{G}}}\!$
homology$_{_{^{\!\!\hbox{\fivebf G}}}}$ $n$-manifold with boundary
then, ${\rm Bd}_{_{{\!\hbox{\fivebf G}}}}\!\Sigma$\nobreak\ is an orientable$_{_{^{\!\!\hbox{\fivebf G}}}}\!\!$
homology$_{_{^{\!\!\hbox{\fivebf G}}}}$ $(n-1)$-manifold without boundary.

\smallskip
\noindent {\bf ii.} Moreover; ${\rm Bd}_{{_{{\!\hbox{\fivebf
G}}}}}\!\!\Sigma$ is Gorenstein$_{{_{{\!\hbox{\fivebf G}}}}}\!.$

\noindent
{\rm(This theorem is a special case of a theorem, attributed to M. Hochster, given in \cite{31} p.~70 Th.~7.3, where it is proven through the use of the {\it canonical module} of the St-R ring of $\Sigma$)}
\end{theorem}

\begin{proof}
{\bf i.}
Induction over the dimension, using
Theorem~\ref{TheoremP29:11}.{\bf i}-{\bf ii},
once the connectedness
of the\nobreak\ boundary is established through
Lemma~\ref{LemmaP32:2}, while orientability$_{_{^{\!\!\hbox{\fivebf
G}}}}\!$ resp.
 ${\rm Bd}_{_{{\!\hbox{\fivebf
G}}}}\!({\rm Bd}_{_{{\!\hbox{\fivebf G}}}}\!\Sigma)\!=\!\emptyset$
follows from
Corollary~\ref{CorP30:2.iii}.{\bf iii}-{\bf iv} p.~\pageref{CorP30:2.iii}.
\hfill$\triangleright$

\noindent{\bf ii.} $\Sigma\ast(\bullet\bullet)$ is a finite
orientable$_{_{^{\!\!\hbox{\fivebf G}}}}
\!$
{\rm CM}$\!_{_{\mathbf{G}}}\!$
homology$_{_{^{\!\!\hbox{\fivebf G}}}}\!\!$ $(n\!+\!1)$-manifold
with boundary by Th.~\ref{TheoremP31:12} + Th.
~\ref{CorP16:iToTh6}.iii p.~\pageref{CorP16:iToTh6}.
$$\ \!{\rm Bd}_{_{{\!\hbox{\fivebf G}}}}\!(\Sigma\ast\bullet\bullet)
\!=\! \big[{^{\hbox{\scriptsize\rm Th.\ \ref{TheoremP31:12}.}\hbox{\scriptsize\bf
2\ or}}
_{ \hbox{\scriptsize\rm Th.~\ref{TheoremP17:7.2}.}\hbox{\scriptsize\bf 2}}}\big] \!=\!
\big(\Sigma\ast{\rm Bd}_{_{{\!\hbox{\fivebf
G}}}}\!(\bullet\bullet)\big) \cup \big(({\rm Bd}_{_{{\!\hbox{\fivebf
G}}}}\!\Sigma)\ast(\bullet\bullet)\big) =$$
$$= \Sigma\ast\emptyset
\cup ({\rm Bd}_{_{{\!\hbox{\fivebf
G}}}}\!\Sigma)\ast(\bullet\bullet) \!=\! \big({\rm
Bd}_{_{{\!\hbox{\fivebf G}}}}\!\Sigma\big)\ast(\bullet\bullet)$$
where the l.h.s.\ is an orientable$_{_{^{\!\!\hbox{\fivebf
G}}\!\!}}$ homology$_{_{{\!\!_{\!}\hbox{\fivebf G}}\!\!}}$
$n$-manifold without boundary by the first part.
So, ${\rm Bd}\!\raise0.7pt\hbox{$_{_{_{\mathbf{G}}}}$}\!\!\Sigma$ is an
orientable$_{_{^{\!\!\hbox{\fivebf G}}}}\!$
{\rm CM}$\!_{_{\mathbf{G}}}\!$
homology$_{_{^{\!\!\hbox{\fivebf G}}}}$ $\!(n-1)$-manifold without boundary by Th.~\ref{TheoremP31:12} i.e. it is Gorenstein$\!_{_{\mathbf{G}}}\!.$
\end{proof}
\begin{varnote}{\bf$\!$1.} \label{NoteP33:1}
$\emptyset\!\neq\!{{\Delta}}$ is a 2-{\rm CM}$\!_{_{\mathbf{G}}}\!$
hm$\!_{_{\mathbf{G}}}\!$ $\Longleftrightarrow {{\Delta}} \!=\!
\hbox{\rm core}{\Delta}\!$ is Gorenstein$\!_{_{\mathbf{G}}}\!\Longleftrightarrow {{\Delta}}$ is a {\rm
homology}$\!_{_{\mathbf{G}}}\!$ sphere.
\end{varnote}

\vskip-3pt\vskip-3pt
\begin{corollary} \label{corollary:CorP33:1}
{\rm(Cf. \cite{22}
p.\ 190.)}
If $\Sigma$ is a finite orientable$_{_{^{\!\!\hbox{\fivebf G}}}}$ homology$_{_{^{\!\!\hbox{\fivebf G}}}}$ $n$-manifold with boundary,
so is ${\bf2}\Sigma$ except that ${\rm Bd}_{_{{\!\hbox{\fivebf
G}}}}\!({\bf2}\Sigma)=\emptyset$.
${\bf2}\Sigma=\!$``the double of $\Sigma\hbox{''}\!:=\Sigma \mbfcupcap{\cup}{\hskip0.05cm\rm Bd}\ddot{\Sigma}$ where $\ddot{\Sigma}$ is a disjoint copy of $\Sigma$ and
``\mbfcupcap{\cup}{\hskip0.05cm\rm Bd}\hskip-0.1cm''
is ``the union through identification of the boundary vertices''.
If $\Sigma$ is ${\rm CM}_{_{{\!\hbox{\fivebf G}}}}\!$ \hbox{then
\hbox{${\bf2}\Sigma$ is $\hbox{\rm2}\hbox{-}{\rm
CM}_{_{{\!\hbox{\fivebf G}}}}\!.$ }}
\end{corollary}

\begin{proof}
Use \cite{16}
p.\ 57 (23.6) Lemma,
i.e., apply the non-relative augmental {\bf M-$_{\!}$Vs}
\nobreak\ to the pair $({\rm Lk}_{{\Sigma}}\!{\bf v}
, {\rm Lk}_{{\ddot\Sigma}}\!{{\bf v}})$
using Prop.~\ref{PropP34:2.a}.{\bf a} p.~\pageref{PropP34:2.a}
and then to $(\Sigma,\ddot{\Sigma})$ for
the ${\rm CM}_{_{{\!\hbox{\fivebf G}}}}\!$ case.
\end{proof}

\goodbreak
\begin{theorem} \label{TheoremP33:14}
If $\Sigma$ is a finite {\rm CM}$_{_{\mathbb{Z}}}\!$-homology$_{_{{\mathbb{Z}}}}\!\!$ $n$-manifold, then $\Sigma$ is
orientable$_{_{{\mathbb{Z}}}}.$
\end{theorem}

\begin{proof}
$\Sigma$  finite {\rm CM}$_{_{{\mathbb{Z}}}}\! $ $ \Longleftrightarrow
\Sigma\ ${\rm CM}${{_{_{\mathbb{Z}}}}_{_{{\!\hbox{\fivebf
p}}}}}\! $ for all prime fields $ {{{{{\mathbb{Z}}}}}_{_{{\!\hbox{\fivebf
p}}}}}$, i.e. for characteristic ${{\mathbf{p}}}>0$,
(M.A. Reisner, 1976, cf. \cite{Ho} p.~181-2.) by induction over $\dim \Sigma$,
Theorem~\ref{TheoremP14:5} p.~\pageref{TheoremP14:5}
and the {Structure Theorem for Finitely Generated Modules over {\bf PID}s}.
So, $\Sigma$ is a finite {\rm CM}$_{{_{{\mathbb{
Z}}}}_{_{^{\!\hbox{\fivebf p}}}}}\!\!\!$
homology$_{{_{{\mathbb{Z}}}}_{_{^{\!\hbox{\fivebf
p}}}}}\!\!\!$ $n$-manifold for any prime number ${\bf P}$ by
Lemma~\ref{LemmaP32:1}.{\bf i}, since ${\rm
Bd}\raise0.0pt\hbox{$_{_{{\mathbb{Z}}}}$}\!\Sigma
=
{\rm Bd}_{{_{{\mathbb{Z}}}}_{_{^{\!\hbox{\fivebf
p}}}}}\!\!\!\Sigma$ by Lemma~\ref{LemmaP32:1}.{\bf ii} above.
In particular, $
{\rm Bd}_{{_{{\mathbb{Z}}}}_{_{^{\!\hbox{\fivebf
2}}}}}\!\!\!\Sigma$ is Gorenstein$_{{_{{\mathbb{
Z}}}}_{_{^{\!\hbox{\fivebf 2}}}}}\!\!$ by
Lemma~\ref{LemmaP28:1.ii}.{\bf ii} p.~\pageref{LemmaP28:1.ii}
and Theorem~\ref{TheoremP33:13}
above.
If ${\rm Bd}_{_{{\mathbb{Z}}}}\!\Sigma
=
\emptyset$
then $\Sigma$ is orientable$_{_{{\mathbb{Z}}}}$.
Now, if ${\rm Bd}_{{_{{\mathbb{Z}}}}_{_{^{\!\hbox{\fivebf
2}}}}}\!\!\!\Sigma\!\ne\!\emptyset$ then $\dim{\rm
Bd}_{{_{{\mathbb{Z}}}}_{_{^{\!\hbox{\fivebf 2}}}}}\!\!\!
\Sigma\!=\!n\!-\!1$ by
Cor.~\ref{CorP30:2.iii}.{\bf iii}+{\bf iv} p.~\pageref{CorP30:2.iii}.
and, in particular, ${\rm Bd}_{_{{\mathbb{Z}}}}\!\Sigma
=
{\rm Bd}_{{_{{\mathbb{Z}}}}_{_{^{\!\hbox{\fivebf
2}}}}}\!\!\!\Sigma$ is a quasi-$(n\!-\!1)$-manifold.
\normalbaselines

\smallskip
${\rm Bd}_{{_{{\mathbb{Z}}}}_{_{^{\!\hbox{\fivebf 2}}}}}
({\rm Bd}_{{_{{\mathbb{Z}}}}_{_{^{\!\hbox{\fivebf
2}}}}}\!\!\!\Sigma)= \emptyset$ since ${\rm Bd}_{{_{{\mathbb{
Z}}}}_{_{^{\!\hbox{\fivebf 2}}}}}\!\!\!\Sigma$ is
Gorenstein$_{{_{{\mathbb{Z}}}}_{_{^{\!\hbox{\fivebf
2}}}}}\!\!$ and so, $\dim{\rm Bd}_{{_{{\mathbb{Z}}}}} ({\rm
Bd}_{{_{{\mathbb{Z}}}}}\!\Sigma)\le n-4 $ by
Cor.~\ref{CorP30:1} p.~\pageref{CorP30:1}.
So if ${\rm Bd}_{{_{{\mathbb{Z}}}}} ({\rm
Bd}_{{_{{\mathbb{Z}}}}}\!\Sigma)\neq\emptyset$ then, by
Cor.~\ref{CorP30:2.ii}.{\bf ii} p.~\pageref{CorP30:2.ii}.
$ {\rm Bd}_{{_{{\mathbb{Z}}}}}\!\Sigma$ is
nonorientable$_{{_{{\mathbb{Z}}}}}$ i.e.
$$
\mhatH\!\!\!_{_{^{\!n-1}}}\! ({\rm Bd}_{{_{{\mathbb{
Z}}}}}\!\Sigma, {\rm Bd}_{{_{{\mathbb{Z}}}}}({\rm
Bd}_{{_{{\mathbb{Z}}}}}\!\Sigma);{\mathbb{Z}})
\ \!= \ \!
\mhatH\!\!\!_{_{^{\!n-1}}}\! ({\rm Bd}_{{_{{\mathbb{
Z}}}}}\!\Sigma;{\mathbb{Z}})\!=\!0 $$
and the torsion submodule of
$$
\mhatH\!\!\!_{_{^{\!n-2}}}\! ({\rm Bd}_{{_{{\mathbb{
Z}}}}}\!\Sigma, {\rm Bd}_{{_{{\mathbb{Z}}}}}({\rm
Bd}_{{_{{\mathbb{Z}}}}}\!\Sigma);{\mathbb{Z}}) =
\big[{^{\hbox{\scriptsize\rm For dimensio-}}_{\hbox{\scriptsize\rm \ nal
reasons.}}}\big]
=\!\mhatH\!\!\!\!_{_{^{\!n-2}}}\! ({\rm
Bd}_{{_{{\mathbb{Z}}}}}\!\Sigma;{\mathbb{Z}}) $$
is isomorphic to
${\mathbb{Z}}_{_{^{\!\hbox{\fivebf 2}}}}\!$ by
Lemma~\ref{LemmaP28:1.i}.{\bf i} p.~\pageref{LemmaP28:1.i}.
In particular $ \mhatH\!\!\!\!_{_{^{\!n-2}}}\! ({\rm
Bd}_{{_{{\mathbb{Z}}}}}\!\Sigma;{\mathbb{Z}}) \otimes {\mathbb{Z}}_{_{^{\!\hbox{\fivebf 2}}}}\!\ne\!0$ implying, by
Th.~\ref{TheoremP14:5} p.~\pageref{TheoremP14:5},
that
$$ \mhatH\!\!\!\!_{_{^{\!n-2}}}\! ({\rm Bd}{_{_{{\mathbb{
Z}}}}}_{_{^{\!\hbox{\fivebf 2}}}}\!\!\!\Sigma;{\mathbb{Z}}_{_{^{\!\hbox{\fivebf 2}}}}) \!= \mhatH\!\!\!\!_{_{^{\!n-2}}}\!
({\rm Bd}{_{_{{\mathbb{Z}}}}}\!\Sigma;{\mathbb{Z}}_{_{^{\!\hbox{\fivebf 2}}}}) \!=
\mhatH\!\!\!\!_{_{^{\!n-2}}}\!({\rm Bd}{_{_{{\mathbb{
Z}}}}}\!\Sigma;{\mathbb{Z}}) \otimes {\mathbb{Z}}_{\bf 2}\! \oplus \hbox{\rm Tor}^{^{_{\mathbb{Z}}}}_{\bf 1}
\big(\mhatH\!\!\!_{_{^{\!n-3}}}\! ({\rm Bd}{_{_{{\mathbb{
Z}}}}}\!\Sigma;{\mathbb{Z}}),{\mathbb{Z}}_{\bf 2}\!
\big) \!\ne 0 $$
contradicting the Gorenstein$_{_{{\!{\mathbb{Z}}}_{\bf 2}}}\!\!$-ness of ${\rm
Bd}\raise0.5pt\hbox{$_{_{{\!{\mathbb{Z}}}_{\bf 2}}}$}\!\!\Sigma.$
\end{proof}

\begin{corollary} \label{corollary:CorP33:2}\label{CorP33:2}
Each simplicial homology$_{_{{\mathbb{Z}}}}\!$ manifold $(\mathrm{hm}_{_{{\mathbb{Z}}}})$ is locally
orientable.
\end{corollary}
\begin{proof} Proposition~\ref{PropP22:1} p.~\pageref{PropP22:1}.\end{proof}
\quad

\begin{varnote}{\bf$\!$2.} \label{NoteP33:2}
Corollary\ \ref{corollary:CorP33:2} confirms G.E. Bredon's
conjecture stated in \cite{1} Remark p. 384 just after the
definition of Borel-Moore cohomology manifolds with boundary, the
reading of which is aiming at Poincar\'{e} duality and therefore
only allows $\hbox{\rm Bd}{\eighti X}\!=\!\emptyset$ or
$\dim\hbox{\rm Bd}{\eighti
X}\!=\!n-1$, cf. our Ex.~3 p.~\pageref{ExampleP32:3}.
Weak homology manifolds over {\bf PID}s are defined in \cite{1} p.\
329 and again we would like to draw the attention to their
connection to Buchsbaum rings, cf. our
p.~\pageref{DefP22},
and that, for polytopes, ``join'' becomes ``tensor product'' under the
St-R ring functor as mentioned in
\hbox{our Note~{\bf iii} p.~\pageref{NoteP21:iii}.}
\end{varnote}

\vfill\break
\part{Appendices and Addenda} \label{Appendices&Addendum}

\section{Appendices} \label{Appendices}

\subsection{The 3$\times$3-lemma (also called ``the 9-lemma'')}  \label{3x3-lemma}

 In this appendix we will make use of the relative
Mayer-Vietoris
sequence ({\bf M-$\!$Vs}) to form a 9-Lemma-grid
{\rm clarifying some} {\rm relations between Products and Joins,
e.g. that $\{X_{_{^{\!1}}} \ \!\!\mhatast Y_{_{^{\!2}}},
X_{_{^{\!2}}} \ \!\!\mhatast Y_{_{^{\!1}}}\}$ is excisive
\underbar{iff} \hbox{$\{X_{_{^{\!1}}}\!\!\times\!\! Y_{_{^{\!2}}},
X_{_{^{\!2}}}\!\!\times\!\! Y_{_{^{\!1}}}\}$ is.}

\scriptsize

\medskip
\noindent
\phantom{.}
\framebox
{\hsize=15.4cm
{
\noindent \hbox{\noindent
\vbox{\hsize=0.82\hsize
{\eightrm\noindent\cite{Rotman:HomAlg} {\eightbf Ex. 6.16} p.~175-6 (3$\times$3 Lemma)\ \
{\eightrm (J.J. Rotman: {\eighti An} {\eighti Intr}. {\eighti to}
{\eighti Homological} {\eighti Algebra}.)}
Consider the commutative diagram of modules (on the right): If the
columns are exact and if the bottom two (or the top two) rows
are exact, then the top row (or the bottom row) is exact.
(Hint: Either use the Snake lemma or proceed as follows: first show
that $\alpha\alpha^{{_{\prime}}}\!=0,$ then regard each row as a
complex and the diagram as a short exact sequence of complexes, and
apply theorem\ 6.3.)
\hrule \vskip 1pt \hrule
\vskip 0.5pt
\noindent
{\eightrm \noindent
{{\eightrm \cite{21}} {\eightbf Ex.~5,} {\eightrm p.~208.}  (S. Mac$\ \!$Lane: {\eighti Categories} {\eighti for} {\eighti the}
{\eighti Working} {\eighti Mathematician}.)}
\hfill\break
\indent A 3$\times$3 diagram is one of the  form (to the right)
(bordered by zeros).\hfill\break
{\eightbf(a)} Give a direct proof of the 3$\times$3 lemma: If a
3$\times$3 diagram is commutative and all three columns and the last
(first) two rows are short exact sequences, then so is the first
(last) row.
\hfill\break
{\eightbf(b)}  Show that this lemma also follows from the ker-coker sequence.\hfill\break
{\eightbf(c)}
Prove the middle 3$\times$3 lemma: If a 3$\times$3 diagram is commutative, and all three columns and the first and third rows are short exact sequences, then so is the middle row.
\hrule \vskip 1pt \hrule
\vskip 0.5pt
\noindent
\cite{17} {\eightbf Ex. 16}\ p.~227 {\eightrm (P.J. Hilton and S. Wylie: {\eighti Homology} {\eighti Theory}.)}
Let (the diagram on the right) be a commutative diagram in which
each row and each column is an exact sequence of differential groups
(Def. p. 99). Then there are defined (see 5.5.1 (p.\ 196))
homomorphism ${\nu}\!_{_{^{\ast}}}\!\!^{{_{C}}}\!\!:\ \! ${\eighti
{\char"48}}({\eighti {\char"43}}\lower2pt\hbox{\fiverm3})
$\rightarrow$ {\eighti {\char"48}}({\eighti
{\char"43}}\lower2pt\hbox{\fiverm1}),
${\nu}\!_{_{^{\ast}}}\!\!^{{_{1}}}\!\!:\ \! ${\eighti
{\char"48}}({\eighti {\char"43}}\lower2pt\hbox{\fiverm1})
$\rightarrow$ {\eighti {\char"48}}({\eighti
{\char"41}}\lower2pt\hbox{\fiverm1}),
${\nu}\!_{_{^{\ast}}}\!\!^{{_{3}}}\!\!:\ \! ${\eighti
{\char"48}}({\eighti {\char"43}}\lower2pt\hbox{\fiverm3})
$\rightarrow$ {\eighti {\char"48}}({\eighti
{\char"41}}\lower2pt\hbox{\fiverm3}),
${\nu}\!_{_{^{\ast}}}\!\!^{{_{A}}}\!\!:\ \! ${\eighti
{\char"48}}({\eighti {\char"41}}\lower2pt\hbox{\fiverm3})
$\rightarrow$ {\eighti {\char"48}}({\eighti
{\char"41}}\lower2pt\hbox{\fiverm1}).
Prove that $ {\nu}\!_{_{^{\ast}}}\!\!^{{_{C}}}
{\nu}\!_{_{^{\ast}}}\!\!^{{_{1}}} = -
{\nu}\!_{_{^{\ast}}}\!\!^{{_{3}}} {\nu}\!_{_{^{\ast}}}\!\!^{{_{A}}}.
$ (Reversed notation.)
\hrule \vskip 1pt \hrule \vskip 0.5pt \noindent\cite{6} {\eightbf Ex.~2} p.~53
{\eightrm (A. Dold: {\eighti Lectures} {\eighti on}
{\eighti Algebraic} {\eighti Topology}.)}
...if (the diagram on the right is)
a commutative diagram of chain maps with exact rows and columns then
the homology sequence of these rows and columns constitute a
2-dimensional lattice of group homomorphisms. It is commutative
except for the ({\eighti \char"40}$\!_{_{^{\ast}}}\!
$\hbox{\eightbf-}{\eighti \char"40}$\!_{_{^{\ast}}}\!)$-squares
which anticommute. } } } } }

\noindent
$\!\!\!\!\!\!$
{ {
\noindent\hskip-0.3cm\vbox{\hsize=0.25\hsize
{
\vbox{\hsize=0.95\hsize
{$\hskip0.3cm 0 \hskip0.6cm 0\
\hskip0.55cm 0
\hfill\break
\vbox{\vskip 9pt}
\noindent\hskip0.4cm
\hskip0.3cm \big\downarrow
\hskip0.6cm  \big\downarrow
\hskip0.6cm \big\downarrow $
}}
\vbox{\hsize=0.8\hsize
{
$
\noindent\hskip-0.2cm 0\ \! {{\hbox to 0.3cm{\rightarrowfill}}}
A_{_{^{\!1}}}\!\!\! \rlap{${\hbox to 0.45cm{\rightarrowfill}}$}
{\raise4pt\hbox{$\ \alpha^{{_{\prime}}}$}}
\hskip0.1cm A_{_{^{2}}}\!\ \!\! \rlap{${\hbox to
0.4cm{\rightarrowfill}}$} {\raise4pt\hbox{$\ \alpha$}}
\hskip0.2cm A_{_{^{3}}}\hskip-0.1cm {{\hbox to 0.4cm{\rightarrowfill}}} \
_{\!}0
\hfill\break%
\vbox{\vskip 9pt}
\noindent\hskip0.4cm%
\hskip0.3cm \big\downarrow
\hskip0.6cm  \big\downarrow
\hskip0.6cm \big\downarrow $ }}
\vbox{\hsize=0.97\hsize { \vbox{\vskip 1pt} $
\noindent\hskip-0.2cm 0\ \! {{\hbox to 0.3cm{\rightarrowfill}}}
B_{_{^{\!1}}}\! {{\hbox to 0.4cm{\rightarrowfill}}} \
_{\!}B{_{_{^{\!2}}}}\!_{\!} {{\hbox to 0.4cm{\rightarrowfill}}} \
\!_{\!}B{_{_{^{\!3}}}} _{\!} {{\hbox to 0.4cm{\rightarrowfill}}}  0
\hfill\break%
\vbox{\vskip 9pt}
\noindent\hskip0.4cm%
\hskip0.3cm \big\downarrow
\hskip0.6cm  \big\downarrow
\hskip0.6cm \big\downarrow $ }}
\vbox{\hsize=0.97\hsize { \vbox{\vskip 1pt} $
\noindent\hskip-0.2cm%
0\ \! {{\hbox to 0.3cm{\rightarrowfill}}} C_{_{^{\!1}}}\! {{\hbox to
0.4cm{\rightarrowfill}}} \ _{\!}C{_{_{^{\!2}}}}\!_{\!} {{\hbox to
0.4cm{\rightarrowfill}}} \ \!_{\!}C{_{_{^{\!3}}}} _{\!} {{\hbox to
0.4cm{\rightarrowfill}}}  0
\hfill\break%
\vbox{\vskip 9pt}
\noindent\hskip0.39cm%
\hskip0.31cm \big\downarrow \phantom{{\hbox to
0.51cm{\rightarrowfill}}\ \ \!} \big\downarrow
\hskip0.6cm %
\big\downarrow
$ }}
\noindent\hskip-0.1cm\vbox{\hsize=0.95\hsize
{ $
\hskip0.4cm  0\
\hskip0.55cm 0
 \hskip0.58cm 0\
 \hskip0.5cm
 $ }} \vskip0cm
\vbox{\vskip 1.5cm}%
}}}}\hskip-0.7cm}}

\normalsize
\smallskip
\noindent%
{\tenrm Put} $\left\{\hskip-0.0cm\begin{array}{cc}
{\bf \hat C}X\!\!:=\!\!X \mhatast
\{(x\!_{_0}\!,1)\}\ {\rm where}\
\underline{x\!_{_0}\!\!\in\!X_{_{^{\!2}}}\!}\ {\rm and}\
(x\!_{_0}\!,1):=\widetilde{(x\!_{_0}\!,y\!_{_0}\!,1)}\in X\mhatast Y\\
{\bf \hat C}Y\!\!:=\!Y \mhatast
\{(y\!_{_0}\!,0)\}\ {\rm where}\
\underline{y\!_{_0}\!\!\in\!Y_{_{^{\!2}}}\!}\ {\rm and}\
(y\!_{_0}\!,0):=\widetilde{(x\!_{_0}\!,y\!_{_0}\!,0)}\in X\mhatast Y%
\end{array}\right.$%
\hskip-0.1cm {\small($\mhatast$-def. p.~\pageref{EquivJoinDef(x,1)}).}

\smallskip
\tiny
\noindent
{\tenrm Put} $\left\{\hskip-0.2cm\begin{array}{cc}
\left\{\hskip-0.2cm\begin{array}{ll} %
A:={\bf \hat C} X_{_{^{\!1}}}\!\times\!Y_{_{^{\!2}}}
\\
B:={\bf \hat C}X_{_{^{\!2}}}\!\!\times\!\! Y_{_{^{\!1}}}%
\end{array}\right.
\\
\left\{\hskip-0.2cm\begin{array}{cc}
C:=X_{_{^{\!1}}}\!\!\times\!{\bf \hat C}Y_{_{^{\!2}}}\\
D:=X_{_{^{\!2}}}\!\!\times\!{\bf \hat C}Y_{_{^{\!1}}}%
\end{array}\right.%
\end{array}\right.$
%
%
\noindent\hskip-0.2cm\hbox{\hsize6cm{\tenrm then}\
$
\left\{\hskip-0.2cm\begin{array}{ccccc}%
\smallskip%
\left\{\begin{array}{cc}A\cap B\cap C\cap D=\\
=(A\cap D)\cup(B\cap C)
\end{array}%
\right\}&%
\hskip-0.40cm=X_{_{^{\!2}}}\!\!\times\! Y_{_{^{\!2}}}\hskip3.5cm\\%
\smallskip%
\left\{\begin{array}{cc}(A\cap B)\cup(C\cap D)=\\
=(A\cup C)\cap(B\cup D)
\end{array}%
\right\}&%
\hskip-0.30cm=X_{_{^{\!2}}} \ \!\!\mtoppast Y_{_{^{\!2}}} \hskip3.6cm\\
\smallskip%
\left\{\begin{array}{cc}(A\cup B)\cap(C\cup D)=\\
=(A\cap C)\cup(B\cap D)
\end{array}%
\right\}&%
\hskip-0.45cm=X_{_{^{\!1}}}\!\!\times\!\! Y_{_{^{\!2}}}\!\!\cup
\!\!X_{_{^{\!2}}}\!\!\times\!\! Y_{_{^{\!1}}}\hskip2.7cm\\
\smallskip%
\left\{\hskip-0.0cm\begin{array}{cc}%
(A\cup B)\cup(C\cup D)=\\
=(A\cup C)\cup(B\cup D)
\end{array}%
\right\}&%
\hskip-0.45cm=X_{_{^{\!1}}} \ \!\!\mtoppast
Y_{_{^{\!2}}}\cup X_{_{^{\!2}}}\mtoppast Y_{_{^{\!1}}}\hskip2.5cm\\
\hskip-0.5cm(A\cup D)\cap(B\cup C)\phantom{\Big\} } &\hskip-1.6cm=
X_{_{^{\!2}}} \ \!\!\mtoppast Y_{_{^{\!2}}}\cup
[X_{_{^{\!1}}}\!\!\times\!\! Y_{_{^{\!2}}}\!\!\cup
\!\!X_{_{^{\!2}}}\!\!\times\!\! Y_{_{^{\!1}}}]\hskip1.5cm
\end{array}
\right.%
$
\hskip-3.0cm%
{\large$\Longrightarrow$}
%
%
%
\raise0.9cm\vbox{
\def\objectstyle{\tiny}
\def\labelstyle{\tiny}
\xymatrix{\ar @{>->} [dr] |{%
\vbox{\hsize2cm\ \ \\ \hbox{$(A\cup D)\cap(B\cup C)$}}
}
\hbox{${{(A\cap C)\cap(B\cap D)=} \atop {=(A\cap B)\cup(C\cap D)}}$}
\ar[r]\ar[d] %
&
\hbox{${ {(A\cup C)\cap(B\cup D)=} \atop {=(A\cap B)\cup(C\cap D)}
}$}
\ar @{>-<} [dl] |{\phantom{\vbox{\hsize2cm\ \ \\ \hbox{$(A\cup D)\cap(B\cup C)$}}}}\ar[d] \\
\hbox{${{(A\cap C)\cup(B\cap D)=} \atop {=(A\cup B)\cap(C\cup D)}}$}
\ar[r] &
\hbox{${{(A\cup C)\cup(B\cup D)=} \atop {=(A\cup B)\cup(C\cup D)}}$}
}%
}%
}%

\normalsize

\smallskip
All operations are assumed to take place within
$X_{1}\mhatast Y_{1}$.

The last commutative square provide us with a ``$9$-lemma''-grid as
that above, constituted by two horisontal and two vertical relative
Mayer-Vietoris sequences,
using the naturality of the {\bf M-Vs} as in Munkres
\cite{26} p.~186-7.
\newpage
\enlargethispage{6.0cm}
\begin{landscape}%
\vbox{
\noindent\hskip-0.7cm
\hbox{
The entries consists of three levels where the upper (lower)
concerns the vertical (horisontal) {\bf M-$\!$Vs}.}
$X_2$,
$Y_2$ are assumed to be non-empty.

%

\noindent\hskip-1cm
\vbox{
$\hskip3.3cm {\downarrow}\delta_{3\ast} \hskip2.0cm(-\#)\hskip2.2cm
{\downarrow}\delta_{1\ast} \hskip2.3cm\#\hskip2.5cm{\downarrow}
\delta_{2\ast}\hskip3.0cm\# $

\medskip
\noindent$\lower4.5pt\hbox{$^{
 ...{ \buildrel \hbox{$0$} \over \longrightarrow}
\hskip0.6cm\phantom{\ }
{ { { {
\hbox{${\mhatH\!\!\!\!\!_{_{_{n+1}}}}$}\!\!\big(
\Xit_{_{^{\!1}}}\!\mtopast\! \Yit_{_{^{\!1}}}\!, {\hbox{\sevenbf((}}{\hat \Cbf} {
\Xit_{_{^{\!1}}}}\!\!\times\!\Yit_{_{^{\!2}}}{\hbox{\sevenbf)}}\cup
{\hbox{\sevenbf(}}\Xit_{_{^{\!1}}}\!\!\times\! {\hat \Cbf}{\Yit_{_{^{\!2}}}}{\hbox{\sevenbf))}}\cap}
\atop {\cap {\hbox{\sevenbf((}} {\hat \Cbf}
{\Xit_{_{^{\!2}}}}\!\!\times\!\Yit_{_{^{\!1}}}{\hbox{\sevenbf)}}\cup
{\hbox{\sevenbf(}}\Xit_{_{^{\!2}}}\!\!\times\! {\hat \Cbf} {\Yit_{_{^{\!1}}}}{\hbox{\sevenbf))}}\big)\ =
} } \atop {=\hbox{${\mhatH\!\!\!\!\!_{_{_{n+1}}}}$}\!\!
\big(\Xit_{_{^{\!1}}}\!\mtopast \Yit_{_{^{\!1}}}\!,
\Xit_{_{^{\!2}}}\!\mtopast \Yit_{_{^{\!2}}}\big) =} } \atop { {
=\hbox{${\mhatH\!\!\!\!\!_{_{_{n+1}}}}$}\!\!\big(
\Xit_{_{^{\!1}}}\!\mtopast\! \Yit_{_{^{\!1}}}\!, {\hbox{\sevenbf((}}{\hat \Cbf}{
\Xit_{_{^{\!1}}}}\!\!\times\!\Yit_{_{^{\!2}}}{\hbox{\sevenbf)}}\cap {\hbox{\sevenbf(}}{\hat \Cbf}
\Xit_{_{^{\!2}}}\!\!\times\!
{\Yit_{_{^{\!1}}}}{\hbox{\sevenbf))}}\cup}
\atop {\cup \!{\hbox{\sevenbf((}}  { \Xit_{_{^{\!1}}}}\!\!\times\!{\hat \Cbf}
\Yit_{_{^{\!2}}}{\hbox{\sevenbf)}}\cap {\hbox{\sevenbf(}} \Xit_{_{^{\!2}}}\!\!\times\! {\hat \Cbf}
{\Yit_{_{^{\!1}}}}{\hbox{\sevenbf))}}\big)
}} }
\hskip0.3cm\!{\buildrel \hbox{$\delta$}_\ast \over {\hbox to
0.7cm{\rightarrowfill}}}\hskip0.3cm
{ { { {
\hbox{${\mhatH\!_{_{n}}}$}\!\!\big( \Xit_{_{^{\!1}}}\!\times\!
\Yit_{_{^{\!1}}}\!, {\hbox{\sevenbf((}}{\hat \Cbf}{
\Xit_{_{^{\!1}}}}\!\!\times\!\Yit_{_{^{\!2}}}{\hbox{\sevenbf)}}\cap
{\hbox{\sevenbf(}}\Xit_{_{^{\!1}}}\!\!\times\! {\hat \Cbf}
{\Yit_{_{^{\!2}}}}{\hbox{\sevenbf))}}\cap}
\atop {\cap {\hbox{\sevenbf((}} {\hat \Cbf}
{\Xit_{_{^{\!2}}}}\!\!\times\!\Yit_{_{^{\!1}}}{\hbox{\sevenbf)}}\cap
{\hbox{\sevenbf(}}\Xit_{_{^{\!2}}}\!\!\times\! {\hat \Cbf}
{\Yit_{_{^{\!1}}}}{\hbox{\sevenbf))}}\big)\ = } } \atop
{=\hbox{${\mhatH\!_{_{n}}}$}\!\! \big(\Xit_{_{^{\!1}}}\!\times
\Yit_{_{^{\!1}}}\!, \Xit_{_{^{\!2}}}\!\times \Yit_{_{^{\!2}}}\big)
=} } \atop { {
=\hbox{${\mhatH\!_{_{n}}}$}\!\!\big( \Xit_{_{^{\!1}}}\!\times\!
\Yit_{_{^{\!1}}}\!, {\hbox{\sevenbf((}}{\hat \Cbf}{
\Xit_{_{^{\!1}}}}\!\!\times\!\Yit_{_{^{\!2}}}{\hbox{\sevenbf)}}\cap
{\hbox{\sevenbf(}}{\hat \Cbf} \Xit_{_{^{\!2}}}\!\!\times\!
{\Yit_{_{^{\!1}}}}{\hbox{\sevenbf))}}\cap}
\atop {\cap \!{\hbox{\sevenbf((}}  {
\Xit_{_{^{\!1}}}}\!\!\times\!{\hat \Cbf}
\Yit_{_{^{\!2}}}{\hbox{\sevenbf)}}\cap
{\hbox{\sevenbf(}}\Xit_{_{^{\!2}}}\!\!\times\! {\hat \Cbf}
{\Yit_{_{^{\!1}}}}{\hbox{\sevenbf))}}\big)
}} }
\indent{\hbox to 0.5cm{\rightarrowfill}}\ \ \
{ { { {
\hbox{${\mhatH\!_{_{n}}}$}\!\!\big( {\hat \Cbf}
\Xit_{_{^{\!1}}}\!\times\! \Yit_{_{^{\!1}}}\!,
{\hbox{\sevenbf(}}{\hat \Cbf}{
\Xit_{_{^{\!1}}}}\!\!\times\!\Yit_{_{^{\!2}}}{\hbox{\sevenbf)}}\cap
{\hbox{\sevenbf(}}{\hat \Cbf} \Xit_{_{^{\!2}}}\!\!\times\!
{\Yit_{_{^{\!1}}}}{\hbox{\sevenbf)}}\big)\oadd}
\atop {\oadd \hbox{${\mhatH\!_{_{n}}}$}\!\!\big( \Xit_{_{^{\!1}}}
\!\times\! {\hat \Cbf}\Yit_{_{^{\!1}}}\!,
{\hbox{\sevenbf(}}{\Xit_{_{^{\!1}}}}\!\!\times\!{\hat
\Cbf}\Yit_{_{^{\!2}}}{\hbox{\sevenbf)}}\cap
{\hbox{\sevenbf(}}\Xit_{_{^{\!2}}}\!\!\times\! {\hat \Cbf}
{\Yit_{_{^{\!1}}}}{\hbox{\sevenbf)}}\big)\ = } } \atop
{=\hbox{${\mhatH\!_{_{n}}}$}\!\! \big(\Yit_{_{^{\!1}}}\!,
\Yit_{_{^{\!2}}}\big) \oadd \hbox{${\mhatH\!_{_{n}}}$}\!\!
\big(\Xit_{_{^{\!1}}}\!, \Xit_{_{^{\!2}}}\big) =} } \atop { {
=\hbox{${\mhatH\!_{_{n}}}$}\!\!\big( {\hat \Cbf}
\Xit_{_{^{\!1}}}\!\times\! \Yit_{_{^{\!1}}}\!,
{\hbox{\sevenbf(}}{\hat \Cbf}{
\Xit_{_{^{\!1}}}}\!\!\times\!\Yit_{_{^{\!2}}}{\hbox{\sevenbf)}}\cap
{\hbox{\sevenbf(}}{\hat \Cbf} \Xit_{_{^{\!2}}}\!\!\times\!
{\Yit_{_{^{\!1}}}}{\hbox{\sevenbf)}}\big)\oadd}
\atop {\oadd \hbox{${\mhatH\!_{_{n}}}$}\!\!\big( \Xit_{_{^{\!1}}}
\!\times\! {\hat \Cbf}\Yit_{_{^{\!1}}}\!, {\hbox{\sevenbf(}}{
\Xit_{_{^{\!1}}}}\!\!\times\!{\hat
\Cbf}\Yit_{_{^{\!2}}}{\hbox{\sevenbf)}}\cap {\hbox{\sevenbf(}}
\Xit_{_{^{\!2}}}\!\!\times\! {\hat \Cbf}
{\Yit_{_{^{\!1}}}}{\hbox{\sevenbf)}}\big) }} }
\hskip0.8cm{\vbox{\hbox{$ {\buildrel  \hbox{$0$} \over {\hbox to
0.7cm{\rightarrowfill}}}
$}}} }$}$

$ \hskip3.3cm{\Bigg\downarrow}(i_{_3},-j_{_3}) \hskip1.6cm\#
\hskip2.3cm\Bigg\downarrow(i_{_1},-j_{_1})\hskip1.3cm
\#\hskip2.3cm\Bigg\downarrow(i_{_2},-j_{_2}) \hskip2.3cm \#$

\noindent$\lower4.5pt\hbox{$^{
\raise5pt\hbox{...$\buildrel \hbox{$0$} \over {\hbox to
0.5cm{\rightarrowfill}}$}
{ { { {
\hbox{${\mhatH\!\!\!\!\!_{_{_{n+1}}}}$}\!\!\big(
\Xit_{_{^{\!1}}}\!\mtopast\! \Yit_{_{^{\!1}}}\!, {\hat \Cbf}{
\Xit_{_{^{\!1}}}}\!\!\times\!\Yit_{_{^{\!2}}}\cup
\Xit_{_{^{\!1}}}\!\!\times\! {\hat \Cbf}
{\Yit_{_{^{\!2}}}}\big)\oadd}
\atop {\oadd \hbox{${\mhatH\!\!\!\!\!_{_{_{n+1}}}}$}\!\!\big(
\Xit_{_{^{\!1}}} \!\mtopast\! \Yit_{_{^{\!1}}}\!, {\hat \Cbf}{
\Xit_{_{^{\!2}}}}\!\!\times\!\Yit_{_{^{\!1}}}\cup
\Xit_{_{^{\!2}}}\!\!\times\! {\hat \Cbf} {\Yit_{_{^{\!1}}}}\big)=} }
\atop {=\lower0pt\hbox{${\mhatH\!\!\!\!\!_{_{_{n+1}}}}$}\!\!
\big( \Xit_{_{^{\!1}}} \!\mtopast\! \Yit_{_{^{\!1}}}\!,
\Xit_{_{^{\!1}}}\!\mtopast\! \Yit_{_{^{\!2}}}\big) \oadd
\lower5pt\hbox{$^{\mhatH\!\!\!\!\!_{_{n+1}}}$}\!\! \big(
\Xit_{_{^{\!1}}} \!\mtopast\! \Yit_{_{^{\!1}}}\!, \Xit_{_{^{\!2}}}
\!\mtopast\! \Yit_{_{^{\!1}}}\big) =} } \atop {{
=\hbox{${\mhatH\!\!\!\!\!_{_{_{n+1}}}}$}\!\!\big(
\Xit_{_{^{\!1}}}\!\mtopast\! \Yit_{_{^{\!1}}}\!, {\hat \Cbf}{
\Xit_{_{^{\!1}}}}\!\!\times\!\Yit_{_{^{\!2}}}\cup
\Xit_{_{^{\!1}}}\!\!\times\! {\hat \Cbf}
{\Yit_{_{^{\!2}}}}\big)\oadd}
\atop {\oadd \hbox{${\mhatH\!\!\!\!\!_{_{_{n+1}}}}$}\!\!\big(
\Xit_{_{^{\!1}}} \!\mtopast\! \Yit_{_{^{\!1}}}\!, {\hat \Cbf}{
\Xit_{_{^{\!2}}}}\!\!\times\!\Yit_{_{^{\!1}}}\cup
\Xit_{_{^{\!2}}}\!\!\times\! {\hat \Cbf} {\Yit_{_{^{\!1}}}}\big)
} } }
\raise5pt\hbox{$\buildrel \hbox{$\delta$}_\ast \over {\hbox to
0.4cm{\rightarrowfill}}$}
{{{{
\hbox{${\mhatH\!_{_{n}}}$}\!\!\big( \Xit_{_{^{\!1}}}\!\times\!
\Yit_{_{^{\!1}}}\!, {\hbox{\sevenbf(}}{\hat \Cbf}{
\Xit_{_{^{\!1}}}}\!\!\times\!\Yit_{_{^{\!2}}}{\hbox{\sevenbf)}}\cap
{\hbox{\sevenbf(}}\Xit_{_{^{\!1}}}\!\!\times\! {\hat \Cbf}
{\Yit_{_{^{\!2}}}}{\hbox{\sevenbf)}}\big)\oadd}
\atop {\oadd \hbox{${\mhatH\!_{_{n}}}$}\!\!\big( \Xit_{_{^{\!1}}}
\!\times\! \Yit_{_{^{\!1}}}\!, {\hbox{\sevenbf(}}{\hat \Cbf}{
\Xit_{_{^{\!2}}}}\!\!\times\!\Yit_{_{^{\!1}}}{\hbox{\sevenbf)}}\cap
{\hbox{\sevenbf(}} \Xit_{_{^{\!2}}}\!\!\times\! {\hat \Cbf}
{\Yit_{_{^{\!1}}}}{\hbox{\sevenbf)}}\big)\ = } } \atop
{=\hbox{${\mhatH\!_{_{n}}}$}\!\! \big(\Xit_{_{^{\!1}}}\!\times\!
\Yit_{_{^{\!1}}}\!, \Xit_{_{^{\!1}}}\!\times\!\Yit_{_{^{\!2}}}\big)
\oadd \hbox{${\mhatH\!_{_{n}}}$}\!\!
\big(\Xit_{_{^{\!1}}}\!\!\times\!\Yit_{_{^{\!1}}}\!,
\Xit_{_{^{\!2}}}\!\times\! \Yit_{_{^{\!1}}}\big) =} } \atop { {
=\hbox{${\mhatH\!_{_{n}}}$}\!\!\big( \Xit_{_{^{\!1}}}\!\times\!
\Yit_{_{^{\!1}}}\!, {\hbox{\sevenbf(}}{\hat \Cbf}{
\Xit_{_{^{\!1}}}}\!\!\times\!\Yit_{_{^{\!2}}}{\hbox{\sevenbf)}}\cap
{\hbox{\sevenbf(}}\Xit_{_{^{\!1}}}\!\!\times\! {\hat \Cbf}
{\Yit_{_{^{\!2}}}}{\hbox{\sevenbf)}}\big)\oadd}
\atop {\oadd \hbox{${\mhatH\!_{_{n}}}$}\!\!\big( \Xit_{_{^{\!1}}}
\!\times\! \Yit_{_{^{\!1}}}\!, {\hbox{\sevenbf(}}{\hat \Cbf}{
\Xit_{_{^{\!2}}}}\!\!\times\!\Yit_{_{^{\!1}}}{\hbox{\sevenbf)}}\cap
{\hbox{\sevenbf(}}\Xit_{_{^{\!2}}}\!\!\times\! {\hat \Cbf}
{\Yit_{_{^{\!1}}}}{\hbox{\sevenbf)}}\big) } } }
\raise5pt\hbox{$\buildrel \hbox{$\alpha$}_{\ast}^{^\prime} \over
{\hbox to 0.4cm{\rightarrowfill}}$}
{{{{
\hbox{${\mhatH\!_{_{n}}}$}\!\! \big({\hat \Cbf}
\Xit_{_{^{\!1}}}\!\times\! \Yit_{_{^{\!1}}}\!, {\hat \Cbf}
\Xit_{_{^{\!1}}}\!\times\!\Yit_{_{^{\!2}}}\big) \oadd
\hbox{${\mhatH\!_{_{n}}}$}\!\! \big({\hat \Cbf}
\Xit_{_{^{\!1}}}\!\!\times\!\Yit_{_{^{\!1}}}\!, {\hat \Cbf}
\Xit_{_{^{\!2}}}\!\times\! \Yit_{_{^{\!1}}}\big)
\oadd}
\atop {\oadd \hbox{${\mhatH\!_{_{n}}}$}\!\!
\big(\Xit_{_{^{\!1}}}\!\times\! {\hat \Cbf}\Yit_{_{^{\!1}}}\!,
\Xit_{_{^{\!1}}}\!\times\!{\hat \Cbf}\Yit_{_{^{\!2}}}\big) \oadd
\hbox{${\mhatH\!_{_{n}}}$}\!\!
\big(\Xit_{_{^{\!1}}}\!\!\times\!{\hat \Cbf}\Yit_{_{^{\!1}}}\!,
\Xit_{_{^{\!2}}}\!\times\! {\hat \Cbf}\Yit_{_{^{\!1}}}\big)= } }
\atop {=
\hbox{${\mhatH\!_{_{n}}}$}\!\!\big( \Yit_{_{^{\!1}}}\!,\Yit_{_{^{\!2}}}\big)
\oadd
\hbox{${\mhatH\!_{_{n}}}$}\!\!\big(\Yit_{_{^{\!1}}}\!,\Yit_{_{^{\!1}}}\big)
\oadd
\hbox{${\mhatH\!_{_{n}}}$}\!\!\big( \Xit_{_{^{\!1}}}\!,\Xit_{_{^{\!1}}}\big)
\oadd
\hbox{${\mhatH\!_{_{n}}}$}\!\!\big(\Xit_{_{^{\!1}}}\!,\Xit_{_{^{\!2}}}\big)
= } } \atop { {
=\hbox{${\mhatH\!_{_{n}}}$}\!\! \big({\hat \Cbf}
\Xit_{_{^{\!1}}}\!\times\! \Yit_{_{^{\!1}}}\!, {\hat \Cbf}
\Xit_{_{^{\!1}}}\!\times\!\Yit_{_{^{\!2}}}\big) \oadd
\hbox{${\mhatH\!_{_{n}}}$}\!\! \big({\hat \Cbf}
\Xit_{_{^{\!1}}}\!\!\times\!\Yit_{_{^{\!1}}}\!, {\hat \Cbf}
\Xit_{_{^{\!2}}}\!\times\! \Yit_{_{^{\!1}}}\big)
\oadd}
\atop {\oadd \hbox{${\mhatH\!_{_{n}}}$}\!\!
\big(\Xit_{_{^{\!1}}}\!\times\! {\hat \Cbf}\Yit_{_{^{\!1}}}\!,
\Xit_{_{^{\!1}}}\!\times\!{\hat \Cbf}\Yit_{_{^{\!2}}}\big) \oadd
\hbox{${\mhatH\!_{_{n}}}$}\!\!
\big(\Xit_{_{^{\!1}}}\!\!\times\!{\hat \Cbf}\Yit_{_{^{\!1}}}\!,
\Xit_{_{^{\!2}}}\!\times\! {\hat \Cbf}\Yit_{_{^{\!1}}}\big)} } }
\raise4pt\hbox{$\buildrel {\hbox{$\alpha$}_{\ast}\hbox{$
(\!=\!0\!\big)$}} \over {\hbox to 0.7cm{\rightarrowfill}}$}
}$}$

$\hskip3.3cm {\Bigg\downarrow}(k_{_3},-l_{_3})\hskip1.6cm
\#\hskip2.3cm \Bigg\downarrow(k_{_1},-l_{_1}) \hskip1.4cm
\#\hskip2.2cm \Bigg\downarrow(k_{_2},-l_{_2})\hskip2.3cm \#$

\noindent$\lower4.5pt\hbox{$^{ ...{\buildrel \hbox{$0$} \over {\hbox
to 0.5cm{\rightarrowfill}}}
\!\!\!\!\!\!
{ { { {
\hbox{${\rm H}\!\!\!\!\!_{_{_{n+1}}}$}\!\!\big(\hbox{${{
\Deight}\!\!^{^{_\wp}}\!}$}{\hbox{\sevenbf(}}{\hat
\Cbf}\Xit_{_{^{\!1}}}\!\!\mtopast\!
\Yit_{_{^{\!1}}}\!{\hbox{\sevenbf)}} {\hbox{\bf/}}
{\hbox{\bf[}}\hbox{${{\bf
\Deight\!\!^{^{_\wp}}}\!\!}$}{\hbox{\sevenbf((}} \!{\hat
\Cbf}\Xit_{_{^{\!1}}}\!\!\times \Yit_{_{^{\!2}}}\!{\hbox{\sevenbf)}}
\cup {\hbox{\sevenbf(}}\Xit_{_{^{\!1}}}\!\!\times {\hat \Cbf}
\Yit_{_{^{\!2}}}\!{\hbox{\sevenbf))}}
{\bf+}}
\atop {{\bf+} \lower1pt\hbox{${{\bf \Deight\!\!^{^{\wp}}}\!}$}
{\hbox{\sevenbf((}} \!{\hat \Cbf}\Xit_{_{^{\!2}}}\!\!\times
\Yit_{_{^{\!1}}}\!{\hbox{\sevenbf)}} \cup
{\hbox{\sevenbf(}}\Xit_{_{^{\!2}}}\!\!\times {\hat \Cbf}
\Yit_{_{^{\!1}}}\!{\hbox{\sevenbf))}}{\hbox{\bf]}}\big)\ (=)
} } \atop
{
(=) \hbox{${\mhatH\!\!\!\!\!_{_{_{n+1}}}}$}\!\! \big(
\Xit_{_{^{\!1}}}\!\mtopast\! \Yit_{_{^{\!1}}},
{\hbox{\sevenbf(}}\Xit_{_{^{\!1}}}\!\mtopast\!
\Yit_{_{^{\!2}}}\!{\hbox{\sevenbf)}} \cup
{\hbox{\sevenbf(}}\Xit_{_{^{\!2}}}\!\mtopast\! \Yit_{_{^{\!1}}}{\hbox{\sevenbf)}}\big)\ =
} } \atop { \bigg\{ { {
\hbox{${=\ \!{\rm H}\!\!\!\!\!_{_{_{n+1}}}}$}\!\!\big(\hbox{${{\bf
\Deight\!\!^{^{_\wp}}}\!}$}\! {\hbox{\sevenbf(}}
\Xit_{_{^{\!1}}}\!\!\mtopast\!
\Yit_{_{^{\!1}}}{\hbox{\sevenbf)}}{\hbox{\bf/}}{\hbox{\bf[}}\hbox{${{\bf
\Deight\!\!^{^{_\wp}}}\!}$} {\hbox{\sevenbf((}}{\hat
\Cbf}\Xit_{_{^{\!1}}}\times \Yit_{_{^{\!2}}}{\hbox{\sevenbf)}}\cup
{\hbox{\sevenbf(}}{\hat \Cbf} \Xit_{_{^{\!2}}}\times \Yit_{_{^{\!1}}}
{\hbox{\sevenbf))}}+}
\atop {+ \hbox{${{\bf \Deight\!\!^{^{_\wp}}}\!}$}
{\hbox{\sevenbf((}}\Xit_{_{^{\!1}}}\times{\hat \Cbf}
\Yit_{_{^{\!2}}}{\hbox{\sevenbf)}}\cup\ \!{\hbox{\sevenbf(}}
\Xit_{_{^{\!2}}}\times {\hat \Cbf}
\Yit_{_{^{\!1}}}{\hbox{\sevenbf))}}{\hbox{\bf]}} \big) } } } \bigg\}
}
\!\!\!\!\!\!
{\buildrel \hbox{$\delta$}_\ast \over {\hbox to
0.7cm{\rightarrowfill}}}\!\!\!\!\!\!
{ { { {
\hbox{${\rm H}\!_{_{n}}$}\!\!\big(\hbox{${{\bf
\Deight\!\!^{^{_\wp}}}\!}$} {\hbox{\sevenbf(}}{\hat
\Cbf}\Xit_{_{^{\!1}}}\!\!\times \Yit_{_{^{\!1}}}\!{\hbox{\sevenbf)}}
{\hbox{\bf/}} {\hbox{\bf[}}\hbox{${{\bf
\Deight\!\!^{^{_\wp}}}\!\!}$}{\hbox{\sevenbf((}} \!{\hat
\Cbf}\Xit_{_{^{\!1}}}\!\!\times \Yit_{_{^{\!2}}}\!{\hbox{\sevenbf)}}
\cap {\hbox{\sevenbf(}}\!\Xit_{_{^{\!1}}}\!\!\times {\hat \Cbf}
\Yit_{_{^{\!2}}}\! {\hbox{\sevenbf))}}
{\bf+}}
\atop {{\bf+}
\lower1pt\hbox{${{\bf \Deight\!\!^{^{\wp}}}\!}$} {\hbox{\sevenbf((}}
\!{\hat \Cbf} \Xit_{_{^{\!2}}}\!\!\times \Yit_{_{^{\!1}}}\!
{\hbox{\sevenbf)}}\cap {\hbox{\sevenbf(}}\Xit_{_{^{\!2}}}\!\!\times
{\hat \Cbf} \Yit_{_{^{\!1}}}\!{\hbox{\sevenbf))}}{\hbox{\bf]}}\big)
\ (=)
} } \atop
{
(=)\ \hbox{${\mhatH\!_{_{n}}}$}\!\! \big(
\Xit_{_{^{\!1}}}\!\!\times \Yit_{_{^{\!1}}},
{\hbox{\sevenbf(}}\Xit_{_{^{\!1}}}\!\!\times
\Yit_{_{^{\!2}}}\!{\hbox{\sevenbf)}} \cup
{\hbox{\sevenbf(}}\Xit_{_{^{\!2}}}\!\!\times\!
\Yit_{_{^{\!1}}}{\hbox{\sevenbf)}}\big)\ =
} } \atop { \bigg\{ { {
\hbox{${=\ \!\mhatH\!_{_{n}}}$}\!\!\big(
\Xit_{_{^{\!1}}}\!\!\times \Yit_{_{^{\!1}}},\
{\hbox{\sevenbf((}}{\hat \Cbf}\Xit_{_{^{\!1}}}\!\!\times
\Yit_{_{^{\!2}}}{\hbox{\sevenbf)}}\cup {\hbox{\sevenbf(}}
{\hat \Cbf}\Xit_{_{^{\!2}}}\!\!\times \Yit_{_{^{\!1}}}
{\hbox{\sevenbf))}}\cap}
\atop {\cap {\hbox{\sevenbf((}}\Xit_{_{^{\!1}}}\times{\hat \Cbf}
\Yit_{_{^{\!2}}}{\hbox{\sevenbf)}}\cup\ \!
{\hbox{\sevenbf(}}\Xit_{_{^{\!2}}}\times {\hat \Cbf}
\Yit_{_{^{\!1}}} {\hbox{\sevenbf))}}\big) } } } \bigg\} } \!\!\!
{\hbox to 0.7cm{\rightarrowfill}}
\!\!\!\!\!\! { { { {
\hbox{${\rm H}\!_{_{n}}$}\!\!\big(\hbox{${{\bf
\Deight\!\!^{^{_\wp}}}\!}{\hbox{\sevenbf(}}$} {\hat
\Cbf}\Xit_{_{^{\!1}}}\!\!\times \Yit_{_{^{\!1}}}\!{\hbox{\sevenbf)}}
{\hbox{\bf/}} {\hbox{\bf[}}\hbox{${{\bf
\Deight\!\!^{^{_\wp}}}\!\!}{\hbox{\sevenbf(}}$} \!{\hat
\Cbf}\Xit_{_{^{\!1}}}\!\!\times \Yit_{_{^{\!2}}}\!{\hbox{\sevenbf)}}
{\bf +} \hbox{${{\bf \Deight\!\!^{^{_\wp}}}\!\!}{\hbox{\sevenbf(}}$}
\!{\hat \Cbf}\Xit_{_{^{\!2}}}\!\!\times
\Yit_{_{^{\!1}}}\!{\hbox{\sevenbf)}}{\hbox{\bf]}}\big) \oadd
}
\atop {\oadd
\hbox{${\rm H}\!_{_{n}}$}\!\!\big(\hbox{${{\bf
\Deight\!\!^{^{_\wp}}}\!}{\hbox{\sevenbf(}}$}
\Xit_{_{^{\!1}}}\!\!\times {\hat \Cbf}
\Yit_{_{^{\!1}}}\!{\hbox{\sevenbf)}} {\hbox{\bf/}}
{\hbox{\bf[}}\hbox{${{\bf
\Deight\!\!^{^{_\wp}}}\!\!}{\hbox{\sevenbf(}}$}
\!\Xit_{_{^{\!1}}}\!\!\times {\hat \Cbf}
\Yit_{_{^{\!2}}}\!{\hbox{\sevenbf)}} {\bf +} \hbox{${{\bf
\Deight\!\!^{^{_\wp}}}\!\!}{\hbox{\sevenbf(}}$}
\!\Xit_{_{^{\!2}}}\!\!\times {\hat \Cbf}
\Yit_{_{^{\!1}}}\!{\hbox{\sevenbf)}}{\hbox{\bf]}}\big) \ %
(=)
}} \atop {(=)
\ \ \ \hbox{${\mhatH\!_{_{n}}}$}\!\!
\big(\Yit_{_{^{\!1}}}\!,\Yit_{_{^{\!1}}}\big) \oadd
\hbox{${\mhatH\!_{_{n}}}$}\!\!
\big(\Xit_{_{^{\!1}}}\!,\Xit_{_{^{\!1}}}\big) \ \
\hbox{\eightbf=0}
} } \atop { \bigg\{ { {
=\ \hbox{${\!\mhatH\!_{_{n}}}$}\!\!\big( {\hat \Cbf}
\Xit_{_{^{\!1}}}\!\!\times \Yit_{_{^{\!1}}},\
{\hbox{\sevenbf(}}{\hat \Cbf} \Xit_{_{^{\!1}}}\times
\Yit_{_{^{\!2}}}{\hbox{\sevenbf)}}\cup {\hbox{\sevenbf(}} {\hat
\Cbf}
{\Xit_{_{^{\!2}}}}\times \Yit_{_{^{\!1}}}{\hbox{\sevenbf)}}
\big) \oadd
}
\atop {\oadd \hbox{${\mhatH\!_{_{n}}}$}\!\!\big(
\Xit_{_{^{\!1}}}\!\!\times {\hat \Cbf}\Yit_{_{^{\!1}}},
{\hbox{\sevenbf(}}\Xit_{_{^{\!1}}}\times {\hat
\Cbf}\Yit_{_{^{\!2}}}{\hbox{\sevenbf)}}\cup
{\hbox{\sevenbf(}}{\Xit_{_{^{\!2}}}}\times {\hat \Cbf}
\Yit_{_{^{\!1}}}{\hbox{\sevenbf)}}
\big) } } \bigg\} } }
\hskip-0.3cm{\buildrel  \hbox{$0$} \over {\hbox to
0.5cm{\rightarrowfill}}}
}$}$

\medskip
$\hskip3.3cm {\downarrow}\delta_{3\ast}\hskip2.1cm (-\#) \hskip2.2cm
{\downarrow}\delta_{1\ast} \hskip2.2cm \# \hskip2.5cm
{\downarrow}\delta_{2\ast}\hskip2.7cm \# $
}

\normalsize
\bigskip
\noindent\hskip-1.5cm
The groups within curly braces in the last row doesn't play any part
in the deduction of the excisivity equivalence but fits the
underlying pattern.
}
\end{landscape}

\textwidth 12.5cm
\normalbaselines

\subsection{Simplicial calculus} \label{SecP34:Appendix}

The complex$_{o}\!\!$ of all subsets of a simplex$_{o}\!$
$\sigma$ is denoted ${\bar{\sigma}}$, while the {\it boundary of}
$\sigma$, ${\dot{\sigma}}$, is the set of all proper subsets.

So, ${ \dot{\sigma}}\!:=\{\tau\mid \tau\varsubsetneqq\sigma\}
={\bar{\sigma}}{\raise0.5pt\hbox{\eightmsbm \char"72}}
\{\sigma\}$,
 ${ \bar{\emptyset}}\!_{o}\!=\{\emptyset\!_{o}\!\}$ and
${ \dot{\emptyset}}\!_{o}\!=\emptyset$.

\smallskip
\noindent ``The {\it closed star} of $\sigma$ with respect to
${\Sigma}\hbox{''}=$ $\overline{\hbox{\rm {st}}}_{_{{\!\Sigma}}}\! \sigma\!
:= \{\tau\in \Sigma\mid\ \!\sigma\cup \tau\! \in \Sigma\}. $

\smallskip
\noindent ``The {\it open $($realized$)$ star} of $\sigma$ with respect to $\Sigma\hbox{''}\!\!=\!$ ${\hbox{\rm
{st}}}_{_{{\!\Sigma}}}\!\sigma\!:=\! \{\alpha\in|\Sigma|\mid [{\bf
v}\in\sigma] \Longrightarrow [\alpha({\bf v})\neq0]\}.$
So, ${\hbox{\rm {st}}}_{_{{\Sigma}}}\! (\sigma)
=\{\alpha\in|\Sigma|\mid \alpha \in
\mbfcupcap{\bigcap}{\!v\!\in\sigma}
{{\hbox{\rm {st}}}_{_{ \!{\Sigma}}}\!\!\{v\}} \}
$,
$\alpha_0\not\in{\hbox{\rm {st}}}_{_{{\!\Sigma}}}\!\sigma$ except for
${\hbox{\rm {st}}}_{_{{\Sigma}}}\! \emptyset\!_{o}\!=|\Sigma|$. $|\Sigma|$ and $\alpha_0$ are defined in p.~\pageref{DefP8:alpha0}.

\noindent
``The {\it closed geometrical simplex} $|\sigma|$ with respect to $\Sigma\hbox{''}\!\!=\!\!$
``The {\it $($realized$)$ closure} of $\sigma$ with respect to
$\Sigma\hbox{''}\!\!\!=\!$
$|\sigma|\!:=\!\{\alpha\in|\Sigma| \mid { [\alpha({\bf v})\!\neq0]
\Longrightarrow [{\bf v}\!\in\! \sigma] \}.}$
So, $|\emptyset\!_{o}\!|\!:=\!\{\alpha_0\}.$

\smallskip
\noindent ``The {\it interior}\label{SimplexInterior} of $\sigma$
with respect to ${\Sigma}\hbox{''}\!\!=\! {\mathcal h}\sigma{\mathcal i}
\!=\! \hbox{\rm Int} (\sigma) \!:=\! \{\alpha\in|\Sigma|\mid[{\bf
v}\in\nobreak\sigma ] \Longleftrightarrow [\alpha({\bf v})\neq0]
\}.$
So,
$\hbox{\rm Int}_{_{\!^o}}\!(\sigma\!_{_{\!^{\ \!\!}}}) $ is an open
subspace of $|\Sigma|$ \underbar{iff} $\sigma$ is a maximal simplex
in $\Sigma$ and
$\hbox{\rm Int}(\emptyset\!_{o}\!)\!:= \{\alpha_0\}$.

\smallskip
\noindent
The {\it barycenter} $\hat\sigma$ of\
$\sigma\neq\emptyset\!_{o}$ is the
$\alpha\!\in\!\hbox{Int}(\sigma)\!\in\!\vert\Sigma\vert$
fulfilling
${\bf v}\!\in\! \sigma
\!\Rightarrow\!
\alpha({\sevenrm v})\!=\!{1 \over{\#\sigma}}\!$ while
$\hat{\emptyset}\!_{o}\!\!:= _{\!}\alpha_0.$

\smallskip\noindent
$\!$``The {\it link} of $\sigma$ with respect to $\Sigma\hbox{''}=$
$\hbox{\rm {Lk}}_{_{\!_{\Sigma}}}\!\!\!\sigma\!:= \{\tau\!\in
\Sigma| [\sigma\cap \tau =\emptyset]\land [\sigma\cup \tau \in
\Sigma]\}.$
So, ${\rm Lk}\!\!_{_{_{\Sigma}}}\!\emptyset\!_{o}\!=\!\Sigma,$
$\sigma\!\in\!{\rm
Lk}\!\!_{_{_{\Sigma}}}\!\tau\Longleftrightarrow\tau\!\in\!{\rm
Lk}\!\!_{_{_{\Sigma}}}\!\sigma$
and
${\rm Lk}\!\!_{_{_{\Sigma}}}\!\sigma\!=\!\emptyset\!_{_{}}$
\underbar{iff} $\sigma\!\not\in\!\Sigma,$
while
${\rm Lk}\!\!_{_{_{\Sigma}}}\!\tau\!=\{\emptyset\!_{o}\!\} $
\underbar{iff} $\tau\!\in\!{\Sigma}$ is maximal.

\smallskip
\noindent The {\it join}
$ \Sigma_{1}\!\ast\! \Sigma_{2}:= \{\sigma\!_{1}\cup
\sigma\!_{2}\vert \sigma\!_{i}\in \Sigma_{i}\ (i\!=\!1,2)\}.$
In\ particular $\Sigma\ast\! \{\emptyset_o\!\}=
\{\emptyset_o\}\ast\! \Sigma = \Sigma.$

\begin{proposition} \label{proposition:PropP34:1} \label{PropP34:1}
If $\Sigma_{1}$ and $\Sigma_{2}$  has mutually exclusive vertex
sets, $\ V_{\Sigma_{1}}\cap V_{\Sigma_{2}}=\emptyset \ $ {then},
$$
[\tau\in \Sigma_{1}\ast_o \Sigma_{2}]
\Longleftrightarrow
[\exists !\ \sigma_{i}\in \Sigma_{i}\ \textrm{so\ that}\ \tau
=\sigma_{1}\cup \sigma_{2}],\leqno{{\bf a})}$$
and:
$$ {\rm Lk}\!\!\!\!\!\!
 \lower1.1pt\hbox{$_{_{\Sigma_{1}\!\ast\Sigma_{2}}}$}\!\!(\sigma_{1}\!
\cup\sigma_{2})
=
({\rm Lk}\!\!
\lower1.1pt\hbox{${_{_{\Sigma_{1}}}}$}\!\!\!\sigma_{1})\ast
({\rm Lk}\!\!
\lower1.1pt\hbox{${_{_{\Sigma_{2}}}}$}\!\!\!\sigma_{2}),\leqno{{\bf
b})}.
$$
\quad
\end{proposition}

\begin{proof}
a) follows immediately from the definition.

\medskip\noindent
b)
{\baselineskip12pt
Suppose that $V_{\Sigma_{1}}\cap V_{\Sigma_{2}}=\emptyset$ and
 that $\sigma_{1}\cup \sigma_{2}\in \Sigma_{1}\ast \Sigma_{2}$ where
 $\sigma_{i}\in \Sigma_{i}$  (i=1,2).
$${\rm Lk}_{\Sigma_{1}\ast \Sigma_{2}}(\sigma_{1}\cup
\sigma_{2})\buildrel Def. \over =\{\mu_o\in \Sigma_{1}\ast
\Sigma_{2}|\mu_o\cap(\sigma_{1}\cup \sigma_{2})=\emptyset\ and\
\mu_o\cup(\sigma_{1}\cup \sigma_{2})\in\Sigma_{1}\ast
\Sigma_{2}\}=$$
=[This equality follows from the first part i.e. from {\bf a}]=
$$\hskip-0.5cm=\{\mu_o=\mu_{1}\cup\mu_{2}|[\mu_{i}\in \Sigma_{i}\ (i=1,2)]\ and\
[\mu_{1}\cup\mu_{2}\cap(\sigma_{1}\cup \sigma_{2}=\emptyset)]\ and\ $$
\noindent$[(\mu_{1}\cup\mu_{2})\cup(\sigma_{1}\cup \sigma_{2})\in
\Sigma_{1}\ast \Sigma_{2}]\}$
= [Here we use that $V_{\Sigma_{1}}\cap V_{\Sigma_{2}}=\emptyset]\ =$
$$\hskip-0.5cm\{\mu_{1}\cup\mu_{2}|[\mu_{i}\in \Sigma_{i}\ (i=1,2)]\ and\
[[\mu_{1}\cap\sigma_{1}=\emptyset)]\land[[\mu_{2}\cap\sigma_{2}=\emptyset)]]\
and\ $$
\noindent$[(\mu_{1}\cup\sigma_{1})\cup(\mu_{2}\cup\sigma_{2})\in
\Sigma_{1}\ast \Sigma_{2}]\}$
= [Uniqueness from {\bf a}] =$\{\mu_{1}\cup\mu_{2}|[\mu_{i}\in \Sigma_{i}\ (i=$
$$
1,2)]\ and\
[[\mu_{1}\cap\sigma_{1}=\emptyset)]\land[[\mu_{2}\cap\sigma_{2}=\emptyset)]]\
and\ [[\mu_{1}\cup\sigma_{1})\in
\Sigma_{1}]\land [\mu_{1}\cup\sigma_{2})\in \Sigma_{2}]]\}$$
$$=\ \hbox{[By the definition of ``${\rm Lk}$'']}\ =
\{\mu_{1}\cup\mu_{2}\mid\mu_{i}\in {\rm Lk}_{\Sigma_{i}}\sigma_{i}\
(i=1,2) \}=$$
\noindent\hskip4cm$=({\rm Lk}_{\Sigma_{1}}\sigma_{1})\ast ({\rm
Lk}_{\Sigma_{2}}\sigma_{2})$}
\end{proof}
\normalbaselines

\noindent
With the convention
${\rm Lk}\!\!_{_{_{\Sigma}}} \!\hbox{v} := {\rm
Lk}\!\!_{_{_{\Sigma}}} \!\{\hbox{v}\} $,
any link is an iterated link of vertices and

${\rm Lk}
\!\!\!\!_{\lower2.4pt\hbox {$_{{{{\hbox{\sevenrm
Lk}}}}\!\!_{_{_{\Sigma}}}\!\!\sigma}$}}
\!\!\!\!\!\tau
\!=\!
\{\emptyset\!_{o}\!\}
\ast {\rm Lk}
\!\!\!\!_{\lower2.4pt\hbox {$_{{{{\hbox{\sevenrm
Lk}}}}\!\!_{_{_{\Sigma}}}\!\!\sigma}$}}
\!\!\!\!\tau
\!=\! {\rm Lk}
\!\!_{\lower2.4pt\hbox {$_{{{\bar{\hbox{\seveni {\char"1B}}}}}}$}}
\!\sigma\ast
{\rm Lk}
\!\!\!\!_{\lower2.4pt\hbox {$_{{{{\hbox{\sevenrm
Lk}}}}\!\!_{_{_{\Sigma}}}\!\!\sigma}$}}
\!\!\!\!\tau
\!=\!
\big[{^{\hbox{\eightrm Prop.~\ref{PropP34:1}}}
_{\hbox{\eightrm\hskip0.3cm above}}}\big]
\!=\!
{\rm Lk}\!\!
\!\!\!\!_{\lower2.4pt\hbox {$_{{{{\bar{\hbox{\seveni
{\char"1B}}}\ast\hbox{\sevenrm
Lk}}}}\!\!_{_{_{\Sigma}}}\!\!\sigma}$}}
\!\!\!\!\!\!(\sigma\cup\tau)
\!=\!
{\rm Lk}\!\!\!\!_{\hbox{$_{_{{\overline{\hbox{\sevenrm
st}}}\!\!_{_{_{\Sigma}}}\!\!\hbox{\seveni {\char"1B}}}}$}}
\!\!\!(\sigma\cup\tau)
\!=\!
{\rm Lk}\!\!\!\!\!\!_{\hbox{$_{_{{ \overline{\hbox{\sevenrm
st}}}\!\!_{_{_{\Sigma}}}\!\!
(\sigma\cup\tau)
}}$}} \!\!\!\!\!\!\!\!\!\!(\sigma\cup\tau)
\!=\!
{\rm Lk}\!_{_{\Sigma}} \!(\sigma\cup\tau).$

\noindent
So, $\tau\!\notin\!{\rm Lk}\!_{_{\Sigma}}\!\sigma\!\Rightarrow\!$
${\rm Lk} \!\!\!_{\lower2.3pt\hbox {$_{{{{\hbox{\sevenrm Lk}}}}\!\!_{_{_{\Sigma}}}\!\!\sigma}$}} \tau
\!=\!
\emptyset$\ while;

$\tau\!\in\!{\rm Lk}\!\!_{_{_{\Sigma}}}\!\sigma\Rightarrow $
${\rm Lk}
\!\!\!\!_{\lower2.4pt\hbox {$_{{{{\hbox{\sevenrm
Lk}}}}\!\!_{_{_{\Sigma}}}\!\!\sigma}$}}
\!\!\!\tau
\!=\!
{\rm Lk}\!\!_{_{_{\Sigma}}} \!(\sigma\cup\!\tau)
\!=\!
{\rm Lk}
\!\!\!\!_{\lower2.4pt\hbox {$_{{{{\hbox{\sevenrm
Lk}}}}\!\!_{_{_{\Sigma}}}\!\!\tau}$}}
\!\!\!\sigma
\ \big(\!\!\subset\! {\rm Lk}\!_{_{\Sigma}}\!\sigma \!\cap\! {\rm
Lk}\!_{_{\Sigma}}\!\tau\big).$ \hfill({\bf I}\label{EqP34:I})

\goodbreak
\noindent
With $\dim\Sigma=2$, the
\htmladdnormallink{{\it edge contraction}
}{http://en.wikipedia.org/wiki/Edge_contraction}
of $\{v,w\}\in\Sigma$ in $\vert\Sigma\vert$ is topology preserving if and only if
$
{\rm Lk}\!\!_{_{_{\Sigma}}} \!\{v,w\}=
{\rm Lk}\!_{_{\Sigma}}\!v\!\cap\! {\rm
Lk}\!_{_{\Sigma}}\!w.
$

\smallskip\noindent
$\tau\!\in\!{\rm Lk}\!\!_{_{_{\Sigma}}}\!\sigma \Longrightarrow
\overline{\hbox{\rm st}}
\!\!\!\!_{\lower2.4pt\hbox {$_{{{{\hbox{\sevenrm
Lk}}}}\!\!_{_{_{\Sigma}}}\!\!\tau}$}}
\!\!\!\!\sigma
=
\bar\sigma\ast
{\rm Lk}
\hskip-0.2cm_{\lower2.4pt\hbox {$_{{{{\hbox{\sevenrm
Lk}}}}\!\!_{_{_{\Sigma}}}\!\!\tau}$}}
\!\!\!\!\sigma
=
[\hbox{{\rm Eq.$\ \!${I}}}]
=
\bar\sigma
\ast {\rm Lk}
\!\!\!\!_{\lower2.4pt\hbox {$_{{{{\hbox{\sevenrm
Lk}}}}\!\!_{_{_{\Sigma}}}\!\!\sigma}$}}
\!\!\!\!\tau
= {\rm Lk}
\hskip-0.2cm_{\lower2.4pt\hbox {$_{{{\bar{\hbox{\seveni {\char"1B}}}}}}$}}
\emptyset\!_{o} \!\ast {\rm Lk}
\hskip-0.2cm_{\lower2.4pt\hbox {$_{{{{\hbox{\sevenrm
Lk}}}}\!\!_{_{_{\Sigma}}}\!\!\sigma}$}}
\!\!\!\!\tau
=\!\!
\big[{^{\hbox{\eightrm Prop.~\ref{PropP34:1}}}
_{\hbox{\eightrm\hskip0.3cm above}}}\big]
\!\!=
{\rm Lk}\!\!
\!\!\!\!_{\lower2.4pt\hbox {$_{{{{\bar{\hbox{\seveni
{\char"1B}}}\ast\hbox{\sevenrm
Lk}}}}\!\!_{_{_{\Sigma}}}\!\!\sigma}$}}
\!\!\!\!\!\!
\tau
=
{\rm Lk}
\!\!\!_{\lower2.4pt\hbox {$_{{{\overline{\hbox{\sevenrm
st}}}}\!\!_{_{_{\Sigma}}}\!\!\sigma}$}}\!\!\!
\tau.$

\smallskip
\noindent
The $p$-{\it skeleton} is defined by,
$
\Delta\!^{^{_{(_{\!}p_{\!})}}}\!\!\!:= \{\delta\!\in\!\Delta\
\!\vert\ \!\#\delta\le p+1\}$. So,
$\Delta\!^{^{_{(\!n)}}}\!=\Delta$, if
$n:=\dim\Delta$.

The following definitions simplifies notations:
$\Delta\!^{^{_{\prime\!}}}\!:=
\Delta\!^{^{_{(\!n{_{_{^{_{\!}}}}}\!-\!1_{\!})}}}\!
$.
$\Delta\!^{^{_{p\!}}}\!\!:= \Delta\!\!^{^{_{(_{\!}p_{\!})}}}\!
\smallsetminus
\Delta\!\!^{^{_{(_{\!}p\!-\!1_{\!})}}}\!\!.$

\medskip
\noindent
$\Gamma\subset\Sigma$ is said so be {\it full} in $\Sigma$ if for
all $\sigma\!\in\!\Sigma$; $\sigma\subset V_{_{^\Gamma}}
\Longrightarrow
\sigma\in \Gamma.$

\begin{proposition} \label{proposition:PropP34:2} \label{PropP34:2}
$a.$
\begin{enumerate}
\label{PropP34:2.a}
\item[i.]
 $\!{{\rm Lk}\!\!\!\!\!\!\!
\lower1.5pt\hbox{$_{_{{ \Delta\!_{_{{1\!\!}}}\cup
{\Delta_{_{{\!{{2}}}}}}}}}$}\!\!\!\delta} =\ \!{{\rm Lk}\!
\lower0.5pt\hbox{$ _{_{{ {\Delta_{_{{\!{{1}}}}}}}}} $}
\!\!\!\!\delta} \cup {{\rm Lk}\! \lower0.5pt\hbox{$ _{_{{
{\Delta_{_{{\!{{2}}}}}}}}} $} \!\!\!\!\delta},$
\item[ii.]
${{\rm Lk}\!\!\!\!\!\!\! \lower1.5pt\hbox{$_{_{{
\Delta\!_{_{{1\!\!}}}\cap {\Delta_{_{{\!{{2}}}}}}}}}$}\!\!\!\delta}
=\  _{\!} \!{{\rm Lk}\! \lower0.5pt\hbox{$ _{_{{
{\Delta_{_{{\!{{1}}}}}}}}} $} \!\!\!\!\delta} \cap {{\rm Lk}\!
\lower0.5pt\hbox{$ _{_{{ {\Delta_{_{{\!{{2}}}}}}}}} $}
\!\!\!\!\delta},$ \ \

\item[iii]
$\big({{\Delta\!^{^{\!_{\ \!}}}}\!_{_{{1}}}\!{{\ast}}}
{{\Delta\!^{^{\!_{\ \!}}}\!_{_{{2}}}}} \big)
\!^{^{{\hbox{$_{\prime}$}\!}}}\!\!=_{\!}
({{\Delta\!^{^{_{{\hbox{$_{\prime}$}\!}}}}_{1}}}\!{{\ast}}
\Delta\!^{^{\!_{\ \!}}}\!_{_{{2}}}) \cup (\Delta\!^{^{\!_{\
\!}}}\!_{_{{1}}}\!{{\ast}}
{{\Delta\!^{^{_{{\hbox{$_{\prime}$}\!}}}}_{2}}}), $
\item[iv]
$\Delta\ast(\ \! {\cap\!\!\!\!\!_{_{_{{{i\in\hbox{\sevenbf
I}}}}}}}\! {\Delta\!^{^{\!_{\ \!}}}}\!_{_{{i}}}\!)\!=
{\cap\!\!\!\!\!_{_{_{{{i\in\hbox{\sevenbf I}}}}}}}
(\Delta\ast{\Delta\!^{^{\!_{\ \!}}}}\!_{_{{i}}}) $,\ \ \ and\ \ \
$ \Delta\ast(\ \! {\cup\!\!\!\!\!_{_{_{{{i\in\hbox{\sevenbf
I}}}}}}}\! {\Delta\!^{^{\!_{\ \!}}}}\!_{_{{i}}})\! =
{\cup\!\!\!\!\!_{_{_{{{i\in\hbox{\sevenbf I}}}}}}}
(\Delta\ast{\Delta\!^{^{\!_{\ \!}}}}\!_{_{{i}}}).$

\hskip-0.5cm{\rm ({iv} also holds for arbitrary topological spaces
under the $\mhatast$-join.)}
\end{enumerate}

\noindent
$b.$
\label{PropP34:2.b}\ $\Delta\ {pure}\ \!\Longleftrightarrow\!
\big({{\rm Lk}\!_{_{\Delta}}\!\!{{\delta}}}\big)
\!^{^{{\hbox{$_{\prime}$}\!}}} \!= {{\rm
Lk}\!\!_{_{\Delta\!^{^{_{\prime\!\!}}}}}{{\delta}}}\ \ \forall\
\emptyset\!_{_{^{o}}}\!\!\ne\! \delta\!_{_{}}\in\!\Delta.$

\noindent
$c.$
\label{PropP34:2.c} $[{\Gamma} \cap {\rm Lk}\!_{_{\Delta}}\!\delta
= {\rm Lk}_{_{\!_{\Gamma}}}\!\!\delta $ $\ \! \forall\ \!
\delta\in{\Gamma}] \!\Longleftrightarrow\! [{\Gamma} $ is full in
$\Delta]$ $ \!\Longleftrightarrow\! [{\Gamma} \cap
\overline{\hbox{\rm st}}\!_{_{\Delta}}\!\delta = \overline{\hbox{\rm
st}}_{_{\!_{\Gamma}}}\!\!\delta $ \hbox{$\ \! \forall
\delta\in{\Gamma}].\ \ $
{\rm({\bf II})}} \label{EqP34:II}
\end{proposition}
\begin{proof}
{\it a.} Associativity and distributivity of the logical
connectives.\hfill{$\triangleright$}

\medskip
\noindent{\it b.}\ $(\Longrightarrow)$ $ \Delta\ \hbox{\rm pure}\
\Longrightarrow {\Delta\!^{^{\!_{\prime}}}\!}\ \hbox{\rm pure},\
\hbox{\rm so;} $ \noindent ${{\rm
Lk}\!_{_{\Delta\!^{^{\!_{\prime}}}\!}}\!{{\delta}}} \!=\! \{
\tau\!\in\!\Delta \mid
\tau\cap\delta= \!\emptyset \land\
\!\!\tau\cup\delta\!\in\!\Delta\!^{^{\!_{\prime}}}\! \} \!=$
\vskip0.1cm\noindent
$=\{ \tau\in\Delta \mid
\tau\cap\delta= \!\emptyset \land\
\tau\cup\delta\in\!\Delta\ \land\
\!\#(\tau\cup\delta)\le n \} = $ $ \{
\tau\!\in\!{{\rm Lk}\!_{_{\Delta}}\!\!{{\delta}}}\ \! {|
|}\#(\!\tau\cup\delta)\!\le\! n\}
\!=$
\vskip0.1cm\noindent
$
=\{ \tau\!\in\!{{\rm Lk}\!_{_{\Delta}}\!\!{{\delta}}} \mid
\#\tau\!\le\! n\!_{_{^{\ }}}\!\!\!-\!\#{\delta} \} \!=\! \{
\tau\!\in\!{{\rm Lk}\!_{_{\Delta}}\!\!{{\delta}}}\mid
\#\tau\!\le\! \dim{{{{{\rm Lk}\!_{_{\Delta}}\!\!{{\delta}}}}}}
\} \!=\! \big({{\rm Lk}\!_{_{\Delta}}{{\delta}}}\big)
^\prime.$

\medskip
\noindent $(\Longleftarrow)$ If $\Delta$ non-pure, then $\exists\
\!\delta_{\rm m}
\in\Delta^\prime$ maximal in both
${\Delta}^\prime$ and ${\Delta}$ i.e.,

\hfill
$$_{\!}\big({{\rm Lk}\!_{_{\Delta}}\!\!{{\delta}_{\rm m}\!}}\big)^\prime =
\big(\{\emptyset\!_{_{o}}\}\big)^\prime =
\emptyset \not= \{\emptyset\!_{_{o}}\} = {{\rm
Lk}\!_{_{\Delta\!^{^{\!_{\prime}}}\!}}\delta_{\rm m}}. \ \eqno{\triangleright}$$

\medskip
\noindent{\it c.} $ {\rm Lk}\!_{_{_{\Gamma}\!\!}}\sigma \!=
\{\tau\vert\ \tau\cap\sigma=\emptyset\ \&\ \tau\cup\sigma\in\Gamma\}
= \!\big[ {{\hbox{\eightrm true}\ \forall\ \!\!\sigma\in\Gamma}
\atop {\underline{\hbox{\eightrm iff}}\ \Gamma\ \hbox{\eightrm
full}}} \big]\! = \{\tau\in\Gamma\vert\ \tau\cap\sigma =\emptyset\
\&\ \tau\cup\sigma\in\Delta\}
= \Gamma\cap{\rm Lk}\!_{_{\Delta\!}}\sigma. $\\
\indent
$ {\Gamma} \cap \overline{\hbox{\rm st}}\!_{_{\Delta}}\!\delta =
{\Gamma} \cap \{\tau\in\Gamma\vert\ \tau\cup\sigma\in\Delta\} \cup
\{\tau\in\Delta \mid \tau\notin\Gamma\ \&\ \tau\cup\sigma\in\Delta\}
=$
\vskip0.2cm\noindent
$
=\{\tau\in\Gamma\vert\ \tau\cup\sigma\in\Delta\}$
$ = \big[{{\hbox{\eightrm True\ for\ all}\ \sigma\in\ \!\Gamma\
\underline{\hbox{\eightrm iff}}\ {\Gamma}}\atop{\hbox{\eightrm is\
a\ full\ subcomplex.}}}\big] = \{\tau\in\Gamma\vert\
\tau\cup\sigma\in\Gamma\} = \overline{\hbox{\rm
st}}_{_{\!_{\Gamma}}}\!\!\delta.
$
\end{proof}

\vfill\break

\noindent
The {\it contrastar} of\ $\sigma$ with respect to $\Sigma\!={\rm
cost}_{_{\!{\Sigma}}}\!\sigma\!:=\! \{\tau\!\in\!\Sigma\mid\
\tau\!\not\supseteq\sigma\}.$

So,
$\hbox{\rm
cost}_{_{{\!\Sigma}}}\!\emptyset\!_{o}\!\!=\!\emptyset\!$
\ and\
${\hbox{\rm cost}_{_{{\!\Sigma}}}\!\sigma\!=\!\Sigma\
\underline{\hbox{\rm iff}}\ \sigma\!\not\in\!\Sigma}.$

\goodbreak

\begin{proposition} \label{proposition:PropP35:3} \label{PropP35:3}
Whether
${{\delta_{1}}}\!,{{\delta_{2}}} \!\!\in\!\Delta$ or not,
the following holds;

\medskip
\noindent {\it a.}
\label{PropP35:3a}
\medskip
$ \hbox{\rm
cost}\!_{_{\Delta}}\!({{\delta_{1}}}\cup{{\delta_{2}}})
= \hbox{\rm cost}\!_{_{\Delta}}\!{{\delta_{1}}} \cup
\hbox{\rm cost}\!_{_{\Delta}}\!{{\delta_{2}}}$\
\  and  \
$\delta=\{{\bf v}\!_{_{^{1}}},...,{\bf v}\!_{_{^{p}}}\}
\Rightarrow \hbox{\rm cost}\!_{_{\Delta}}\!{{\delta\!_{_{ }}}}
=\nolinebreak
\mbfcupcap{\bigcup}{{\!\!i=1,p}}
\hbox{\rm cost}\!_{_{\Delta}}\!\!{{\bf v}\!_{i}}. $

\medskip
\noindent {\it b.}
\label{PropP35:3b}
\ $ \hbox{\rm cost}\!\!\!\!\!\!\! _{\lower0.8pt\hbox{$
{_{\hbox{\eightrm cost}\!\!_{_{\Delta}}\!\!{{\delta_{1} }}  }
}$}}\! \!{{\delta_{2}}}\! = \hbox{\rm
cost}\!_{_{\Delta}}\!{{\delta_{1}}} \!\cap \hbox{\rm
cost}\!_{_{\Delta}}\!{{\delta_{2}}}\! = \hbox{\rm
cost}\!\!\!\!\!\!\! _{\lower0.8pt\hbox{$ {_{\hbox{\eightrm
cost}\!\!_{_{\Delta}}\!\!{{\delta_{2} }} } }$}}\!
\!{{\delta_{1}}} \indent \ \ \hbox{\rm and} \indent  \ \
\mbfcupcap{\bigcap}{\!\!\!i=1,q}
\hbox{\rm cost}\!_{_{\Delta}}\!\!{{\delta}\!_{i}} =
\hbox{\rm cost} \hskip-0.3cm_{_{{{\rm cost}\hskip0.05cm
\delta\!_{_{^{2}}}}\hskip-0.3cm_{_{ {^{{{\bf .}}}\!_{\!}{{\bf
.}}_{_{\!{\bf .}}}} }}\!}}\hskip0.0cm{{{\delta}}\!_{_{^{1}}}}
\hskip-0.6cm{\lower12pt\hbox{\sevenrm cost$_{\!_{^\Delta}\!\!^{\delta_q}}$}}
$

\medskip\noindent
\label{PropP35:3c}
$c. \
\left\{\begin{array}{ll}
{\bf i.}\ \ {\delta}\!\notin\!{\Delta} \Leftrightarrow
\hbox{\tenbf[} {\rm Lk}\!\!\!\!\!\! {\lower0.4pt\hbox{${
_{_{\hbox{\eightrm cost}\!\!_{_{\Delta}}\!\!\!{\hbox{\eightbf v}}
}}\!\!\!\!{{\delta\!_{_{ }}}} }$}}
=
\hbox{\rm cost}\!\!\!\!\! \lower1pt\hbox{$_{_{{\!\!\hbox{\eightrm
Lk}\!\!_{_{\Delta}} }}}\!\!\lower3pt\hbox{$_{\delta\!}$} $}\! {{\bf
v}}
=
\emptyset \hbox{\tenbf]}. \\
{\bf ii.}\ {\delta}\!\in\!{\Delta} \Leftrightarrow \hbox{\tenbf[}
{\bf v}\notin{\delta} \Leftrightarrow {\delta}\in\hbox{\rm
cost}\!_{_{\Delta}}\!\!{{\bf v}} \Leftrightarrow {\rm
Lk}\!\!\!\!\!\! {\lower0.4pt\hbox{${ _{_{\hbox{\eightrm
cost}\!\!_{_{\Delta}}\!\!\!{\hbox{\eightbf v}}
}}\!\!\!\!{{\delta\!_{_{ }}}} }$}}
=
\hbox{\rm cost}\!\!\!\!\! \lower1pt\hbox{$_{_{{\!\!\hbox{\eightrm
Lk}\!\!_{_{\Delta}} }}}\!\!\lower3pt\hbox{$_{\delta\!}$} $}\! {{\bf
v}}
\!\supset\!
\{\emptyset_o\}\ (\neq\emptyset) \hbox{\tenbf]}.
\end{array}\right.
$

\hskip0.5cm
$
\big({\delta}\!\in\!{\Delta} \Leftrightarrow \hbox{\tenbf[} {\bf
v}\in{\delta} \Leftrightarrow {\delta}\notin\hbox{\rm
cost}\!_{_{\Delta}}\!\!{{\bf v}} \Leftrightarrow {\rm
Lk}\!\!\!\!\!\! {\lower0.4pt\hbox{${ _{_{\hbox{\eightrm
cost}\!\!_{_{\Delta}}\!\!\!{\hbox{\eightbf v}}
}}\!\!\!\!{{\delta\!_{_{ }}}} }$}} = \emptyset \neq \{\emptyset_o\}
\subset {\rm Lk}\!_{_{\Delta}}\!\!{{\delta\!_{_{ }} }} = \hbox{\rm
cost}\!\!\!\! \lower1pt\hbox{$_{_{{\!\!\hbox{\eightrm
Lk}\!\!_{_{\Delta}}\!\!{{{\hbox{\eighti {\char"0E}}} }} } }}$}\!
{{\bf v}} \hbox{\tenbf]}. \big)
$

\medskip
\noindent
\label{PropP35:3d}
$d.$
\ $ {\delta},\tau\!\in\!{\Delta} \!\Longrightarrow\! \hbox{\tenbf[}
{\delta}\in{\hbox{\rm cost}\!_{_{\Delta}}\!\!{{\tau}}}
\!\Longleftrightarrow\!
\emptyset\not={\rm Lk}\!\!\!\!\! {\lower0.6pt\hbox{${ _{_{\!{\rm
cost}\!\!_{_{\Delta}}\!\!\!\tau }} }$}} \!{{\delta\!_{_{ }}}} = {\rm
Lk}\!\!\!\!\!\!\!\!\!\! {\lower1.4pt\hbox{$_{_{\!\hbox{\eightrm
cost}\!\!_{_{\Delta}}\!\!\!\!\! {\lower2.4pt\hbox{
${^{_{{{\hbox{{\sixrm (}{\eighti {\char"1C}}}}\
\!\!{{{{\setminus}}}}\ {\!}\!{\hbox{{\eighti {\char"0E}}{\sixrm
)}}}}}}}$} } }}$}}\! \!\!\!\!\!{{\delta\!_{_{ }}}}
= \hbox{\rm cost}\!\!\!\!
{\lower1.0pt\hbox{${ _{_{{\!\!\hbox{\eightrm
Lk}\!\!_{_{\Delta}}\!\!\!{{{{\hbox{\eighti {\char"0E}}}}} }}}}\!
}$}}
(\tau\lower1pt\hbox{$\ _{\!}\!^{_{_{\setminus}}}\
_{\!}\!\delta$})\hbox{\tenbf]}. $

\ \ \ {\tenbf(}
{\it d}
$
\Rightarrow {\rm Lk}\!\!\!\!\!\!\! {\lower0.5pt\hbox{${
_{_{\!\hbox{\eightrm
cost}\!\!_{_{\Delta}}\!\!\!{\lower2.4pt\hbox{${^{_{{{\hbox{\eighti
{\char"1C}}}}}}}$}} }}\! }$}} \!{{\delta\!_{_{ }}}} =
{\rm Lk}\!_{_{\Delta}}\!\!{{{\delta}} }$
\ \underbar{iff} \
$\delta\cup\tau \notin \Delta $
\ \ or equivalently, \underbar{iff} \
$\tau\lower1pt\hbox{$\ _{\!}\!^{_{_{\setminus}}}\
_{\!}\!\delta$}\notin\ \!\overline{\hbox{\rm
{st}}}\!_{_{\Delta}}\!\delta.{\bf)}$

\medskip\noindent
\label{PropP35:3e}
$e.$
\ {If} $\delta\ne\emptyset$ {then}; \
$ \left\{\begin{array}{ll}
{\bf 1.}\
[\big({\hbox{\rm cost}\!_{_{\Delta\!\!}}{{\delta}}}\big)
\!^{^{{\hbox{$_{\prime}$}\!}}}
=
{\hbox{\rm cost}\!\!_{_{\Delta\!^{^{_{\prime\!\!}}}}}{{\delta}}}]
\Longleftrightarrow [n\!_{_{{\delta}}}\!\!
=
\!n\!_{_{\Delta}}\ \!\!]
\\
{\bf 2.}\
[{\hbox{\rm cost}\!_{_{\Delta\!\!}}{{\delta}}} \ =\ {\hbox{\rm
cost}\!\!_{_{\Delta\!^{^{_{\prime\!\!}}}}}{{\delta}}}]
\Longleftrightarrow [n\!_{_{{\delta}}}\!\! =\!n\!_{_{\Delta}}\
\!\!\!\!-1]
\\
\end{array}\right.
$
,\vskip0.1cm \text{where}
$ n\!_{_{\lower1pt\hbox{${^{_{{\varphi}}}}$}} }
\!\!:=\!\dim\hbox{\rm
cost}_{_{\!\Delta}}\!{\raise1pt\hbox{$\varphi$}} $
 \quad \text{and} \quad
$\Delta^\prime:=\Delta^{(n_\Delta-1)}$,
\quad \text{with} \quad
$n\!_{_{\Delta}}\!\!:=\!\dim\!\Delta<\infty.$\\
\medskip
{\rm(It is \underbar{always} true that;
$\tau\!\subset\!\delta\Rightarrow n\!_{_{\Delta}}\!\!-\!1\ \!\le
n_{_{\!\raise1pt\hbox{${_{{{\tau}}}}$}} } \!\le
n\!_{_{\raise1pt\hbox{${_{{{\delta}}}}$}} } \!\le n\!_{_{\Delta}}$,\
\text{if}\ $\emptyset\!_{o}\ne\tau,\delta$.)}

\noindent
\label{PropP35:3f}
$f.$
\ $\ \!\overline{\hbox{\rm {st}}}\!_{_{\Delta}}\!\delta \cap
{\hbox{\rm cost}\!_{_{\Delta\!\!}}{{\delta}}}
=
\hbox{\rm cost}\!\!\!\!\lower3.5pt\hbox{${\phantom{\ }
\over{^{\hbox{\eightrm st}}}}$} {\lower2.4pt\hbox{${
_{_{^{{\!\!\!_{_{\Delta}}\!\!\delta }}}}\! }$}} \!\!\delta = {\dot{\delta}}\ast {{\rm Lk}\!_{_{\Delta\!\!}}{{\delta}}}. $
\indent\hskip0cm{\tenbf(}{\rm So, $\hbox{\rm
cost}\!\!\!\!\lower3.5pt\hbox{${\phantom{\ } \over{^{\hbox{\eightrm
st}}}}$} {\lower2.4pt\hbox{${
_{_{^{{\!\!\!_{_{\Delta}}\!\!\delta }}}}\! }$}} \!\!\delta $} {\rm is a quasi-/homology manifold if $\Delta$ is, by Th.\ref{LemmaP28:2.i} p.~\pageref{LemmaP28:2.i}, Th.\ref{TheoremP31:12} p.~\pageref{TheoremP31:12} and \ref{PropP22:1} p.~\pageref{PropP22:1}.$)$}
\end{proposition}

\begin{proof}
If $ \Gamma\!\subset\Delta $ then $ \hbox{\rm
cost}_{_{\Gamma}}\!{{\gamma}} \!=\! \hbox{\rm
cost}\!_{_{\Delta}}\!{{\gamma}} \cap {\Gamma}\ \forall\ \!\gamma
$
and $ \hbox{\rm
cost}\!_{_{\Delta}}\!({{\delta_{1}}}\cup{{\delta_{2}}})
\!=\! \{\tau\!\in\!\Delta \mid
\neg[{{\delta_{1}}}\cup{{\delta_{2}}}\!\!\subset\!\!\tau]\}
\!=\! \{\tau\!\in\!{{\Delta}} \mid
\neg[\delta_{1}\!\!\subset\!\!\tau]\lor
\!\neg[\delta_{2}\!\!\subset\!\!\tau] \} \!\!= $ $ \hbox{\rm
cost}\!_{_{\Delta}}\!{{\delta_{1}}} \cup \hbox{\rm
cost}\!_{_{\Delta}}\!{{\delta_{2}}}, $ giving {\it a} and {\it
b}.

\medskip\noindent
{\it c}.\ A ``brute force''-check gives {\it c}, which is the
$``\tau\!\!=\!\!$ $\{$v$\}\hbox{''}$-case of {\it d}.

\medskip\noindent
{\it d}.\ $ {\tau}\!\in\!{{\rm Lk}\!_{_{\Delta}}\!{{\delta}}}\
\!\!\Longrightarrow\!\!{{\{{\bf v}\}}}\!\in\!{{\rm
Lk}\!_{_{\Delta}}\!{{\delta}}} \ \ \!\forall\ \!\hbox{\rm
\{v\}}\!\in\!{\tau}\in\Delta \!\Longleftrightarrow
{{\delta}}\!\in\!{{\rm Lk}\!_{_{\Delta}}\!\!{{\bf v}}}\ \!
(\subset\hbox{\rm cost}\!_{_{\Delta}}\!\!{{\bf v}})\ \ \!\! \forall\
\!{\bf v}\!\in\!{\tau}\in\Delta $
gives {\it d} from {\it a}, {\it c} and Proposition~\ref{PropP34:2}.a.{\bf
i} above.

\medskip\noindent
{\it e}.\ ${\hbox{\rm cost}\!_{_{\Delta}}\!\!{{\delta}}}\!$ =
[\underbar{iff} $n\!_{_{{\delta}}}\! \!=\!n\!_{_{\Delta}}\!\!-1\
\!\!]$ = $\underline{({\hbox{\rm cost}\!_{_{\Delta}}\!\!{{\delta}}})
\cap{\Delta\!^{^{\hbox{$_{\prime}$}\!}}}} = {\hbox{\rm
cost}\!\!_{_{\Delta\!^{^{_{\prime}}}\!}}\!{{\delta}}}$ $\ \!
\rlap{$^{^{_{\hbox{\sevenrm?}}}}$}{\!=}\ ({\hbox{\rm
cost}\!_{_{\Delta\!\!}}{{\delta}}})^{^{{\hbox{$_{\prime}$}\!}}}=$\
[\underbar{iff} $n\!_{_{{\delta}}}\! \!=\!n\!_{_{\Delta}}] =
\underline{({\hbox{\rm cost}\!_{_{\Delta}}\!\!{{\delta}}})
\cap{\Delta\!^{^{\hbox{$_{\prime}$}\!}}}}.$

\medskip\noindent
{\it f}.\ \ $\overline{\hbox{\rm {st}}}\!_{_{\Delta}}\!\delta\!=\! {
\bar{\delta}}\!\ast\!\hbox{\rm {Lk}}\!_{_{\Delta\!\!}}\delta=
{\dot{\delta}}\!\ast\!\hbox{\rm {Lk}}\!_{_{\Delta\!\!}}\delta
\ \cup
\{\tau\!\in\!\Delta \mid
\delta\!\subset\!\tau \}. $
\end{proof}

$ \hbox{\rm With}\ \delta\!_{_{^{\ \!\!}}}\!, \tau\!\!_{_{^{\
\!\!}}}\!\in\!\Sigma;$
$\delta\!_{_{^{\ \!\!}}}\cup \tau\!\!_{_{^{\
\!\!}}}\notin\!\Sigma\!_{o}\!
\!\Longleftrightarrow\!
\delta\!_{_{^{\ \!\!}}}\!\notin \overline{\hbox{\rm {st}}}\!\!
\lower1.1pt\hbox{${_{_{\Sigma\!_{_{ }}}}}$}\!\tau
\!\Longleftrightarrow\!
{\hbox{\rm {st}}}\!\! \lower1.1pt\hbox{${_{_{\Sigma\!_{_{
}}}}}$}\!\delta\ \!\cap\ \! |\overline{\hbox{\rm {st}}}\!\!
\lower1.1pt\hbox{${_{_{\Sigma\!_{_{ }}}}}$}\!\tau|
\!=\!\{\alpha\!_{_{^0}}{_{\!}}\}
\!\Longleftrightarrow
{\hbox{\rm {st}}}\!\! \lower1.1pt\hbox{${_{_{\Sigma\!_{_{
}}}}}$}\!\tau
\ \!\cap\ \!
{\hbox{\rm {st}}}\!\! \lower1.1pt\hbox{${_{_{\Sigma\!_{_{
}}}}}$}\!\delta
=\!
\{\alpha\!_{_{^0}}{_{\!}}\}
$
\indent and \indent
$ |\overline{\hbox{\rm {st}}}\!\!
\lower1.1pt\hbox{${_{_{\Sigma\!_{_{ }}}}}$}\!\sigma\!_{_{\!^{\
\!\!}}}|
\smallsetminus
{\hbox{\rm {st}}}\!\! \lower1.1pt\hbox{${_{_{\Sigma\!_{_{
}}}}}$}\!\sigma\!
=
|{\dot{\sigma}}\!_{_{\!^{\ \!\!}}}\ast
\hbox{\rm {Lk}}\!\! \lower1.1pt\hbox{${_{_{\Sigma\!_{_{
}}}}}$}\!\sigma|
=
|\overline{\hbox{\rm {st}}}\!\! \lower1.1pt\hbox{${_{_{\Sigma\!_{_{
}}}}}$}\!\sigma\cap
{\hbox{\rm {cost}}}\!\! \lower1.1pt\hbox{${_{_{\Sigma\!_{_{
}}}}}$}\!\sigma|$,
{\rm by\ \cite{26}
p.\ 372,\ 62.6\
and Proposition~\ref{PropP35:3}.{\it f} \ above.
}

\medskip\noindent
$$ \ \!\overline{\hbox{\rm {st}}}\!\!
\lower1.1pt\hbox{${_{_{\Sigma\!_{_{ }}}}}$}\!\sigma\!_{_{\!^{\
\!\!}}}
=
\{\tau\!_{_{\!^{\ \!\!}}}\in\Sigma\vert\ \sigma\cup\tau\in\Sigma\} =
{\bar{\sigma}}\!_{_{\!^{\ \!\!}}}\ast
\hbox{\rm {Lk}}\!\! \lower1.1pt\hbox{${_{_{\Sigma\!_{_{
}}}}}$}\!\sigma.
\eqno{\rm({\bf III})}\label{EqP34:III}$$

Identifying $|{\hbox{\rm {cost}}}\!\!
\lower1.1pt\hbox{${_{_{\Sigma}}}$}\!\sigma|$ with its homeomorphic
image in $|{\Sigma}|$ we get;
$
|{\hbox{\rm {cost}}}\!\! \lower1.1pt\hbox{${_{_{\Sigma}}}$}\!\sigma|
\simeq
|{\Sigma}| \smallsetminus
{\hbox{\rm {st}}}\!\! \lower1.1pt\hbox{${_{_{\Sigma\!_{_{
}}}}}$}\!\sigma
$
and
$
{\hbox{\rm {st}}}\!\! \lower1.1pt\hbox{${_{_{\Sigma\!_{_{
}}}}}$}\!\sigma
\!=\!
|{\Sigma}| \smallsetminus
|{\hbox{\rm {cost}}}\!\!
\lower1.1pt\hbox{${_{_{\Sigma}}}$}\!\sigma|.
$

\subsection{Homology groups for arbitrary simplicial joins}

\medskip
As said in p.~\pageref{EasyIso}, %
it is easily seen what Whitehead's $\widetilde{\bf S} (X\!\ast Y)$
need to fulfill to make
$\widetilde{\bf S} (X\!\ast Y) \approx\ \widetilde{\bf S}(X) \otimes
\widetilde{\bf S}(Y)$
true, since the right hand side is well known as soon as
$\widetilde{\bf S}(X)$ and $\widetilde{\bf S}(Y)$ are known.
But, a priori, we don't know what $\widetilde{\bf S} ((X, X\!_{_2})\ast
(Y, Y\!\!_{_2}))$ actually looks like.
Now, %
Theorem~\ref{TheoremP13:3} p.~\pageref{TheoremP13:3} %
gives the following complete answer:
\vskip-0.3cm
$${ \Delta\!\!^{^{\wp}}}\! {\bf (}(\!X,X\!_{_2}) \
\!\mhatast \ \! (Y,Y\!\!_{_2}){\bf )}
\msmallRingHom{Z}
 {\bf s}{\bf (}
{ \Delta\!\!^{^{\wp}}} \!(\!X\!_{_{\
\!}}\!,X\!_{_2})\nobreak\otimes\nobreak{\Delta \!\!^{^{\wp}}}
\!(Y_{_{\!\!\ \!}}\!,Y\!\!_{_2}){\bf )}.$$
\vskip-0.1cm

Below, we will express
$\mhatH_{_{k}}((\!\Gamma\!_{_1},\Delta\!_{_1})\ast
(\!\Gamma\!_{_2},\Delta\!_{_2}); {\bf G}{\rlap{$_{_1}$}{\
\otimes}}{\bf G}_{_2})$ in terms of
$\mhatH_{_{l}}(\!\Gamma\!_{_i},\Delta\!_{_i};{\bf G}_{_i})$ for
arbitrary simplicial complexes.
By Prop.~\ref{PropP34:1} p.~\pageref{PropP34:1} %
every simplex $\gamma\in \Gamma\!_{_1}\!\!\ast\Gamma\!_{_2}$ splits
\underbar{uniquely} into $\gamma_{_1}\cup \gamma_{_2}$ where
$\gamma_{_i}\in \Gamma_{_i}\ (i=1,2).$  Let the simplices in
$\Gamma_{_i}$ retain their original ordering and  put every vertex
in $\Gamma\!_{_2}$ ``after'' the vertices in $\Gamma\!_{_1}$. Let
$[\gamma_{_1}]\uplus[\gamma_{_2}]:=[\gamma]$ stand for the chosen
generator representing $\gamma$ in
$C^{o}\!\!\!\!_{_{^{v\!_{_1}\!+\!v\!_{_2}\!+\!1}}}
(\Gamma\!_{_1}\ast\Gamma\!_{_2})$ and,  extend the ``$\uplus$''
distributively and linearly, c.f. \cite{30}~p.~228,
\vskip-0.3cm
$$C^{o}\!\!\!\!_{_{^{v_{_1}\!\!+v_{_2}\!+\!1}}}\!\!
(\Gamma\!_{_1}\ast\Gamma\!_{_2})= \!\!\!\!\!\!
{\rlap{$_{_{_{_{_{p+q={v_{_1}+v_{_2}}}}}}}$} {\ \ \raise0pt\hbox{$\
\ \bigoplus$}}} (C^{o}\!\!\!_{_{^{p}}}(\Gamma\!_{_1})\uplus
C^{o}\!\!\!_{_{^{q}}}(\Gamma\!_{_2})).$$
\vskip-0.1cm

The following construction is taken from \cite{Fors2} p.~14.

The boundary function operates on generators as;
\vskip-0.5cm
$$\!\!\delta{\rlap{$\raise2pt\hbox{$_{_{_{^{v+1}}}}$}$}
{\lower5pt\hbox{$^{^{^{_{_{\Gamma_{_1}\ast\Gamma_{_2}}}}}}$}}[\gamma]
=
\delta^{\raise3pt\hbox{$_{_{\Gamma_{_1}\ast\Gamma_{_2}}}$}}
_{\hskip-0.4cm_{_{\dim\gamma_{_1}+\dim\gamma_{_2}+1}}}
\hskip-0.4cm[\gamma]\!} %
={\bigl((\delta{\rlap{$\lower5pt\hbox{$^{^{^{_{_{\Gamma_{_1}}}}}}$}$}
{\raise1pt\hbox{${_{_{\dim\gamma_{_1}}}}$}}}[\gamma_{_1}])\uplus
[\gamma_{_2}]\bigr)+} %
{\bigl([\gamma_{_1}]\uplus(\!-1\!)\!\!^{^{\dim\gamma_{_1}\!+\!1}}
\!\!(\delta{\rlap{$\lower-2pt\hbox{$^{{{_{_{\Gamma_{_2}}}}}}$}$}
{\raise1pt\hbox{${_{_{\dim\gamma_{_2}}}}$}}}
[\gamma_{_2}])\bigr).}$$
\vskip-0.1cm
Define a graded morphism $f$ through its effect on the generators of\\
$C^{o}\!\!_{_{v+1}}\!({\rlap{$\Gamma$}{\ _{_1}}}\!\ast\!
{\rlap{$\Gamma$}{\ _{_2}}})\!\!=\!\!
C^{o}\!\!_{_{v\!_{_1}\!+\!v\!_{_2}\!+\!1}}({\rlap{$\Gamma$}{\
_{_1}}} \!\ast\!{\rlap{$\Gamma$}{\ _{_2}}})$  by %
\vskip-0.2cm
$$f\!\!_{_{v+1}}\!([\gamma])\!\!=\!\!f_{_{v+1}}([\gamma_{_1}]\uplus[\gamma_{_2}])\!:=\!
[\gamma_{_1}]\otimes(\!-1\!)\!^{^{v+1}}[\gamma_{_2}],$$
\vskip-0.4cm
into
\vskip-0.4cm
$$[{\bf C^{o}}({\rlap{$\Gamma$}{\ _{_1}}})\otimes{\bf C^{o}}
({\rlap{$\Gamma$}{\
_{_2}}})]_{_{{v\!_{_1}\!\!+v\!_{_2}}}}\!\!\!=\!\!\!\!\!
{\rlap{$_{_{_{_{_{p+q={v_{_1}\!\!+v_{_2}}}}}}}$} {\ \
\raise0pt\hbox{$\ \bigoplus$}}} (C^{o}_{_{^{p}}}({\rlap{$\Gamma$}{\
_{_1}}})\!\otimes\! C^{o}_{_{^{q}}} (\Gamma\!\!{\ _{_2}})),$$
\vskip-0.1cm
\noindent
The boundary function is given through
its effect on\nobreak\ generators:
\vskip-0.3cm
$$\noindent\phantom{I}\hskip-1.75cm\delta{\rlap{$\raise1pt\hbox{$_{_{_{v}}}$}$}
{\lower4pt\hbox{$^{^{^{{\bf C\!^{o}}\!\!(\Gamma\!_{_1})\otimes {\bf
C\!^{o}}(\Gamma\!_{_2})}}}$}}}\![\tau]\!=\!
\delta{{\rlap{$\raise2pt\hbox{$_{_{_{\dim\gamma_{_1}+\dim\gamma_{_2}}}}$}$}
{\lower3pt\hbox{$^{^{^{{\bf C\!^{o}}(\Gamma\!_{_1})\otimes {\bf
C\!^{o}}(\Gamma\!_{_2})}}}$}}}}
\big([\gamma_{_1}]\otimes[\gamma_{_2}]\big)=$$
$$\noindent\phantom{II}\hskip1.75cm=\bigl((\delta{\rlap{$\lower4pt\hbox
{$^{^{^{_{_{\Gamma\!\!_{_1}\!\!}}}}}$}$}
{\raise0pt\hbox{${_{_{\dim\gamma_{_1}}}}$}}}[\gamma_{_1}])\otimes
[\gamma_{_2}]\bigr)+
\bigl([\gamma_{_1}]\otimes(\!-1\!)\!^{^{\dim\!\gamma_{_1}}}
\!(\delta{\rlap{$\lower4pt\hbox
{$^{^{^{_{_{\Gamma\!\!_{_2}\!\!}}}}}$}$}
{\raise0pt\hbox{${_{_{\dim\!\gamma_{_2}}}}$}}}
[\gamma_{_2}])\bigr).$$
\vskip-0.1cm\noindent
$f_{_{\ast}}$ is obviously a chain isomorphism of degree -1, and we
can conclude that,
\vskip-0.3cm
$${{\bf C}^{o}\bigl({\rlap{$\Gamma$}{\ _{_1}}}\ast
{\rlap{$\Gamma$}{\ _{_2}}}\bigr)\ {\rm and}\ {\bf s}\bigl({\bf
C}^{o}({{\rlap{$\Gamma$}{\ _{_1}}})\!\otimes\! {\bf
C}^{o}({\rlap{$\Gamma$}{\ _{_2}}}})\bigr)}\ {\rm are\ isomorphic\
as\ \underline{chains}},%
$$
\vskip-0.1cm
\noindent where the
``{\bf s}'' stands for ``suspension'' meaning that the suspended
chain  equals the original except that the dimension $i$ in the
original is   dimension $i+1$ in the suspended chain.

The argument now motivating the formula for \underbar{relative} simplicial homology
is the same as that for the ``relative singular
homology''-case in p.~\pageref{RelSingHom}ff. %

\medskip
\noindent{\bf Lemma.} \label{SimplEilenbergZilber} (The
Eilenberg-Zilber theorem for simplicial join.)
{\sl On the category of ordered pairs of arbitrary simplicial pairs
$(\Gamma\!\!_{_{1}},\Delta_{_{1}}\!)\!$ and
$(\Gamma\!\!_{_{2}},\!\Delta_{_{2}}\!)$
there is a natural chain equivalence of} \hskip0.5cm %
\vskip-0.4cm
$${\bf C}^{o}(\Gamma\!\!_{_{1}}\!\!\ast\!\Gamma\!\!_{_{2}})/ {\bf
C}^{o}\big((\Gamma\!\!_{_{1}}\ast\Delta_{_{2}})\cup
(\Delta_{_{1}}\ast\Gamma\!\!_{_{2}})\big)$$
\vskip-0.4cm
{\sl with}
\vskip-0.4cm
$${\bf s}\big\{\big[{\bf
C}^{o}(\Gamma\!\!_{_{1}})/ {\bf C}^{o}(\Delta_{_{1}})\big]
\otimes\big[{\bf C}^{o}(\Gamma\!\!_{_{2}})/{\bf C}^{o}
(\Delta_{_{2}})\big]\big\}\eqno{\square}$$

\begin{varprop}{\bf 1.}\ \
(The Künneth Theorem for Arbitrary Simplicial Joins.) With
{\bf R} a {\bf PID}, {\bf G} and ${\bf G^\prime}$ {\bf R}-modules
and ${\rm Tor}_1^{\bf R}({\bf G,G^\prime})=0$, {\rm\cite{30}~p.~231 Cor.~4} gives;
\vskip-0.4cm
$$
\begin{array}{l}
\noindent\hskip-0.5cm\mhatH_{_{q+1}}
        ((\!\Gamma\!_{_1},\Delta\!_{_1})\ast
(\Gamma\!_{_2},\Delta\!_{_2}); {\bf G}\otimes{\bf G^{\prime}})
\widetilde=\nonumber\\
\noindent\hskip-0.0cm\widetilde={\rlap{$_{_{_{i+j=q}}}$} {\ \raise2pt\hbox{$\bigoplus$}}}
[\mhatH_{i}
           ({\rlap{$\Gamma$}{\ _{_1}}},\Delta\!_{_1};{\bf G})
\otimes
        \mhatH_{j}
({\rlap{$\Gamma$}{\ _{_2}}},\Delta\!_{_2};{\bf G^{\prime}})]
\oplus{\rlap{$_{_{_{i+j=q-1}}}$}{\ \ \raise2pt\hbox{$\bigoplus$}}}
{\rm Tor}_1^{\bf R}\bigl(\mhatH_{i}
           ({\rlap{$\Gamma$}{\ _{_1}}},\Delta\!_{_1};{\bf G}),
 \mhatH_{j}
      ({\rlap{$\Gamma$}{\ _{_2}}},\Delta\!_{_2};
{\bf G^{\prime}})\bigr)
\hskip0.1cm\square\hskip-0.2cm %
\end{array}%
\eqno{^{\hbox{({$\diamondsuit$})\hskip-0.0cm}}}
$$
\end{varprop}

\break
\subsection{%
Local homology for products and joins of arbitrary simplicial complexes}

\medskip
\noindent
{\bf Corollary}~\ref{CorP19:ToTh6}. (from p.~\pageref{CorP19:ToTh6})\label{LinkBevisS88}
Let ${\bf G}$ and ${\bf G^{\prime}}$
be arbitrary modules over a  principal ideal domain {\bf R} such that \
${\rm Tor}_1^{\bf R}({\bf G},{\bf G}^\prime)=0$,
then, for any
$\emptyset_o\ne\sigma\in\Sigma_{1}\times\Sigma_{2}$ with
$\eta\big({\rm Int}(\sigma_{})\big)\subset
{\rm Int}(\sigma_{1})\times{\rm Int}(\sigma_{2})$ and
$c_{\sigma}\!:=\!\dim\sigma_{_{\!1}}\!+
\dim\sigma_{_{\!2}}\!-\dim\sigma$;

\smallskip

\indent
${\underline{\underline{\mhatH_{i+c_{\!\sigma}+1}
({\rm Lk}_{_{\Sigma_{1}\times\Sigma_{2}}}\sigma;
{\bf G}\otimes_{_{\bf R}}{\bf G^{\prime}})}}}
\mringHom{R}
{\rlap{$_{_{_{{{{{p+q=i}\atop{p,q\ge -1}}}}}}}$}{\ \ \raise2pt\hbox{$\ \bigoplus$}}}\ [\mhatH_{p}
           ({\rm Lk}_{\Sigma_{1}}\sigma_{1};{\bf G})
\otimes_{_{\bf R}}
        \mhatH_{q}
({\rm Lk}_{\Sigma_{2}}\sigma_{2});{\bf G^{\prime}})]\oplus
$

\smallskip
\rightline{${\rlap{$\oplus$}{\ {\rlap{$_{_{_{_{_{{{p+q=i-1}\atop{p,q\ge
-1}}}}}}}$} {\ \ \raise2pt\hbox{$\ \bigoplus$}}}}}
{\rm Tor}_1^{\bf R}
\bigl(\mhatH_{p}
           ({\rm Lk}_{\Sigma_{1}}\sigma_{1};{\bf G}),\
 \mhatH_{q}
      ({\rm Lk}_{\Sigma_{2}}\sigma_{2};{\bf G^{\prime}})\bigr)
\mringHom{R}\ \!
{\underline{\underline{\mhatH_{i+1}
        ({\rm Lk}_{_{\Sigma_{1}\ast\Sigma_{2}}}(\sigma_1\cup\sigma_2); {\bf
G}\otimes_{_{\bf R}}{\bf G^{\prime}})}}}.$}

\smallskip
\noindent
If $\dim\sigma\!=\!\dim\sigma_{_{\!1}}\!+
\dim\sigma_{_{\!2}}$ then
$c_{\!\sigma}\!=0$ and so:
$\mhatH_{i}
({\rm Lk}\hskip-0.3cm\lower0.05cm\hbox{$_{_{\Sigma_{1}\times\Sigma_{2}}}$} \hskip-0.3cm\sigma;\ {\bf G})
\mringHom{R}
\mhatH_{i}
     ({\rm Lk}\hskip-0.3cm\lower0.04cm\hbox{$_{_{\Sigma_{1}\times\Sigma_{2}}}$}\hskip-0.2cm(\sigma_{_{\!1}}\!
\cup\sigma_{_{\!2}})
                ;{\bf G})$.

\noindent
\hskip-0.0cmSince ${\rm Tor}_1^{\bf R}({\bf G,G^\prime})\!=\!0,$
there is no torsion in degree $-1$.
So, if $\sigma\!\ne\!\emptyset\!_{_{^{o}}}\!$ and $c_{\!\sigma}\!=\!0$;

\smallskip
\noindent%
$\mhatH_{_{0}}
 ({\rm Lk}\hskip-0.3cm\lower0.04cm\hbox{$_{_{\Sigma_{1}\times\Sigma_{2}}}$}
\hskip-0.3cm\sigma;{\bf G}\otimes{\bf G^{\prime}})
\mringHom{R}
{
\mhatH_{_{0}}({\rm Lk}_{_{\Sigma\!_{_{1}}}}\!\!\!\sigma\!_{_{1}};
{\bf G})
\otimes\mhatH\!\!_{_{-\!1}}\!
({\rm Lk}_{_{\Sigma\!_{_{2}}}}\!\!\!\sigma\!_{_{2}};
{\bf G^{\prime}}) \!\oplus\!
\mhatH\!\!_{_{-\!1}}\!
({\rm Lk}_{_{\Sigma\!_{_{1}}}}\!\!\!\sigma\!_{_{1}};{\bf G})}
\!\otimes\!
        \mhatH_{_{0}}({\rm Lk}_{_{\Sigma\!_{_{2}}}}\!\!\!\sigma\!_{_{2}};
       {\bf G^{\prime}}).$

\smallskip
\noindent{\bf Proof.} The one-vertex complex $\bullet:=\{\{v\},\emptyset\!_o\}$
is a kind of unit-element with respect to ``$\times$'', so, our claims
are true if either  $\Sigma_{1}$ or $\Sigma_{2}$ equals $\bullet$.
Therefore, suppose that  $\Sigma_{1}\!\ne\!\bullet\!\ne\!\Sigma_{2}$.
Lemma~\ref{LemmaP18} p.~\pageref{LemmaP18} and %
\cite{7} %
Th.~12.4 p.\ 89 gives $(\ast)$,

\medskip

\noindent%
{\bf $(\ast)$}\indent
$\eta:(|\Sigma_{_{\!1}}\nabla \Sigma_{_{\!2}}\!|,
|\Sigma_{_{\!1}}\!\nabla \Sigma_{_{\!2}}\!|\setminus
{(\widetilde{\alpha_{_{\!1}}\!,\alpha_{_{\!2}}\!})}\})\
\rlap{$\longrightarrow$}{\raise4pt\hbox{$\ \simeq$}}\ \ \
(|\Sigma_{_{\!1}}\!|{\bf {\bar\nabla}}
|\Sigma_{_{\!2}}\!|,  |\Sigma_{_{\!1}}\!|{\bf {\bar\nabla}}
|\Sigma_{_{\!2}}\!|
\setminus \{(\alpha_{_{\!1}}\!,\alpha_{_{\!2}}\!)\}),$ with\\
\indent\indent$\eta(\widetilde{\alpha_{_{\!1}}\!,\alpha_{_{\!2}}\!})=\nobreak
(\alpha_{_{\!1}}\!,\alpha_{_{\!2}}\!).$

\medskip
Note that
$\sigma\ne\emptyset_o\Rightarrow\sigma_{_{\!\!j}}\ne\emptyset_o,\ j=1,2.$
Put $v_{_{\!j}}\!:=\!{\rm dim}\sigma_{_{\!\!j}},\ j=1,2\ {\rm and}\
v\!:=\!v_{_{\!1}}\!+\!v_{_{\!2}}.$

\smallskip
\noindent $\mhatH\hskip-0.4cm_{_{_{s\lower1.1pt\hbox{-}
(v\lower1.1pt\hbox{-}c_{\!\sigma})\lower1.1pt\hbox{-}1}}}
\hskip-0.2cm
({\rm Lk}\hskip-0.3cm_{_{_{\Sigma_{1}\!\times\!\Sigma_{2}}}}\hskip-0.3cm\sigma;
{\bf G}\otimes{\bf G^{\prime}})
\mringHom{R}
\bigl[{\small\rm{Prop.~\ref{PropP15:1}\ p.~\pageref{PropP15:1}}\atop{{\rm dim}\sigma=v-c_{\!\sigma}}}\bigr]\!
\mringHom{R}
\mhatH_{s}(|\Sigma_{1}{\bf{\times}}\Sigma_{2}|,
                        |\Sigma_{1}{\bf {\times}}\Sigma_{2}|\setminus_o
                       {(\widetilde{\alpha_1,\alpha_2})};
                   {\bf G}\otimes{\bf G^{\prime}})
\mringHom{R}
$

\smallskip
\noindent
$
\mringHom{R}\
\bigl[{{\bf (\ast)}\atop{\rm above}}\bigr]
\mringHom{R}
\!$
$\mhatH_{s}
(|\Sigma_{1}|{\bf {\times}}|\Sigma_{2}|,
|\Sigma_{1}|{\bf {\times}}|\Sigma_{2}|\!\setminus_o\!
(\alpha_1,\alpha_2);{\bf G}\otimes{\bf G^{\prime}})
\ \!\mringHom{R}\!
$
$\bigl[{{\rm \text{Def. pair-product p.}~\pageref{DefP11:PairDef}}
\atop{\text{+ simple calculus}}}\bigr]\!
\mringHom{R}$

\smallskip
\noindent$\mringHom{R}\
\mhatH_{s}
          (|\Sigma_{1}|,|\Sigma_{1}|\setminus_o\alpha_1)
           {\bf {\times}}(|\Sigma_{2}|,|\Sigma_{2}|\setminus_o\alpha_2)
           ;{\bf G}\otimes{\bf G^{\prime}})
\mringHom{R}
$
$\bigl[{\small\rm{Eq.~(3) p.~\pageref{EqP14:3}}
\atop{\text{ +~Lemma p.~\pageref{lemmaP15}.}}}\bigr]
\mringHom{R}$

\smallskip
\noindent
$\mringHom{R}$
$[\mhatH_{p_1}
     (|\Sigma_{1}|,|\Sigma_{1}|\setminus_o \alpha_1;
\noindent{\bf G})\ \otimes
$
$
\mhatH_{p_2}(|\Sigma_{2}|,|\Sigma_{2}|\setminus_o\alpha_2;
                 {\bf G^{\prime}})]
\!\!\!\!_{{p\!_{_{^{j}}}\ge0}\atop{{ p}_1+{ p}_2={ s}}}
\!\!\oplus$

\smallskip
\noindent
$\oplus[{\rm Tor}_1^{\bf R}\bigl(\mhatH_{p_1}
           (|\Sigma_{1}|,|\Sigma_{1}|\setminus_o \alpha_1;{\bf G}),\
 \mhatH_{p_2}
(|\Sigma_{2}|,|\Sigma_{2}|\setminus_o \alpha_2;
{\bf G^{\prime}})\bigr)]\!\!\!\!\!_{{p\!_{_{^{j}}}\ge0}\atop{p_1+p_2=s-1}}
\!\!\!\mringHom{R}
$

\smallskip
\noindent
$
\mringHom{R}\
[\text{Prop.~\ref{PropP15:1}\ p.~\pageref{PropP15:1}}]
\mringHom{R}
$
$
[\mhatH\hskip-0.3cm_{_{p_1-v_1-1}}
                    ({\rm Lk}_{_{\Sigma_{1}}}\!\!\sigma_{1};{\bf G})
    \otimes \mhatH\hskip-0.3cm_{_{p_2-v_2-1}}
           ({\rm Lk}_{_{\Sigma_{2}}}\!\!\sigma_{2};{\bf G^{\prime}})]
\hskip-0.6cm
_{{p\!_{_{^{j}}}\!-v\!_{_{^{j}}}\ge0}\atop{(p_1-v_1)+(p_2-v_2)=s-v}}
\hskip-0.3cm\oplus
$

\noindent
$
\oplus\ [{\rm Tor}_1^{\bf R}\bigl(
\mhatH\hskip-0.3cm_{_{p_1\!-\!v_1\!-\!1}}
                    ({\rm Lk}_{_{\Sigma_{1}}}\!\!\sigma_{1};
{\rlap{\bf G}{\ \ \ \!)}}
    ,\mhatH\hskip-0.3cm_{_{p_2\!-\!v_2\!-\!1}}
           ({\rm Lk}_{_{\Sigma_{2}}}\!\!\sigma_{2};
{\bf G^{\prime}})\bigr)]\hskip-0.7cm
_{{{p\!_{_{^{j}}}\!-\!v\!_{_{^{j}}}\ge0}\atop
{(p_1\!-\!v_1)+(p_2\!-\!v_2)=s\!-\!v-\!1}}}$
$\hskip-0.8cm =[\text{Put}\ {{q_{_{^{j}}}:=p\!_{_{^{j}}}-v\!_{_{^{j}}}-1}}]=$\\
$
=\big[\mhatH_{_{\!q_{_1}\!}}\!\!
                    ({\rm Lk}_{_{\Sigma_{1}}}\!\!\sigma_{1};{\bf G})\otimes\ \mhatH_{_{\!q_{_2}\!}}\!\!
           ({\rm Lk}_{_{\Sigma_{2}}}\!\!\sigma_{2};{\bf G^{\prime}})]
\hskip-0.4cm_{{q_{_{^{j}}}\ge-1}\atop{q_1+q_2=s-v-2}}$
$\hskip-0.4cm\oplus
[{\rm {\rm Tor}_1^{\bf R}}\bigl(
\mhatH_{_{\!q_{_1}}}\!\!
                    ({\rm Lk}_{_{\Sigma_{1}}}\!\!\sigma_{1};{\bf G})
,\mhatH_{_{\!q_{_2}}}\!\!
           ({\rm Lk}_{_{\Sigma_{2}}}\!\!\sigma_{2};{\bf G^{\prime}})\bigr)]
\hskip-0.1cm{_{{q_{_{^{j}}}\ge-1}\atop{\hskip-0.6cm q_1+q_2=s-v-3}}}\hskip-0.15cm.
$\\
Now, put $i:=s-v-2.\ \triangleright$

\smallskip
Proposition~\ref{PropP34:1} p.~\pageref{PropP34:1} and Theorem~\ref{TheoremP14:4} p.~\pageref{TheoremP14:4} gives the
second isomorphism even for $\sigma\!_{_{1}}\!\!=\!\emptyset\!_{_{^{o}}}\
{\rm and}\ \!\!/\ \!\!{\rm or}\ \sigma\!_{_{2}}\!\!=\!\emptyset\!_{_{^{o\!}}}.$
\hfill$\square$

\break
\subsection{%
On Relative Homeomorphisms and Quotients}
\vskip-0.01cm
\vskip-0.01cm

An excision map $(X\setminus P, Y\setminus P)\mapsto (X, Y)$
induces an isomorphism in singular homology if ${\bf\bar{\it
\!P}}\subset {\rm Int}(Y)$
while ``$P$ open in $X$" is not needed, cf. \cite{26}
p.\ 180-1.
This will be used in the next proposition. Note that $\mhatH_{_{i\!}}({X}, \bullet;{\bf G})
\equiv\mhatH_{_{i\!}}({X};{\bf G})$.

\begin{varprop} \label{PropP88}
Let $ f\!\!:\!(X,A)\!\!\rightarrow\! \!(Y,B) $ be a
relative homeomorphism, i.e., $f\!:X\!\rightarrow\!Y$
is continuous, and $f\!:X\setminus A\!\rightarrow\!
Y\setminus B$ is a homeomorphism. If $A\ne \emptyset,\{\wp\}$ is a strong
deformation retract of a
neighborhood $N$ and $B$ and $f(N)$ are closed subsets of $N^{^{_{\prime}}}\!\!:=f(N\setminus A)\cup B=f(N)\cup B$, since $f(A)\subset B$, then;

\noindent {\bf i.}\label{PropP88i} {\rm(from p.~\pageref{LemmaP12})} $B$ is a strong deformation retract of $N^{\prime}$.

\indent If also ${\bf\bar{\it \!A}}\!\subset {\rm
Int}\raise1pt\hbox{$_{_{{\!X}}}$}\!(N)$,\ \ ${\bf\bar{\it
\!B}}\!\subset {\rm
Int}\raise1pt\hbox{$_{_{{Y}}}$}\!(N^{^{_{\prime}}})$ and
$f\!_{{{|}}}:(X\setminus A,N\setminus A)\longrightarrow(Y\setminus B,N^{^{_{\prime}}}\setminus B)$
 is a pair homeomorphism {\rm(which it is if $X,Y$ is Hausdorff and $N$ compact)} then;

\smallskip
\noindent
{\bf ii.} \label{PropP88ii}{\rm (cp. \cite{36}
p.~66 Lemma~7.3)}
$ \mhatH\raise1pt\hbox{$_{_{^{i}}}$} (X, A;{\bf G})\!\
\cong\ \mhatH\raise1pt\hbox{$_{_{^{i}}}$} (Y\!, B;{\bf G})
$ {\sl for} $i\in {\bf Z}$ {\sl and}
${\bf \Delta{\raise1.5pt\hbox{$\!\!^{\wp}$}}}\! \hbox{\bf (}X,
A\hbox{\bf )} \!\approx\! {\bf
\Delta{\raise1.5pt\hbox{$\!\!^{\wp}$}}}\! \hbox{\bf (}Y\!,
B\hbox{\bf )}. $
{\bf G} is an {\bf A}-module where {\bf A} is a commutative ring with unit.

\smallskip
\noindent {\bf iii.} \label{PropP88iii} $A$ closed \hbox{$\Rightarrow
\mhatH_{_{i\!}}({X},{{ A}};{\bf G}) \cong
\mhatH_{_{i\!}}({X}/{{ A}}, \bullet;{\bf G})
\cong\mhatH_{_{i\!}}({X}/{{ A}};{\bf G})
$.
${\bf \Delta{\raise1.5pt\hbox{$\!\!^{\wp}$}}}\! \hbox{\bf
(}X, A\hbox{\bf )} \!\approx\! $
${
\Delta{\raise1.5pt\hbox{$\!\!^{\wp}$}}}\! \hbox{\bf (}{X}/{{
A}},\bullet\hbox{\bf )}. $}
\end{varprop}

\vskip-0.2cm
\noindent {\bf Proof.} {\rm (cp.~ \cite{30}
p.~202 Th.~9, proof)} Let
$F:N\!\times{\bf I}\longrightarrow N$
be the postulated strong neighborhood deformation
retraction of $N$ down onto $A$ and define;
\vskip-0.6cm
$$\hskip-0.2cmF^\prime : N^\prime\!\!\times{\bf
I}\!\rightarrow N^\prime\ ;\
\left\{\begin{array}{ll}
(y,t) \mapsto\
y & \text{if}\
y \in B,\ t \in {\bf I} \\
(y,t)\mapsto
{f\circ F (f
^{^{_{-1}}}\! (y), t) }  & \text{if}\ y \in
f(N)\setminus B = f(N\setminus{A}),\ t
\in {\bf I}. %
\end{array}\right.
$$
\vskip-0.2cm
\noindent
{\bf i.}
$F^{\prime}$ is well-defined and it is continuous as being so when restricted to any one of the closed subspaces
$f(N)\times{\bf I}$ resp. $B\times{\bf I}$ that together
cover $N^{\prime}\times{\bf I}$,
cf.~\cite{30} ~{\bf 4} p.~5.
({\bf N.B.} Changing the domain $f(N)\setminus B$ to $f(N)$ doesn't alter $F^{\prime}$.)\\
{\bf ii.} The following standard diagram gives ii.

\FFrame{0.0pt}{0.0pt}{\hsize=1.00\hsize
\noindent
\raise3pt\hbox{\FFrame{0.01cm}{0.01pt}{\hsize=0.76\hsize\lineskip-2pt
\noindent $
\ \mhatH\raise2pt\hbox{$_{_{{\!\star}}}$}\!(X,A;{\bf G})\ {\rlap{\ \raise3.5pt\hbox{{\ $\cong$}}}
{\vbox{\moveleft0pt\hbox{{$\longrightarrow$}}}}}\
\mhatH\raise2pt\hbox{$_{_{{\!\star}}}$}\!(X,N\!;{\bf G})\ {\rlap{\ \raise3.5pt\hbox{$\cong$}} {\vbox{\moveleft2pt\hbox{{$\longleftarrow$}}}}}\
\mhatH\raise1pt\hbox{$_{_{^{\!\star}}}$}\! (X\setminus A,
N\setminus A;{\bf G}) $\\
$ {\hskip-0.1cm
\ \ \ \msmallhatH\hbox{{\eightrm(}{\eighti f}{\eightrm)}}
{\Big\downarrow}
\hskip2.5cm
\Big\downarrow
\hskip2.5cm{\Big\downarrow}\!\cong
\hbox{\hsize4cm\Big(${{{\fivei f}\!_{^{{|}}}\
\text{is a pair ho-} }%
\atop{{\rm meomorphism} }}$\Big)}\!\! } $\\
$\indent\hskip-0.3cm \mhatH\raise2pt\hbox{$_{_{{\!\star}}}$}\!(Y,B;
{\bf G})\ \!
{\rlap{\ \raise3.5pt\hbox{{\ $\cong$}}}
{\vbox{\moveleft0pt\hbox{{$\longrightarrow$}}}}}\ \!
\mhatH\raise2pt\hbox{$_{_{{\!\star}}}$}\!(Y,N^{^{_{\prime\!}}}\!;{\bf
G})\ \!
{\rlap{\ \raise3.5pt\hbox{{\ $\cong$}}}
{\vbox{\moveleft1pt\hbox{{$\longleftarrow$}}}}}\
\!\mhatH\raise1pt\hbox{$_{_{^{\!\star}}}$}\! (X\!\setminus\!B,
N^{^{_{\prime}}}\!\setminus\!B;{\bf G})$}}}
\vskip-0.1cm

\noindent
{\bf iii.} (\cite{7}
p.~125 Ex.~1) $f_{(|)}:= (\text{restriction of})\ {\rm projection}$. Use as
\noindent second diagram row;
$
\mhatH\raise1pt\hbox{$_{_{^{\!\star}}}$}\!
(X\!/\!A,\bullet;{\bf G})\ {\rlap{\ \raise2.5pt\hbox{{\ $\cong$}}}
{\vbox{\moveleft0pt\hbox{{$\longrightarrow$}}}}}\
\mhatH\raise2pt\hbox{$_{_{{\!\star}}}$}\!(X\!/\!A,N\!/\!A;{\bf
G})\ {\rlap{\ \raise2.5pt\hbox{{\ $\cong$}}}
{\vbox{\moveleft1pt\hbox{{$\longleftarrow$}}}}}\
\mhatH\raise1pt\hbox{$_{_{^{\!\star}}}$}\!
(X\!/\!A\!\setminus\!(A\!/\!A),
N\!/\!A\!\setminus\!(A\!/\!A);{\bf G}).\hfill\triangleright $

\normalbaselines
\noindent
\cite{26}
Th.~46.2 p.~279, see p.~\pageref{TheoremP13:Munkres}, gives the chain equivalences in {\bf ii} and {\bf iii}.$\hfill\square$

\medskip
\normalbaselines \noindent
\FFrame{0.0pt}{0.0pt}{\hsize=0.54\hsize
{ {{
\vbox{
\noindent
The quotient maps $p\!_{_{1}}, p\!_{_{2}}$ w.r.t. $\ast$-join and sus\-pen\-sion, makes
$p\!_{_{2}}\!\circ\!p_{_{1}}\!\!^{_{
\raise-1pt\hbox{\fiverm$\ \!$\hbox{\sevenbf-}1}}}$
con\-tinuous and $(p\!_{_{4\!\!}}\circ_{_{\!}}p\!_{_{2\!}})\circ
(p_{_{\!3\!}}\circ_{_{\!}}
p\!_{_{1\!\!}})^{\hbox{\sevenbf-}\!1}$
a homeomorphism,
cf. \cite{7}
p.~134 Ex.~5.\\
${\bf v}\!_{_{^1\!}},{\bf v}\!_{_{^2\!}}$ are the
suspension points.
\hskip-0.1cm}}}}}\hskip-0.1cm
\hbox{\FFrame{0.03cm}{0.01pt}{\hsize=0.42\hsize
\lineskip=-3pt
\vbox{\moveright1.9cm \hbox{$\hbox{\eightbf S}
\lower1pt\hbox{$ _{_{^{\!{\bf v}\!_{_{^1\!}}\!,\!{\bf
v}\!_{_{^2}}\!}}}\!\!( {\hbox{\eighti
X}_{\!\wp}\!\!\times\!\hbox{\eighti Y}_{\!\wp}} $} ) $} \noindent
$\hbox{{\eighti X}$_{\!\wp_{\!}}$}\times\hbox{{\eighti
Y}$_{\!\wp_{\!}}$} \times\hbox{{\eightbf I}$_{_{\!}\wp}$}\ \!
\rlap{\raise6.5pt\hbox{$\rlap{\ {$\nearrow$}} {^{_{\hbox{{\eighti p}}_{_{^{\!2}}}}}}$}}
{\lower6.5pt\hbox{$\rlap{\ {$\searrow$}}{_{^{\hbox{{\eighti p}}_{_{^{\!1}}}}}}$}} $\ \
$\hbox{{\eighti p}}_{_{^{\!2}}}\! {\lower3.5pt\hbox{$^\circ$}}
\hbox{{\eighti
p}}_{_{^{\!1}}}\!\!\!^{\lower0.5pt\hbox{\sevenbf-{\fivebf1}}}\!\!
\Big\uparrow \ \ \rlap{\raise6.5pt\hbox{$\rlap{$\searrow$} {\ \
^{\hbox{{\eighti p}}_{_{^{\!4}}}}}$}}
{\lower6.5pt\hbox{$\rlap{$\nearrow$} {\ \ _{^{\hbox{{\eighti
p}}_{_{^{\!3}}}}}}$}}\! {{\lower0.0pt\hbox{{\eighti
X}$_{\!\wp}\ast${\eighti Y}$_{\!\wp}$}} \over
{\lower2.0pt\hbox{{\eighti X}$_{\!\wp}$+{\eighti Y}$_{\!\wp}$}}}
\!\simeq\! { {\bf S}\!_{_{^{{\bf v}\!_{_{^1}}\!\!,\!{\bf
v}\!_{_{^2}}\!\!\!}}} ( {{X}_{\!\wp}\!\times\!{Y}_{\!\wp}\!} )
\over \!\! \lower2.5pt\hbox{$\{{\bf v}\!_{_{^1}}\!,\!{\bf
v}\!_{_{^2}}\!,\!\wp\}$} }\\
$ {\lower0.0pt\hbox{ \vbox{\moveright2.1 cm\hbox{
{\lower3.0pt\hbox{{\eighti X}$_{\!\wp}\ast${\eighti Y}$_{\!\wp}$\hskip-0.1cm
}} }} }}}}} \normalbaselines
\vskip-0.1cm\noindent
\vskip-0.1cm
\vskip-0.1cm

\begin{varlemma}
{\rm(\cite{6}
p.\ 225)}
$ f\!{\!}:{\!}Z\!{\!}\rightarrow\!{\!} W {\!}$ is null-homotopic
\underbar{iff} W is a retract of the mapping cone $
({\bf C}\!Z)\cup W.\ ({\bf C}\!Z\!:=\!Z\!\ast\bullet.)\hfill\square$
\end{varlemma}
\vskip-0.1cm
\vskip-0.1cm

\noindent
$Z\!\subset W\!$,
$i:={\rm id}_{{{Z}}}$ and excising the cone-top followed by a retraction gives;

\normalbaselines
\vskip-0.06cm
\noindent\FFrame{0.0cm}{0.0pt}{\hbox{\hsize=1.05\hsize\hskip-0.05cm
\vbox{\FFrame{0.01cm}{0.0pt}{\hsize=0.58\hsize
\noindent
$
\mhatH\!\raise1pt\hbox{$_{_{_{q}}}$}\! ({\bf
C}\!Z\cup\!_{_{^{i}}}{_{\!}}{_{\!}}W) \ \!\cong\ \!
\mhatH\!\raise1pt\hbox{$_{_{_{q}}}$}\! ({\bf
C}\!Z\cup\!_{_{^{i}}}{_{\!}}{_{\!}}W,{\bf C}Z) \ \!\cong\ \!
\mhatH\!\raise1pt\hbox{$_{_{_{q}}}$}\!(W,Z).$\\
So, if $i\!_{_{\ }}\!$ is null-homotopic, $i\!_{_{\ast}}\!$ is trivial
and $l\!_{_{\ast}}\!$ splits, i.e:\\
\indent\hskip2cm$\mhatH\!\!\lower0.5pt\hbox{${_{_{^{\star+1}}}}$}\!(W,Z)\
\!\widetilde=\
\mhatH\!\!\lower0.5pt\hbox{${_{_{^{\star+1}}}}$}\!(W)
\oplus \mhatH\!\lower0pt\hbox{$_{_{{\star}}}$}\!(Z).$\\
Now, ${
i{_{_{^{\!^X}}}}\!\!+i{_{_{^{\!^Y}}}}\!\!:{X}\!+ {Y}
\!\hookrightarrow\!
{X}\!\ast\!{Y} } $
is null-homotopic, so:
\hskip-0.2cm
\vskip-0.1cm}}\hskip-0.05cm
\vbox{\FFrame{0.08cm}{0.1pt}{\hsize=0.325\hsize
\lineskip=-2pt
$ \indent\hskip0.22cm
\!\!{_{\!}}{_{\!}}{\lower4.6pt\hbox{$|$}}
\!{\buildrel{r\!_{_{\ast}}}\over\longleftarrow}
\mhatH\!\raise1pt\hbox{$_{_{_{\star}}}$}\! (({\bf
C}Z)\!\cup\!\!_{_{^{i}}}{_{\!}}{_{\!}}W)\\
\indent\hskip0.5cm
\raise1.5pt\hbox{{\vbox{\moveleft3.4pt\hbox{$\downarrow$}}}}
\hskip0.1cm
^{j{_{_{\!\ast}}}}\!\!\!\lower0pt\hbox{$\nearrow$} \hskip0.6cm \downarrow\cong\\
\noindent
\lower2.55pt\hbox{$
\buildrel{i\!_{_{\ast}}}\over\rightarrow$}
\mhatH\raise0pt\hbox{$_{_{{\!\star}}}$}\!(W)
{\buildrel{l\!_{_{\ast}}}\over\rightarrow}
\mhatH\raise0pt\hbox{$_{_{{\!\star}}}$}\!(W,Z)
{\buildrel{\partial\!_{_{\ast}}}\over\rightarrow}
\mhatH\raise1pt\hbox{$_{_{^{\!\star\hbox{\fivebf-}1}}}$}\!(Z)
{\buildrel{i\!_{_{\ast}}}\over\rightarrow
}$\hskip-0.2cm \vskip-0.0cm}}}}
\vskip-0.05cm
\noindent
$
{{\mhatH}}\!\!_{{{q\!+\!1}}\!}(X \ast Y ,X\!+Y;{\bf
G})\mringHom{A}$
\vskip-0.0cm
$\mringHom{A}\mhatH_{q+1}\!(X\!\ast\!Y\!;{\bf
G})\oplus \mhatH_{q}(X+ Y;{\bf G})
\mringHom{A}
\mhatH_{q}\!(X\!\ast\!Y\!;{\bf
G}) \oplus
\mhatH_{q}(X;{\bf G}) \oplus
\mhatH_{q} (Y,\{\wp\};{\bf G}).
$

The l.h.s. equals $
\mhatH_{q}(
{{X}\times{Y}},\{\wp\};{\bf G})$. Deleting
$\{\wp\}$ 
gives Eq.\ 2 p.~\pageref{EqP12:2} first line.

\section{Addenda}\label{Addendium:A.1}
\vskip-0.3cm

\subsection{The importance of simplicial complexes}
In \cite{9}
p.\ 54 Eilenberg and Steenrod wrote:
\begin{quotation}
Although triangulable spaces appear to form a rather narrow
class, a major portion of the spaces occurring in applications
within topology, geometry and analysis are of this type.
Furthermore, it is shown in Chapter\ X that any compact space can be expressed as a limit of triangulable spaces in a reasonable sense.
In this sense, triangulable spaces are dense in the family of
compact spaces.
\end{quotation}
Moreover, not only are the major portion of the spaces occurring in applications triangulable, but the triangulable spaces plays an even more significant theoretical role than that, when it comes to determining the homology and homotopy groups of arbitrary topological spaces, or as Sze-Tsen Hu, in his {\it Homotopy Theory} (\cite{18} 1959) p.\ 171, states, after giving a description of the Milnor realization of the Singular Complex {\it S}({\it X}) of an
\hbox{arbitrary space {\it X}:}

\begin{quotation}
The
significance of this is that, in computing the homotopy groups of a
space {\it X}, we may assume without loss of generality that {\it X}
is triangulable and hence locally contractible. In fact, we may
replace {\it X} by {\it S}({\it X}).
\end{quotation}

Recall that the Milnor realization of
{\it S}({\it X}) of any space {\it X} is weakly homotopy equivalent
to {\it X} and triangulable, as is any Milnor
realization\label{MilnorRealization} of any simplicial set, cf.
\cite{11}
p.\ 209 Cor.\ 4.6.12.
Triangulable spaces are frequently encountered and their singular
homology is often most easily determined
by calculating the simplicial homology of a triangulation of the
space.
Well, the usefulness also goes the other way: E.G. Sklyarenko: {\it
Homology and Cohomology Theories of General Spaces}
(\cite{Sklyarenko} p.\ 132):
\begin{quotation}
The singular theory is irreplaceable in problems of homotopic
topology, in the study of spaces of continuous maps and in the
theory of fibrations. Its importance, however, is not limited to the
fact that it is applicable beyond the realm of the category of
polyhedra. The singular theory is necessary for the in depth
understanding of the homology and cohomology theory of polyhedra
themselves. \dots
In particular it is used in the problems of homotopic classification
of continuous maps of polyhedra and in the description of homology
and cohomology of cell complexes. The fact that the simplicial
homology and cohomology share the same properties is proved by
showing that they are isomorphic to the singular one.
\end{quotation}

\goodbreak
\subsection{The simplicial category}\label{UnderavdP4:TheSimplicialCategory}

The simplicial category
is a somewhat controversial concept among algebraic topologists.
Indeed, this becomes apparent as Fritsch and Piccinini are quite explicit about it when defining their version of the
\htmladdnormallink{{\it simplicial category}}
{http://planetmath.org/encyclopedia/SimplicialCategory.html},
which they call the {\it finite ordinals} and we quote \hbox{from \cite{11}
p.\ $\!$220-1:}

\begin{quote}
A short but systematic treatment of the category of finite ordinals
can be found in MacLane (1971) under the name ``the simplicial
category''. Our category of finite ordinals is not exactly the same
as MacLane's; however, it is isomorphic to the subcategory obtained
from MacLane's category by removing its initial object.
\end{quote}

The last sentence pinpoints the disguised ``{\it initial
simplex}''-deletion described in the following quotation
from Mac\ $\!$Lane \cite{21}
p.\ 178 (=2:nd ed. of ``Mac\ $\!$Lane (1971)''):
\begin{quote} \label{TopologistDelta}
By $\hbox{\tenib\char"01}^+\!$ we denote the full subcategory of
$\hbox{\tenib\char"01}$ with objects all the positive
ordinals$^{\!}$
$\{1,2,3,\dots\}$ $($omit only zero$)$. Topologists use this
category, call it $\hbox{\tenib\char"01}$, and rewrite its objects $($using the geometric dimension$)$ as $\{0,1,2,\dots\}$.
Here we stick to our $\hbox{\tenib\char"01}$, which contains the real $0$, an object which is necessary if all face and degeneracy operations are to be expressed as in $(3)$, in terms of binary product $\mu$ and unit $\eta$.
\end{quote}

This mentioned Eq.\ (3) is part of
Mac\ ${\!}$Lane's description in \cite{21}
p.\ 175ff, of the {\it Simplicial Category} with objects being the
generic ordered abstractions of the
{\it simplices} in Definition~1  p.~\pageref{DefP1:1},
and
with the empty set $\emptyset$, denoted $0$ by Mac\ $\!$Lane, as an
initial object which induces a unit operation $\eta$ in a
``universal monoid''-structure imposed by {\it ordinal
addition}\label{P4:OrdinalAddition}. \noindent

``Ordinal addition''\label{OrdinalAddition:Addendum} is very much
equivalent of the union ``$\sigma_1\cup\sigma_2$'' involved in the
join-definition in Definition\ 1 p.~\pageref{DefP1:1}, once a
compatible strict vertex order is imposed on each maximal simplex in
the join-factors together with the agrement that all the vertices of
the left factor is strictly less than all those of the right one.
Since this {\it vertex-ordering} is needed also in the definition
p.~\pageref{DefP18:OrderedSimpCartProd} of
the binary simplex-operation called {\it ordered simplicial
cartesian product}, we have used these {\it ordered simplicial
complexes} all through this paper. In particular we used it when we
defined the action of the boundary function related to the chain
complexes in our homology definitions.

The empty simplex plays a central role in modern mathematical
physics, as hinted at in \cite{koch:frobalgand2dimtqft} p.~188 where
J. Kock describe the classical use of the ``topologist $\Delta$''
within category theory, there named ``topologist's delta'' and
denoted $\triangle$, in the following words:
\begin{quote}
In this way, the topologist's delta is a sort of bridge between category theory and topology. In these contexts the empty ordinal (empty simplex) is not used, but in our context it is important to keep it, because without it we could not have a
\htmladdnormallink{{\it monoidal}}
{http://en.wikipedia.org/wiki/Monoidal_category},\
 structure.
\end{quote}
\enlargethispage{\baselineskip}

In his series
\htmladdnormallink{{\it This Week's Finds in Mathematical
Physics}}{http://math.ucr.edu/home/baez/TWF.html},\ John Baez gives a compact but very lucid and accessible presentation (Week~115ff) of the modern category theory based homotopy theory, which, on his suggestion, could be called {\it homotopical algebra}. He ends section {\bf B} in Week~115 with the following, after a description of {\it The Category of Simplices}, which he calls Delta, the
objects of which are the nonempty simplices, ``$1$'' consisting of a single vertex, ``the simplex with two vertices'' $:=$ an intervall, ``the simplex with three vertices'' $:=$ a triangle, and so on:
\begin{quote}
We can be slicker if we are willing to work with a category
*equivalent* to Delta (in the technical sense described in
``week76''), namely, the category of *all* nonempty totally ordered sets, with order-preserving maps as morphisms. This has a lot more objects than just $\{0\}, \{0,1\}, \{0,1,2\}$, etc., but all of its objects are isomorphic to one of these. In category theory, equivalent categories are the same for all practical purposes - so we brazenly call this category Delta, too. If we do so, we have following *incredibly* slick description of the category of simplices: it's just the category of finite nonempty totally ordered sets!

If you are a true mathematician, you will wonder ``why not use the empty set, too?'' Generally it's bad to leave out the empty set. It may seem like "nothing", but "nothing" is usually very important. Here it corresponds to the ``empty simplex'', with no vertices! Topologists often leave this one out, but sometimes regret it later and put it back in (the buzzword is ``augmentation''). True category theorists, like Mac Lane, never leave it out. They define Delta to be the category of *all* totally ordered finite sets. For a beautiful introduction to this approach, try: \cite{21} Saunders Mac Lane, Categories for the Working Mathematician, Springer, Berlin, 1988.
\end{quote}

\subsection{Simplicial sets}\label{UnderavdP4:RelCombToLogics}

A
\htmladdnormallink{{\it simplicial set}}
{http://en.wikipedia.org/wiki/Simplicial_set},
is a
\htmladdnormallink{{\it contravariant functor}}
{http://en.wikipedia.org/wiki/Functor}
from the
\htmladdnormallink{{\it simplicial cat-}}
{http://planetmath.org/encyclopedia/SimplicialCategory.html}
\htmladdnormallink{{\it egory}}
{http://planetmath.org/encyclopedia/SimplicialCategory.html},
to the
\htmladdnormallink{{\it category of sets}}
{http://planetmath.org/encyclopedia/SimplicialCategory.html}.
The
\htmladdnormallink{{\it natural transformations}}
{http://en.wikipedia.org/wiki/Natural_transformation}
constitute the morphisms.

So, the classical simplicial sets are objects in the following functor category:
$$(\text{\rm Category of sets})^{\text{(Simplicial Category)}^{\text{op}}}.$$

A simplicial set is usually perceived through a representation in its target as a dimension-indexed collection of its non-degenerate simplices, which itself usually is envisioned through its realization, i.e. as a geometric simplicial complex, cf. p.~\pageref{MilnorRealization} middle and Theorem~\ref{Theorem:CountableRealization}
p.~\pageref{Theorem:CountableRealization}.
The morphisms in the category of simplicial sets, i.e. the natural transformations, now materializes as realisations of simplicial maps.
By mainly just well-ordering \label{well-ordering} the vertices in each simplex of a simplicial complex -- in an overall consistent way, cf. \cite{11} p.\ 152 Ex.\ 4,
each ``classical simplicial complex''
defined as in Definition~2 p.~\pageref{DefP2:2},
can be regarded as a ``classical simplicial set'' in the functor category of classical {\it Simplicial Sets}
\label{DefP4:SimpSet}.

The augmented simplicial complexes,
resulting from Definition~1 p.~\pageref{DefP1:1},
can, equivalently, be associated with an ``augmented simplicial set'' in the following functor category:
$$(\hbox{\rm Category of sets})^{\text{\rm (Augmented Simplicial Category)}\!^{{\rm op}}}\label{DefP4:AugmSimpSet},$$
as S. Mac$\ {\!}$Lane does in \cite{21} p.~179.
There Mac$\ {\!}$Lane, in p.~178, \underbar{\it chooses} the  empty topological space $\emptyset$ as the standard $(-1)$-dimensional affine simplex, implying that the empty topological space $\emptyset$ would have a cyclic chain generator in degree $-1$, which isn't compatible with the ``Eilenberg-Steenrod''-formalism for pair-spaces, where $\emptyset$ has a fixt privileged role, including nothing but trivial homology groups in all dimensions.

In
\htmladdnormallink{``Week 116~{\bf
E}''}{http://math.ucr.edu/home/baez/week116.html}
J. Baez gives an alternative specification of a simplicial set as: -``... a presheaf on the category Delta'' and then he also pays attention to the {\it realization functor}.

In
\htmladdnormallink{``Week 117~{\bf
J}''}{http://math.ucr.edu/home/baez/week117.html}
John Baez describe the
\htmladdnormallink{{\it Nerv of a
Category}}{http://planetmath.org/?op=getobj&from=objects&name=Nerve},
as a functor from the Category of Small Categories (Cat) to the Category of Simplicial Sets, and ends the paragraph with the following:
\begin{quotation}
I should also point out that topologists usually do this stuff
with the topologist's version of Delta, which does not include the ``empty simplex''.
\end{quotation}
This exclusion of the ``empty simplex'' within general topology  is probably a consequence of the fact that there is no ``natural'' $(-1)$-dimensional topological space in classical general topology that can serve as the realisation of the simplicial complex $\{\emptyset_o\}$.

\break
\subsection{A category theoretical realization functor}\label{CatDefRealization}

Our plea for the use of an \underbar{augmented} topological category rests on the otherwise lost {\it realization functor}. In \cite{KanAdjFunc}, D.M. Kan introduced the {\it realization functor} and the {\it simplicial singular functor} as an example of an adjoint pair of functors. In \cite{11} p.~303, with its suggestive notations, Kan's example is given a somewhat generalized category theoretical formulation as follows:

\smallskip
 Let $D$ be a small category and let $DSets$ denote the category of contravariant functors $D\longrightarrow Sets$. One might view the objects of $DSets$ as sets graded by the objects of $D$ with the morphisms of $D$ operating on the right. Let $\Phi: D\longrightarrow Sets$ be an arbitrarily given covariant functor; associated to it, construct a pair of adjoint functors
\begin{eqnarray}
\Gamma_\Phi: DSets\longrightarrow Sets,\hskip0.0cm \nonumber\\
 S_\Phi: Sets\longrightarrow DSets, \nonumber
\end{eqnarray}
\vskip-0.4cm
as follows:
\begin{itemize}
\item[(l)]\ The left adjoint functor $\Gamma_\Phi$ - called {\it realization functor} associates to each object $X$ of $DSets$ the set of equivalence classes of pairs $(x,t)\sqcup X(d) \times \Phi(d)$ modulo the relation\hskip0.4cm
$(x\alpha,t)\sim(x,\Phi(\alpha)(t)).$\\
Given a morphism $f : Y \longrightarrow X$, i.e., a natural transformation, one has the function\hskip0.6cm
$\Gamma_\Phi f : \Gamma_\Phi Y \longrightarrow\Gamma_\Phi X, \ [x,t]\longrightarrow[f(x),t]$\hskip0.3cm
where $[x,t]$ denotes the equivalence class represented by the pair $(x,t)$.
\item[(2)]\ The right adjoint functor $S_\Phi$ - called {\it singular functor} - associates to each set $Z$ the object of $DSets$ given by\hskip0.6cm
$(S_\Phi Z)(d)=Z^{\Phi(d)},$\hskip0.3cm
for each object $d$ of $D$, and\hskip0.3cm
$\alpha^\ast : Z^{\Phi(d)}\longrightarrow Z^{\Phi(d^\prime)},\hskip0.2cm
x\longrightarrow x\circ\Phi(a),$\hskip0.4cm
for each morphism $\alpha:d^\prime\longrightarrow d$ of $D$.
\end{itemize}
\vskip-0.2cm
\vskip-0.2cm



\begin{thebibliography}{10}
%
\providecommand{\bysame}{\leavevmode\hbox
to3em{\hrulefill}\thinspace}

\bibitem{Baclawski}
K.~Baclawski, \emph{Cohen-{M}acaulay {C}onnectivity and {G}eometric
  {L}attices}, Europ. J. Combinatorics \textbf{3} (1982), 293--305.

\bibitem{Baez}
J.~Baez, 1998, John Baez' web site,
\htmladdnormallink{http://math.ucr.edu/home/baez/}
  {http://math.ucr.edu/home/baez/}.

\bibitem{1}
G.E. Bredon, \emph{Sheaf {T}heory}, Springer Verlag, N.Y., 1997.

\bibitem{Brown:TenProductTopologies}
R.~Brown, \emph{Ten {T}opologies for ${X}\times {Y}$}, Quart. J.
Math. Oxford,
  (2) \textbf{14} (1963), 303--319.

\bibitem{2}
\bysame, \emph{Topology}, Ellis Horwood Limited, 1988.

\bibitem{3}
W.~Bruns and J.~Herzog, \emph{Cohen-{M}acaulay {R}ings}, Cambr.
Univ. Press,
  1998.

\bibitem{CartanandEilenberg}
H.~Cartan and S.~Eilenberg, \emph{Homological {A}lgebra}, Princeton
University
  Press, New Jersey, 1956.

\bibitem{4}
D.E. Cohen, \emph{Products and {C}arrier {T}heory}, Proc. London
Math. Soc.
  \textbf{VII} (1957), 219--248.

\bibitem{5}
G.E. Cook and R.L. Finney, \emph{Homology of {C}ell {C}omplexes},
Princeton
  University Press, New Jersey, 1967.

\bibitem{6}
A.~Dold, \emph{Lectures on {A}lgebraic {T}opology}, Springer-Verlag,
1972.

\bibitem{7}
J.~Dugundji, \emph{Topology}, Allyn \& Bacon Inc., Boston, 1966.

\bibitem{Ehlers&Porter}
P.J. Ehlers and T.~Porter, \emph{Joins for ({A}ugmented)
{S}implicial {S}ets},
  J. Pure Appl. Alg. (1) \textbf{245} (2000), 37--44, 
  \htmladdnormallink{http://front.math.ucdavis.edu/math.CT/9904039}
  {http://front.math.ucdavis.edu/math.CT/9904039}.

\bibitem{8}
S.~Eilenberg and S.~MacLane, \emph{Group {E}xtensions and
{H}omology}, Ann. of
  Math. (4) \textbf{43} (1942).

\bibitem{9}
S.~Eilenberg and N.~Steenrod,\ \nobreak\emph{Foundations of {A}lgebraic {T}opology},
  Princeton Univ. Press,
  \nobreak1952.\nobreak{ \htmladdnormallink{http://www.archive.org/details/foundationsofalg033540mbp}
  {http://www.archive.org/details/foundationsofalg033540mbp}}

\bibitem{10}
G.~Fors, \emph{Algebraic {T}opological {R}esults on
{S}tanley-{R}eisner
  {R}ings}, Commutative Algebra (Trieste, Italy 1992) (A.~Simis, N.V. Trung,
  and G.~Valla, eds.), World Scientific Publ. Co. Pte. Ltd., 1994, pp.~69--88.

\bibitem{Fors2}
\bysame, \emph{A {H}omology {T}heory {B}ased on the {E}xistence of a
  $(-1)$-simplex}, Tech. Report~3, Department of Mathematics, University of
  Stockholm, Sweden, 1994.

\bibitem{Fritsch&Golasinski}
R.~Fritsch and M.~Golasi$\acute{\rm n}$ski, \emph{Topological, simplicial and
  categorical joins}, Archiv der Mathematik \textbf{82} (2004), 468--480. \ \ \ ps-preprint at:\ \
\htmladdnormallink{http://www.mathematik.uni-muenchen.de/$\sim$fritsch/joins.ps}
  {http://www.mathematik.uni-muenchen.de/~fritsch/joins.ps}.

\bibitem{11}
R.~Fritsch and R.A. Piccinini, \emph{Cellular {S}tructures in
{T}opology},
  Cambr. Univ. Press., 1990.

\bibitem{12}
G.~Gräbe H, \emph{$\ddot{U}$ber den {S}tanley-{R}eisner-{R}ing von
  {Q}uasimannigfaltigkeiten}, Math. Nachr. \textbf{117} (1984), 161--174.

\bibitem{13}
\bysame, \emph{$\ddot{U}$ber den {R}and von
{H}omologiemannigfaltigkeiten},
  Beiträge Algebra. Geom. \textbf{22} (1986), 29--37.

\bibitem{14}
T.~Hibi, \emph{Union and gluing of a family of
{C}ohen-{M}acaulay partially ordered sets}, Nagoya Math. J. \textbf{107} (1987), 91--119.

\bibitem{15}
\bysame, \emph{Level rings and algebras with straightening
laws}, J. Algebra \textbf{117} (1988), 343--362.

\bibitem{16}
\bysame, \emph{Combinatorics on {C}onvex {P}olytopes}, Carslaw
Publ., 1992.

\bibitem{17}
P.J. Hilton and S.~Wylie, \emph{Homology {T}heory}, Cambr. Univ.
Press., 1960.

\bibitem{Ho}
M. Hochster; \emph{Cohen-{M}acaulay {R}ings, {C}ombinatorics, and {S}implicial {C}omplexes}
Second Oklahoma Ring Theory Conf. March, 1976 (Dekker, N.Y., 1977).

\bibitem{18}
S.-T. Hu, \emph{Homotopy {T}heory}, Academic Press, 1959.

\bibitem{JandT}
A.~Joyal and M.~Tierney, \emph{An introduction to simplicial
homotopy theory}, Aug.\ 5 1999, \ \ \ pdf-preprint of chap. 01 at:
  \htmladdnormallink{http://hopf.math.purdue.edu//Joyal-Tierney/JT-chap-01.pdf}
{http://hopf.math.purdue.edu//Joyal-Tierney/JT-chap-01.pdf}.

\bibitem{KanAdjFunc}
D.M. Kan, \emph{Adjoint functors}, Trans. Amer. Math. Soc. \textbf{87} (1958),
  294--329.

\bibitem{19}
J.L. Kelley, \emph{General {T}opology}, Van. Nostrand, 1955.

\bibitem{koch:frobalgand2dimtqft}
J.~Kock, \emph{Frobenius {A}lgebras and 2{D} {T}opological {Q}uantum {F}ield
  {T}heories}, Cambridge Univ. Press, 2003. Frontmatter and index
    at:\hfill\break \htmladdnormallink{http://www.cambridge.org/uk/catalogue/catalogue.asp?isbn=9780521540315}
{http://www.cambridge.org/uk/catalogue/catalogue.asp?isbn=9780521540315}.

\bibitem{20}
J.~Lawson and B.~Madison, \emph{Comparisons of notions of weak hausdorffness
  topology}, Lecture Notes in Pure and Appl. Math. \textbf{24} (1976),
  207--215, Marcel Dekker Inc. New York.

\bibitem{21}
\htmladdnormallink{S.~MacLane,}
{http://en.wikipedia.org/wiki/Saunders_Mac_Lane}
\emph{Categories for the {W}orking {M}athematician},
  Springer, 1988.

\bibitem{MacLane}
S.~MacLane, \emph{The development of mathematical ideas by collision: {T}he
  case of categories and topos theory}, Categorical {T}opology (S.~MacLane and
  J.~Adámek, eds.), World Scient. Publ., 1-9 {A}ugust 1988, {P}rague,
  {C}zechoslovakia, 1989.

\bibitem{22}
C.R.F. Maunder, \emph{Algebraic {T}opology}, Van Nostrand Reinhold, 1970.

\bibitem{23}
J.W. Milnor, \emph{Construction of universal boundles. {II}}, Ann. of Math. (2)
  \textbf{63} (1956), 430--436.

\bibitem{24}
\bysame, \emph{The geometric realization of a semi-simplicial complex}, Ann. of  Math. (2) \textbf{65} (1957), 357--362.

\bibitem{Mi}
M. Miyazaki, \emph{On 2-{B}uchsbaum complexes}, J. Math. Kyoto Univ. 30-3 (1990), 367--392.
\htmladdnormallink{http://projecteuclid.org/euclid.kjm/1250520019}
{http://projecteuclid.org/euclid.kjm/1250520019}.

\bibitem{25}
J.R. Munkres, \emph{Elementary {D}ifferential {T}opology}, Ann. of Math. Stud., vol.~54, Princeton University Press, New Jersey, 1966.

\bibitem{26}
\bysame, \emph{Elements of {A}lgebraic {T}opology}, Benjamin/Cummings, 1984.

\bibitem{27}
\bysame, \emph{Topological results in combinatorics}, Michigan Math. J. \textbf{31} (\nolinebreak\ 1984), 113--128.

\bibitem{28}
A.A. Ranicki, \emph{On the {H}auptvermutung}, The {H}auptvermutung {B}ook (A.A. Ranicki, ed.), Dordrecht, 1996, pp.~3--31.

\bibitem{Rotman:HomAlg}
J.J. Rotman, \emph{Introd. to {H}omological {A}lgebra}, Academic Press, Inc, \nolinebreak\ 1979.

\bibitem{29}
M.E. Rudin, \emph{An unshellable triangulation of a tetrahedron}, Bull. Am. Math. Soc. \textbf{64} (1958), 90--91. Free pdf-version from AMS at
\htmladdnormallink{http://www.ams.org/journals/bull/1958-64-03/home.html}
{http://www.ams.org/journals/bull/1958-64-03/home.html}.

\bibitem{Sklyarenko}
E.G. Sklyarenko, \emph{Homology and {C}ohomology {T}heories of {G}eneral
  {S}paces}, {E}nc. {M}ath. Sc.~50, Springer, 1996, pp.~124--256.

\bibitem{30}
E.H. Spanier, \emph{Algebraic {T}opology}, McGraw-Hill, 1966.

\bibitem{31}
R.~Stanley, \emph{Combinatorics and {C}ommutative {A}lgebra}, 2nd. ed
  Birkhäuser, 1996.

\bibitem{32}
B.~St$\ddot{\rm u}$ckrad and W.~Vogel, \emph{Buchsbaum {R}ings and
  {A}pplications}, VEB/Springer, 1986.

\bibitem{33}
R.M. Vogt, \emph{Convenient categories of topological spaces for homotopy
  theory}, Arch. Math. \textbf{XXII} (1971), 545--555.

\bibitem{34}
J.E. Walker, \emph{Canonical homeomorphisms of posets}, European J. Combin.
  (1988), 97--107.

\bibitem{36}
G.W. Whitehead, \emph{Homotopy groups of joins and unions}, Trans. Amer. Math.
  Soc. \textbf{83} (1956), 55--69.

\bibitem{35}
\bysame, \emph{Elements of {H}omotopy {T}heory}, GTM\ 61, Springer,
  \nolinebreak\ 1978.

\bibitem{ZariskiAndSamuel} O. Zariski and P. Samuel;
{\it Commutative {A}lgebra}; Vol. {\bf I}, Springer -V., 1986.

\end{thebibliography}
\end{document}